%% file: sums_of_fft.tex
\documentclass[11pt]{article}
\usepackage{amsmath, amssymb,theorem,apacite}
\usepackage{graphicx}
\usepackage{caption}
\usepackage{subcaption}
\usepackage{multirow}
\newcommand{\cov}{\mathrm{cov}}
\newcommand{\var}{\mathrm{var}}
\newcommand{\cum}{\mathrm{cum}}
\newcommand{\Ex}{\mathrm{E}}

\newtheorem{theorem}{Theorem}[section]
\newtheorem{lemma}{Lemma}[section] 
\newtheorem{assumption}{Assumption}[section] 
 
\newtheorem{corollary}{Corollary}[section] 
\newtheorem{example}{Example}[section]
\newtheorem{remark}{Remark}[section] 

\newcommand{\Pcon}{\stackrel{\mathcal{P}}{\rightarrow}}
\newcommand{\Dcon}{\stackrel{\mathcal{D}}{\rightarrow}}

\newcommand{\sinc}{\mathrm{sinc}}
\newcommand{\Sinc}{\mathrm{Sinc}}
\newcommand{\ub}{\boldsymbol{s}}
\newcommand{\ubb}{\boldsymbol{u}}
\newcommand{\vbb}{\boldsymbol{v}}
\newcommand{\jb}{\boldsymbol{j}}
\newcommand{\rb}{\boldsymbol{r}}
\newcommand{\ob}{\boldsymbol{\omega}}

\newcommand{\kb}{\boldsymbol{k}}
\newcommand{\sbb}{\boldsymbol{s}}

\newcommand{\mb}{\boldsymbol{m}}

\title{{\bf Statistical inference for spatial statistics defined in
    the Fourier domain}}
\author{Suhasini Subba Rao\\
Department of Statistics, \\
Texas A\&M University, \\
College Station, U.S.A.\\
{\tt suhasini@stat.tamu.edu}} 
\date{\today}

\begin{document}
\maketitle
\begin{abstract}
A class of Fourier based statistics for irregular spaced spatial data
is introduced, examples include, the Whittle likelihood, a
parametric estimator of the covariance function based on the 
$L_{2}$-contrast function and  a simple nonparametric estimator of
the spatial autocovariance which is a non-negative function.
The Fourier based statistic is a quadratic form of 
a discrete Fourier-type transform of the spatial data.
Evaluation of the statistic is computationally
tractable, requiring $O(nb^{})$ operations, where 
$b$ are the number Fourier
frequencies used in the definition of the statistic 
 and $n$ is the sample size. The asymptotic sampling properties of
the statistic are derived using  both increasing domain and fixed domain spatial asymptotics.  
These results are used to construct a statistic which is asymptotically pivotal.

{\it Keywords and phrases:} Fixed and increasing domain
asymptotics, irregular spaced locations, quadratic forms, 
spectral density function,
stationary spatial random fields.
\end{abstract}

\input{1_introduction_4}

\input{2_statisticR2}

\input{3_summary}

\input{4_uniform}

\input{5_student}

\subsection*{Acknowledgments}

This work has been partially supported by the National Science
Foundation, DMS-1106518 and DMS-1513647. The author gratefully acknowledges  
Gregory Berkolaiko and Mikyoung Jun for many very useful discussions. Furthermore, the
author gratefully acknowledges the editor, associate editor and three anonymous
referees who made substantial improvements to every part of the manuscript.

\appendix

\section{Appendix}\label{appendix:A}

\input{A_1_mean}

\input{A_2_variance}

\bibliographystyle{apacite}
\bibliography{spatial}

\newpage

\section*{Supplementary material}

\input{A_3_nonuniform}

\input{B_1_spatial}

\input{B_2_cumulant}

\input{B_3}

\input{B_5_grids}

\input{B_6_fixed_domainv2}

\input{B_4_Q}

\input{B_7_simulations}

\input{tables}

\end{document}

%% file: 1_introduction_4.tex
\section{Introduction}

In recent years irregular spaced spatial 
data has become ubiquitous in several disciplines as varied as the geosciences to econometrics. 
The analysis of such data poses several challenges which do not arise in data which 
is sampled on a regular lattice. A major obstacle is the 
computational costs when dealing with large irregular sampled
data sets. If spatial data are sampled on a regular lattice then
algorithms such as the Fast Fourier transform can be employed to reduce
the computational burden (see, for example, \citeA{p:che-hur-06}). 
Unfortunately, such algorithms have little benefit if the spatial data are
irregularly sampled. To address this issue, within the spatial domain, several authors, including, 
\citeA{p:vec-88},
\citeA{p:cre-hua-99}, \citeA{p:ste-04}, have proposed estimation
methods which are designed to reduce the computational burden.

In contrast to the above references, \citeA{p:fue-07} and
\citeA{p:mat-yaj-09} argue that
working within the frequency domain often simplifies the computational
burden. Both authors focus on parametric estimation using a
Whittle-type likelihood. 
\citeA{p:fue-07} assumes that the irregular spaced data can be embedded on a grid and the 
missing mechanism is deterministic and `locally smooth'. 
A possible drawback of this construction, is that the local smooth
assumption will not hold if the locations are extremely
irregular. Therefore \citeA{p:mat-yaj-09} propose a Whittle likelihood
approach to parameter estimation which takes into account the
irregular nature of the locations. 
The focus of most Fourier domain estimators have been on the Whittle
likelihood (the exception being the recent paper by
\citeA{p:ban-lah-nor-13}, which we discuss later). 
In this paper we argue that several estimators, both parametric and
nonparametric can be defined within the Fourier domain. 
 For example, within the Fourier domain, we propose a nonparametric, non-negative definite estimator of the spatial
covariance. Nonparametric estimators of the spatial autocovariance are
often defined using kernel smoothing  methods (defined in \citeA{p:hal-fis-hof-94}) or
the empirical variogram (defined in \citeA{b:cre-93}). However, these
"raw" covariance estimators may not be
non-negative functions, and a second step is required, which involves 
taking the Fourier transform of a finite discretisation of the sample
autocovariance, setting negative values to zero and
inverting back, to ensure that the resulting estimator
is a non-negative function. In contrast, by 
defining the covariance estimator within the Fourier domain the estimator is
guaranteed to be a non-negative definite function, thus rendering the
need to make the second nonlinear step unnecessary. The 
  purpose of this paper is two fold. The first is to demonstrate that
  several parameters can be estimated within the Fourier
  domain. The second is to obtain a comprehensive understanding of
  quadratic forms of irregular sampled spatial processes.

In order to define estimators within the Fourier domain we adopt
the approach pioneered by \citeA{p:mat-yaj-09} and
\citeA{p:ban-lah-10} who assume that
the irregular locations are independent, identically distributed
random variables (thus allowing the data to be
extremely irregular) and define the irregular sampled discrete Fourier transform (DFT) as 
\begin{eqnarray}
\label{eq:dft}
J_{n}(\ob) = \frac{\lambda^{d/2}}{n}\sum_{j=1}^{n}Z(\ub_{j})\exp(i\ub_{j}^{\prime}\ob),
\end{eqnarray}
where $\ub_{j}\in [-\lambda/2,\lambda/2]^{d}$ denotes the spatial locations observed in the 
space $[-\lambda/2,\lambda/2]^{d}$ and $\{Z(\ub_{j})\}$ denotes the spatial random field at these 
locations. It's worth mentioning a similar transformation on irregular sampled data goes back to
\citeA{p:mas-78}, who defines the discrete Fourier transform of 
Poisson sampled continuous time series. Using this definition, 
\citeA{p:mat-yaj-09} define the Whittle likelihood by taking the
weighted integral of the periodogram, $|J_{n}(\ob)|^{2}$. 
Of course, in practice the weighted integral needs to be approximated
by a Riemann sum. Indeed in Remark 2, \citeA{p:mat-yaj-09}, suggest using
the frequency grid $\{\ob_{\kb} = 2\pi\kb/\lambda;\kb\in
\mathbb{Z}^{d}\}$ when constructing the Whittle likelihood, though no
justification is given for this discretisation. 

%In this paper we show that that approximating an integral by a sum can
%change the sampling properties of the proposed estimator.  
%A clear advantage of this approach is that  it avoids the inversion of a large matrix. However, one still needs to 
%evaluate the integral. 

Motivated by the integrated Whittle likelihood, our aim is to consider
estimators with the form $\int g_{\theta}(\ob)|J_{n}(\ob)|^{2}d\ob$
(such quantities have been widely studied in time series, dating 
back to at least as far as \citeA{p:par-61}, but has received very little
attention in the spatial literature). 
In practice this integral cannot be evaluated, and needs to be approximated
by a Riemann sum
\begin{eqnarray}
\label{eq:QQQ}
Q_{a,\Omega,\lambda}(g_{\theta};0) =  \frac{1}{\Omega^{d}}\sum_{k_{1},\ldots,k_{d}=-a}^{a}
g_{\theta}(\ob_{\Omega,\kb})|J_{n}(\ob_{\Omega,\kb})|^{2} 
\end{eqnarray}
where $\{\ob_{\Omega,\kb} = 2\pi \kb/\Omega, \kb =
(k_{1},\ldots,k_{d}), -a\leq  k_{i}\leq a\}$ is the frequency grid
over which the sum is evaluated. 
%Note that this class ofstatistics includes the Whittle likelihood (with $g_{\theta}(\ob_{\Omega,\kb})
%= [f_{\theta}(\ob_{\Omega,kb})+\eta^{2}]^{-1}$, where
%$\{f_{\theta}(\cdot);\theta\in \Theta\}$ denotes the known parametric family) and a
%nonparametric estimator of the covariance.  
In terms of computation, evaluation of
$\{J_{n}(\ob_{\Omega,\kb});\kb=(k_{1},\ldots,k_{d}),k_{j}=-a,\ldots,a\}$
requires $O(a^{d}n)$ operations. However, once
$\{J_{n}(\ob_{\Omega,\kb})\}$ has been
evaluated the evaluation of $Q_{a,\Omega,\lambda}(g_{\theta};0)$ only requires
$O(a^{d})$ operations. 

As far as we are aware there exists no results on the sampling
properties of the general quadratic form defined in (\ref{eq:QQQ}). 
To derive the asymptotic sampling properties of 
$Q_{a,\Omega,\lambda}(g_{\theta};0)$ we will work under two
asymptotic frameworks that are commonly used in spatial statistics. 
Our main focus will be the increasing domain
framework, introduced in 
\citeA{p:hal-pat-94} (see also \citeA{p:hal-fis-hof-94} and used in, for example, \citeA{p:lah-03}, 
\citeA{p:mat-yaj-09}, \citeA{p:ban-lah-10}, \citeA{p:ban-lah-nor-13} and 
\citeA{p:ban-sub-15}). This is where the number of  observed locations 
$n\rightarrow \infty$ as the size of the spatial domain
$\lambda \rightarrow \infty$ (we usually assume
$\lambda^{d}/n \rightarrow 0$). We also analysis the sampling
properties of $Q_{a,\Omega,\lambda}(g_{\theta};0)$ within the 
fixed-domain asymptotic framework
(where $\lambda$ is kept fixed but the number of locations, $n$ grows)
considered in Stein (1994, 1999), \nocite{b:ste-99}\nocite{p:ste-94}
\citeA{p:zha-04} and \citeA{p:zha-05}. Not too surprisingly, the sampling properties of
$Q_{a,\Omega,\lambda}(g_{\theta};0)$ differ according to the framework
used. 
%in particular the estimator $Q_{a,\Omega,\lambda}(g_{\theta};0)$
%consistently estimates the parameter of interest in the increasing
%domain framework but not in the fixed domain framework. 
%However, we show that the fixed domain analysis does allow us to correct for finite
%domain bias in some situations and understand the role that the number
%of locations $n$ play in the choice of $a$.  

We show in Section \ref{sec:summary} and \ref{sec:uniform} that
$Q_{a,\Omega,\lambda}(g_{\theta};0)$ is a consistent estimator of 
the functional $I(g_{\theta};\frac{a}{\Omega})$ as $\lambda\rightarrow
\infty$ and $\Omega \rightarrow \infty$,
where
\begin{eqnarray}
\label{eq:I}
I\left(g_{\theta};\frac{a}{\Omega} \right) =
\frac{1}{(2\pi)^{d}}\int_{[-2\pi a/\Omega,2\pi a/\Omega]^{d}}g_{\theta}(\ob)f(\ob)d\ob.
\end{eqnarray}
However, the frequency grid $\ob_{\Omega,\kb}$ plays a vital role in the
rate of convergence. In particular, we show that 
\begin{eqnarray*}
\Ex[Q_{a,\Omega,\lambda}(g_{\theta};0)] =
I\left(g_{\theta};\frac{a}{\Omega} \right) + O\left(\frac{\log
    \lambda}{\lambda} + \frac{1}{\Omega} + \frac{1}{n}\right)
\end{eqnarray*}
and 
\begin{eqnarray}
\label{eq:varQQ}
\var[Q_{a,\Omega,\lambda}(g_{\theta};0)] = 
\left\{
\begin{array}{cc}
O\left(\lambda^{-d}\right) & \Omega \geq \lambda \\
O\left(\Omega^{-d}\right) & \Omega \leq \lambda \\
\end{array}
\right..
\end{eqnarray}
In other words, using a frequency grid which is coarser than $\{2\pi
\kb/\lambda;\kb \in \mathbb{Z}^{d}\}$ leads to an estimator with a larger bias and variance,
whereas using a frequency grid finer than $\{2\pi \kb/\lambda;\kb \in \mathbb{Z}^{d}\}$ is
computationally more cumbersome but with no real improvement in mean
squared error. Thus balancing efficiency with computational burden, in general, it is optimal to use the frequency
grid $\{\ob_{\lambda,\kb} = 2\pi \kb/\lambda\}$ and in Section \ref{sec:uniform}
we focus on deriving the sampling properties on this frequency grid. 
This result justifies the insightful suggestion of Remark 2, \citeA{p:mat-yaj-09}.
As mentioned above \citeA{p:ban-lah-nor-13} also use the Fourier domain for spatial
inference, however their objectives are very different to those in
this paper. \citeA{p:ban-lah-nor-13} 
 show that $\{J_{n}(\ob_{\Omega,\kb})\}_{\kb}$ are
asymptotically independent if $\Omega/\lambda \rightarrow 0$ as
$\Omega \rightarrow \infty$ and $\lambda \rightarrow \infty$
(this corresponds to a very coarse
frequency grid). Based on this
property they use $Q_{a,\Omega,\lambda}(g;0)$, where $\Omega$ is such that
$\Omega/\lambda \rightarrow 0$, to construct the empirical
likelihood. The justification for their construction is that the
distribution of the resulting empirical likelihood is asymptotically
pivotal as $\lambda \rightarrow \infty$. 
The sampling properties of $Q_{a,\Omega,\lambda}(g;0)$ are
not derived in \citeA{p:ban-lah-nor-13}, however, it is clear from
(\ref{eq:varQQ}) for $\Omega<<\lambda$ the resulting estimator
is far from efficient. Therefore a drawback of using this
construction for estimation is a
substantial increase in mean square error.  

Since $Q_{a,\Omega,\lambda}(g;0)$ is optimal when using the frequency grid
$\Omega = \lambda$, in Section \ref{sec:uniform} we focus on deriving
the sampling properties of
$Q_{a,\lambda,\lambda}(g_{\theta};0)$. We consider the slightly more general
statistic
\begin{eqnarray}
\label{eq:QQQr}
Q_{a,\lambda,\lambda}(g_{\theta};\rb) =  \frac{1}{\lambda^{d}}\sum_{k_{1},\ldots,k_{d}=-a}^{a}
g_{\theta}(\ob_{\lambda,\kb})J_{n}(\ob_{\lambda,\kb})
\overline{J_{n}(\ob_{\lambda,\kb+\rb})}, \qquad \rb \in \mathbb{Z}^{d}
\end{eqnarray}
and show asymptotic normality of
$Q_{a,\lambda,\lambda}(g_{\theta};\rb)$ when the random field is stationary and
Gaussian and obtain the second order properties of 
$Q_{a,\lambda,\lambda}(g_{\theta};\rb)$ when the random field is
stationary (but not necessarily Gaussian). The sampling properties of
$Q_{a,\lambda,\lambda}(g_{\theta};0)$ when the domain is kept fixed are
considered in Section \ref{sec:fixed-asy}.
The variance of $Q_{a,\lambda,\lambda}(g_{\theta};0)$ is usually
difficult to directly estimate. However, in Section \ref{sec:variance2} we
show that if the locations
are independent, uniformly distributed random variables, then  
$\{Q_{a,\lambda,\lambda}(g_{\theta};\rb)\}$ forms a `near
uncorrelated' sequence whose variance is asymptotically equivalent to 
$Q_{a,\lambda,\lambda}(g_{\theta};\rb)$. 
More precisely, if $Q_{a,\lambda,\lambda}(g_{\theta};0)$ is real we
define the the studentized statistic
\begin{eqnarray*}
T_{\mathcal{S}}=\frac{Q_{a,\lambda,\lambda}(g_{\theta};0) - 
I(g_{\theta};\frac{a}{\lambda})}{\sqrt{\frac{1}{|\mathcal{S}|}\sum_{\rb
  \in \mathcal{S}}|Q_{a,\lambda,\lambda}(g_{\theta};\rb)|^{2}}},
\end{eqnarray*}
for some fixed set $\mathcal{S}\subset \mathbb{Z}^{d}/\{0\}$.
We show that $T_{\mathcal{S}}\Dcon t_{2|\mathcal{S}|}$ as $\lambda
\rightarrow \infty$, where
$t_{2|\mathcal{S}|}$ denotes a $t$-distribution with $2|\mathcal{S}|$
degrees of freedom and $|\mathcal{S}|$ denotes the cardinality of $\mathcal{S}$. 

We now summarize the paper. 
In Section \ref{sec:estimator} we state the assumptions and sampling
properties of $\{J_{n}(\ob_{\Omega,\kb})\}$ which are used to motivate
the examples in Section \ref{sec:Qexamples}. In Section
\ref{sec:summary} we summarize the sampling properties of
$\{Q_{a,\Omega,\lambda}(g_{\theta};0)\}$. In Section \ref{sec:uniform}
we focus on $Q_{a,\lambda,\lambda}(g;\rb)$ and these results are used
to study $T_{2|\mathcal{S}|}$
in Section \ref{sec:variance2}. 
$Q_{a,\Omega,\lambda}(g_{\theta};0)$ is a quadratic form for
irregular sampled spatial processes. 
An outline of the proofs in the case of uniform sampling
can be found in Appendix \ref{appendix:A}. 
However, the main proofs can be found in the supplementary
material,  \citeA{p:sub-14}; many of these results build on the work
of \citeA{p:kaw-59} and may be of independent interest. A simulation
study  to illustrate the performance of the nonparametric non-negative definite estimator of the
  spatial covariance is given in Section \ref{sec:simulations}, \citeA{p:sub-14}.

%% file: 2_statisticR2.tex
\section{Assumptions and Examples}\label{sec:estimator}

\subsection{Assumptions and notation}

In this section we state the required assumptions and notation. This
section can be skipped on first reading. 

We observe the spatial random field 
$\{Z(\ub);\ub \in \mathbb{R}^{d}\}$ at the locations
$\{\ub_{j}\}_{j=1}^{n}$ where $\ub_{j}\in [-\lambda/2,\lambda/2]^{d}$.
Throughout this paper we will use the following assumptions on the spatial random field. 

\begin{assumption}[Spatial random field]\label{assum:S}
\begin{itemize}
\item[(i)]$\{Z(\ub);\ub \in \mathbb{R}^{d}\}$ is a second order stationary random field
with mean zero and covariance function $c(\ub_{1}-\ub_{2}) = \cov(Z(\ub_{1}),Z(\ub_{2})|\ub_{1},\ub_{2})$. 
We define the spectral density function as $f(\ob) = \int_{\mathbb{R}^{d}}c(\ub)\exp(-i\ub^{\prime}\ob)d\ub$. 
\item[(ii)] $\{Z(\ub);\ub\in \mathbb{R}^{d}\}$ is a stationary Gaussian random field.
\end{itemize}
\end{assumption}
%\item[(i)] $|cov(\ub_{1};\ub_{2}|
%\ub_{1},\ub_{2})| \leq {\boldsymbol \alpha}(\ub_{1}-\ub_{2})$. 

We require the following definitions. For some
finite $0<C<\infty$ and $\delta>0$, let 
\begin{eqnarray}
\label{eq:def-alpha}
%\alpha(s) =
%\left\{
%\begin{array}{cc}
%C & |s|\in [-1,1] \\
%C|s|^{-2} & |s| > 1 \\
%\end{array}
%\right.\textrm{ and }
\beta_{\delta}(s) =
\left\{
\begin{array}{cc}
C & |s|\in [-1,1] \\
C|s|^{-(1+\delta)} & |s| > 1 \\
\end{array}
\right..
\end{eqnarray}
Let $\beta_{\delta}(\ub) = \prod_{j=1}^{d}\beta_{\delta}(s_{j})$. Let 
$\xi_{\eta}(j) = C[I(j=0)+I(j\neq 0)|j|^{-\eta}]$ (for some finite
constant $C$). 
To minimise notation we will often use 
$\sum_{\kb=-a}^{a}$ to denote the multiple sum
$\sum_{k_{1}=-a}^{a}\ldots\sum_{k_{d}=-a}^{a}$. Let $\|\cdot\|_{1}$
denote the $\ell_{1}$-norm of a vector and $\|\cdot\|_{2}$ denote the $\ell_{2}$-norm.
Let $\Re X$ and $\Im X$ denote the real and imaginary parts of  $X$.
We make heavy use of the sinc function which is defined as 
\begin{eqnarray}
\label{eq:sinc-def}
\sinc(\omega) = \frac{\sin(\omega)}{\omega} \textrm{ and }\Sinc(\ob) = \prod_{j=1}^{d}\sinc(\omega_{j}).
\end{eqnarray}
We use the notation $\{\ob_{\Omega,\kb} = 2\pi\kb/\Omega;\kb\in
\mathbb{Z}^{d}\}$ for a general frequency grid. Further, as mentioned in the
introduction using $\lambda = \Omega$ is optimal, therefore to
reduce notation we let $\{\ob_{\kb} = 2\pi\kb/\lambda;\kb\in
\mathbb{Z}^{d}\}$. Closely related to the sinc function is the
triangle kernel, $T:\mathbb{R}\rightarrow \mathbb{R}$ where $T(u) = 1-|u|$
for $u\in [-1,1]$ and zero elsewhere.

We adopt the assumptions of \citeA{p:hal-pat-94}, \citeA{p:mat-yaj-09} and \citeA{p:ban-lah-10} and 
assume that $\{\ub_{j}\}$ are iid random variables with density
$\frac{1}{\lambda^{d}}h(\frac{\cdot}{\lambda})$, where
$h:[-1/2,-1/2]^{d}\rightarrow \mathbb{R}$.
%We use the following assumptions on the sampling density $h$.
\begin{assumption}[Non-uniform sampling]\label{assum:nonuniform}
The locations $\{\ub_{j}\}$ are independent distributed random variables on 
$[-\lambda/2,\lambda/2]^{d}$, where the density of $\{\ub_{j}\}$ is $\frac{1}{\lambda^{d}}h(\frac{\cdot}{\lambda})$, 
and $h(\cdot)$ admits the Fourier representation 
\begin{eqnarray*}
h(\ubb) = \sum_{\jb\in \mathbb{Z}^{d}}\gamma_{\jb}\exp(i2\pi \jb^{\prime}\ubb),
\end{eqnarray*}
where $\sum_{\jb\in \mathbb{Z}^{d}}|\gamma_{\jb}|<\infty$ such
that $|\gamma_{\jb}| \leq
C\prod_{i=1}^{d}\xi_{1+\delta}(j_{i})$ (for some $\delta>0$).
%\begin{itemize}
%\item[(i)] 
%If $\jb \neq 0$, then 
%the Fourier coefficients satisfy  $|\gamma_{\jb}| \leq
%C\prod_{i=1}^{d}[|j_{i}|^{-2}I(j_{i}\neq 0) +
%I(j_{i=0})]$\footnote{Need to include the fact that
%  $\sum_{\jb}|\gamma_{\jb}|^{2}<\infty$ and not one. Also see whether
%  this assumption can be relaxed to $|j|^{-(1+\delta)}$}. 
This assumption is satisfied if the second derivative of $h$ is bounded on
the $d$-dimensional torus $[-1/2,1/2]^{d}$. 
\end{assumption}

\begin{remark}
If $h$ is such that 
$\sup_{\ub\in [-1/2,1/2]^{d}}|\frac{\partial^{m_{1}+\ldots,m_{d}}h(s_{1},\ldots,s_{d})}{\partial
  s_{1}^{m_{1}}\ldots \partial s_{d}^{m_{d}}}|<\infty$ ($0\leq
m_{i}\leq 2$) but $h$ is not
continuous on the $d$-dimensional torus $[-1/2,1/2]^{d}$ then 
$|\gamma_{\jb}| \leq C\prod_{i=1}^{d}\xi_{1}(j_{i})$ and
the above condition will not be satisfied. However, this assumption
can be induced by tapering the observations such that $Z(\ub_{j})$ is
replaced with $\widetilde{Z}(\ub_{j})$, where
$\widetilde{Z}(\ub_{j}) = t(\ub_{j})Z(\ub_{j})$, $t(\ub) =
\prod_{i=1}^{d}t(s_{i})$ and $t$ is a weight function which has a
bounded second derivative, $t(-1/2)=t(1/2)=0$ and $t^{\prime}(1/2)=t^{\prime}(-1/2)=0$. 
By using  $\widetilde{Z}(\ub_{j})$ instead of $Z(\ub_{j})$, in all the
derivations below we replace the density $h(\ub)$ with $t(\ub)h(\ub)$. This means the results now rely on
the Fourier coefficients of $t(\ub)h(\ub)$, which decay at the rate
$|\int_{[-1/2,1/2]^{d}}t(\ub)h(\ub)\exp(i2\pi \jb^{\prime}\ub)d\ub| \leq C\prod_{i=1}^{d}\xi_{2}(j_{i})$,
and thus the above condition is satisfied. Note that
\citeA{p:mat-yaj-09}, Definition 2, uses a similar data-tapering
scheme to induce a similar condition. 
\end{remark}

The case that random locations follow a uniform distribution is an
example of a distribution which satisfies Assumption
\ref{assum:nonuniform}. It gives rise to several elegant
simplifications. Thus we state the uniform case as a separate
assumption. 
\begin{assumption}[Uniform sampling]\label{assum:uniform}
The locations $\{\ub_{j}\}$ are independent uniformly distributed random variables on 
$[-\lambda/2,\lambda/2]^{d}$.
\end{assumption}

Many of the results in this paper use that the location follow a
random design.  This helps in understand the sampling properties of
these complex estimators. However, it can "mask" the approximation
errors when replacing sums by integrals and the role that the sample
size $n$ plays in these approximations. To get some idea of these
approximations when the domain $\lambda$ is kept fixed but
$n\rightarrow \infty$ we will, on occassion, treat the locations as
deterministic and make  the following assumption.
\begin{assumption}\label{assum:fixedDFT}
Let $\{\sbb_{n,j}^{};j=1,\ldots,n\}$ denote the $n$ locations where
$s_{n,j}^{}=(s_{n,1,j},s_{n,2,j},\ldots,s_{n,d,j})$ the number of
locations $n\rightarrow \infty$ in such a way that for every $1\leq
\ell \leq d$
\begin{eqnarray}
\label{eq:sj-sj-1}
\sum_{j=1}^{n-1}\left|\frac{\lambda}{n^{1/d}} - \left(s_{n,\ell,(j+1)}^{} -
    s_{n,\ell,(j)}^{}\right) \right| =
O\left(\frac{\lambda}{n^{1/d}}\right), \quad \sum_{j=1}^{n-1} \left(s_{n,\ell,(j+1)}^{} -
    s_{n,\ell,(j)}^{}\right)^{2} =
O\left(\frac{\lambda}{n^{1/d}}\right)
\end{eqnarray}
where $s_{n,\ell,(j)}^{}$ denotes the  order statistics
corresponding to $s_{n,\ell,j}^{}$.
\end{assumption}

The integrated periodogram $Q_{a,\Omega,\lambda}(g;0)$ resembles the
integrated periodogram estimator commonly used in time series 
(see for example, \citeA{p:wal-64}, \citeA{p:han-71}, \citeA{p:dun-79}, \citeA{p:dah-jan-96},
\citeA{p:mik-10} and \citeA{p:kre-14}). However, there are some
fundamental differences, between time series estimators and
$Q_{a,\Omega,\lambda}(g;0)$ which makes the analysis very different. 
%We observe that the
%user-chosen frequency grid $\{\ob_{\Omega,\kb};\kb=(k_{1},\ldots,k_{d}),k_{j}\in [-a,a]\}$
%used to define $Q_{a,\Omega,\lambda}(g_{\theta};\rb)$ influences the limits
%of the functional $I(g;\frac{a}{\Omega})$. This is different to
%regular sampled locations where the frequency grid can only be defined over $[0,2\pi]^{d}$ (outside this domain aliasing
%occurs). Thus if the spatial process were observed only on a grid then
%aliasing means that it is impossible to estimate
%$I(g_{\theta};\infty)$. On the other hand, 
Unlike regularly spaced or
near regularly spaced data, `truely' irregular sampling means that the DFT can estimate high frequencies, without the curse
of aliasing (a phenomena which was noticed as early as
\citeA{p:sha-sil-59} and \citeA{p:beu-70}). In this case, 
if the function $g_{\theta}$, 
in the definition of $Q_{a,\Omega,\lambda}(g_{\theta};\rb)$ is bounded, there's no need for
the frequency grid to be bounded, and $a$ can be
magnitudes larger than $\lambda$. Below we state assumptions on the function $g$ 
and the frequency grid.

\begin{assumption}[Assumptions on $g_{\theta}(\cdot)$ and the size of
  frequency grid]\label{assum:G}
\begin{itemize}
\item[(i)] If $g_{\theta}$ is not a bounded function over
  $\mathbb{R}^{d}$ but $\sup_{\ob \in [-C,C]^{d}}|g_{\theta}(\ob)|
  <\infty$, then we must
  restrict the frequency grid to $[-C,C]^{d}$ (thus
  $a=C\Omega$). Further we assume
for all $1\leq j \leq d$, $\sup_{\ob\in [-C,C]^{d}}|\frac{\partial
  g_{\theta}(\ob)}{\partial \omega_{j}}|<\infty$.  
\item[(ii)] If $\sup_{\ob \in
    \mathbb{R}^{d}}|g_{\theta}(\ob)| <\infty$, then the frequency grid  can be 
unbounded (in the sense that $a/\lambda \rightarrow \infty$ as $a$ and
$\lambda\rightarrow \infty$). Further we assume for all $1\leq j \leq d$, 
$\sup_{\ob\in \mathbb{R}^{d}}|\frac{\partial
  g_{\theta}(\ob)}{\partial \omega_{j}}|<\infty$. 
\end{itemize}
\end{assumption}

 \begin{assumption}[Conditions on the spatial process]\label{assum:GG}
\begin{itemize}
\item[(a)] $|c(\ub)|\leq \beta_{1+\delta}(\ub)$.

Required for the bounded frequency grid, to obtain the DFT calculations.
\item[(b)] For some $\delta > 0$, $f(\ob)\leq \beta_{\delta}(\ob)$ 

Required for the unbounded
  frequency grid - using this assumption instead of (a) in the case of
  an bounded frequency grid leads to slightly larger errors bounds 
in the derivation of the first and second order moments. This assumption is
  also used to obtain the CLT result for both the bounded and
  unbounded frequency grids. 
%(however, it is unlikely this assumption is really  required for the bounded frequency grid).
\item[(c)] For all $1 \leq j\leq d$ and some $\delta>0$,
the partial derivatives satisfy 
$|\frac{\partial f(\ob)}{\partial \omega_{j}}|\leq
\beta_{\delta}(\ob)$.

We use this condition to approximate sums with integral for both the
bounded and unbounded frequency grids. It is also used to make a
series of approximations to derive the limiting variance in the case
that the frequency grid is unbounded. 
\item[(d)] For some $\delta > 0$, $|\frac{\partial^{d}f(\ob)}{\partial
    \omega_{1},\ldots,\partial \omega_{d}}|\leq
  \beta_{\delta}(\ob)$. 

Required only in the proof of Theorem \ref{lemma:mean-stat}(ii)(b). 

\item[(e)] For some $\delta >0$, $|f(\ob)|\leq \beta_{1+\delta}(\ob)$.

Required only for the fixed domain asymptotics. 
\end{itemize} 
\end{assumption}

\begin{remark}\label{remark:covariance}
Assumption \ref{assum:GG}(a) satisfied by a wide range of covariance
functions.  Examples include:
\begin{itemize}
\item[(i)] The Wendland
covariance, since its covariance is bounded and has a compact support.
\item[(ii)] The Matern covariance, which for $\nu > 0$ is defined as 
$c_{\nu}(\|\ub\|_{2})=\|\ub\|_{2}^{\nu}
K_{\nu}(\|\ub\|_{2})$ ($K_{\nu}$ is the modified Bessel function of the second
kind); see \citeA{b:ste-99}. 
To see why, we note that  if $\nu >0$ then $c_{\nu}(\ub)$ is a bounded
function. Furthermore, for large $\|\ub\|_{2}$ we
note that $c_{\nu}(\|\ub\|_{2})\sim
C_{\nu}\|\ub\|_{2}^{\nu-0.5}\exp(-\|\ub\|_{2})$ as
$\|\ub\|_{2}\rightarrow \infty$ (where $C_{\nu}$ is a finite constant). Thus by using the 
inequality
\begin{eqnarray*}
d^{-1/2}(|s_{1}|+|s_{2}|+\ldots +|s_{d}|)\leq
\sqrt{s_{1}^{2}+s_{2}^{2}+\ldots + s_{d}^{2}}\leq (|s_{1}|+|s_{2}|+\ldots +|s_{d}|) 
\end{eqnarray*}
we can show $|c_{\nu}(\ub)|\leq
\beta_{1+\delta}(\ub)$ for any $\delta >0$.
\end{itemize}
% recall that $s_{1}^{2}+s_{2}^{2}+2s_{1}s_{2}\leq s_{1}^{2}+s_{2}^{2}+2(s_{1}^{2}+s_{2})$
%Using this unequality
%it is straightforward to show that the exponential and Whittle covariance, defined by $c(\|\ub\|_{2})=\phi\exp(-\|\ub\|_{2})$ and
%$c(\|\ub\|_{2}) = (1+\sqrt{3}\|\ub\|_{2})\exp(-\sqrt{3}\|\ub\|_{2})$ respectively
%(where $\|\cdot\|_{2}$ denotes the Euclidean distance) satisfies this
%condition.  
\end{remark}

\begin{remark}\label{remark:exp}
Assumption \ref{assum:GG}(b,c,d) appears quite technical, but it is 
satisfied by a wide range of spatial covariance functions. For
example, the spectral density of the Matern covariance defined in
Remark \ref{remark:covariance} is $f_{\nu}(\ob)
=\frac{2^{\nu-1}\Gamma(\nu+\frac{d}{2})}{\pi^{d/2}(1+\|\ob\|_{2}^{2})^{(\nu+d/2)}}$
(see \citeA{b:ste-99}, page 49).
It is straightforward to show that this spectral density satisfies
Assumption \ref{assum:GG}(b,c,d), noting that the $\delta$ used to define
$\beta_{\delta}$ will vary
 with $\nu$, dimension $d$ and the order of the derivative $f_{\nu}(\cdot)$. 
%we now show that the exponential covariance (which belongs to the
%Mat\'ern class),  defined by $c(\|\ub\|_{2})=\phi\exp(-\|\ub\|_{2})$
%(where $\|\cdot\|_{2}$ denotes the Euclidean distance) satisfies these
%assumptions. Consider the case $d=1$ and $d=2$. For 
%$d=1$ the spectral density function of the exponential covariance is 
%\begin{eqnarray*}
%f(\omega) = \frac{\phi}{1+\omega^{2}},
%\end{eqnarray*}   
%whereas for $d=2$ the exponential covariance has the spectral density 
%\begin{eqnarray*}
%f(\omega_{1},\omega_{2}) = \frac{2\pi\phi}{(1+\omega_{1}^{2}+\omega_{2}^{2})^{3/2}}.
%\end{eqnarray*}   
%Both these model satisfy Assumption \ref{assum:G}, however the $\delta$
% which is used to define $\beta_{\delta}$ in (\ref{eq:def-alpha}) will vary
% with each model and condition. 
%For $d=1$,  Assumption \ref{assum:G}(ii)(a) uses $\delta = 1$, 
% Assumption \ref{assum:G}(ii)(b) uses $\delta=2$. For $d=2$, Assumption
%\ref{assum:G}(ii)(a) uses $\delta=1/2$, Assumption
%\ref{assum:G}(ii)(b) uses $\delta = 1/2$, whereas  Assumption \ref{assum:G}(ii)(c)
%uses $\delta = 3/2$.
%It is straightforward to show that these spectral density functions
%satisfy Assumption \ref{assum:G}(ii)(a,b). 
\end{remark}

If the spatial random field is non-Gaussian we require the following
assumptions on the higher order cumulants. 
\begin{assumption}[Non-Gaussian random fields]\label{assum:nonGaussian}
$\{Z(\ub);\ub \in\mathbb{R}^{d}\}$ is a fourth order stationary spatial 
random field, in the sense that $\Ex[Z(\ub)]=\mu$,
$\cov[Z(\ub_{1}),Z(\ub_{2})]=c(\ub_{1}-\ub_{2})$,  $\cum[Z(\ub_{1}),Z(\ub_{2}),Z(\ub_{3})] = 
\kappa_{2}(\ub_{1}-\ub_{2},\ub_{1}-\ub_{3})$
and 
$\cum[Z(\ub_{1}),Z(\ub_{2}),Z(\ub_{3}),Z(\ub_{4})] = 
\kappa_{4}(\ub_{1}-\ub_{2},\ub_{1}-\ub_{3},\ub_{1}-\ub_{4})$, for some
functions $\kappa_{3}(\cdot)$ and 
$\kappa_{4}(\cdot)$ and all $\ub_{1},\ldots,\ub_{4}\in \mathbb{R}^{d}$. 
We define the fourth order spectral density as 
$f_{4}(\ob_{1},\ob_{2},\ob_{3}) = \int_{\mathbb{R}^{3d}}\kappa_{4}(\ub_{1},\ub_{2},\ub_{3})
\exp(-i\sum_{j=1}^{3}\ub_{j}^{\prime}\ob_{j})d\ob_{1}d\ob_{2}d\ob_{3}$. We
assume that for some $\delta>0$
the spatial  tri-spectral density function is such that
$|f_{4}(\ob_{1},\ob_{2},\ob_{3})|\leq
\beta_{\delta}(\ob_{1})\beta_{\delta}(\ob_{2})\beta_{\delta}(\ob_{3})$
and $|\frac{\partial
  f_{4}(\omega_{1},\dots,\omega_{3d})}{\partial \omega_{j}}|\leq \beta_{\delta}(\ob_{1})\beta_{\delta}(\ob_{2})\beta_{\delta}(\ob_{3})$.
\end{assumption}

\subsection{Properties of Fourier transforms}\label{sec:DFT}

In this section we briefly summarize some of the characteristics of 
the Fourier transforms $J_{n}(\ob_{\Omega,\kb})$. These results will be
used as the basis of several estimator defined in the following
section. 

\begin{theorem}[Increasing domain asymptotic]\label{lemma:nonuniformdft}
Let us suppose that $\{Z(\ub);\ub\in \mathbb{R}^{d}\}$ is a stationary spatial
random field whose covariance function (defined in Assumption
\ref{assum:S}(i)) satisfies Assumption \ref{assum:GG}(a) for some $\delta >0$.
Furthermore, the locations $\{\ub_{j}\}$ satisfy Assumption \ref{assum:nonuniform}. Then we have 
\begin{eqnarray}
&&\cov\left[J_{n}(\ob_{\kb_{1}}),J_{n}(\ob_{\kb_{2}})\right]
= \langle \gamma,\gamma_{(\kb_{2}-\kb_{1})}\rangle f\left(\ob_{\kb_{1}}\right)
+\frac{c(0)\gamma_{\kb_{2}-\kb_{1}}\lambda^{d}}{n} +
O\left(\frac{1}{\lambda}\right), \label{eq:dft-nonuniform}
\end{eqnarray}
where the bounds are uniform in $\kb_{1},\kb_{2}\in \mathbb{Z}^{d}$
and  
$\langle \gamma,\gamma_{\rb}\rangle =
\sum_{\jb\in\mathbb{Z}^{d}}\gamma_{\jb}\gamma_{\rb-\jb}$.

Further, suppose Assumption \ref{assum:GG}(c) holds 
and that $|\gamma_{\jb}|\leq
\prod_{i=1}^{d}\xi_{2+\delta}(j_{i})$
for some $\delta > 0$. If 
$\Omega > \lambda$, then 
\begin{eqnarray}
\cov\left[J_{n}(\ob_{\Omega,\kb_{1}}),J_{n}(\ob_{\Omega,\kb_{2}})\right]
&=&
f\left(\ob_{\Omega,\kb_{1}}\right)\sum_{\jb_1,\jb_2}\gamma_{\jb_1}\gamma_{\jb_2}\Sinc\left(
  \pi \left[(\jb_1+\jb_2)-\frac{\lambda}{\Omega}(\kb_1 - \kb_2) 
    \right]\right) \nonumber\\
&&+ 
O\left(\frac{\lambda^{d}}{n}+\frac{\log \lambda}{\lambda}\right). \label{eq:dft-nonuniformOmega}
\end{eqnarray}
\end{theorem}
{\bf PROOF} See Section \ref{sec:proof-nonuniform}, \citeA{p:sub-14}. \hfill $\Box$

\vspace{1mm}
Comparing (\ref{eq:dft-nonuniform}) with
(\ref{eq:dft-nonuniformOmega}) we see that if $\kb_{1}$ and $\kb_{2}$
are such that $\lambda\max|\ob_{\Omega,\kb_1} - \ob_{\Omega,\kb_2}|<1$
then there is a high amount of correlation between the Fourier
transforms. On the other hand if $\lambda\max|\ob_{\Omega,\kb_1} -
\ob_{\Omega,\kb_2}|>1$ then there is a decay in correlation. This means if 
the frequency grid is very coarse, $\Omega <<
\lambda$ (as considered in \citeA{p:ban-lah-10} and
\citeA{p:ban-lah-nor-13}) the DFTs are almost uncorrelated. On the
other hand, if the frequency grid is very fine, $\Omega > \lambda$ (see (\ref{eq:dft-nonuniformOmega}))
frequencies which are close to each other are highly correlated. These observations
suggest that estimators based on $J_{n}(\ob_{\Omega,\kb})$ don't
gain in efficiency when $\Omega>\lambda$ but there will be a loss in
efficiency when $\Omega <\lambda$. We show this heuristic to be true
in Section \ref{sec:summary}.

\vspace{3mm} The results in the above theorem give the limit within the
increasing domain framework. In Theorem \ref{theorem:fixedDFT1},
\citeA{p:sub-14}, we obtain the properties of
$J_{n}(\ob_{\Omega,\kb})$ within the fixed domain framework, where
$\lambda$ is kept fixed but $n\rightarrow \infty$. The expressions are
long, but we summarize the most relevant parts in the remark below.

\begin{remark}\label{remark:fixed}
\begin{itemize}
\item[(i)] 
%Using that
%\begin{eqnarray*}
%&&\frac{(2\lambda)^{d}}{(\pi)^{d}}
%\int_{\mathbb{R}^{d}}f(\ob)\Sinc\left(\frac{\lambda}{2}\left[\ob+\frac{2\pi \kb_{1}}{\Omega}\right]\right)
%\Sinc\left(\frac{\lambda}{2}
%\left[\ob + \frac{2\pi
%    \kb_{2}}{\Omega}\right]\right)d\ob \\
%&=& \int_{[-\lambda/2,\lambda/2]^{2d}}
%c(\ubb_1-\ubb_2)
%\exp\left(\frac{2\pi i\ubb_{1}^{\prime} \kb_1}{\Omega}\right)\exp\left(-
%\frac{2\pi i\ubb_{2}^{\prime} \kb_2}{\Omega}\right)  d\ubb_1d\ubb_2
%\end{eqnarray*}
Let 
\begin{eqnarray}
\label{eq:Akb}
A_{\lambda}\left(\frac{\kb}{\Omega}\right) = 
\int_{[-\lambda,\lambda]^{d}} T\left(\frac{\ubb}{\lambda}\right)c(\ubb)\exp\left(\frac{2i\pi \kb^{\prime}\ubb}{\Omega}\right)d\ubb.
\end{eqnarray}
Using Theorem \ref{theorem:fixedDFT1},
\citeA{p:sub-14}, under Assumption \ref{assum:uniform} we have 
$
\var\left[J_{n}\left(\frac{2\pi \kb}{\Omega}\right)\right] =
  A_{\lambda}(\frac{\kb}{\Omega}) + O\left(\frac{\lambda^{d}}{n}\right)$.
Under Assumption \ref{assum:fixedDFT} we obtain a similar result
$\var\left[J_{n}\left(\frac{2\pi \kb}{\Omega}\right)\right] =
  A_{\lambda}(\frac{\kb}{\Omega}) +
  O\left(\frac{\lambda}{n^{1/d}}(\frac{|\kb|_{1}}{\Omega}+1)\right)$.

Therefore if the sampling frequency, $\Omega$, is chosen such that $\Omega \geq 2\lambda$, then
$\frac{1}{\Omega^{d}}A_{\lambda}(\frac{\kb}{\Omega})$ are the Fourier coefficients
of $T\left(\frac{\ubb}{\lambda}\right)c(\ubb)$ defined on the domain 
$[-\Omega/2,\Omega/2]^{d}$. In Section \ref{sec:Qexamples} we show
that this fixed domain approximation can, in some cases, be used to
obtain unbiased estimators and also quantify the bias when it cannot
be avoided.

%\item[(ii)] In the case that $d=1$ we have have the simplification 
%\begin{eqnarray*}
%\cov\left[J_{n}(\omega_{k_1}),J_{n}(\omega_{k_2}) \right]
%= 
%\left\{
%\begin{array}{cc}
%A_{\lambda}\left(\frac{k}{\lambda}\right) +O\left(\frac{\lambda}{n}(\frac{|k|_{1}}{\Omega}+1)\right  & k_1 = k_2 (=k) \\
%\frac{\Omega}{(k_1-k_2)\pi\lambda} \Im\left[
%e^{-\lambda i\pi (k_1-k_2)/\Omega}\left(\int^{\lambda}_{0}c(v) e^{2\pi i k_2 v/\lambda}dv
% - \int^{\lambda}_{0}c(v) e^{2\pi i k_1 v/\lambda}dv
%\right)\right] + O\left(\frac{\lambda}{n}(\frac{|k_{1}|+|k_{2}|}{\Omega}+1)\right & k_1\neq k_2 \\
%\end{array}
%\right.
%\end{eqnarray*}
%where
%\begin{eqnarray}
%\label{eq:Bkb}
%B_{\lambda}\left(\frac{k}{\lambda}\right)  = \int_{0}^{\lambda}c(u)\sin\left(\frac{2\pi
%k u}{\lambda}\right)du,
%\end{eqnarray}
%noting that a similar result also holds for $d>1$. 
%We note that for large $k$ $A_{\lambda}(\frac{k}{\lambda})\sim
%|k|^{-2}$, however in general $B_{\lambda}(\frac{k}{\lambda})\sim
%|k|^{-1}$ (since on the torus defined by $[0,1]$ $c(u)$ is
%discontiuous at the end points). This implies even when the domain $\lambda$ is
%fixed, $|\cov[J_{n}(\omega_{k_1},J_{n}(\omega_{k_2}))]|\sim
%|k_{1}-k_{2}|^{-2}$ \emph{however} within the fixed domain  
%correlations between the Fourier transforms
%will not decay to zero as $|k_{1}-k_{2}|\rightarrow \infty$. 
\item[(ii)] Using Theorem \ref{theorem:fixedDFT1},
\citeA{p:sub-14} and under Assumption \ref{assum:uniform} we have 
Assumption \ref{assum:nonuniform}
we have 
\begin{eqnarray}
&&\var\left[J_{n}\left(\frac{2\pi\kb}{\Omega}\right)\right] \label{eq:Jfixed}\\
&\approx& 
\frac{\lambda^{d}}{(2\pi)^{d}}\sum_{\jb_{1},\jb_{2}\in \mathbb{Z}^{d}}\gamma_{\jb_1}\gamma_{\jb_2}
\int_{\mathbb{R}^{d}}f(\ob)\Sinc\left(\frac{\lambda}{2}\left[\ob + \frac{2\pi \jb_{1}}{\lambda}
+\frac{2\pi \kb_{}}{\Omega}\right]\right)\Sinc\left(\frac{\lambda}{2}
\left[\ob - \frac{2\pi \jb_{2}}{\lambda} + \frac{2\pi \kb_{1}}{\Omega}\right]\right)d\ob.\nonumber
\end{eqnarray}
Comparing (\ref{eq:Jfixed}) with (\ref{eq:dft-nonuniform}) we observe
that when the locations do not satisfy either Assumption
\ref{assum:uniform} or \ref{assum:fixedDFT}, then within the fixed
domain framework  the sampling design and the spectral
density function are convolved. However, under increasing domain
asymptotics, see Theorem 
\ref{lemma:nonuniformdft}, the location density is separated from the
spectral density. Using (\ref{eq:Jfixed}) to construct an estimator
can be difficult because of the influence of design density. In which
case it is preferable to use the approximation in
(\ref{eq:dft-nonuniform}) as the basis of any estimation scheme.
%but if
%the range of the covariance is sufficiently small with respect to the
%spatial domain $\lambda$, then the approximation in
%(\ref{eq:dft-nonuniform}) will hold. 
%It is easier to define spatial estimators when the
%spatial design and spectrum are separated. For this reason, if neither 
%Assumption
%\ref{assum:uniform} or \ref{assum:fixedDFT}

\end{itemize}
%Furthermore, we observe that when the spatial domain $\lambda$ is
%allowed to grow, the Fourier transforms $J_{n}(\ob_{\kb})$ are either 
%asymptotically uncorrelated at different frequencies or are
%asymptotically uncorrelated as $\kb$ grows.
% Surprisingly, in
%the case that the domain is fixed the correlation between the Fourier transforms $J_{n}(\omega \kb_{1})$ and
%$J_{n}(\omega \kb_{2})$ decreases as the distance between $\kb_{1}$ and $\kb_{2}$ grows. 
%However, it should be noted that the bounds are based on approximating summands
%by integrals. This means that for very large $\kb$
%(of order close to $n$) $\var[J_{n}(\pi \kb)]$ is unlikely to be
%negligible despite $A(\kb) = O(\prod_{j=1}^{d}|k_{j}|^{-2})$. 
\end{remark}

\begin{remark}\label{remark:reconcile}
We observe from Theorems \ref{lemma:nonuniformdft} and
\ref{theorem:fixedDFT1} that in the increasing domain framework
$|J_{n}(\ob_{\Omega,\kb})|^{2}$
is an estimator of the spectral density function, $f(\ob_{\Omega,\kb})$
whereas within the fixed domain framework
$|J_{n}(\ob_{\Omega,\kb})|^{2}$ is an estimator of the Fourier
coefficient $A_{\lambda}(\kb/\Omega)$. 
However, if 
$\int_{\mathbb{R}^{d}}|f(\ob)|d\ob <\infty$ then 
$A_{\lambda}(\frac{\kb}{\Omega})\rightarrow f(\omega_{\kb,\Omega})$ as
$\lambda \rightarrow \infty$ and if the stronger condition
$\int_{\mathbb{R}^{d}}(1+\|\ubb\|_{1})|c(\ubb)|d\ubb <\infty$ holds, then
$|A_{\lambda}(\frac{\kb}{\Omega})- f(\omega_{\kb,\Omega})| = O(\lambda^{-1})$.
\end{remark}

\subsection{Examples of estimators defined within the Fourier domain}
\label{sec:Qexamples}

Many parameters or quantities of interest
can be written as a linear functional involving the spectral density
function $f$. In Theorem \ref{lemma:nonuniformdft} and Remark
\ref{remark:reconcile} we showed that if $\lambda$ is large 
and the design of locations uniform then $\Ex[J_{n}\left(\ob_{\Omega,\kb}\right)]\approx
f(\ob_{\Omega,\kb})$ (if the design is not uniform then there will be
a multiplicative constant).  Motivated by this observation, in this
section we consider estimators (or criterions) which take the form
stated in (\ref{eq:QQQ}).
If the locations follow a uniform distribution then for some of the
examples below it is possible to reduce the (fixed domain) bias in the estimator.
In such cases the function $g(\cdot)$ is replaced with a fixed domain
approximation $g_{\lambda}(\cdot)$. Note as $\lambda \rightarrow
\infty$ these functions are asymptotically equivalent.

\subsubsection{The Whittle likelihood}\label{sec:whittle}

Suppose the stationary spatial process $\{Z(\sbb);\sbb\in
\mathbb{R}^{d}\}$ has spectral density $f_{\theta}(\ob)$ where
$\theta$ is unknown but belongs to the compact parameter space $\Theta$. 
\citeA{p:mat-yaj-09} propose using the integrated Whittle likelihood to estimate
$\theta$. More precisely, they define the Whittle likelihood as
\begin{eqnarray*}
\mathcal{L}_{I,n}(\theta,\eta^{2}) = \int_{\Omega}\left(\log
  \left[f_{\theta}(\ob) + \eta^{2}\right] + \frac{|J_{n}(\ob)|^{2}}{\left[f_{\theta}(\ob)+\eta^{2}\right]} \right)d\ob,
\end{eqnarray*}
and use $(\widehat{\theta},\widehat{\eta})\in \arg\min_{\theta,\eta}  \mathcal{L}_{n}(\theta,\eta)$ as an
estimator of $\theta$ and $\eta$ (where $\eta$ is an estimator of
the "ridge effect"). Of course, this integral cannot be evaluated
in practice and a Riemann sum approximation is
necessary. 
%In Theorem \ref{lemma:nonuniformdft} 
%we showed that at the frequencies $\ob_{\kb} = \ob_{\lambda,\kb}$ the
%covariance between between the Fourier transforms decay, 
Using $\{\ob_{\kb}  = 2\pi\kb/\lambda; \kb  =
(k_{1},\ldots,k_{d}), -C\lambda \leq k_{j}\leq  C\lambda\}$ we
approximate the integral with the sum
\begin{eqnarray*}
\mathcal{L}_{S,n}(\theta,\eta^{2}) = \frac{1}{\lambda^{d}}\sum_{\kb=-a}^{a}\bigg(\log 
\left[f_{\theta}(\ob_{\kb}) + \eta^{2}\right] +
\frac{|J_{n}(\ob_{\kb})|^{2}}{\left[ f_{\theta}(\ob_{\kb})+ \eta^{2}\right]} \bigg).
\end{eqnarray*}
A heuristic motivation for the above likelihood is that in the case
the locations are uniformly distribution then
$\{J_{n}(\ob_{\kb})\}$ are near uncorrelated random
variables with asymptotic variance $f_{\theta_0}(\ob_{\kb})+\eta_0^{2}$.  If $J_{n}(\ob_{\kb})$ were
Gaussian, uncorrelated random variables with variance
$f_{\theta_0}(\ob_{\kb})+\eta_0^{2}$ then
$\mathcal{L}_{S,n}(\theta,\eta^{2})$ would be the true likelihood. 
The choice of $a=C\lambda$ is necessary (where $C$ does not depend on $\lambda$), since 
$|f_{\theta}(\ob)|\rightarrow 0$ as  $\|\ob\|_{}\rightarrow \infty$
(for any norm $\|\cdot\|$), thus the 
discretized Whittle likelihood is only well defined over a bounded
frequency grid. The choice of $C$ is tied to how fast the tails in 
the parametric class of spectral density functions
$\{f_{\theta};\theta \in \Theta\}$ decay to zero, therefore Assumption \ref{assum:G}(i)
is satisfied.

If either Assumption \ref{assum:uniform} or \ref{assum:fixedDFT} is
satisfied, from Theorem \ref{theorem:fixedDFT1} and Remark
\ref{remark:fixed} it is clear that $A_{\lambda}(\frac{\kb}{\lambda};\theta)$ is
a better approximation of $\var[J_{n}(\ob)]$ (where
$A_{\lambda}(\cdot)$ is defined in (\ref{eq:Akb})). 
Therefore if Assumption \ref{assum:uniform} or \ref{assum:fixedDFT} is
satisfied, a better finite sample approximation can be obtained by 
using $\widehat{\theta} = \arg\min
\mathcal{L}_{S,n}(\theta)$ as an estimator of $\theta$, where
\begin{eqnarray}
\label{eq:fixedLike}
\mathcal{L}_{S,n}(\theta,\eta^{2}) = \frac{1}{\lambda^{d}}\sum_{\kb=-a}^{a}\bigg(\log 
\left[A_{\lambda}\left(\frac{\kb}{\lambda};\theta\right) + \eta^{2}\right] +
\frac{|J_{n}(\ob_{\kb})|^{2}}{\left[A_{\lambda}\left(\frac{\kb}{\lambda};\theta\right) + \eta^{2}\right]} \bigg).
\end{eqnarray}
%Note that in the case that $\var[J_{n}(\ob_{\kb})] =
%\sigma^{2}A_{\lambda}(\frac{\kb}{\lambda};\theta)$ where $\theta$ is
%known then the Whittle likelihood estimator of $\sigma^{2}$ is 
%\begin{eqnarray*}
%\widehat{\sigma}^{2} = \frac{1}{a}\sum_{\kb=-a}^{a}\frac{|J_{n}(\ob_{\kb})|^{2}}{A_{\lambda}(\frac{\kb}{\lambda};\theta)}.
%\end{eqnarray*}

\subsubsection{The spectral density estimator}

We recall from Theorem \ref{lemma:nonuniformdft} that
$\var[J_{n}(\ob_{\kb})] \approx \langle \gamma,\gamma_{0} \rangle
f(\ob_{\kb})$. Since $f(\cdot)$ is locally constant in a neighbourhood
of $\ob$ and motivated by spectral methods in time series we use
$\widehat{f}_{\lambda,n}(\ob)$ as a nonparametric estimator  of $f$
(or a constant multiple of it), where
\begin{eqnarray*}
\widehat{f}_{\lambda,n}(\ob_{}) =
\sum_{\kb=-\lambda/2}^{\lambda/2}W_{b}(\ob_{}-\ob_{\kb})|J_{n}(\ob_{\kb})|^{2}  
%\frac{1}{n}\sum_{\kb=-a}^{a}W_{b}(\ob-\ob_{\kb}) 
= \frac{1}{b^{d}}Q_{a,\lambda,\lambda}(W_{b},0),
%\frac{1}{n}\sum_{j=1}^{n}Z(\sbb_{j})^{2} = \frac{1}{b^{d}}\widetilde{Q}_{a,\lambda}(W_{b},0),
\end{eqnarray*}
$W_{b}(\ob) = b^{-d}\prod_{j=1}^{d}W(\frac{\omega_{j}}{b})$ and
$W:[-1/2,-1/2]\rightarrow \mathbb{R}$ is a spectral window. In this
case we set the number of frequencies  $a=\lambda/2$, 
and Assumption \ref{assum:G}(i)
is satisfied.
%It is worth mentioning that the spectral density estimator given in
%Example \ref{example:fixed}(ii) does not exactly satisfy the
%Assumptions \ref{assum:G}(i)(a,b) (since $\sup_{\ob}|W_{b}(\ob)|=O(b^{-d})$), however it is straightforward to
%adapt the sampling properties derived in the later sections to this case.

\subsubsection{A  nonparametric non-negative definite estimator of the
  spatial covariance}\label{sec:nonparametric}

In this section we propose a nonparametric estimator of the
covariance. The estimator is based on the representation
\begin{eqnarray*}
c(\ubb) = \frac{1}{(2\pi)^{d}}\int_{\mathbb{R}^{d}}f(\ob)\exp(i\ob^{\prime}\ubb)d\ob.
\end{eqnarray*}
Since the expectation of $|J_{n}(\ob_{\kb})|^{2}$ is approximately
$f(\ob_{\kb})$, to estimate the spatial covariance we propose
approximating the above integral with a sum and the spectral density
with the absolute square of the Fourier transform. However, using the
frequency grid $\ob_{\kb} = \frac{2\pi \kb}{\lambda}$ is problematic
outside the region $[-\lambda/2,\lambda/2]^{d}$. % due to reflections.
Instead,  we propose using a finer grid 
to estimate the covariance, namely
\begin{eqnarray}
\label{eq:nonparcov2}
\widetilde{c}_{\Omega, n}(\ubb) = 
\frac{1}{\Omega^{d}}\sum_{\kb  = -a}^{a}
\left|J_{n}\left(\ob_{\kb,\Omega}\right)\right|^{2}\exp(i\ubb^{\prime}\ob_{\kb,\Omega})
\end{eqnarray}
where $\Omega \geq 2\lambda$. In this case, $a=a(\lambda)$
can be chosen such that
$a/\lambda \rightarrow \infty$ as $a\rightarrow\infty$ and
$\lambda\rightarrow \infty$ (thus Assumption \ref{assum:G}(ii)
is satisfied). We observe that when $d=1$ then for all $u\in
[0,\Omega/2]$ we have
$\widetilde{c}_{\Omega,n}(u) = \widetilde{c}_{\Omega,n}(\Omega - u)$,
however as we are only interested in estimating the autocovariance
within $[-\lambda,\lambda]^{d}$ and $\Omega/2 \geq \lambda$ this is not an issue.

A disadvantage with the above ``raw estimator'' of the covariance is
that there is no guarantee that it yields a non-negative definite
spatial auto-covariance function.
However, this can easily be remedied by multiplication with the triangle kernel. More
precisely, define the estimator 
\begin{eqnarray}
\widehat{c}_{\Omega,n}(\ubb) &=&  
T\left(\frac{\ubb}{\widehat{\Omega}}\right)
\widetilde{c}_{\Omega,n}(\ubb)\label{eq:cnhat} 
\end{eqnarray}
where $T({\boldsymbol u})=\prod_{j=1}^{d}T(u_{j})$. This 
covariance estimator has the advantage that it is zero outside
the region $[-\widehat{\Omega},\widehat{\Omega}]^{d}$. Moreover,
$\widehat{c}_{\Omega,n}(\ubb)$ is a non-negative definite sequence. To
show this result, we use that the  Fourier transform of the 
triangle kernel, $T(u)$ is  $\sinc^{2}(\frac{\omega}{2})$. Thus the Fourier transform of 
$\widehat{c}_{\Omega,n}(\ubb) $ is 
\begin{eqnarray*}
\widehat{f}_{\Omega,\lambda}(\ob) &=&
\int_{[-\widehat{\Omega},\widehat{\Omega}]^{d}}\widehat{c}_{\Omega,n}(
\ubb)\exp(-i\ob^{\prime}\ubb)d\ubb
  =
  \frac{\widehat{\Omega}^{d}}{\Omega^{d}}\sum_{\kb=-a}^{a}|J_{n}(\ob_{\kb,\Omega})|^{2}
\Sinc^{2}\left[\frac{\widehat{\Omega}}{2}(\ob_{\kb,\Omega}-\ob)\right].
\end{eqnarray*}
Clearly, $\widehat{f}_{\Omega,\lambda}(\ob)\geq 0$, therefore, the
estimator $\{\widehat{c}_{\Omega,n}(\ubb)\}$ is a non-negative definite
function and thus a valid covariance function. 

In Section \ref{sec:simulations}, \citeA{p:sub-14}, 
we illustrate the performance of the nonparametric non-negative definite estimator of the
  spatial covariance through some simulations. 

\subsubsection{A nonlinear least squares estimator of a parametric
  covariance function}

 We recall that Whittle likelihood can only be defined on
  a bounded frequency grid. This
  can be an issue if the observed locations are dense on
  the spatial domain and thus contain a large amount of high frequency
  information which would be missed by the Whittle likelihood. An alternative method for parameter estimation
  of a spatial process is to use a different loss function.
Motivated by \citeA{p:ric-79}, the discussion  on the Whittle
estimator in Secton \ref{sec:whittle} and Theorem \ref{lemma:nonuniformdft}
we define the quadratic loss function
\begin{eqnarray}
\label{eq:L2estimator}
L_{n}(\theta) = \frac{1}{\lambda^{d}}\sum_{k_{1},\ldots,k_{d} = -a}^{a}\left(|J_{n}(\ob_{\kb})|^{2}
  - \langle \gamma, \gamma_{0}\rangle f_{\theta}(\ob_{\kb}) \right)^{2},
\end{eqnarray} 
and let $\widehat{\theta}_{n} = \arg\min_{\theta \in \Theta} L_{n}(\theta)$ or equivalently
solve $\nabla_{\theta}L_{n}(\theta)=0$, where
\begin{eqnarray}
\label{eq:DLn}
\nabla_{\theta}L_{n}(\theta) &=& -\frac{2 \langle \gamma, \gamma_{0}\rangle}{\lambda^{d}}\sum_{k_{1},\ldots,k_{d} = -a}^{a}
 \nabla_{\theta}f_{\theta}(\ob_{\kb})\left\{|J_{n}(\ob_{\kb})|^{2}
  - \langle \gamma, \gamma_{0}\rangle f_{\theta}(\ob_{\kb}) \right\}  \nonumber\\
&=& -\langle \gamma, \gamma_{0}\rangle\left[Q_{a,\lambda,\lambda}(2 
\nabla_{\theta}f_{\theta}(\cdot);0) -\langle \gamma, \gamma_{0}\rangle
\frac{2}{\lambda^{d}}\sum_{\kb = -a}^{a}f_{\theta}(\ob_{\kb})
\nabla_{\theta}f_{\theta}(\ob_{\kb})\right].
\end{eqnarray}
It is
well known that the distributional properties of a quadratic loss
function are determined by its first derivative.
In particular,  the asymptotic sampling properties of
$\widehat{\theta}_{n}$ are determined by 
$Q_{a,\lambda,\lambda}(2 
\nabla_{\theta}f(\cdot;\theta);0)$ . In this case $a$
can be such that
$a/\lambda \rightarrow \infty$ as $\lambda\rightarrow \infty$ and Assumption \ref{assum:G}(ii)
is satisfied. An estimator of $\langle \gamma, \gamma_{0}\rangle$ is
given in Remark \ref{remark:gamma}. Note that in the definition of $L_{n}(\theta)$,
$\langle \gamma, \gamma_{0}\rangle$ can be replaced with  $\sigma^{2}$, in
which case one is estimating a multiple of $f_{\theta_0}(\ob)$.

If either Assumption \ref{assum:uniform} or \ref{assum:fixedDFT} is
satisfied then we can replace
$f_{\theta}(\ob_{\kb})$ with
$A_{\lambda}(\frac{\kb}{\lambda};\theta)$, to obtain a better fixed
domain approximation.

%% file: 3_summary.tex
\section{A summary of the sampling properties of $Q_{a,\Omega,\lambda}(g;0)$}\label{sec:summary}

In this section we consider the sampling properties of 
$Q_{a,\Omega,\lambda}(g;0)$ for the general frequency grid
$\{\ob_{\Omega,\kb} = \frac{2\pi \kb}{\Omega}\}$. The proof and more
general results can be found in Section \ref{sec:fine-grid},
\citeA{p:sub-14}. To simplify notation in this section we mainly consider the case that the
locations are uniformly distributed.  We show that using the
 frequency grid $\ob_{\lambda,\kb}$ is optimal and in
 Section \ref{sec:uniform} we focus on this frequency grid.
% deriving the sampling properties of $Q_{a,\lambda,\lambda}(g;0)$ the the locations 
%satisfy the general assumptions in Assumption \ref{assum:nonuniform}. 

\begin{lemma}\label{lemma:summaryE}
Assumptions \ref{assum:S}(i), \ref{assum:uniform} and
\ref{assum:GG}(a,c) or (b,c) holds. Let $A_{\lambda}(\cdot)$ be
defined as in (\ref{eq:Akb}).
\begin{itemize}
\item[(i)] Then 
\begin{eqnarray*}
\Ex\left[Q_{a,\Omega,\lambda}(g;0)\right] 
 &=& \frac{c_{2}}{\Omega^{d}}\sum_{\kb = -a}^{a}g\left(\frac{2\pi \kb}{\Omega}\right)A_{\lambda}\left(\frac{\kb}{\Omega}\right) + 
 \frac{c(0)\lambda^{d}}{n^{}\Omega^{d}}\sum_{\kb =  -a}^{a}g(\omega_{\Omega,\kb}) \\
 &=& 
 \frac{c_{2}}{(2\pi)^{d}}\int_{[-2\pi a/\Omega,2\pi a/\Omega]^{d}}g(\ob)f(\ob)d\ob
 + \frac{c(0)\lambda^{d}}{n\Omega^{d}}\sum_{\kb = -a}^{a}g(\omega_{\Omega,\kb})+
O\left(\frac{\log \lambda}{\lambda}+\frac{1}{\Omega}\right).
\end{eqnarray*}
\item[(ii)] If Assumption \ref{assum:fixedDFT} holds instead of   \ref{assum:uniform} then
\begin{eqnarray*}
\Ex\left[Q_{a,\Omega,\lambda}(g;0)\right] 
&=& \frac{1}{\Omega^{d}}
\sum_{\kb = -a}^{a}g\left(\frac{2\pi \kb}{\Omega}\right)A_{\lambda}\left(\frac{\kb}{\Omega}\right) + 
O\left(\frac{1}{\Omega^{d}}\sum_{\kb =
    -a}^{a}|g(\omega_{\Omega,\kb})|\left\{\frac{\lambda}{n^{1/d}}\left(1+\frac{\|\kb\|_1}{\Omega}\right) \right\}\right)
\end{eqnarray*} 
\item[(iii)] If Assumption \ref{assum:fixedDFT} holds instead of
  \ref{assum:uniform} 
and $\sup|g(\cdot)|<\infty$, then 
\begin{eqnarray*}
\Ex\left[Q_{a,\Omega,\lambda}(g;0)\right] 
&=& \frac{1}{\Omega^{d}}
\sum_{\kb = -\infty}^{\infty}g\left(\frac{2\pi \kb}{\Omega}\right)A_{\lambda}\left(\frac{\kb}{\Omega}\right) + 
 O\left(\frac{\lambda a^{d}}{\Omega^{d}n^{1/d}} + \frac{\lambda a^{d+1}}{\Omega^{d+1}n^{1/d}} +\frac{1}{a}\right).
\end{eqnarray*} 
\end{itemize}
\end{lemma}
We observe that in the increasing domain framework
$Q_{a,\Omega,\lambda}(g;0)$ is estimating an integral of the spectral
density function. On the other hand, within the fixed domain framework
$Q_{a,\Omega,\lambda}(g;0)$ is estimating a weighted sum of the
Fourier coefficients $A_{\lambda}(\kb/\Omega)$. 

We apply the above results 
to some of the examples considered in the previous section. The
results are given in the general case that the 
locations may be not uniformly distributed (see Theorem
\ref{theorem:nonuniformmean} and Lemma \ref{lemma:meangrid}, \citeA{p:sub-14}). 
\begin{example}
\begin{itemize}
\item[(i)] \underline{The Whittle likelihood} Under Assumption
  \ref{assum:nonuniform}, increasing domain asymptotics and by
  using Theorems \ref{lemma:nonuniformdft} and  \ref{theorem:nonuniformmean} we have 
\begin{eqnarray*}
\Ex\left[\mathcal{L}_{S,n}(\theta,\eta^{2})\right] = \frac{1}{(2\pi)^{d}}\int_{[-a/\lambda,a/\lambda]^{d}}\left(
\log\left[f_{\theta}(\ob) + \eta^{2}\right] +
\frac{\langle \gamma,\gamma_{0}\rangle f_{\theta_0}(\ob) + \gamma_{0}\eta_{0}^{2}}{f_{\theta}(\ob) + \eta^{2}}\right) + O\left(\frac{1}{\lambda}\right),
\end{eqnarray*}
where $f_{\theta_0}(\cdot)$ denotes the true spectral density,
$c(\ub;\theta_0)$ the corresponding spatial covariance and
$\eta_{0}^{2}=\lambda^{d}n^{-1}c(0;\theta_0)$
(note that $\eta_0^{2} = O(\lambda^{d}/n)$ hence for large $n$ it is
close to zero). The $\theta$ which minimises $\Ex\left[\mathcal{L}_{S,n}(\theta,\eta^{2})\right] $ is such that
$f_{\theta}(\ob) = \langle \gamma,\gamma_{0}\rangle
f_{\theta_0}(\ob)$ for all $\ob\in \mathbb{R}^{d}$.

However, under fixed domain asymptotics and assuming
Assumption \ref{assum:fixedDFT} holds then the expectation of the
fixed domain likelihood defined in (\ref{eq:fixedLike}) is 
\begin{eqnarray*}
\Ex[\mathcal{L}_{S,n}(\theta,\eta^{2})] = \frac{1}{\lambda^{d}}\sum_{\kb=-a}^{a}\bigg(\log 
\left[A_{\lambda}\left(\frac{\kb}{\lambda};\theta\right) + \eta^{2}\right] +
\frac{A_{\lambda}(\frac{\kb}{\lambda};\theta_0)}{\left[A_{\lambda}\left(\frac{\kb}{\lambda};\theta\right)
    + \eta^{2}\right]} \bigg) +
O\left(\frac{a^{3d+1}}{n^{1/d}\lambda^{d}} + \frac{a^{3d}}{\lambda^{d-1}n^{1/d}}\right),
\end{eqnarray*}
where in the above error bound we use that $|A_{\lambda}(\kb/\lambda)|\leq \prod_{i=1}^{d}\xi_{2}(k_{i})$.
\item[(ii)] \underline{The nonparametric covariance}
Under increasing domain asymptotics and using 
Lemma \ref{lemma:meangrid}(ii), \citeA{p:sub-14}, for $\ubb \in [-\min(\lambda,\Omega/2),\min(\lambda,\Omega/2)]^{d}$ we have 
\begin{eqnarray*}
\Ex[\widetilde{c}_{\Omega, n}(\ubb)] = \langle \gamma,
\gamma_{0}\rangle c(\ubb) + O\left(\frac{\log \lambda}{\lambda} + \frac{1}{\Omega}\right),
\end{eqnarray*}
and $\Ex[\widehat{c}_{\Omega, n}(\ubb)] = \langle \gamma,
\gamma_{0}\rangle T(\frac{\ubb}{\widehat{\Omega}})c(\ubb) + O\left(T(\frac{\ubb}{\widehat{\Omega}})[\frac{\log \lambda}{\lambda} + \frac{1}{\Omega}]\right)$
where $\langle \gamma,
\gamma_{0}\rangle$ is defined in Theorem \ref{lemma:nonuniformdft}.

In order to understand the properties of $\widetilde{c}_{\Omega,
  n}(\ubb)$ under fixed domain asymptotics we assume that Assumption
\ref{assum:fixedDFT} holds (though a similar result holds under Assumption
\ref{assum:uniform}). If $\Omega \geq 2\lambda$ and $\vbb \in
[-\lambda,\lambda]$ then by using Lemma \ref{lemma:summaryE}(iii) we have
\begin{eqnarray*}
\Ex\left[ \widetilde{c}_{\Omega, n}(\ubb) \right] &=& 
\frac{1}{\Omega^{d}}
\sum_{\kb = -\infty}^{\infty}A_{\lambda}\left(\frac{\kb}{\Omega}\right)\exp\left(i\ubb^{\prime}\ob_{\Omega,\kb}\right) + 
 O\left(\frac{\lambda a^{d}}{\Omega^{d}n^{1/d}} + \frac{\lambda
     a^{d+1}}{\Omega^{d+1}n^{1/d}} +\frac{1}{a}\right) \\
&=& T\left(\frac{\ubb}{\lambda}\right)c(\ubb) + O\left(\frac{\lambda a^{d}}{\Omega^{d}n^{1/d}} + \frac{\lambda
     a^{d+1}}{\Omega^{d+1}n^{1/d}} +\frac{1}{a}\right),
\end{eqnarray*}
where we recall $T(\cdot)$ denotes the triangle kernel.
The above follows from the representation that for any $\Omega \geq
2\lambda$ and $\vbb \in [-\Omega/2,\Omega/2]^{d}$ then 
\begin{eqnarray*}
T\left(\frac{\vbb}{\lambda}\right)c(\ubb) = \frac{1}{\Omega^{d}}
\sum_{\kb \in \mathbb{Z}^{d}}A_{\lambda}\left(\frac{\kb}{\Omega}\right)\exp\left(\frac{2\pi
      i \kb^{\prime}\ubb}{\Omega} \right).
\end{eqnarray*}
Thus we observe that as $\lambda \rightarrow \infty$ and $\Omega \rightarrow \infty$, 
$\widetilde{c}_{\Omega,  n}(\ubb)$ is asymptotically a consistent estimator of 
$c(\ubb)$.  However, for finite $\lambda$,
$\widetilde{c}_{\Omega,n}(\ubb)$ will have the bias
$[1-T(\frac{\ubb}{\lambda})]c(\ubb)$ (and correcting this bias will
lead to an increase in variance). 
%which is large close to the boundary of $[-\lambda,\lambda]^{d}$. 
\end{itemize}
\end{example}

In Lemma \ref{lemma:summaryE} we observe that the term
$\frac{c(0)\lambda^{d}}{n\Omega^{d}}\sum_{\kb =
  -a}^{a}g(\ob_{\Omega,\kb})$ arises.
 It can be removed by
using a bias corrected version of $Q_{a,\Omega,\lambda}(g;0)$
\begin{eqnarray}
\widetilde{Q}_{a,\Omega,\lambda}(g;0) = \frac{1}{\Omega^{d}}\sum_{k_{1},\ldots,k_{d}=-a}^{a}
g(\ob_{\Omega,\kb})\left|J_{n}(\ob_{\Omega,\kb})\right|^{2}
- \frac{\lambda^{d}}{\Omega^{d}n}\sum_{\kb=-a}^{a}g(\ob_{\Omega,\kb})
\frac{1}{n}\sum_{j=1}^{n}Z(\sbb_{j})^{2}.
\end{eqnarray}
For the remainder of this section we focus on the 
bias corrected estimator
$\widetilde{Q}_{a,\Omega,\lambda}(g;\rb)$. However, similar results to
those established for $\widetilde{Q}_{a,\Omega,\lambda}(g;\rb)$ also
hold for $Q_{a,\Omega,\lambda}(g;\rb)$ (see Section
\ref{sec:nonparametric-prop}, \citeA{p:sub-14}).
%However, in some situations
%some situations 
%using $Q_{a,\Omega,\lambda}(g;\rb)$ is preferable and
%the sampling properties of this estimator are given in Section \ref{}, \citeA{p:sub-14}.

We show in Section \ref{sec:fine-grid}, \citeA{p:sub-14} that in 
 the case that $\{Z(\ubb);\ubb\in \mathbb{R}^{d}\}$ is a Gaussian
stationary spatial process
\begin{eqnarray*}
\var\left[\widetilde{Q}_{a,\Omega,\lambda}(g;0)\right]
 =  C_{1}\left(\frac{a}{\Omega}\right)
\frac{1}{\Omega^{d}}\sum_{k_{1},\ldots,k_{d}=-2a}^{2a}\Sinc^{2}\left(
  \frac{\lambda}{\Omega}\kb 
  \pi\right) +
O\left(\widetilde{\ell}_{a,\Omega,\lambda}\left[\frac{1}{\lambda^{d}}I_{\frac{\lambda}{\Omega}<1}
  + \frac{1}{\Omega^{d}}I_{\frac{\lambda}{\Omega}\geq 1}\right]\right)\label{eq:varQtilde}
\end{eqnarray*}
where $\widetilde{\ell}_{a,\lambda.\Omega}$ is defined in
(\ref{eq:widetildeell}), Section \ref{sec:fine-grid}, \citeA{p:sub-14} and 
\begin{eqnarray}
\label{eq:Calambda}
C_{1}\left(\frac{a}{\Omega}\right) = 
\frac{1}{(2\pi)^{d}}\int_{-[2\pi a/\Omega,2\pi a/\Omega]^{d}} 
f(\ob_{})^{2}\left[g(\ob)|^{2} + g(\ob)\overline{g(-\ob)}\right]d\ob. 
\end{eqnarray}  
Observe that the rate of convergence of
$\var\left[\widetilde{Q}_{a,\Omega,\lambda}(g;0)\right]$
is determined by the term 
\begin{eqnarray*}
\sum_{k_{1},\ldots,k_{d}=-2a}^{2a}\Sinc^{2}\left(
  \frac{\lambda}{\Omega}\kb 
  \pi\right) = \prod_{i=1}^{d}\sum_{k_i=-2a}^{2a}\sinc^{2}\left(
  \frac{\lambda}{\Omega}k_{i}
  \pi\right).
\end{eqnarray*}
It is this term along with the following result which 
gives the crucial insight into the rate of convergence for
different  frequency grids $\{\ob_{\Omega,\kb}\}$. 
If  $a\rightarrow \infty$ then 
\begin{eqnarray*}
\frac{1}{\Omega}\sum_{k=-\infty}^{\infty}\sinc^{2}\left( \frac{\lambda}{\Omega}k
  \pi\right) = 
\left\{
\begin{array}{cc}
\frac{1}{\lambda} & \frac{\lambda}{\Omega} <1 \\
\frac{1}{\Omega} & \frac{\lambda}{\Omega} \in \mathbb{Z} \\
O\left(\frac{1}{\Omega}\right) & \frac{\lambda}{\Omega} > 1 \textrm{
  and }\frac{\lambda}{\Omega} \notin \mathbb{Z},
\end{array}
\right.
\end{eqnarray*}
further, if $\lambda/\Omega \rightarrow \infty$ then 
$\sum_{k=-\infty}^{\infty}\sinc^{2}\left( \frac{\lambda}{\Omega}k
  \pi\right)\rightarrow 1$ (see  Section \ref{sec:fine-grid},
\citeA{p:sub-14} for the proof). This result implies that 
\begin{eqnarray*}
\var[\widetilde{Q}_{a,\Omega,\lambda}(g;0)] = 
\left\{
\begin{array}{cc}
O\left(\frac{1}{\lambda^{d}}\right)& \lambda < \Omega\\
O\left(\frac{1}{\Omega^{d}}\right)& \lambda \geq \Omega\\
\end{array}
\right. 
\end{eqnarray*}
In other words, a frequency grid of $\ob_{\lambda,\kb} = \frac{2\pi
  \kb}{\lambda}$ or finer will yield a rate of convergence of
$O(\lambda^{-d})$ and $\var[\widetilde{Q}_{a,\Omega,\lambda}(g;0)]
\approx \lambda^{-d}C_{1}(\frac{a}{\Omega})$. However a coarse frequency grid $\ob_{\Omega,\kb} = \frac{2\pi
  \kb}{\Omega}$ where $\Omega > \lambda$ will yield a slower rate of
convergence of $O(\Omega^{-d})$. In the theorem below we make this
precise.

For the following theorem we consider general stationary spatial
random fields, this requires the following definition
\begin{eqnarray}
\label{eq:Dalambda}
D_{1}\left(\frac{a}{\Omega}\right) = 
\frac{1}{(2\pi)^{2d}}
\int_{[-2\pi a/\Omega,2\pi a/\Omega]^{2d}}g(\ob_{1})\overline{g(\ob_{2})}
f_{4}\left(-\ob_{1},-\ob_{2},\ob_{2}\right)d\ob_{1}d\ob_{2}. 
\end{eqnarray}

\begin{theorem}\label{theorem:summaryVar}
Suppose Assumptions \ref{assum:S}(i), \ref{assum:nonuniform},
\ref{assum:G}(i) or (ii), \ref{assum:GG}(b,c) and
\ref{assum:nonGaussian} hold. Let $C_{1}(\cdot)$ and $D_{1}(\cdot)$ 
be defined as in (\ref{eq:Calambda}) and (\ref{eq:Dalambda}) respectively. 
\begin{itemize}
\item[(i)] If a fine frequency grid is used ($\frac{\lambda}{\Omega}<1$)  then 
\begin{eqnarray*}
\lambda^{d}\var\left[\widetilde{Q}_{a,\Omega,\lambda}(g;0)\right] = 
C_{1}\left(\frac{a}{\Omega}\right)\left[
\frac{\lambda}{\Omega}\sum_{\kb=-2a}^{2a}\Sinc^{2}\left( \frac{\lambda}{\Omega}\kb
  \pi\right)\right] + D_{1}\left(\frac{a}{\Omega}\right)
+ O\left(\ell_{a,\Omega,\lambda}^{(2)}\right)
\end{eqnarray*}
\item[(ii)] If a course frequency grid is used
  ($\frac{\lambda}{\Omega}\geq 1$) then
\begin{eqnarray*}
\Omega^{d}\var[\widetilde{Q}_{a,\Omega,\lambda}(g;0)] = C_{1}\left(\frac{a}{\Omega}\right)
\left[\sum_{\kb=-2a}^{2a}\Sinc^{2}\left( \frac{\lambda}{\Omega}\kb 
  \pi\right)\right] + \left(\frac{\Omega}{\lambda}\right)^{d}D_{1}\left( \frac{a}{\Omega}\right)
+ O\left(\widetilde{\ell}_{a,\Omega,\lambda}^{(2)}\right)
\end{eqnarray*}
\end{itemize}
where $\widetilde{\ell}_{a,\Omega,\lambda}^{(2)}$ is defined in
(\ref{eq:widetildeell2}), Section \ref{sec:fine-grid}, \citeA{p:sub-14}.
\end{theorem}
As mentioned above one important implication of the above result is
that the rate of convergence depends on whether the sampling frequency
on the frequency grid is coarser or finer than $1/\lambda$, where
$\lambda$ is the length of the spatial domain. In terms of 
the asymptotic sampling properties
(see Lemma \ref{lemma:summaryE} and Theorem \ref{theorem:summaryVar})
there seems to be little benefit using a very fine frequency grid, as
it does not reduce the bias or variance (but is computationally
costly). 

We observe that if the spatial process is non-Gaussian then an
additional term, $D_{1}(a/\Omega)$, arises.
However, if the frequency grid is extremely coarse in the sense that
$\lambda/\Omega \rightarrow \infty$ as $\lambda \rightarrow \infty$
and $\Omega \rightarrow \infty$, then the fourth order cumulant
$D_{1}(a/\Omega)$ is asymptotic negligible  compared with the leading
term which is a function of the spectral density. For example, if $\lambda/\Omega \in \mathbb{Z}^{+}$ then 
\begin{eqnarray*}
\Omega^{d}\var[\widetilde{Q}_{a,\Omega,\lambda}(g;0)] = C_{1}\left(\frac{a}{\Omega}\right)
+ O\left(\ell_{a,\Omega,\lambda}^{(2)} + \frac{\Omega^{d}}{\lambda^{d}}\right).
\end{eqnarray*}
Thus if $\Omega/\lambda \rightarrow \infty$ as $\lambda\rightarrow
\infty$ and $\Omega\rightarrow \infty$ (and $\overline{g(\ob)} = g(-\ob)$)
we have verified condition (C.4) in  \citeA{p:ban-lah-nor-13};
\begin{eqnarray*}
\frac{\var\left[\sum_{\kb =
      -a}^{a}g(\ob_{\Omega,\kb})|J_{n}(\ob_{\Omega,\kb})|^{2}\right]}{2\sum_{\kb
      = -a}^{a}|g(\ob_{\Omega,\kb})|^{2}f(\ob_{\Omega,\kb})^{2}}
\Pcon 1,
\end{eqnarray*}
which is required for their proposed spatial spectral empirical
likelihood methodology. Therefore a very coarse grid has the advantage that the term $D_{1}(\cdot)$ is
negligible and thus does not need to estimated. However, we see from 
Lemma \ref{lemma:summaryE} and Theorem \ref{theorem:summaryVar}
that the disadvantage is that there is a substantial increase in both
variance and bias. 

Since the grid, $\ob_{\kb} = 2\pi \kb/\lambda$, yields optimal samping properties 
in Section \ref{sec:uniform} we focus on deriving 
sampling properties of $\widetilde{Q}_{a,\lambda}(g;\rb)$, where 
\begin{eqnarray}
\label{eq:tildeQQQ}
\widetilde{Q}_{a,\lambda}(g;\rb) = \frac{1}{\lambda^{d}}\sum_{k_{1},\ldots,k_{d}=-a}^{a}
g(\ob_{\kb})J_{n}(\ob_{\kb})
\overline{J_{n}(\ob_{\kb+\rb})}
- \frac{1}{n}\sum_{\kb=-a}^{a}g(\ob_{\kb})
\frac{1}{n}\sum_{j=1}^{n}Z(\sbb_{j})^{2}e^{-i\sbb_{j}^{\prime}\ob_{\rb}}.
\end{eqnarray}
We consider the case that the locations come from a
general sampling scheme (not just uniform) and $\rb \neq 0$. As mentioned above on this
sampling grid the fourth order cumulant term $D_{1}(\cdot)$ is not negligible and for
inference needs to be estimated. However, under the condition that the
sampling locations are uniformly distributed the variance can be
estimated using $\{\widetilde{Q}_{a,\lambda}(g;\rb);\rb \neq
0\}$ and a studentized statistic constructed. 

Note that $\widetilde{Q}_{a,\lambda}(g;\rb)$ is the bias corrected
version of  $Q_{a,\lambda,\lambda}(g;\rb)$ defined in
(\ref{eq:QQQr}). In some
situations $Q_{a,\lambda}(g;\rb) =Q_{a,\lambda,\lambda}(g;\rb)$
(estimator with no bias correction) may be a more suitable estimator
than $\widetilde{Q}_{a,\lambda}(g;\rb)$, and the sampling properties of this estimator can be found in Section
\ref{appendix:Q}, \citeA{p:sub-14}. 
%This is one of the main
%motivations for our considering the sampling properties of
%$\widetilde{Q}_{a,\lambda}(g;\rb)$ when $\rb \neq 0$. 

%% file: 4_uniform.tex
\section{Sampling properties of $\widetilde{Q}_{a,\lambda}(g;\rb)$}\label{sec:uniform}

In this section we consider the sampling properties of
$\widetilde{Q}_{a,\lambda}(g;\rb)$ (defined (\ref{eq:tildeQQQ})).

We show that under the increasing domain framework
$\widetilde{Q}_{a,\lambda}(g;\rb)$ is a consistent estimator of
$I(g;\frac{a}{\lambda})$ (or some multiple of it), where
$I(g;\frac{a}{\lambda})$ is defined in (\ref{eq:I}). 
The sampling properties in the fixed domain
framework are given in Section \ref{sec:fixed-asy}.

\subsection{The expectation of $\widetilde{Q}_{a,\lambda}(g;\rb)$}\label{sec:3mean}

We start with the expectation of
$\widetilde{Q}_{a,\lambda}(g;\rb)$. We show 
if $\sup_{\ob\in \mathbb{R}^{d}}|g(\ob)|<\infty$, the choice
of $a$ does not play a significant role in the asymptotic properties
of $\widetilde{Q}_{a,\lambda}(g;\rb)$. However, if $a>>\lambda$,
the analysis of $\widetilde{Q}_{a,\lambda}(g;\rb)$ requires more delicate
techniques than those used to prove Theorem
\ref{lemma:nonuniformdft}. We start by stating some
pertinent features in the analysis of
$\widetilde{Q}_{a,\lambda}(g;\rb)$, which gives a flavour of our approach. 
By writing $\widetilde{Q}_{a,\lambda}(g;\rb)$ as a quadratic form it is straightforward to show that 
\begin{equation}
\label{eq:Eastat}
 \Ex\left[\widetilde{Q}_{a,\lambda}(g;\rb)\right] \nonumber\\
= c_{2}\sum_{\kb=-a}^{a}g(\ob_{\kb})
\frac{1}{\lambda^{d}}\int_{[-\lambda/2,\lambda/2]^{d}}c(\ub_{1}-\ub_{2})\exp(i\ob_{k}^{\prime}(\ub_{1}-\ub_{2})-i\ub_{2}^{\prime}\ob_{\rb})
h\left(\frac{\ub_{1}}{\lambda}\right)h\left(\frac{\ub_{2}}{\lambda}\right)d\ub_{1}d\ub_{2},
\end{equation}
where $c_{2}=n(n-1)/n^{2}$. The proof of 
Theorem \ref{lemma:nonuniformdft}  is
 based on making a change of variables $v = \ub_{1}-\ub_{2}$ and then systematically changing
the limits of the integral. 
This method works if the frequency grid
$[-a/\lambda,a/\lambda]^{d}$ is fixed for all $\lambda$.  
However, if the frequency grid $[-a/\lambda,a/\lambda]^{d}$ is allowed
to grow with $\lambda$, applying this brute force method to $\Ex\left[\widetilde{Q}_{a,\lambda}(g;\rb)\right]$ has 
the disadvantage that it aggregrates the errors  within
the sum of $\Ex\left[\widetilde{Q}_{a,\lambda}(g;\rb)\right]$.  Instead, 
to further the analysis, we
replace $c(\ub_{1}-\ub_{2})$ by its Fourier transform
$c(\ub_{1}-\ub_{2}) =
\frac{1}{(2\pi)^{d}}\int_{\mathbb{R}^{d}}f(\ob)\exp(i\ob^{\prime}(\ub_{1}-\ub_{2}))d\ob$
and focus on the case  that the sampling design is 
uniform; $h(\ub/\lambda) =
\lambda^{-d}I_{[-\lambda/2,\lambda/2]}(\ub)$ (later we consider
general densities). This reduces 
the first term in $\Ex\left[\widetilde{Q}_{a,\lambda}(g;\rb)\right]$ to the Fourier transforms of step functions,  
%($\lambda^{-d}I(\ub \in [-\lambda/2,\lambda/2]^{d})$), 
which is the product of sinc functions. Specifically, we obtain 
\begin{eqnarray*}
\Ex\left[\widetilde{Q}_{a,\lambda}(g;\rb)\right] &=& \frac{c_{2}}{(2\pi)^{d}}\sum_{\kb=-a}^{a}g(\ob_{\kb})\int_{\mathbb{R}^{d}}f(\ob)
\Sinc\left(\frac{\lambda \ob}{2}+\kb\pi\right)\Sinc\left(\frac{\lambda \ob}{2}+(\kb+\rb)\pi\right)d\ob \\
 &=&  \frac{c_{2}}{\pi^{d}}\int_{\mathbb{R}^{d}}
\Sinc({\boldsymbol y})\Sinc({\boldsymbol y}+\rb\pi)\left[\frac{1}{\lambda^{d}}
\sum_{\kb=-a}^{a}g(\ob_{\kb})f(\frac{2{\boldsymbol y}}{\lambda}-\ob_{\kb})\right]
d{\boldsymbol y},
\end{eqnarray*}
where the last line above is due to a change of variables ${\boldsymbol y}=\frac{\lambda \ob}{2}+\kb\pi$.
Since the spectral density function is absolutely integrable it is clear that $\left[\frac{1}{\lambda^{d}}
\sum_{\kb=-a}^{a}g(\ob_{\kb})f(\frac{2{\boldsymbol y}}{\lambda}-\ob_{\kb})\right]$ is uniformly bounded over 
${\boldsymbol y}$ and that $\Ex\left[\widetilde{Q}_{a,\lambda}(g;\rb)\right]$ is finite for all $\lambda$.
Furthermore, if $f(\frac{2{\boldsymbol y}}{\lambda}-\ob_{\kb})$ were replaced with 
$f(-\ob_{\kb})$, then what remains in the integral are two shifted $\sinc$ functions, which is zero if 
$\rb\in \mathbb{Z}^{d}/\{0\}$, i.e.  
\begin{eqnarray*}
\Ex\left[\widetilde{Q}_{a,\lambda}(g;\rb)\right] 
 &=&  \frac{c_{2}}{\pi^{d}}\int_{\mathbb{R}^{d}}
\Sinc({\boldsymbol y})\Sinc({\boldsymbol y}+\rb\pi)\left[\frac{1}{\lambda^{d}}
\sum_{\kb=-a}^{a}g(\ob_{\kb})f(-\ob_{\kb})\right]
d{\boldsymbol y} + R_{},
\end{eqnarray*}
where
\begin{eqnarray*}
R_{} &=& \frac{c_{2}}{\pi^{d}}\int_{\mathbb{R}^{d}}
\Sinc({\boldsymbol y})\Sinc({\boldsymbol y}+\rb\pi)\left[\frac{1}{\lambda^{d}}
\sum_{\kb=-a}^{a}g(\ob_{\kb})\left(f(\frac{2{\boldsymbol y}}{\lambda}-\ob_{\kb})-f(-\ob_{\kb})\right)\right]
d{\boldsymbol y}. 
\end{eqnarray*}
In the following theorem we show
that under certain conditions on $f$, $R$ is asymptotically negligible. 
Let $b=b(\rb)$ denote the number of zero elements in the vector $\rb\in
\mathbb{Z}^{d}$. 
%In Theorem \ref{lemma:mean-stat} we show that under certain conditions on $f$ and the rate of 
%growth of $a$ with respect to $\lambda$ that this difference is asymptotically negligible. 

\begin{theorem}\label{lemma:mean-stat}
Let $I(g;\cdot)$ be defined as in (\ref{eq:I}). 
Suppose Assumptions \ref{assum:S}(i) and \ref{assum:uniform} hold. 
\begin{itemize}
\item[(i)] If Assumptions \ref{assum:G}(i) and \ref{assum:GG}(a,c) hold,
then we have 
\begin{eqnarray}
\label{eq:QaiT}
\Ex\left[\widetilde{Q}_{a,\lambda}(g;\rb)\right] 
= 
\left\{
\begin{array}{cl}
O(\frac{1}{\lambda^{d-b}}) & \rb\in \mathbf{Z}^{d}/\{0\} 
 \\
I\left(g;C\right)+ O(\frac{1}{\lambda}) & \rb=\bf{0} \\
\end{array}
\right. \qquad
\end{eqnarray}
\item[(ii)] Suppose Assumptions \ref{assum:G}(ii) holds and  
\begin{itemize}
\item[(a)] Assumption\ref{assum:GG}(b) holds,
then 
$\sup_{a}\left|\Ex\left[\widetilde{Q}_{a,\lambda}(g;\rb)\right]\right|<\infty$.
\item[(b)] Assumption \ref{assum:GG}(b,c,d) holds and 
$\{m_{1},\ldots,m_{d-b}\}$ is the subset of
non-zero values in $\rb =(r_{1},\ldots,r_{d})$, then we have  
\begin{eqnarray}
\label{eq:QaiiT}
\Ex\left[\widetilde{Q}_{a,\lambda}(g;\rb)\right] 
=\left\{ 
\begin{array}{cl}
 O\left(\frac{1}{\lambda^{d-b}}
 \prod_{j=1}^{d-b}\left(\log\lambda + \log|m_{j}|\right)\right) & \rb\in \mathbf{Z}^{d}/\{0\} \\
I\left(g;\frac{a}{\lambda}\right)+ O\big(\frac{\log\lambda }{\lambda}+\frac{1}{n}\big)
& \rb = \bf{0}
\end{array}
\right..
\end{eqnarray}
\item[(c)]If only  Assumption \ref{assum:GG}(b,c) holds, then the
$O\left( \frac{1}{\lambda^{d-b}}
\prod_{j=1}^{d-b}\left(\log\lambda + \log|m_{j}|\right)\right)$ term
in (b)
is replaced with the slower rate $O\left( \frac{1}{\lambda}\left(\log\lambda + \log\|\rb_{}\|_{1}\right)\right)$.
\end{itemize}
Note that the above bounds for (b) and (c) are uniform in $a$.
\end{itemize}
\end{theorem}
{\bf PROOF} See Section \ref{sec:mean}. \hfill $\Box$

\vspace{3mm}
We observe that if $\rb\neq \bf{0}$, then
$\widetilde{Q}_{a,\lambda}(g;\rb)$ is estimating zero. It would appear that these terms don't contain
any useful information, however in Section \ref{sec:variance2} we show how
these terms can be used to  estimate nuisance parameters. 

In order to analyze $\Ex\left[\widetilde{Q}_{a,\lambda}(g;\rb)\right]$
in the case that the locations are not from a uniform distribution we
return to (\ref{eq:Eastat}) and replace $c(\ub_{1}-\ub_{2})$ and $h(\cdot)$ by their
Fourier representations 
\begin{eqnarray*}
&&\Ex\left[\widetilde{Q}_{a,\lambda}(g;\rb)\right] \\
&=& \frac{c_{2}}{\pi^{d}}\sum_{\jb_{1},\jb_{2}\in \mathbb{Z}}\gamma_{\jb_{1}}\gamma_{\jb_{2}}
\frac{1}{\lambda^{d}}\sum_{\kb=-a}^{a}g(\omega_{k})
\int_{-\infty}^{\infty}f\left(\frac{2{\boldsymbol
      y}}{\lambda}-\ob_{\kb}\right)\Sinc({\boldsymbol
  y})\Sinc({\boldsymbol y}+(\rb-\jb_{1}-\jb_{2})\pi)d{\boldsymbol y}.
\end{eqnarray*}
This representation allows us to use similar techniques to those used
in the uniform sampling case to prove the following result.

\begin{theorem}\label{theorem:nonuniformmean}
Let $I(g;\cdot)$ be defined as in (\ref{eq:I}). 
Suppose Assumptions \ref{assum:S}(i) and \ref{assum:nonuniform} hold. 
\begin{itemize}
\item[(i)] If in addition Assumptions \ref{assum:G}(i) and \ref{assum:GG}(a,c) hold, then we have 
\begin{eqnarray*}
\Ex\big[ \widetilde{Q}_{a,\lambda}(g;\rb)\big] =\langle
\gamma,\gamma_{\rb}\rangle I\left(g;\frac{a}{\lambda}\right)
+ O(\lambda^{-1}),
\end{eqnarray*}
$O(\lambda^{-1})$ is uniform over $\rb\in \mathbb{Z}^{d}$.
\item[(ii)] If in addition Assumptions  \ref{assum:G}(ii) and \ref{assum:GG}(b,c) hold, then we have 
\begin{eqnarray*}
\Ex\big[ \widetilde{Q}_{a,\lambda}(g;\rb)\big] =\langle \gamma,\gamma_{-\rb}\rangle
I\left(g;\frac{a}{\lambda}\right)
+ O\left(\frac{\log\lambda +I(\rb \neq {\bf 0})\log \|\rb\|_{1}}{\lambda}\right).
\end{eqnarray*}
\end{itemize}
\end{theorem}
{\bf PROOF} See Section \ref{sec:proof-nonuniform} in \cite{p:sub-14}. \hfill $\Box$

\vspace{3mm}
We observe that by applying Theorem \ref{theorem:nonuniformmean} to the case
that $h$ is uniform (using that $\gamma_{0}=1$ else $\gamma_{\jb}=0$)
gives $\Ex\big[ \widetilde{Q}_{a,\lambda}(g;\rb)\big]
=O(\lambda^{-1})$ for $\rb\neq 0$. Hence, in the case that the
sampling is uniform, Theorems
\ref{lemma:mean-stat}  and
\ref{theorem:nonuniformmean} give similar results, though the bounds
in Theorem \ref{lemma:mean-stat} are sharper.

%Furthermore, by comparing bounds we observe that in the case the frequency grid is bounded
%Theorems \ref{lemma:mean-stat}(i) and \ref{theorem:nonuniformmean}(i) 
%can be proved under Assumption \ref{assum:GG}(b,c,d) rather than \ref{assum:GG}(a,c),
% however the error bounds would be slightly larger.

\begin{remark}[Estimation of $\sum_{\jb\in
    \mathbb{Z}^{d}}|\gamma_{\jb}|^{2}$]\label{remark:gamma}
The above lemma implies that $\Ex\big[
\widetilde{Q}_{a,\lambda}(g;0)\big] =\langle \gamma,\gamma_{{\bf 0}}\rangle
I\left(g;\frac{a}{\lambda}\right)$. 
%\begin{eqnarray*}
%\Ex\big[ \widetilde{Q}_{a,\lambda}(g;0)\big] \approx \frac{\langle
%  \gamma,\gamma_{{\bf 0}} \rangle}{(2\pi)^{d}}
%\int_{2\pi[-a/\lambda,a/\lambda]^{d}}g(\ob)f(\ob)d\ob 
%\end{eqnarray*}
Therefore, to estimate
$I\left(g;\frac{a}{\lambda}\right)$ we require an estimator of $\langle
  \gamma,\gamma_{{\bf 0}} \rangle$. To do this, we recall that 
\begin{eqnarray*}
\langle \gamma,\gamma_{{\bf 0}} \rangle = \sum_{\jb\in
    \mathbb{Z}}|\gamma_{\jb}|^{2} = 
\frac{1}{\lambda^{d}}\int_{[-\lambda/2,\lambda/2]^{d}}h_{\lambda}(\ob)^{2}d\ob. 
\end{eqnarray*}
Therefore one method for estimating the above integral is to 
define a grid on $[-\lambda/2,\lambda/2]^{d}$ and estimate 
$h_{\lambda}$ at each point, then to take the average squared
over the grid (see Remark 1,
\citeA{p:mat-yaj-09}). An alternative, computationally simpler method,
is to use the method proposed in \citeA{p:gin-nic-08}, that is 
\begin{eqnarray*}
\widehat{\langle \gamma,\gamma_{{\bf 0}} \rangle}= 
\frac{2}{n(n-1)b}\sum_{1\leq j_{1}<j_{2}\leq n}K\left(\frac{\ub_{j_{1}}-\ub_{j_{2}}}{b}\right)^{2},
\end{eqnarray*}
as an estimator of $\langle \gamma,\gamma_{{\bf 0}} \rangle$, 
where $K:[-1/2,1/2]^{d}\rightarrow \mathbb{R}$ is a kernel
function. Note that multiplying the above kernel with
$\exp(-i\ob_{\rb}^{\prime}\ub_{j_{2}})$ results in an estimator of $\langle \gamma,\gamma_{\rb} \rangle$.
In the case $d=1$ and under certain regularity conditions, \citeA{p:gin-nic-08} show if the bandwidth $b$ is
selected in an appropriate way then $\widehat{\langle \gamma,\gamma_{{\bf 0}}
\rangle}$ attains the classical $O(n^{-1/2})$ rate under suitable
regularlity conditions (see, also, \citeA{p:bic-rit-88} and
\citeA{p:lau-96}). It seems plausible a similar result holds for
$d>1$ (though we do not prove it here). Therefore, an estimator of $I\left(g;\frac{a}{\lambda}\right)$
is $\widetilde{Q}_{a,\lambda}(g;\rb)/\widehat{\langle \gamma,\gamma_{{\bf 0}} \rangle}$.
\end{remark}

\subsection{The covariance and asymptotic normality}\label{sec:3var}

In the previous section we showed that the expectation of 
$\widetilde{Q}_{a,\lambda}(g;\rb)$ depends only on the number of
frequencies $a$ through the limit of the integral 
$I(g;\frac{a}{\lambda})$ (if $\sup_{\ob
  \mathbb{R}^{d}}|g(\ob)|<\infty$). In this section, we show that $a$ plays
a mild role in the higher order properties of
$\widetilde{Q}_{a,\lambda}(g;\rb)$. We focus on the case that the
random field is Gaussian and later describe how the results differ in
the case that the random field is non-Gaussian.

\begin{theorem}\label{theorem:nonuniformvar}
Suppose Assumptions \ref{assum:S}, \ref{assum:nonuniform} and  
\begin{itemize}
\item[(i)] 
Assumption \ref{assum:G}(i) and \ref{assum:GG}(a,c) hold. Then uniformly for all
  $0\leq \|\rb_{1}\|_{1},\|\rb_{2}\|_{1}\leq C|a|$ (for some finite
  constant $C$) we have 
\begin{eqnarray*}
\lambda^{d}\cov\left[\widetilde{Q}_{a,\lambda}(g;\rb_{1}),\widetilde{Q}_{a,\lambda}(g;\rb_{2})\right]=    
U_{1}(\rb_{1},\rb_{2};\ob_{\rb_{1}},\ob_{\rb_{2}}) + O\left(\frac{1}{\lambda^{}}+\frac{\lambda^{d}}{n}\right) 
\end{eqnarray*}
and 
\begin{eqnarray*}
 \lambda^{d}\cov\left[\widetilde{Q}_{a,\lambda}(g;\rb_{1}),\overline{\widetilde{Q}_{a,\lambda}(g;\rb_{2})}\right]=  
U_{2}(\rb_{1},\rb_{2};\ob_{\rb_{1}},\ob_{\rb_{2}})+ O\left(\frac{1}{\lambda^{}}+\frac{\lambda^{d}}{n}\right) 
\end{eqnarray*}
\item[(ii)]   Assumption \ref{assum:G}(ii) and \ref{assum:GG}(b) hold. Then
  we have 
\begin{eqnarray*}
\lambda^{d}\sup_{a,\rb_{}}\var\left[\widetilde{Q}_{a,\lambda}(g;\rb_{})\right]
<\infty \quad\textrm{and}\quad
\lambda^{d}\sup_{a,\rb_{}}\left|\cov\left[\widetilde{Q}_{a,\lambda}(g;\rb_{}),
\overline{\widetilde{Q}_{a,\lambda}(g;-\rb_{})}\right]\right|<\infty,
\end{eqnarray*}
if $\lambda^{d}/n\rightarrow c$ (where $0\leq c<\infty$) as $\lambda\rightarrow \infty$ and
$n\rightarrow \infty$. 
\item[(iii)] Assumption \ref{assum:G}(ii) and \ref{assum:GG}(b,c)
  hold. Then uniformly for all
  $0\leq \|\rb_{1}\|_{1},\|\rb_{2}\|_{1}\leq C|a|$ (for some finite constant $C$) we have 
\begin{eqnarray*}
\lambda^{d}\cov\left[\widetilde{Q}_{a,\lambda}(g;\rb_{1}),\widetilde{Q}_{a,\lambda}(g;\rb_{2})\right]=
U_{1}(\rb_{1},\rb_{2};\ob_{\rb_{1}},\ob_{\rb_{2}})  + O(\ell_{\lambda,a,n}) 
\end{eqnarray*}
\begin{eqnarray*}
\lambda^{d}\cov\left[\widetilde{Q}_{a,\lambda}(g;\rb_{1}),\overline{\widetilde{Q}_{a,\lambda}(g;\rb_{2})}\right]=
U_{2}(\rb_{1},\rb_{2};\ob_{\rb_{1}},\ob_{\rb_{2}})
+ O(\ell_{\lambda,a,n}), 
\end{eqnarray*}
where
\begin{eqnarray}
\label{eq:ell-def}
\ell_{\lambda,a,n} = 
\log^{2}(a)\bigg[
\frac{\log a+\log \lambda}{\lambda}\bigg] +
\frac{\lambda^{d}}{n}.
\end{eqnarray}
The expressions 
$U_{1}(\rb_{1},\rb_{2};\ob_{\rb_{1}},\ob_{\rb_{2}})$ and
$U_{2}(\rb_{1},\rb_{2};\ob_{\rb_{1}},\ob_{\rb_{2}})$ are rather
cumbersome and are defined in (\ref{eq:Ur}), \citeA{p:sub-14}.
Note, if we drop the restriction on $\rb_{1}$ and $\rb_{2}$ and simply let
$\rb_{1},\rb_{2}\in \mathbb{Z}^{d}$ then the bound for
$\ell_{\lambda,a,n}$ needs to include the additional term $\log^{2}(a)(\log
\|\rb_{1}\|_{1} + \log \|\rb_{2}\|_{1})/\lambda$.
\end{itemize}

\end{theorem}
{\bf PROOF} See Section \ref{sec:proof-nonuniform} in \citeA{p:sub-14}. \hfill $\Box$

\vspace{3mm}

We now briefly discuss the above results. 
From Theorem \ref{theorem:nonuniformvar}(iii)  we see that
$\widetilde{Q}_{a,\lambda}(g;\rb)$ is a mean squared consistent estimator of $\langle \gamma, \gamma_{\rb} \rangle
I(g;\frac{a}{\lambda})$, i.e. 
$\Ex[\widetilde{Q}_{a,\lambda}(g;\rb) - \langle \gamma, \gamma_{\rb}
\rangle I(g;\frac{a}{\lambda})]^{2} = 
O(\lambda^{-d}+(\frac{\log \lambda}{\lambda}+\frac{1}{n})^{2})$
as $a\rightarrow \infty$ and $\lambda \rightarrow \infty$. 
%Indeed, if $\sup_{\ob\in \mathbb{R}^{d}}|g(\ob)|<\infty$, then the rate that $a/\lambda \rightarrow \infty$ plays no role.
%\begin{eqnarray*}
%\Ex\left|Q_{a,\lambda}(g_{};0) - \langle {\boldsymbol
%    \gamma},{\boldsymbol \gamma}_{0} \rangle I(g;\infty) \right|^{2} =
%O\left(\frac{1}{\lambda^{d}} +
%  \left(\frac{\lambda}{a}\right)^{\beta-1}  \right),
%\end{eqnarray*}

In order to obtain an explicit expression for the variance
additional conditions are required. 
%In the case that the frequency grid is bounded we obtain an expression for the
%variance of $\widetilde{Q}_{a,\lambda}(g;\rb)$. On the other hand, we
In particular, Theorem \ref{theorem:nonuniformvar}(iii) states 
that if the frequency grid is unbounded we require 
some additional conditions on the spectral density function and 
some mild constraints on the rate of grow of the frequency
domain $a$.  More precisely, $a$ should be such that $a = O(\lambda^{k})$
for some $1\leq k < \infty$.
If these conditions are fulfilled, then the  asymptotic  variance of
$\widetilde{Q}_{a,\lambda}(g;\rb)$ (up to the
limits of an integral) are equivalent for both the bounded and unbounded  
frequency grid. Furthermore by choosing $a=\lambda^{k}$, 
Theorems \ref{theorem:nonuniformmean} and \ref{theorem:nonuniformvar}(i) yields 
 the following bound for the mean squared error when estimating $I(g;\infty)$
\begin{eqnarray*}
\Ex\left(\widetilde{Q}_{a,\lambda}(g;0) - \langle \gamma, \gamma_{\rb}
\rangle I(g;\infty) \right)^{2} = 
O\left(\frac{1}{\lambda^{d}} + \left[\frac{\log \lambda }{\lambda} + \frac{1}{n} + \left(\frac{\lambda}{\lambda^{k}}\right)^{\delta}\right]^{2}\right)
\end{eqnarray*}
where $f(\ob)\leq \beta_{\delta}(\ob)$. This bound guides
the choice of $k$. Since the rate of convergence cannot better
$O(\lambda^{-1})$ we choose $k$ such that
$(\lambda/\lambda^{k})^{\delta} = \lambda^{-1}$, 
thus $k = \delta^{-1}+1$ and we choose $a=\lambda^{1+\delta^{-1}}$.
%where the rate of decay, $\delta$, of the tails of the spectral
%density can be  estimated by plotting $\|\ob\|_{1}$ against $\log
%\widehat{f}_{\lambda,n}(\ob)$ (where the spectral density estimator
%$\widehat{f}_{\lambda,n}$ is defined in Example
%\ref{example:fixed}(ii)) and using the gradient of the line of best
%it as an estimator of $\delta+1$.

\begin{remark}[Selecting $a$ in practice]
The above  gives theoretical guidelines. In practice, if
$\sup_{\ob\in \mathbb{R}^{d}}|g(\ob)|<\infty$ we suggest
plotting $|J_{n}(\ob_{\kb})|^{2}$ against
$\ob_{\kb}$. $|J_{n}(\ob_{\kb})|^{2}$ will drop  close to zero for
large $|\ob_{\kb}|$ (see Figure \ref{fig:1}, Section
\ref{sec:simulations}, \citeA{p:sub-14}). Thus $a$ should be selected such that it lies after
this point. The precise value does not matter too much as the results
are not too sensitive to the choice of $a$. 
\end{remark}

By comparing Theorem \ref{theorem:nonuniformvar}(i)  and (ii)
 we observe that in the case the frequency grid is bounded, the same
 result can be proved under
Assumption \ref{assum:GG}(b,c) rather than \ref{assum:GG}(a,c),
 however the error bounds are slightly larger.

\vspace{3mm}

The expressions for
$\cov[\lambda^{d}\cov\left[\widetilde{Q}_{a,\lambda}(g;\rb_{1}),\overline{\widetilde{Q}_{a,\lambda}(g;\rb_{2})}\right]]$
(see (\ref{eq:Ur})) are unwieldy, however, some
simplifications can be made if $\|\rb_{1}\|_{1}<<\lambda$ and $\|\rb_{2}\|_{1}<<\lambda$.
\begin{corollary}\label{cor:C}
Suppose Assumptions \ref{assum:nonuniform}, \ref{assum:G} and 
\ref{assum:GG}(a,c) or \ref{assum:GG}(b,c) hold.  
Then we have 
\begin{eqnarray}
U_{1}(\rb_{1},\rb_{2};\ob_{\rb_{1}},\omega_{\rb_{2}}) &=& \Gamma_{\rb_{1}-\rb_{2}}C_{1}
+ O\left(\frac{\|\rb_{1}\|_{1}+\|\rb_{2}\|_{1}+1}{\lambda}\right) \nonumber\\
U_{2}(\rb_{1},\rb_{2};\ob_{\rb_{1}},\omega_{\rb_{2}}) &=& \Gamma_{\rb_{1}+\rb_{2}}C_{2}
+ O\left(\frac{\|\rb_{1}\|_{1}+\|\rb_{2}\|_{1}+1}{\lambda}\right)\label{eq:C1C2approx}
\end{eqnarray}
where $\Gamma_{\rb} =
\sum_{\jb_{1}+\jb_{2}+\jb_{3}+\jb_{4}=\rb}\gamma_{\jb_{1}}\gamma_{\jb_{2}}\gamma_{\jb_{3}}\gamma_{\jb_{4}}$
and 
\begin{eqnarray*}
C_{1}
&=&\frac{1}{(2\pi)^{d}}\int_{2\pi[-a/\lambda,a/\lambda]^{d}}f(\ob)^{2}\left[|g(\ob)|^{2}+g(\ob)\overline{g(-\ob)}\right]d\ob \\
C_{2}
&=&\frac{1}{(2\pi)^{d}}\int_{2\pi[-a/\lambda,a/\lambda]^{d}}f(\ob)^{2}\left[g(\ob)g(-\ob)+g(\ob)g(\ob)\right]d\ob. 
\end{eqnarray*}
Recall that $C_{1} = C_{1}(a/\lambda)$ 
(where $C_{1}(\cdot)$ is defined in (\ref{eq:Calambda})).
\end{corollary}
Using the above expressions for
$\lambda^{d}\cov\left[\widetilde{Q}_{a,\lambda}(g;\rb_{1}),\widetilde{Q}_{a,\lambda}(g;\rb_{2})\right]$
and
$\lambda^{d}\cov\left[\widetilde{Q}_{a,\lambda}(g;\rb_{1}),\overline{\widetilde{Q}_{a,\lambda}(g;\rb_{2})}\right]$ 
the variance of the real and imaginary parts of
$\widetilde{Q}_{a,\lambda}(g;\rb_{})$ can easily be deduced. 

In the following theorem we derive bounds
for the cumulants of $\widetilde{Q}_{a,\lambda}(g;\rb)$,
which are subsequently used to show asymptotical normality of
$\widetilde{Q}_{a,\lambda}(g;\rb)$.

\begin{theorem}\label{theorem:cumulantsnonuniform}
Suppose Assumptions \ref{assum:S}, \ref{assum:nonuniform}, 
\ref{assum:G} and \ref{assum:GG}(b) hold. Then for all
$q\geq 3$ and uniform in $\rb_{1},\ldots,\rb_{q}\in \mathbb{Z}^{d}$ we have 
\begin{eqnarray}
\label{eq:cumqa2}
\cum_{q}\left[\widetilde{Q}_{a,\lambda}(g,\rb_{1}),\ldots,\widetilde{Q}_{a,\lambda}(g,\rb_{q})\right]
&=&  O\bigg(\frac{\log^{2d(q-2)}(a)}{\lambda^{d(q-1)}}\bigg)
%  + 
% \frac{\log^{2}(n)}{\lambda n } + \frac{1}{n^{2}}\bigg).
\end{eqnarray}
if $\frac{\lambda^{d}}{n\log^{2d}(a)}\rightarrow 0$ as $n\rightarrow \infty$, 
$a\rightarrow \infty$ and $\lambda \rightarrow \infty$.
\end{theorem}
{\bf PROOF} See Section \ref{sec:proof-nonuniform} in \citeA{p:sub-14}. \hfill $\Box$

\vspace{2mm}
From the above theorem we see that if
$\frac{\lambda^{d}}{n\log^{2d}(a)}\rightarrow 0$ and
$\log^{2}(a)/\lambda^{1/2}\rightarrow 0$ as $\lambda\rightarrow
\infty$, $n\rightarrow \infty$ and $a\rightarrow \infty$, then we have 
$\lambda^{dq/2}\cum_{q}(\widetilde{Q}_{a,\lambda}(g,\rb_{1}),\ldots,\widetilde{Q}_{a,\lambda}(g,\rb_{q}))\rightarrow
0$ for all $q\geq 3$. Using this result we show
asymptotic Gaussianity of $\widetilde{Q}_{a,\lambda}(g,\rb_{})$.

\begin{theorem}\label{theorem:CLTnu}
Suppose Assumptions \ref{assum:S}, \ref{assum:nonuniform},  
\ref{assum:G} and \ref{assum:GG}(b,c) hold. Let $C_{1}$ and  
$C_{2}$, be defined as in Corollary \ref{cor:C}.
%We define the $m$-dimension complex random vectors $\widetilde{\boldsymbol{Q}}_{m} =
%(\widetilde{Q}_{a,\lambda}(g,\rb_{1}),\ldots,\widetilde{Q}_{a,\lambda}(g,\rb_{m}))$,
%where $\rb_{1},\ldots,\rb_{m}$ are such that $\rb_{i}\neq -\rb_{j}$ and $\rb_{i}\neq
%0$. %then for any $i,j\in \{1,\ldots,m\}$ 
Under these conditions we have 
\begin{eqnarray*}
\lambda^{d/2}\Delta^{-1/2}
\left(
\begin{array}{c}
\Re\left( \widetilde{Q}_{a,\lambda}(g,\rb_{1}) -  \langle
  \gamma,\gamma_{-\rb}\rangle I(g;\frac{a}{\lambda})\right)
\\
\Im\left( \widetilde{Q}_{a,\lambda}(g,\rb_{1}) - \langle \gamma,\gamma_{-\rb}\rangle
I(g;\frac{a}{\lambda})
\right)
\end{array}
\right) 
\Dcon \mathcal{N}(0,I_{2}),
\end{eqnarray*}
where $\Re X$ and $\Im X$ denote the real and imaginary parts of the
random variable $X$, and 
\begin{eqnarray*}
\Delta = 
\frac{1}{2}\left(
\begin{array}{cc}
\Re(\Gamma_{0}C_{1}+\Gamma_{2\rb}C_{2}) & -\Im(\Gamma_{0}C_{1}+\Gamma_{2\rb}C_{2})\\
-\Im(\Gamma_{0}C_{1}+\Gamma_{2\rb}C_{2}) & \Re(\Gamma_{0}C_{1}-\Gamma_{2\rb}C_{2})\\
\end{array}
\right)
\end{eqnarray*}
with $\frac{\log^{2}(a)}{\lambda^{1/2}}\rightarrow
0$ and $\lambda^{d}/n\rightarrow 0$  as $\lambda \rightarrow \infty$, $n\rightarrow \infty$ and
$a\rightarrow \infty$.
\end{theorem}
{\bf PROOF} See Section \ref{sec:cumulants}, \citeA{p:sub-14}. 
\hfill $\Box$

\vspace{2mm}
It is likely that the above result also holds when the
assumption of Gaussianity of the spatial random field is relaxed and
replaced with  the conditions stated in 
Theorem \ref{theorem:nonGaussiannonUniform} (below) together with some mixing-type 
assumptions. We leave this for future work.   However, in 
the following theorem, we obtain an expression for the variance
of $\widetilde{Q}_{a,\lambda}(g;\rb)$ for non-Gaussian
random fields. 

\begin{theorem}\label{theorem:nonGaussiannonUniform}
Let us suppose that  $\{Z(\ub);\ub \in\mathbb{R}^{d}\}$ is a fourth order stationary spatial 
random field that satisfies Assumption \ref{assum:S}(i),  
\ref{assum:nonuniform}, \ref{assum:G}, \ref{assum:nonGaussian} and \ref{assum:GG}(a,c) or \ref{assum:GG}(b,c)  are satisfied. 

If $\|\rb\|_{1}, \|\rb\|_{2}<<\lambda$, then we have 
\begin{eqnarray}
\label{eq:nonGaussian1nu}
\lambda^{d}\cov\left[\widetilde{Q}_{a,\lambda}(g;\rb_{1}),\widetilde{Q}_{a,\lambda}(g;\rb_{2})\right]=
\Gamma_{\rb_{1}-\rb_{2}}(C_{1}  + D_{1})+
O\left(\ell_{\lambda,a,n}+ \frac{(a\lambda)^{d}}{n^{2}}+\frac{\|\rb_{1}\|_{1}+\|\rb_{2}\|_{1}}{\lambda}\right) 
\end{eqnarray}
\begin{eqnarray}
\label{eq:nonGaussian2nu}
\lambda^{d}\cov\left[\widetilde{Q}_{a,\lambda}(g;\rb_{1}),\overline{\widetilde{Q}_{a,\lambda}(g;\rb_{2})}\right]=
\Gamma_{\rb_{1}+\rb_{2}}(C_{2} + D_{2}) + O\left(\ell_{\lambda,a,n}+\frac{(a\lambda)^{d}}{n^{2}}+\frac{\|\rb_{1}\|_{1}+\|\rb_{1}\|_{1}}{\lambda}\right),
\end{eqnarray}
where $C_{1}$ and $C_{2}$ are defined as in Corollary \ref{cor:C} and 
\begin{eqnarray*}
D_{1}&=& \frac{1}{(2\pi)^{2d}}\int_{2\pi[-a/\lambda,a/\lambda]^{2d}}g(\ob_{1})\overline{g(\ob_{2})}
f_{4}(-\ob_{1},-\ob_{2},\ob_{2})d\ob_{1}d\ob_{2} \\
D_{2} &=& \frac{1}{(2\pi)^{2d}}\int_{2\pi[-a/\lambda,a/\lambda]^{2d}}g(\ob_{1})g(\ob_{2})
f_{4}(-\ob_{1},-\ob_{2},\ob_{2})d\ob_{1}d\ob_{2}.
\end{eqnarray*}
Note that we can drop the restriction that $\|\rb\|_{1},
\|\rb\|_{2}<<\lambda$ and let $\rb_{1},\rb_{2}\in \mathbb{Z}^{d}$,
under this more general assumption we obtain expressions that are
similar to those for $U_{1}$ and $U_{2}$ (given in Theorem \ref{theorem:nonuniformvar})
and we need to include $\log^{2}(a)[\|\rb\|_{1}+\|\rb_{2}\|_{1}]/\lambda$
in the bounds. 
\end{theorem}
{\bf PROOF} See Section \ref{sec:proof-nonuniform}, \citeA{p:sub-14}.\hfill $\Box$

\vspace{2mm}
We observe that to ensure the term
$\frac{(a\lambda)^{d}}{n^{2}}\rightarrow 0$ we need to choose $a$ such
that $a^{d} = o(n^{2}/\lambda^{d})$. In contrast for Gaussian
random fields $a = O(\lambda^{k})$ for some $1\leq k <\infty$ was
sufficient for obtaining an expression for the variance and asymptotic
normality. 

\subsubsection*{Mixed Domain verses Pure Increasing Domain asymptotics}\label{sec:alt-asymptotics}

The asymptotics in this paper are mainly done using mixed domain asymptotics,
that is, as the domain $\lambda\rightarrow \infty$, the number of
locations observed grows at a faster rate than $\lambda$, in other
words $\lambda^{d}/n\rightarrow 0$ as $n\rightarrow
\infty$. However, as rightly pointed out by a referee, for a given
application it may be difficult to disambiguate Mixed Domain (MD) from
the Pure Increasing Domain (PID), where $\lambda^{d}/n\rightarrow c$ ($0<
c< \infty$) set-up. We briefly discuss how the results change under
PID asymptotics and the implications of this. We find that the results
point to a rather intriguing difference for spatial processes that are 
Gaussian and non-Gaussian. 

In the case that spatial process is Gaussian, using 
both MD and PID asymptotics we have 
$\lambda^{d}\var[\widetilde{Q}_{a,\lambda}(g;\rb_{})] = O(1)$ 
(see Theorem  \ref{theorem:nonuniformvar}(i)). Furthermore,  an asymptotic
expression for the variance is 
\begin{eqnarray*}
\lambda^{d}\var[\widetilde{Q}_{a,\lambda}(g;\rb_{})]=\Gamma_{0}[C_{1}+E_{1}]
+ O\left(\log^{2}(a)\bigg[
\frac{\log a+\log \lambda}{\lambda}\bigg] +\frac{\|\rb\|_{1}+1}{\lambda}\right)
\end{eqnarray*}
where $E_{1}=O(\lambda^{d}/n)$ and is a function of the spectral
density, this term is not asymptotically 
negligible under PID asymptotics.  From the above we see
that if we choose $a$ such
that $a=O(\lambda^{k})$ for some $1< k < \infty$ then the 
frequency grid is unbounded and similar results as those stated in
Sections \ref{sec:3mean} and \ref{sec:3var} hold under PID
asymptotics.  In the case that the process is non-Gaussian, using 
Theorem \ref{theorem:nonGaussiannonUniform} we have 
\begin{eqnarray*}
\lambda^{d}\var[\widetilde{Q}_{a,\lambda}(g;\rb_{})]=\Gamma_{0}[C_{1}+D_{1}+E_{1}+F_{1}]
+ O\left(\log^{2}(a)\bigg[
\frac{\log a+\log \lambda}{\lambda}\bigg] + \frac{(a\lambda)^{d}}{n^{2}}
+\frac{\|\rb\|_{1}+1}{\lambda}\right),
\end{eqnarray*}
where $F_{1}=O(\lambda^{d}/n)$ is a function of the fourth order spectral density function.
However, there arises an additional term $O((a\lambda)^{d}/n^{2})$. 
It can be see from the proof of Theorem
\ref{theorem:nonGaussiannonUniform},  if 
$\frac{(a\lambda)^{d}}{n^{2}}\rightarrow \infty$ as
$a,\lambda,n\rightarrow \infty$, then
$\lambda^{d}\var[\widetilde{Q}_{a,\lambda}(g;\rb_{})]$ is not
bounded. Thus, the number of frequencies, $a$, should be such that  $(a\lambda)^{d}/n^{2}\rightarrow 0$.
In the case of  MD asymptotics,  we
choose $a$ such that $(a\lambda)^{d}/n^{2}\rightarrow 0$ and 
$\log^{3}(a)/\lambda \rightarrow 0$. Under these two
conditions the frequency grid can be unbounded  and grow at the rate $a/\lambda$ as 
$\lambda \rightarrow \infty$. 
However, under PID asymptotics (where $\lambda = O(n^{1/d})$)
in order to ensure that $(a\lambda)^{d}/n^{2}=O(1)$
we require $a=O(n^{1/d})=O(\lambda)$. This constrains the
frequency grid to be bounded. To summarize, in the case that the
spatial process is non-Gaussian and $n$ is not that large in order that
$\var[\widetilde{Q}_{a,\lambda}(g;\rb_{})] =
O(\lambda^{-d})$, the frequency grid must be bounded or a coarser
frequency grid $\ob_{\Omega,\kb}$ (where $\lambda<\Omega$) used.

\subsection{Fixed domain asymptotics}\label{sec:fixed-asy}

We now turn our attention to asymptotic sampling properties of 
$Q_{a,\lambda}(g;0)$ when the domain
$\lambda$ is kept fixed but the number of sampling locations
$n\rightarrow \infty$. In order to simplify notation we consider the
case $d=1$. We will work under the near uniform design of locations, where $\{u_{n,j}\}$
satisfies Assumption \ref{assum:fixedDFT}. 
We first observe that the DFT of 
$\{Z(s_{n,j});j=1,\ldots,n\}$ approximates the Fourier transform of
 $\{Z(s);s\in [-\lambda/2,\lambda/2]\}$ 
 i.e. given that the spatial  autocovariance  is Lipschitz
continuous we have 
\begin{eqnarray*}
%\frac{\lambda^{1/2}}{n}\sum_{j=1}^{n}Z(s_{n,j})\exp(i2\pis_{n,j}\omega_{k}) 
J_{n}(\omega_{k})= \mathcal{J}_{\lambda}\left(\frac{k}{\lambda}\right) +
O_{p}\left(\frac{1}{\lambda^{1/2}}\sum_{j=2}^{n}\left\{\frac{\lambda}{n}(|\omega_{k}|+1)|s_{n,(j)}
    - s_{n,(j-1)}|^{1/2} + \left|\frac{\lambda}{n}-(s_{n,(j)}-s_{n,(j-1)})\right|\right\}\right)
\end{eqnarray*}
where
\begin{eqnarray*} 
\mathcal{J}_{\lambda}\left(\frac{k}{\lambda}\right) = \frac{1}{\lambda^{1/2}}\int_{-\lambda/2}^{\lambda/2}Z(s)\exp\left(\frac{2\pi i k
  s}{\lambda}\right)ds.
\end{eqnarray*}
It is the approximation of $J_{n}(\omega_{k})$ by $\mathcal{J}_{\lambda}(k/\lambda)$
that determine the sampling properties of $Q_{a,\lambda}(g;0)$. 
More precisely, in Section \ref{sec:appendixfixed}, \citeA{p:sub-14} we show that 
\begin{eqnarray*}
\var\left[\mathcal{J}_{\lambda}\left(\frac{k}{\lambda}\right)\right] =
A_{\lambda}\left(\frac{k}{\lambda}\right) \textrm{ and }
\cov\left[\mathcal{J}_{\lambda}\left(\frac{k_1}{\lambda}\right),\mathcal{J}_{\lambda}\left(\frac{k_2}{\lambda}\right)\right] =
\frac{(-1)^{k_1-k_2+1}}{\pi(k_{1}-k_{2})}\left[B_{\lambda}\left(\frac{k_1}{\lambda}\right) -
B_{\lambda}\left(\frac{k_2}{\lambda}\right)\right]
\end{eqnarray*}
where 
\begin{eqnarray}
\label{eq:ABkb}
A_{\lambda}\left(\frac{k}{\lambda}\right) =
\int_{-\lambda}^{\lambda}\left(1-\frac{|u|}{\lambda}\right)c(u)e^{2\pi i
    u/\lambda}du \textrm{ and }
B_{\lambda}\left(\frac{k}{\lambda}\right)  = \int_{0}^{\lambda}c(u)\sin\left(\frac{2\pi
k u}{\lambda}\right)du.
\end{eqnarray}
This result and Theorem \ref{theorem:fixedDFT1} are used to prove the following theorem.
%\begin{eqnarray*}
%Q_{a,\lambda}(g;0) = \frac{1}{\lambda}\sum_{k=1}^{a}g(\omega_{k})|J_{n}(\omega_{k})|^{2}
%= \sum_{k=1}^{a}g(\omega_{k})|\mathcal{J}_{\lambda}(k)|^{2} + 
%O_{p}\left( \frac{1}{n}\left[\frac{1}{a}\sum_{k=1}^{a}(|k|+1)|g(\omega_{k})|\right] \right)
%\end{eqnarray*}

\begin{theorem}
Suppose Assumptions \ref{assum:S}, \ref{assum:fixedDFT} and
\ref{assum:GG}(e) hold. Then 
\begin{eqnarray*}
\Ex\left[Q_{a,\lambda}(g;0)\right] 
&=& 
\sum_{\kb = -a}^{a}g\left(\omega_{k}\right)A_{\lambda}\left(\frac{k}{\lambda}\right) + 
O\left( \frac{1}{n}\left[\frac{1}{a}\sum_{k=1}^{a}(|k|+1)|g(\omega_{k})|\right] \right)
\end{eqnarray*}
and  if $\sup_{k}|g(\omega_{k})|<\infty$
\begin{eqnarray}
\label{eq:fixedvar}
\var\left[ Q_{a,\lambda}(g;0)\right] 
&=& 
\sum_{k_{1},k_{2} =-\infty}^{\infty}g\left(\omega_{k}\right) 
\overline{g\left(\omega_{k_1}\right)}B_{\lambda}(k_{1},k_{2}) +
O\left(\frac{a^{4}}{n^{2}} + \frac{a\log a}{n} + \frac{1}{a}\right)
\end{eqnarray}
 where 
\begin{eqnarray*}
B_{\lambda}(k_1,k_2) = \left\{
\begin{array}{cc}
I(k\neq 0)[A_{\lambda}\left(\frac{k}{\lambda}\right)^{2}  + 
\frac{1}{\pi^{2}k^{2}}B_{\lambda}(\frac{k}{\lambda})^{2}] + 2I(k=0)A_{\lambda}(0)^{2}
& k_1 = k_2 (=k) \\
\frac{1}{\pi^{2}(k_1 -
  k_2)^{2}}\left[B_{\lambda}(\frac{k_{1}}{\lambda}) -
  B_{\lambda}(\frac{k_{2}}{\lambda}) \right]^{2}+
\frac{1}{\pi^{2}(k_1 +
  k_2)^{2}}\left[B_{\lambda}(\frac{k_{1}}{\lambda}) +
  B_{\lambda}(\frac{k_{2}}{\lambda}) \right]^{2}& k_1\neq k_2 \\ 
\end{array}
\right.
\end{eqnarray*}
and $A_{\lambda}(\cdot)$ and $B_{\lambda}(\cdot)$ are defined in
(\ref{eq:ABkb}).
\end{theorem}
From the above result we see that if $|g(\omega_{k})|<\infty$ and the
number of Fourier transforms, $a$, in $Q_{a,\lambda}(g;0)$,
grows at a sufficiently slow rate as $n\rightarrow \infty$ then
$Q_{a,\lambda}(g;0)$ is asymptotically an unbiased estimator of
$\sum_{k\in \mathbb{Z}}g(\omega_{k})A_{\lambda}(\frac{k}{\lambda})$
and the variance bounded (which fits with the conclusions of Theorem
\ref{theorem:nonuniformvar}(ii)). Moreover, by replacing the
$|J_{n}(\omega_{k})|^{2}$ in $Q_{a,\lambda}(g;0)$ with
$|\mathcal{J}_{\lambda}\left(\frac{k}{\lambda}\right)|^{2}$ and 
using infill asymptotics we see that 
the distribution of
$Q_{a,\lambda}(g;0)$ will not be normal but will be a mix of dependent
chi-squares.

We now consider the case that $g(\omega_{k})$ is not bounded. In particular, we focus on
the case that the covariance of a spatial Gaussian process is 
$\sigma^{2}c(u)$ where $c(u)$ is known but $\sigma^{2}$ is
unknown and our aim is to estimate $\sigma^{2}$. Let
$A_{\lambda}(\frac{k}{\lambda}) =
\int_{-\lambda}^{\lambda}c(u)(1-|u|/\lambda)\exp(2\pi i k u/\lambda)du$. Since 
$\Ex[|J_{n}(\omega_{k})|^{2}] =
\sigma^{2}A_{\lambda}(\frac{k}{\lambda}) + O(n^{-1}(\lambda + |k|))$ 
it seems natural to use $\widehat{\sigma}^{2}$ as an estimator of
$\sigma^{2}$, where 
\begin{eqnarray*}
\widehat{\sigma}^{2} = \frac{1}{a}\sum_{k=1}^{a}\frac{|J_{n}(\omega_{k})|^{2}}{A_{\lambda}(\frac{k}{\lambda})}.
\end{eqnarray*}
Note $\widehat{\sigma}^{2}$ corresponds to the Whittle likelihood
estimator of $\sigma^{2}$ when $c(\cdot)$ is known.
%as an estimator of $\sigma^{2}$ (this estimator also
%corresponds to the variance estimator using the Whittle likelihood when $f$ is known). 
Using the above theorem we have 
\begin{eqnarray*}
\Ex[\widehat{\sigma}^{2}] =
\sigma^{2} + O\left( \frac{1}{n}\left[\frac{1}{a}\sum_{k=1}^{a}(|k|+1)|k|^{2}\right] \right).
\end{eqnarray*}
Hence, is $a$ is chosen such that the error above goes to zero as
$n\rightarrow \infty$ then $\widehat{\sigma}^{2}$ is asymptotically
and unbiased estimator of $\sigma^{2}$. Furthermore, it can be shown that 
\begin{eqnarray*}
\var[\widehat{\sigma}^{2}] 
&=& \frac{1}{a^{2}}\sum_{k=1}^{a}\frac{B_{\lambda}(k,k)}{A_{\lambda}(\frac{k}{\lambda})^{2}}
+\frac{2}{a^{2}}\sum_{k_{1}=1}^{a}\sum_{k_2 = k_1
  +1}^{a}\frac{B_{\lambda}(k_1,k_2)}{A_{\lambda}(\frac{k_1}{\lambda})
  A_{\lambda}(\frac{k_2}{\lambda})} + O\left(\frac{a^{6}}{n^{2}} + \frac{a^{3}}{n}\right)
\end{eqnarray*}
But here we run into a problem. One would expect that if $a$ is
chosen such that $a^{4}<<n$ that
$\var[\widehat{\sigma}^{2}] = O(a^{-1})$. This is true for
the first term on the right hand side of the above, \emph{but} it is not
necessarily true for the
second term, which for most covariances will remain of order $O(1)$
for all $a$. The problem is
that the denominator is comprised of
$|A_{\lambda}(\frac{k_1}{\lambda}) A_{\lambda}(\frac{k_2}{\lambda})| \sim |k_1k_2|^{-2}$ whereas the numerator is
comprised of $B_{\lambda}(\frac{k}{\lambda}) =
\frac{1}{\lambda}\int_{0}^{\lambda}c(u)\exp(2\pi i k u/\lambda)du$. If
$B_{\lambda}(k/\lambda)\sim |k|^{-(1+\delta)}$ (for some $\delta>0$),
then the second term would converge to zero as $a\rightarrow \infty$. 
However, we see that when $c(u)$ is defined on the torus $[0,\lambda]$, then
for most covariances there will be a discontinuity at zero
since $c(0)\neq c(\lambda)$, this means that $B_{\lambda}(k/\lambda)\sim |k|^{-1}$. 
Consequently, in general $\var[\widehat{\sigma}^{2}] = O(1)$. 
One exception is when $c(u) = c(\lambda-u)$ for all $u \in
[0,\lambda/2]$, in this case $B_{\lambda}(\frac{k}{\lambda})=0$ for
all $k$ and $\var[\widehat{\sigma}^{2}]\approx a^{-1}$. However, such
a covariance is unlikely to arise in most spatial
applications. 

Therefore, in general, it seems that we cannot
consistently estimate $\sigma^{2}$ using a Fourier domain approach.
 We conjecture that the only
transformation of the data that will consistently estimate
$\sigma^{2}$ is a transformation with the eigenfunctions and values associated with the
covariance operator $c$.
In contrast, \citeA{p:zha-04} and \citeA{p:zha-05} showed that if the maximum
likelihood were used to estimate $\sigma^{2}$ in a Gaussian random
field with covariance $\sigma^{2}c(\cdot)$ where $c(\cdot)$ is a known
and Matern covariance function, then even within the fixed domain framework $\sigma^{2}$
can be consistently estimated. This demonstrates that there exists
situations where there are clear gains by working within the likelihood
framework (if the correct distribution is specified). However, if the true covariance is $c(u) =
c(u;\theta)$ and $\theta$ is also unknown, then even within the Gaussian
likelihood framework (with correct specification of the distribution)
one cannot consistently estimate $\sigma^{2}$ and $\theta$.
%(the maximum likelihood estimator is
%equivalent to the estimator generated with the eigenfunctions and
%values of $c$).

\subsection{Uniformly sampled locations}\label{sec:uniform1}

In many situations it is reasonable to suppose that the locations are
uniformly distributed over a region. In this case, many of the
results stated above can be simplified. These simplifications allow us
to develop a simple method for estimating nuisance parameters 
in Section \ref{sec:variance2}.

We first recall that if the locations are uniformly sampled then
$h_{\lambda}(\ub) = \lambda^{-d}I_{[-\lambda/2,\lambda/2]^{d}}(\ub)$,
therefore the Fourier coefficients of $h_{\lambda}(\cdot)$ are simply
$\gamma_{{\bf 0}}=1$ and $\gamma_{\jb}=0$ if $\jb\neq {\bf 0}$.
This implies that $\Gamma_{\rb}$ (defined in Corollary \ref{cor:C}) is such that
$\Gamma_{\rb}=0$ if $\rb\neq {\bf 0}$ and $\Gamma_{{\bf
    0}}=1$. Therefore if $\{Z(\ub);\ub\in \mathbb{R}^{d}\}$ is a stationary
spatial random process that satisfies the assumptions in Theorem
\ref{theorem:nonGaussiannonUniform} (with uniform sampling of the locations) and $\|\rb\|_{1}<<\lambda$ 
\begin{eqnarray}
\lambda^{d}\cov\left[\Re \widetilde{Q}_{a,\lambda}(g;\rb_{1}),\Re \widetilde{Q}_{a,\lambda}(g;\rb_{2})\right]=
\left\{
\begin{array}{cl}
\frac{1}{2}[C_{1}+D_{1}]  + O\left(\ell_{\lambda,a,n}+\frac{(a\lambda)^{2}I_{\textrm{S=NG}}}{n^{2}}+\frac{\|\rb\|_{1}}{\lambda}\right) & \rb_{1}=\rb_{2}(=\rb)\\
O(\ell_{\lambda,a,n}+\frac{(a\lambda)^{2}I_{\textrm{S=NG}}}{n^{2}}) & \rb_{1}\neq\rb_{2} \textrm{ or }-\rb_{2}\\
\end{array}
\right., \nonumber\\
\lambda^{d}\cov\left[\Im \widetilde{Q}_{a,\lambda}(g;\rb_{1}),\Im \widetilde{Q}_{a,\lambda}(g;\rb_{2})\right]=
\left\{
\begin{array}{cl}
\frac{1}{2}[C_{1}+D_{1}]  + O\left(\ell_{\lambda,a,n}+\frac{(a\lambda)^{2}I_{\textrm{S=NG}}}{n^{2}}+\frac{\|\rb\|_{1}}{\lambda}\right) & \rb_{1}=\rb_{2}(=\rb) \\
O(\ell_{\lambda,a,n}+\frac{(a\lambda)^{2}I_{\textrm{S=NG}}}{n^{2}}) & \rb_{1}\neq\rb_{2} \textrm{ or }-\rb_{2}\\
\end{array}
\right. \label{eq:asymptoticUniform}
\end{eqnarray}
and $\lambda^{d}\cov\left[\Re \widetilde{Q}_{a,\lambda}(g;\rb_{1}),\Im
  \widetilde{Q}_{a,\lambda}(g;\rb_{2})\right]=O(\ell_{\lambda,a,n}+\frac{(a\lambda)^{2}I_{\textrm{S=NG}}}{n^{2}})$,
where $I_{\textrm{S=NG}}$ denotes the indicator variable and is one if
the random field is Gaussian, else it is zero (see Section \ref{sec:variance}). From the above we
observe that $\{\Re \widetilde{Q}_{a,\lambda}(g;\rb_{}),\Im
\widetilde{Q}_{a,\lambda}(g;\rb_{});\rb\in \mathbb{Z}^{d}\}$ is a
`near uncorrelated' sequence where in the case $\rb\neq {\bf 0}$,
 $\Re \widetilde{Q}_{a,\lambda}(g;\rb_{})$ and $\Im
\widetilde{Q}_{a,\lambda}(g;\rb_{})$ are consistent estimators of
zero. It is this property we exploit in Section \ref{sec:variance2}
when studentizing $\widetilde{Q}_{a,\lambda}(g;0)$. Below we show asymptotic
normality in the case of uniformly distributed locations. 

\begin{theorem}\label{theorem:CLT}[CLT on real and imaginary parts]
Suppose Assumptions \ref{assum:S}, \ref{assum:uniform}, 
\ref{assum:GG}(b,c) and \ref{assum:G}(i) or \ref{assum:G}(ii)  hold. Let $C_{1}$ and  
$C_{2}$, be defined as in Corollary \ref{cor:C}.
We define the $m$-dimension complex random vectors $\widetilde{\boldsymbol{Q}}_{m} =
(\widetilde{Q}_{a,\lambda}(g,\rb_{1}),\ldots,\widetilde{Q}_{a,\lambda}(g,\rb_{m}))$,
where $\rb_{1},\ldots,\rb_{m}$ are such that $\rb_{i}\neq -\rb_{j}$ and $\rb_{i}\neq
0$. %then for any $i,j\in \{1,\ldots,m\}$ 
Under these conditions we have 
\begin{eqnarray}
\label{eq:cltQtilde}
\frac{2\lambda^{d/2}}{C_{1}}\left(
%\begin{array}{ccc}
\frac{C_{1}}{C_{1}+\Re C_{2}}\Re  \widetilde{Q}_{a,\lambda}(g,0),  
\Re \widetilde{\boldsymbol{Q}}_{m},
\Im \widetilde{\boldsymbol{Q}}_{m}
%\end{array}
\right)\Pcon
 \mathcal{N}\big(0,I_{2m+1}\big)
\end{eqnarray}
with $\frac{\log^{2}(a)}{\lambda^{1/2}}\rightarrow
0$, $\lambda^{d}/n\rightarrow 0$  as $\lambda \rightarrow \infty$, $n\rightarrow \infty$ and
$a\rightarrow \infty$.
\end{theorem}
{\bf PROOF} See Section \ref{sec:cumulants}, \citeA{p:sub-14}.  
\hfill $\Box$

%% file: 5_student.tex
\section{A studentized $\widetilde{Q}_{a,\lambda}(g;0)$-statistic}\label{sec:variance2}

The expression for the variance 
$\widetilde{Q}_{a,\lambda}(g;0)$ given in the examples
above, is rather unwieldy and difficult to estimate directly.
In this section we describe simple method for estimating
the variance of $\widetilde{Q}_{a,\lambda}(g;0)$ under the assumption
the locations are uniformly distributed. This estimator is used to
obtain a simple studentised statistic for $\widetilde{Q}_{a,\lambda}(g;0)$.
Our approach is motivated by the method of
orthogonal samples for time series proposed in \citeA{p:sub-15}, where the
idea is to define a  sample which by construction shares some
of the properties as the estimator of interest. In this section we
show that $\{\widetilde{Q}_{a,\lambda}(g;\rb);\rb\neq 0\}$ is an
orthogonal sample associated with $\widetilde{Q}_{a,\lambda}(g;\rb)$.

%We note that the Gaussianity assumption can be relaxed by
%imposing some additional conditions on the higher order cumulants. On
%the other hand, in order to relax the uniform design assumption we
%have to  use asymptotically far frequencies. 

We first focus on estimating the variance 
of $\widetilde{Q}_{a,\lambda}(g;0)$. Assuming for ease of presentation
that the spatial random field is Gaussian and using 
Theorem \ref{lemma:mean-stat} we have 
$\Ex[\widetilde{Q}_{a,\lambda}(g;0)] = I\left(g;\frac{a}{\lambda}
\right) +o(1)$ and 
\begin{eqnarray*}
\lambda^{d}\var[\widetilde{Q}_{a,\lambda}(g;0)] &=& C_{1}  + O(\ell_{\lambda,a,n}). %+ 2C_{3}G_{\lambda} +
%2f_{2}G_{\lambda}^{2} + O(\ell_{\lambda,a,n}),
%\lambda^{d}\var[\widetilde{Q}_{a,\lambda}(g;0)] &=& C_{1} + O(\ell_{\lambda,a,n})
\end{eqnarray*}
%where to simplify notation we have assumed for all $\ob \in
%\mathbb{R}^{d}$ that $\overline{g(\ob)} = g(-\ob)$, which implies 
%$\Re[C_{3}G_{\lambda}]=C_{3}G_{\lambda}$ and $G_{\lambda} =
%\frac{1}{n}\sum_{\kb = -a}^{a}g(\ob_{\kb})$ is real. 
In contrast, we observe that if no elements of the vector $\rb$ are
zero, then by Theorem \ref{lemma:mean-stat}
$\Ex[\widetilde{Q}_{a,\lambda}(g;\rb)]  = O(\prod_{i=1}^{d}[\log\lambda + \log
|r_{i}|]/\lambda^{d})$ (slightly slower rates are obtained when $\rb$
contains zeros). In other words, $\widetilde{Q}_{a,\lambda}(g;\rb)$ is
estimating zero. On the other hand, by (\ref{eq:asymptoticUniform}) we
observe that 
\begin{eqnarray*}
\lambda^{d}\cov\left[\Re \widetilde{Q}_{a,\lambda}(g;\rb_{1}),\Re \widetilde{Q}_{a,\lambda}(g;\rb_{2})\right]&=&
\left\{
\begin{array}{cl}
\frac{1}{2}C_{1} 
+  O\left(\ell_{\lambda,a,n}+ \frac{\|\rb\|_{1}}{\lambda}\right) &
\rb_{1}=\rb_{2} (=\rb)\\
O\left(\ell_{\lambda,a,n}\right) &
\rb_{1}\neq\rb_{2}, \rb_{1}\neq - \rb_{2}\\
\end{array}
\right..
\end{eqnarray*}
A similar result holds for $\{\Im \widetilde{Q}_{a,\lambda}(g;\rb)\}$,
furthermore we have 
$\lambda^{d}\cov[\Re \widetilde{Q}_{a,\lambda}(g;\rb_{1}),\Im \widetilde{Q}_{a,\lambda}(g;\rb_{2})]=O(\ell_{\lambda,a,n})$. 

In summary, if $\|\rb\|_{1}$ is not too large, then 
$\{\Re \widetilde{Q}_{a,\lambda}(g;\rb), \Im \widetilde{Q}_{a,\lambda}(g;\rb)\}$ are `near
uncorrelated' random variables whose variance is approximately the
same as $\widetilde{Q}_{a,\lambda}(g;0)/\sqrt{2}$. 
This suggests that we can use $\{\Re \widetilde{Q}_{a,\lambda}(g;\rb), \Im
\widetilde{Q}_{a,\lambda}(g;\rb);\rb \in \mathcal{S}\}$ to
estimate $\var[\widetilde{Q}_{a,\lambda}(g;0)]$, where the set
$\mathcal{S}$ is defined as 
\begin{eqnarray}
\label{eq:defS}
\mathcal{S}=\left\{\rb;\|\rb\|_{1}\leq M, \rb_{1}\neq \rb_{2}
  \textrm{ and all elements of $\rb$ are non-zero}\right\}. 
\end{eqnarray}
This leads to the following estimator
\begin{eqnarray}
\label{eq:estV}
\widetilde{V} = \frac{\lambda^{d}}{2|\mathcal{S}|}\sum_{\rb\in \mathcal{S}}^{}\left(2|\Re
  \widetilde{Q}_{a,\lambda}(g;\rb)|^{2} + 2|\Im \widetilde{Q}_{a,\lambda}(g;\rb)|^{2}\right) =
\frac{\lambda^{d}}{|\mathcal{S}|}\sum_{\rb\in \mathcal{S}}^{}|\widetilde{Q}_{a,\lambda}(g;\rb)|^{2},
\end{eqnarray} 
where $|\mathcal{S}|$ denotes the cardinality of the set
$\mathcal{S}$. Note that we specifically select the set $\mathcal{S}$
such that no element $\rb$ contains zero, this is to
ensure that $\Ex[\widetilde{Q}_{a,\lambda}(g;\rb)]$ is small and does not bias
the variance estimator. 

In the following theorem we obtain a mean squared bound for
$\widetilde{V}$.
\begin{theorem}\label{theorem:var-est}
Let $\widetilde{V}$ be defined as
in (\ref{eq:estV}), where
$\mathcal{S}$ is defined in (\ref{eq:defS}). 
Suppose Assumptions \ref{assum:S}, \ref{assum:uniform} and
\ref{assum:GG}(a,b,c) hold and either Assumption \ref{assum:G}(i) or
(ii) holds. Then we have  
\begin{eqnarray*}
\Ex\left(\widetilde{V} -
    \lambda^{d}\var[\widetilde{Q}_{a,\lambda}(g;0)]\right)^{2} =
  O(|\mathcal{S}|^{-1}+|M|\lambda^{-1}+
  \ell_{\lambda,a,n}+\lambda^{-d}\log^{4d}(a))
\end{eqnarray*}
as
$\lambda \rightarrow \infty$, $a\rightarrow \infty$ and $n\rightarrow
\infty$ (where $\ell_{a,\lambda,n}$ is defined in
(\ref{eq:ell-def})).
%\item[(ii)] $\Ex\left(\widehat{V} -
%    \lambda^{d}\var[Q_{a,\lambda}(g;0)]\right)^{2} =
%  O(|\mathcal{S}|^{-1}+|M|\lambda^{-1}+\ell_{\lambda,a,n}+\lambda^{-d}\log^{4d}(a))$ (if $a^{d} = O(n)$).
%\end{itemize}
\end{theorem}
PROOF. See \citeA{p:sub-14}, Section \ref{sec:b3}. \hfill $\Box$

\vspace{3mm}
Thus it follows from the above result that if the set $\mathcal{S}$
grows at a rate such that $|M|\lambda^{-1}\rightarrow 0$ as $\lambda
\rightarrow \infty$, then 
 $\widetilde{V}$ is a mean square consistent
estimators of
$\lambda^{d}\var[\widetilde{Q}_{a,\lambda}(g;0)]$. We use this result
to define an asymptotically pivotal statistic. Let
\begin{eqnarray*}
T_{\mathcal{S}}=\frac{\lambda^{d/2}[\widetilde{Q}_{a,\lambda}(g;0) -
  I(g;\frac{a}{\lambda})]}{\sqrt{\widetilde{V}}}.
\end{eqnarray*}
By using Theorem \ref{theorem:CLT} we can immediately show
that for fixed $\mathcal{S}$,  $T_{\mathcal{S}}\Dcon t_{2|\mathcal{S}|}$ as $\lambda
\rightarrow \infty$. Therefore $T_{\mathcal{S}}$ is asymptotically
pivotal and can be used to construct confidence intervals and test on
the parameter $I(g;\frac{a}{\lambda})$.
%where $t_{2|\mathcal{S}|}$ denotes a $t$-distribution with $2|\mathcal{S}|$
%degrees of freedom and $|\mathcal{S}|$ denotes the cardinality of $\mathcal{S}$. 
%Note that for finite samples the statistic may be right skewed and we
%recommend using  power transform, similar to that defined in
%\citeA{p:che-04} to make it the statistic more normal. 

We note that the same approach and studentisation can be used in the case that the
random field is non-Gaussian and also for the non-bias corrected
statistic $Q_{a,\lambda}(g;0)$.  However, it is a trickier to relax the assumption that
the locations are uniformly distributed. This is because in the case
of a non-uniform design $\Ex[\widetilde{Q}_{a,\lambda}(g;\rb)]$
($\rb \neq 0$) will not, necessarily, be estimating zero. 
%We note that the same method can be used to estimate the variance of
%the non-bias corrected statistic
%$Q_{a,\lambda}(g;0)$.
%However, we know from Theorem \ref{lemma:nonuniformdft}
%that frequencies which are far in distance are asymptotically
%uncorrelated. Therefore, it may be possible to use the set
%$\mathcal{S}$ where the elements $\rb$ are sufficiently far apart to
%estimate $\var[Q_{a,\lambda}(g;\rb)]$, 
%In the case that Assumption
%\ref{assum:nonuniform} holds it is possible one could use the 

%% file: A_1_mean.tex
\subsection{Proof of Theorem \ref{lemma:mean-stat}}\label{sec:mean}

It is clear from the motivation at the start of Section \ref{sec:3mean} that the sinc function plays 
an important role in the analysis of
$Q_{a,\lambda}(g;\rb)$ and 
$\widetilde{Q}_{a,\lambda}(g;\rb)$. Therefore we now summarise some of its
properties. It is well known that
$\frac{1}{\lambda}\int_{-\lambda/2}^{\lambda/2}\exp(ix\omega)dx =
\sinc(\frac{\lambda \omega}{2})$ and 
\begin{eqnarray}
\label{eq:orthogonal}
\int_{-\infty}^{\infty} \sinc(u)du = \pi \textrm{ and } 
\int_{-\infty}^{\infty} \sinc^{2}(u)du = \pi. 
%\textrm{ and }
%\int_{-\infty}^{\infty} \frac{sin^{2}(\lambda u)}{u^{2}}du = \lambda \cdot \pi.
\end{eqnarray}
We next state a well know result that is an important component of the proofs in this paper. 
%The sequence of sinc functions $\{\sinc(x+s\pi);s\in \mathbb{Z}\}$ are orthogonal. We
%will summarise the proof below, noting that similar techniques will be
%used through out the appendix. 
\begin{lemma}\label{lemma:orthogonal}[Orthogonality of the $\sinc$ function]
\begin{eqnarray}
\label{eq:sincA}
\int_{-\infty}^{\infty} \sinc(u) \sinc(u+x)du =  \pi \sinc(x)
\end{eqnarray}
and if $s\in \mathbb{Z}/\{0\}$ then 
\begin{eqnarray}
\label{eq:sincB}
\int_{-\infty}^{\infty} \sinc(u) \sinc(u+s\pi)du = 0.
\end{eqnarray}
\end{lemma}
{\bf PROOF} In Section \ref{appendix:tech-proofs}, \citeA{p:sub-14}. \hfill $\Box$

\vspace{2mm}

We use the above results in the proof below.
\vspace{2mm}

{\bf PROOF of Theorem \ref{lemma:mean-stat}}. 
We first prove (i). By using the same proof used to prove 
\citeA{p:ban-sub-15}, Theorem 1 we can
show that 
\begin{eqnarray*}
\Ex\left[J_{n}(\ob_{\kb_{1}})\overline{J_{n}(\ob_{\kb_{2}})}-\frac{\lambda^{d}}{n}\sum_{j=1}^{n}Z(\ub_{j})^{2}e^{i\ub_{j}^{\prime}(\ob_{\kb_{1}}-\ob_{\kb_{2}})} \right]
 =\left\{
\begin{array}{cc}
f(\ob_{k}) + O(\frac{1}{\lambda}) & \kb_{1} = \kb_{2} (= \kb), \\
O(\frac{1}{\lambda^{d-b}}) & 
  \kb_{1}-\kb_{2}\neq 0.
\end{array}
\right.,
\end{eqnarray*}
where $b=b(\kb_{1}-\kb_{2})$ denotes the number of zero elements in the vector $\rb$.
Therefore, since $a=O(\lambda)$,
by taking expectations of $\widetilde{Q}_{a,\lambda}(g;\rb)$ we  use
the above to give 
\begin{eqnarray*}
\Ex\left[\widetilde{Q}_{a,\lambda}(g;\rb)\right] 
&=&
\left\{
\begin{array}{cc}
\frac{1}{\lambda^{d}}\sum_{\kb=-a}^{a}g(\ob_{\kb})f(\ob_{\kb}) +
O(\frac{1}{\lambda^{}}) & \rb = 0 \\
O(\frac{1}{\lambda^{d-b}})& \rb  \neq 0 \\ 
\end{array}
\right..
\end{eqnarray*} 
Therefore, by replacing the summand with the integral we obtain
(\ref{eq:QaiT}). 
%The proof of (\ref{eq:QaiT}) uses the proof of 
%Theorem \ref{theorem:dft} in \citeA{p:ban-sub-14} and is very similar to the above. 

The above method cannot be used to prove 
(ii) since $a/\lambda\rightarrow \infty$, this leads to bounds 
which may not converge. Therefore, as discussed in Section
\ref{sec:3mean} we consider an alternative approach. To do this we 
expand $\widetilde{Q}_{a,\lambda}(g;\rb)$ as a quadratic form to give 
\begin{eqnarray*}
\widetilde{Q}_{a,\lambda}(g;\rb) 
 &=& \frac{1}{n^{2}}\sum_{\stackrel{j_{1}, j_{2}=1}{j_{1}\neq j_{2}}}^{n}
\sum_{\kb=-a}^{a}g(\ob_{\kb})Z(\ub_{j_{1}})Z(\ub_{j_{2}}) \exp(i\ob_{k}^{\prime}(\ub_{j_{1}}-\ub_{j_{2}}))
\exp(-i\ob_{\rb}^{\prime}\ub_{j_{2}}). 
\end{eqnarray*}
Taking expectations gives
 \begin{eqnarray*}
 \Ex\left[\widetilde{Q}_{a,\lambda}(g;\rb)\right] &=& c_{2}\sum_{\kb=-a}^{a}g(\ob_{\kb})
\Ex\left[c(\ub_{1}-\ub_{2})\exp(i\ob_{k}^{\prime}(\ub_{1}-\ub_{2})-i\ub_{2}^{\prime}\ob_{\rb})\right] 
\end{eqnarray*}
where $c_{2} = n(n-1)/2$.
In the case that $d=1$ the above reduces to 
\begin{eqnarray}
 \Ex\left[\widetilde{Q}_{a,\lambda}(g;r)\right] &=& c_{2}\sum_{k=-a}^{a}g(\omega_{k})
\Ex\left[c(s_{1}-s_{2})\exp(i\omega_{k}(s_{1}-s_{2})-is_{2}\omega_{r})\right] \nonumber\\
 &=& \frac{c_{2}}{\lambda^{2}}\sum_{k=-a}^{a}g(\omega_{k})
\int_{-\lambda/2}^{\lambda/2}\int_{-\lambda/2}^{\lambda/2}c(s_{1}-s_{2})
\exp(i\omega_{k}(s_{1}-s_{2})-is_{2}\omega_{r})ds_{1}ds_{2}. \quad \label{eq:EC}
\end{eqnarray}
Replacing $c(s_{1}-s_{2})$ with the Fourier representation of the covariance function 
gives 
\begin{eqnarray*}
\Ex\left[\widetilde{Q}_{a,\lambda}(g;r)\right] = \frac{c_{2}}{2\pi}
\sum_{k=-a}^{a}g(\omega_{k})
\int_{-\infty}^{\infty}f(\omega)\sinc\left(\frac{\lambda \omega}{2}+k\pi\right)
\sinc\left(\frac{\lambda \omega}{2}+(k+r)\pi\right)d\omega.
\end{eqnarray*}
By a change of variables $y = \lambda \omega/2 + k\pi$ and replacing
the sum with an integral we have 
\begin{eqnarray*}
\Ex\left[\widetilde{Q}_{a,\lambda}(g;r)\right] &=&  \frac{c_{2}}{\pi\lambda}
\sum_{k=-a}^{a}g(\omega_{k})
\int_{-\infty}^{\infty}f(\frac{2y}{\lambda}-\omega_{k})\sinc(y)\sinc(y+r\pi)dy \\ 
&=&  \frac{c_{2}}{\pi}
\int_{-\infty}^{\infty}\sinc(y)\sinc(y+r\pi)
\bigg(\frac{1}{\lambda}\sum_{k=-a}^{a}g(\omega_{k})
f(\frac{2y}{\lambda}-\omega_{k})\bigg)dy \\
&=& \frac{c_{2}}{2\pi^{2}}\int_{-\infty}^{\infty}\sinc(y)\sinc(y+r\pi)\int_{-2\pi
  a/\lambda}^{2\pi a/\lambda}g(u)
f(\frac{2y }{\lambda} - u)du dy  +O\left(\frac{1}{\lambda}\right),
\end{eqnarray*}
where the $O(\lambda^{-1})$ comes from Lemma
\ref{lemma:sum-integral}(ii) in \citeA{p:sub-14}. 
Replacing $f(\frac{2y }{\lambda} -  u)$ with $f(- u)$ gives 
\begin{eqnarray*}
\Ex\left[\widetilde{Q}_{a,\lambda}(g;r)\right] 
&=&
\frac{c_{2}}{2\pi^{2}}\int_{-\infty}^{\infty}\sinc(y)\sinc(y+r\pi)\int_{-2\pi
  a/\lambda}^{2\pi a/\lambda}g(u)
f(- u)du dy + R_{n} +  O(\lambda^{-1}),
\end{eqnarray*}
where
\begin{eqnarray*}
R_{n}& =& 
\frac{c_{2}}{2\pi^{2}}\int_{-\infty}^{\infty}\sinc(y)\sinc(y+r\pi)\bigg(\int_{-2\pi
  a/\lambda}^{2\pi a/\lambda}g(u)\left(
f(\frac{2y }{\lambda} - u) - f(- u)\right)du dy.
\end{eqnarray*}
By using Lemma \ref{lemma:1star} in \citeA{p:sub-14}, we have 
$|R_{n}| = O(\frac{\log\lambda +I(r\neq 0)\log |r|}{\lambda})$. Therefore, by using Lemma
\ref{lemma:orthogonal}, replacing $c_{2}$ with one (which leads to the
error $O(n^{-1})$) and (\ref{eq:orthogonal}) we have 
\begin{eqnarray*}
&& \Ex\left[\widetilde{Q}_{a,\lambda}(g;r)\right] \\
&=&
\frac{c_{2}}{2\pi^{2}}\int_{-\infty}^{\infty}\sinc(y)\sinc(y+r\pi)\int_{-2\pi
  a/\lambda}^{2\pi a/\lambda}g(u)
f(- u)du dy +O\left(
\frac{\log\lambda +I(r\neq 0)\log|r|}{\lambda}\right) \\
 &=& \frac{I(r=0)}{2\pi}\int_{-2\pi a/\lambda}^{2\pi a/\lambda}g(u)
f(u)du dy + 
O\left(\frac{1}{\lambda} + \frac{\log\lambda + I(r\neq 0)\log|r|}{\lambda}
\right),
\end{eqnarray*}
which gives (\ref{eq:QaiiT}) in the case $d=1$. 

To prove the result for $d>1$, we will only consider the case $d=2$,
as it highlights the difference from the $d=1$ case. 
By substituting the spectral representation into (\ref{eq:Eastat}) we have  
\begin{eqnarray*}
&& \Ex\left[\widetilde{Q}_{a,\lambda}(g;(r_{1},r_{2}))\right] \\
&=&  \frac{c_{2}}{\pi^{2}\lambda^{2}}
\sum_{k_{1},k_{2}=-a}^{a}g(\omega_{k_1},\omega_{k_2})
\int_{\mathbb{R}^{2}}f\left(\frac{2u_{1}}{\lambda} -
  \omega_{k_1},\frac{2u_{2}}{\lambda} - \omega_{k_2}\right)\sinc(u_{1})\sinc(u_{1}+r_{1}\pi)\sinc(u_{2})\sinc(u_{2}+r_{2}\pi)du_{1}du_{2}\\
&&.
\end{eqnarray*}
In the case that either $r_{1}=0$ or $r_{2}=0$ we can use the same
proof given in the case $d=1$ to give 
\begin{eqnarray*}
\label{eq:Qalambda}
\Ex[\widetilde{Q}_{a,\lambda}(g;\rb)] 
&=&  \frac{1}{(2\pi)^{2}\pi^{2}\lambda^{2}}
\int_{[-a/\lambda,a/\lambda]^{2}}^{}g(\omega_{1},\omega_{2}) f\left(\omega_{1},\omega_{2}\right) \\
&&\int_{\mathbb{R}^{d}}\sinc(u_{1})\sinc(u_{1}+r_{1}\pi)\sinc(u_{2})\sinc(u_{2}+r_{2}\pi)du_{1}du_{2}
 +R_{n},
\end{eqnarray*}
where $|R_{n}| = O(\frac{\log\lambda +I(\|\rb\|_{1}\neq
  0)\log\|\rb\|_{1}}{\lambda} + n^{-1})$, which gives the desired result. 
However, in the case that both
$r_{1}\neq 0$ and $r_{2}\neq 0$, we can use  Lemma \ref{lemma:1star}, equation
(\ref{eq:lemma1star1d}),  \citeA{p:sub-14}  to obtain 
\begin{eqnarray*}
&&\Ex\left[\widetilde{Q}_{a,\lambda}(g;(r_{1},r_{2}))\right] =  \frac{c_{2}}{\pi^{2}\lambda^{2}}
\sum_{k_{1},k_{2}=-a}^{a}g(\omega_{k_1},\omega_{k_2})\int_{\mathbb{R}^{d}}f\left(\frac{2u_{1}}{\lambda} -
  \omega_{k_1},\frac{2u_{2}}{\lambda} -  \omega_{k_2}\right)\times\\
&&\sinc(u_{1})\sinc(u_{1}+r_{1}\pi)\sinc(u_{2})\sinc(u_{2}+r_{2}\pi)du_{1}du_{2}
= O\left(\frac{\prod_{i=1}^{2}[\log \lambda + \log|r_{i}|]}{\lambda^{2}}\right)
\end{eqnarray*}
thus we obtain the faster rate of convergence.  

It is straightforward to  generalize these arguments to $d>2$.
 \hfill $\Box$

%% file: A_2_variance.tex
\subsection{Variance calculations in the case of uniformly sampled locations}\label{sec:variance}

We start by proving the results in Section \ref{sec:uniform} for the
case that the locations are uniformly distributed. The proof here forms
the building blocks for the case of non-uniform sampling of the locations,
which can found in Section \ref{sec:proof-nonuniform}, \citeA{p:sub-14}.

In Lemma \ref{lemma:covariance}, Appendix \ref{sec:b3},
\citeA{p:sub-14} we show that if the spatial process is Gaussian and
the locations are uniformly sampled then for a bounded frequency grid
we have 
\begin{eqnarray*}
&&\lambda^{d}\cov\left[\widetilde{Q}_{a,\lambda}(g;\rb_{1}),\widetilde{Q}_{a,\lambda}(g;\rb_{2})\right]=    
\left\{
\begin{array}{cl}
C_{1}(\ob_{\rb})
+ O(\frac{1}{\lambda^{}}+\frac{\lambda^{d}}{n}) & \rb_{1} = \rb_{2} (=\rb) \\
O(\frac{1}{\lambda^{}}+\frac{\lambda^{d}}{n}) & \rb_{1} \neq \rb_{2} 
\end{array}
\right.
\end{eqnarray*}
and 
\begin{eqnarray*}
&& \lambda^{d}\cov\left[\widetilde{Q}_{a,\lambda}(g;\rb_{1}),\overline{\widetilde{Q}_{a,\lambda}(g;\rb_{2})}\right]=  
\left\{
\begin{array}{cl}
C_{2}(\ob_{\rb})
+ O(\frac{1}{\lambda^{}}+\frac{\lambda^{d}}{n}) & \rb_{1} = -\rb_{2} (= \rb) \\
O(\frac{1}{\lambda^{}}+\frac{\lambda^{d}}{n}) & \rb_{1} \neq -\rb_{2} 
\end{array}
\right.
\end{eqnarray*}
where
\begin{eqnarray}
\label{eq:Cr}
C_{1}(\ob_{\rb})  = C_{1,1}(\ob_{\rb}) + C_{1,2}(\ob_{\rb})
\textrm{ and } C_{2}(\ob_{\rb})  = C_{2,1}(\ob_{\rb}) + C_{2,2}(\ob_{\rb}),
\end{eqnarray}
with 
\begin{eqnarray}
\label{eq:Cr11}
C_{1,1}(\ob_{\rb}) 
&=&\frac{1}{(2\pi)^{d}}\int_{2\pi[-a/\lambda,a/\lambda]^{d}}f(\ob)f(\ob+\ob_{\rb}) 
|g(\ob)|^{2} d\ob, \nonumber\\
C_{1,2}(\ob_{\rb}) 
&=&\frac{1}{(2\pi)^{d}} \int_{\mathcal{D}_{\rb}}f(\ob)f(\ob+\ob_{\rb}) 
g(\ob)\overline{g(-\ob-\ob_{\rb})} d\ob, \nonumber\\
C_{2,1}(\ob_{\rb}) 
&=&\frac{1}{(2\pi)^{d}}\int_{2\pi[-a/\lambda,a/\lambda]^{d}}f(\ob)f(\ob+\ob_{\rb}) 
g(\ob)g(-\ob) d\ob,\nonumber\\ 
 C_{2,2}(\ob_{\rb}) 
&=&\frac{1}{(2\pi)^{d}}\int_{\mathcal{D}_{\rb}}f(\ob)f(\ob+\ob_{\rb}) 
g(\ob)g(\ob+\ob_{\rb}) d\ob.
\end{eqnarray}
On the other other hand, if the frequency grid is unbounded, then
under  Assumptions \ref{assum:G}(ii) and \ref{assum:GG}(b) we have 
\begin{eqnarray*}
\lambda^{d}\cov\left[\widetilde{Q}_{a,\lambda}(g;\rb_{1}),\widetilde{Q}_{a,\lambda}(g;\rb_{2})\right] =
A_{1}(\rb_{1},\rb_{2})+A_{2}(\rb_{1},\rb_{2}) + O(\frac{\lambda^{d}}{n}),
\end{eqnarray*}
\begin{eqnarray*}
\lambda^{d}\cov\left[\widetilde{Q}_{a,\lambda}(g;\rb_{1}),\overline{\widetilde{Q}_{a,\lambda}(g;\rb_{2})}\right]
= 
A_{3}(\rb_{1},\rb_{2})+A_{4}(\rb_{1},\rb_{2}) + O(\frac{\lambda^{d}}{n}),
\end{eqnarray*}
where 
\begin{eqnarray*}
A_{1}(\rb_{1},\rb_{2}) &=& \frac{1}{\pi^{2d}\lambda^{d}}
\sum_{\mb=-2a}^{2a}\sum_{\kb=\max(-a,-a+\mb)}^{\min(a,a+\mb)} 
\int_{\mathbb{R}^{2d}} f(\frac{2\ubb}{\lambda} -\ob_{\kb_{}})
f(\frac{2\vbb}{\lambda} + \ob_{\kb_{}}+\ob_{\rb_{1}}) \\
&& g(\ob_{\kb})\overline{g(\ob_{\kb}-\ob_{\mb})}
\Sinc(\ubb-\mb \pi)\Sinc(\vbb+(\mb+\rb_{1}-\rb_{2})\pi)
\Sinc(\ubb)\Sinc(\vbb)d\ubb d\vbb
\end{eqnarray*}
and expressions for
$A_{2}(\rb_{1},\rb_{2}),\ldots,A_{4}(\rb_{1},\rb_{4})$ can be found in
Lemma \ref{lemma:covariance}, \citeA{p:sub-14}. Interestingly
$\sup_{a}|A_{1}(\rb_{1},\rb_{2})|<\infty$,
$\sup_{a}|A_{2}(\rb_{1},\rb_{2})|<\infty$,
$\sup_{a}|A_{3}(\rb_{1},\rb_{2})|<\infty$ and
$\sup_{a}|A_{4}(\rb_{1},\rb_{2})|<\infty$. Hence if
$\lambda^{d}/n\rightarrow c$ ($c<\infty$)
we have 
\begin{eqnarray*}
\lambda^{d}\sup_{a}\left|\cov\left[\widetilde{Q}_{a,\lambda}(g;\rb_{1}),\widetilde{Q}_{a,\lambda}(g;\rb_{2})\right]\right|
<\infty \quad\textrm{and}\quad
\lambda^{d}\sup_{a}\left|\cov\left[\widetilde{Q}_{a,\lambda}(g;\rb_{1}),
\overline{\widetilde{Q}_{a,\lambda}(g;\rb_{2})}\right]\right|<\infty,
\end{eqnarray*}

\vspace{3mm}
We now obtain simplified expressions for the terms
$A_{1}(r_{1},r_{2}),\ldots,A_{4}(r_{1},r_{2})$ under the slightly
stronger condition that Assumption \ref{assum:GG}(c) also holds. 
\begin{lemma}\label{lemma:var-asymp}
Suppose Assumptions \ref{assum:G}(ii) and \ref{assum:GG}(b,c)
hold. Then
for $0 \leq |r_{1}|, |r_{2}|\leq C|a|$ (where $C$ is some finite constant)  we have 
\begin{eqnarray*}
A_{1}(r_{1},r_{2})  
&=&
\left\{
\begin{array}{cc}
C_{1,1}(\omega_{r}) + O(\ell_{\lambda,a,n}) & r_{1}=r_{2}(=r)\\
O(\ell_{\lambda,a,n}) & r_{1}\neq r_{2}
\end{array}
\right.
\end{eqnarray*}
\begin{eqnarray*}
A_{2}(r_{1},r_{2}) &=&
\left\{
\begin{array}{cc}
C_{1,2}(\omega_{r}) + O(\ell_{\lambda,a,n}) & r_{1}=r_{2}(=r)\\
O(\ell_{\lambda,a,n}) & r_{1}\neq r_{2}
\end{array}
\right.
\end{eqnarray*}
\begin{eqnarray*}
A_{3}(r_{1},r_{2})&=& 
\left\{
\begin{array}{cc}
C_{2,1}(\omega_{r}) + O(\ell_{\lambda,a,n}) & r_{1}=-r_{2}(=r)\\
O(\ell_{\lambda,a,n}) & r_{1}\neq -r_{2}
\end{array}
\right. 
\end{eqnarray*}
and 
 \begin{eqnarray*}
A_{4}(r_{1},r_{2})&=& 
\left\{
\begin{array}{cc}
C_{2,2}(\omega_{r}) + O(\ell_{\lambda,a,n}) & r_{1}=-r_{2}(=r)\\
O(\ell_{\lambda,a,n}) & r_{1}\neq -r_{2}
\end{array}
\right.,
\end{eqnarray*}
where $C_{1,1}(\omega_{r}),\ldots,C_{2,2}(\omega_{r})$ (using $d=1$) and $\ell_{\lambda,a,n}$ are defined in 
Lemma \ref{lemma:covariance}.
%\begin{eqnarray*}
%\ell_{\lambda,a,n} =  
%O\bigg(\log^{2}(a)\bigg[\frac{\Gamma(a/\lambda)(\log a+\log \lambda)}{\lambda}\bigg]\bigg). 
%\end{eqnarray*}
\end{lemma}
PROOF. We first consider $A_{1}(r_{1},r_{2})$ and write it as 
\begin{eqnarray*}
&&A_{1}(r_{1},r_{2}) = \\
&&\frac{1}{\pi^{2}}\int_{-\infty}^{\infty}\int_{-\infty}^{\infty}
\sum_{m=-2a}^{2a}\sinc(u) \sinc(v)\sinc(u-m\pi)\sinc(v+(m+r_{1}-r_{2})\pi)
H_{m,\lambda}\left(\frac{2u}{\lambda},\frac{2v}{\lambda};r_{1}\right)  
dudv,
\end{eqnarray*}
where 
\begin{eqnarray*}
H_{m,\lambda}\left(\frac{2u}{\lambda},\frac{2v}{\lambda};r_{1}\right) &=& 
\frac{1}{\lambda} 
\sum_{k=\max(-a,-a+m)}^{\min(a,a+m)}f\left(-\frac{2u}{\lambda} + \omega_{k}\right) 
 f\left(\frac{2v}{\lambda} + \omega_{k}+\omega_{r}\right)g(\omega_{k})\overline{g(\omega_{k}-\omega_{m})}
\end{eqnarray*}
noting that $f(\frac{2u}{\lambda}-\omega) = f(\omega -
\frac{2u}{\lambda})$. 

If $f(-\frac{2u}{\lambda} + \omega_{k})$  and 
$ f(\frac{2v}{\lambda} + \omega_{k}+\omega_{r})$ are replaced with 
$f(\omega_{k})$ and $ f(\omega_{k}+\omega_{r})$ respectively, then  we can exploit 
the orthogonality property of the sinc functions. This requires the following series of approximations. 
\begin{itemize}
\item[(i)]
We start by defining a similar version of $A_{1}(r_{1},r_{2})$ but with the
sum replaced with an integral. Let 
\begin{eqnarray*}
&&B_{1} (r_{1}-r_{2};r_{1})\\
&=& \frac{1}{\pi^{2}}\int_{-\infty}^{\infty}\int_{-\infty}^{\infty}
\sum_{m=-2a}^{2a}\sinc(u)\sinc(u-m\pi)
\sinc(v)\sinc(v+(m+r_{1}-r_{2})\pi)H_{m}\left(\frac{2u}{\lambda},\frac{2v}{\lambda};r_{1}\right)dudv,
\end{eqnarray*} 
where 
\begin{eqnarray*}
&&H_{m}\left(\frac{2u}{\lambda},\frac{2v}{\lambda};r_{1}\right)  =\frac{1}{2\pi}\int_{2\pi\max(-a,-a+m)/\lambda}^{2\pi\min(a,a+m)/\lambda}
f(\omega-\frac{2u}{\lambda})f(\frac{2v}{\lambda} + \omega+\omega_{r_{1}})
g(\omega)\overline{g(\omega-\omega_{m})} d\omega.
\end{eqnarray*}
By using Lemma \ref{lemma:A8}, \citeA{p:sub-14} we 
have
\begin{eqnarray*}
|A_{1} (r_{1},r_{2})-B_{1}(r_{1}-r_{2};r_{1})| = O\left(\frac{\log^{2}a}{\lambda}\right).
\end{eqnarray*}
\item[(ii)] Define the quantity 
\begin{eqnarray*}
&&C_{1} (r_{1}-r_{2};r_{1}) \\
&=& \frac{1}{\pi^{2}}\int_{-\infty}^{\infty}\int_{-\infty}^{\infty}
\sum_{m=-2a}^{2a}\sinc(u)\sinc(u-m\pi)
\sinc(v)\sinc(v+(m+r_{1}-r_{2})\pi) H_{m}(0,0;r_{1})dudv.
\end{eqnarray*}
By using Lemma \ref{lemma:A9}, \citeA{p:sub-14}, we can replace 
$B_{1}(r_{1}-r_{2};r_{1})$ with
$C_{1} (r_{1}-r_{2};r_{1})$ to give the replacement error
  \begin{eqnarray*}
|B_{1}(r_{1}-r_{2};r_{1})-C_{1} (r_{1}-r_{2};r_{1})| =O\bigg(\log^{2}(a)\bigg[
\frac{(\log a+\log \lambda)}{\lambda}\bigg]\bigg). 
\end{eqnarray*}
\item[(iii)] Finally, we analyze $C_{1}(r_{1}-r_{2};r_{1})$. Since $H_{m}(0,0;r_{1})$ does not depend on 
$u$ or $v$ we take it out of the integral to give 
\begin{eqnarray*}
&&C_{1} (r_{1}-r_{2};r_{1}) \\
&=& \frac{1}{\pi^{2}}\sum_{m=-2a}^{2a}H_{m}(0,0;r_{1})
\left(\int_{\mathbb{R}}\sinc(u)\sinc(u-m\pi)du\right)\left(\int_{\mathbb{R}}
\sinc(v)\sinc(v+(m+r_{1}-r_{2})\pi)dv\right) \\
&=&  \frac{1}{\pi^{2}}H_{0}(0,0;r_{1})
\left(\int_{\mathbb{R}}\sinc^{2}(u)du\right)\left(\int_{\mathbb{R}}
\sinc(v)\sinc(v+(r_{1}-r_{2})\pi)dv\right),
\end{eqnarray*} 
where the last line of the above is due to orthogonality of
the sinc function (see Lemma \ref{lemma:orthogonal}). If $r_{1}\neq
r_{2}$, then by orthogonality of the sinc function we have 
$C_{1}(r_{1}-r_{2};r_{1})=0$. On the other hand if $r_{1}=r_{2}$ we
have 
\begin{eqnarray*}
C_{1} (0;r) &=& \frac{1}{2\pi}\int_{-2\pi a/\lambda}^{2\pi a/\lambda}
f(\omega)f(\omega+\omega_{r_{}})|g(\omega)|^{2} d\omega = C_{1,1}(\omega_{r}).
\end{eqnarray*}
\end{itemize}
The proof for the remaining terms 
$A_{2}(r_{1},r_{2}),A_{3}(r_{1},r_{2})$ and $A_{4}(r_{1},r_{2})$ is
identical, thus we omit the details. \hfill $\Box$

\vspace{3mm}

\begin{theorem}\label{theorem:asymptotic}
Suppose Assumptions \ref{assum:S}, \ref{assum:uniform},
\ref{assum:G}(ii) and \ref{assum:GG}(b,c) hold. Then
for $0 \leq |r_{1}|, |r_{2}|\leq C|a|$ (where $C$ is some finite
constant)  we have 
\begin{eqnarray*}
\lambda^{d}\cov\left[\widetilde{Q}_{a,\lambda}(g;\rb_{1}),\widetilde{Q}_{a,\lambda}(g;\rb_{2})\right]=
\left\{
\begin{array}{cl}
C_{1}(\ob_{\rb})  + O(\ell_{\lambda,a,n}) & \rb_{1}=\rb_{2}(=\rb) \\
O(\ell_{\lambda,a,n}) & \rb_{1}\neq\rb_{2}\\
\end{array}
\right.
\end{eqnarray*}
\begin{eqnarray*}
\lambda^{d}\cov\left[\widetilde{Q}_{a,\lambda}(g;\rb_{1}),\overline{\widetilde{Q}_{a,\lambda}(g;\rb_{2})}\right]=
\left\{
\begin{array}{cl}
C_{2}(\ob_{\rb}) + O(\ell_{\lambda,a,n}) & \rb_{1}=-\rb_{2}(=\rb) \\
O(\ell_{\lambda,a,n}) & \rb_{1}\neq-\rb_{2}\\
\end{array}
\right.,
\end{eqnarray*}
where $C_{1}(\ob_{\rb})$ and $C_{2}(\ob_{\rb})$ are defined in
(\ref{eq:Cr}) and 
\begin{eqnarray}
\label{eq:ell-def}
\ell_{\lambda,a,n} = 
\log^{2}(a)\bigg[
\frac{\log a+\log \lambda}{\lambda}\bigg] +
\frac{\lambda^{d}}{n}.
\end{eqnarray}
\end{theorem}
{\bf PROOF}. By using Lemmas
\ref{lemma:covariance} and \ref{lemma:var-asymp}
we immediately obtain (in the case $d=1$)
\begin{eqnarray*}
\cov\left[\widetilde{Q}_{a,\lambda}(g;r_{1}),\widetilde{Q}_{a,\lambda}(g;r_{2})\right] =
\left\{
\begin{array}{cc}
C_{1}(\omega_{r}) + O(\ell_{\lambda,a,n}) & r_{1}=r_{2}(=r) \\
O(\ell_{\lambda,a,n}) & r_{1}\neq r_{2}
\end{array}
\right.
\end{eqnarray*}
and 
\begin{eqnarray*}
\cov\left[\widetilde{Q}_{a,\lambda}(g;r_{1}),\overline{\widetilde{Q}_{\lambda,a,n}(g;r_{2})}\right] =
\left\{
\begin{array}{cc}
C_{2}(\omega_{r}) + O(\ell_{\lambda,a,n}) & r_{1}=-r_{2}(=r)\\
O(\ell_{\lambda,a,n}) & r_{1}\neq -r_{2}
\end{array}
\right..
\end{eqnarray*}
This gives the result for $d=1$. To prove the result for $d>1$ we use
the same procedure outlined in the proof of Lemma \ref{lemma:var-asymp}
and the above.  \hfill $\Box$

\vspace{2mm} The above theorem means that the variance for both the
bounded and unbounded frequency grid are equivalent (up to the limits
of an integral). 

\vspace{3mm} Using the above results and the
Lipschitz continuity of $g(\cdot)$ and $f(\cdot)$ we can show that 
\begin{eqnarray*}
C_{j}(\ob_{\rb}) = C_{j} + O\left(\frac{\|\rb\|_{1}}{\lambda}\right),
\end{eqnarray*}
where $C_{1}$ and $C_{2}$ are defined in Corollary \ref{cor:C}.

We now derive an expression for the variance for the non-Gaussian
case. 
\begin{theorem}\label{theorem:variance-nongaussian}
Let us suppose that  $\{Z(\ub);\ub \in\mathbb{R}^{d}\}$ is a fourth order stationary spatial 
random field that satisfies Assumption \ref{assum:S}(i). Suppose
all the assumptions in  Theorem \ref{theorem:nonGaussiannonUniform}
hold with the exception of Assumption \ref{assum:nonuniform} which is
replaced with Assumption \ref{assum:uniform} (i.e. we assume the
locations are uniformly sampled). Then for
$\|\rb\|_{1},\|\rb_{2}\|_{1} <<\lambda$ we have 
\begin{eqnarray}
\label{eq:nonGaussian1}
\lambda^{d}\cov\left[\widetilde{Q}_{a,\lambda}(g;\rb_{1}),\widetilde{Q}_{a,\lambda}(g;\rb_{2})\right]=
\left\{
\begin{array}{cl}
C_{1}  + D_{1}+
O(\ell_{\lambda,a,n}+\frac{(a\lambda)^{d}}{n^{2}} +\frac{\|\rb\|_{1}}{\lambda}) & \rb_{1}=\rb_{2}(=\rb) \\
O(\ell_{\lambda,a,n}+\frac{(a\lambda)^{d}}{n^{2}} ) & \rb_{1}\neq\rb_{2}\\
\end{array}
\right.
\end{eqnarray}
\begin{eqnarray}
\label{eq:nonGaussian2}
\lambda^{d}\cov\left[\widetilde{Q}_{a,\lambda}(g;\rb_{1}),\overline{\widetilde{Q}_{a,\lambda}(g;\rb_{2})}\right]=
\left\{
\begin{array}{cc}
C_{2} + D_{2}+
O(\ell_{\lambda,a,n}^{}+
\frac{(a\lambda)^{d}}{n^{2}}+\frac{\|\rb\|_{1}}{\lambda}) & \rb_{1}=-\rb_{2}(=\rb) \\
O(\ell_{\lambda,a,n}^{}+
\frac{(a\lambda)^{d}}{n^{2}}) & \rb_{1}\neq-\rb_{2}\\
\end{array}
\right.,
\end{eqnarray}
where $C_{1}$ and $C_{2}$ are defined in Corollary
\ref{cor:C} and 
\begin{eqnarray*}
D_{1}&=& \frac{1}{(2\pi)^{2d}}\int_{2\pi[-a/\lambda,a/\lambda]^{2d}}g(\ob_{1})\overline{g(\ob_{2})}
f_{4}(-\ob_{1},-\ob_{2},\ob_{2})d\ob_{1}d\ob_{2} \\
D_{2} &=& \frac{1}{(2\pi)^{2d}}\int_{2\pi[-a/\lambda,a/\lambda]^{2d}}g(\ob_{1})g(\ob_{2})
f_{4}(-\ob_{1},-\ob_{2},\ob_{2})d\ob_{1}d\ob_{2}.
\end{eqnarray*}
\end{theorem}
{\bf PROOF} We prove the result for the 
notationally simple case $d=1$. By using 
indecomposable partitions, conditional cumulants
and (\ref{eq:ahatcov}) in \citeA{p:sub-14} we have 
\begin{eqnarray}
\label{eq:ahatcov2}
\lambda\cov\left[\widetilde{Q}_{a,\lambda}(g;r_{1}),\widetilde{Q}_{a,\lambda}(g;r_{2})\right] 
= A_{1}(r_{1},r_{2}) + A_{2}(r_{1},r_{2}) + B_{1}(r_{1},r_{2}) + B_{2}(r_{1},r_{2})
+ O\left(\frac{\lambda}{n}\right), \quad
\end{eqnarray}
where $A_{1}(r_{1},r_{2})$ and $A_{2}(r_{1},r_{2})$ are defined below
equation (\ref{eq:ahatcov}) and 
\begin{eqnarray*}
B_{1}(r_{1},r_{2}) &=& \lambda c_{4}\sum_{k_{1},k_{2}=-a}^{a}g(\omega_{k_{1}})\overline{g(\omega_{k_{2}})}
\Ex\bigg[\kappa_{4}(s_{2}-s_{1},s_{3}-s_{1},s_{4}-s_{1})e^{is_{1}\omega_{k_{1}}}e^{-is_{2}\omega_{k_{1}+r_{1}}}
e^{-is_{3}\omega_{k_{2}}}e^{is_{4}\omega_{k_{2}+r_{2}}} \bigg] \\
B_{2}(r_{1},r_{2}) &=& \frac{\lambda}{n^{4}}\sum_{j_{1},\ldots,j_{4}\in \mathcal{D}_{3}}\sum_{k_{1},k_{2}=-a}^{a}
g(\omega_{k_{1}})\overline{g(\omega_{k_{2}})}\times \\
 &&\Ex\bigg[\kappa_{4}(s_{j_2}-s_{j_1},s_{j_3}-s_{j_1},s_{j_4}-s_{j_1})e^{is_{j_1}\omega_{k_{1}}}e^{-is_{j_2}\omega_{k_{1}+r_{1}}}
e^{-is_{j_3}\omega_{k_{2}}}e^{is_{j_4}\omega_{k_{2}+r_{2}}} \bigg]
\end{eqnarray*}
with $c_{4}=n(n-1)(n-2)(n-3)/n^{2}$ 
\\*
$\mathcal{D}_{3}=\{j_{1},\ldots,j_{4};j_{1}\neq j_{2}\textrm{ and }
j_{3}\neq j_{4}\textrm{ but some $j$'s are in common}\}$.
The limits of $A_{1}(r_{1},r_{2})$ and $A_{2}(r_{1},r_{2})$
are given in Lemma \ref{lemma:var-asymp}, therefore, all that
remains is to derive bounds for $B_{1}(r_{1},r_{2})$ and 
$B_{2}(r_{1},r_{2})$. We will show that $B_{1}(r_{1},r_{2})$ is the
dominating term, whereas by placing sufficient conditions on the rate of growth of $a$, we will show that
$B_{2}(r_{1},r_{2})\rightarrow 0$. 

In order to analyze $B_{1}(r_{1},r_{2})$ we will use the result
\begin{eqnarray}
\int_{\mathbb{R}^{3}}\sinc(u_{1}+u_{2}+u_{3}+(r_{2}-r_{1})\pi)
\sinc(u_{2})\sinc(u_{3})du_{1}du_{2}du_{3} = 
\left\{
\begin{array}{cc}
\pi^{3} & r_{1}=r_{2} \\
 0 &r_{1}\neq r_{2}.
\end{array}
\right.,\label{eq:sinc3}
\end{eqnarray}
which follows from Lemma \ref{lemma:orthogonal}. In the following
steps we will make a series of approximations which will allow us to
apply (\ref{eq:sinc3}).

We start by substituting  the Fourier
representation of the cumulant function
\begin{eqnarray*}
\kappa_{4}(s_{1}-s_{2},s_{1}-s_{3},s_{1}-s_{4}) = \frac{1}{(2\pi)^{3}}\int_{\mathbb{R}^{3}}
f_{4}(\omega_{1},\omega_{2},\omega_{3})e^{i(s_{1}-s_{2})\omega_{1}}e^{i(s_{1}-s_{3})\omega_{2}}e^{i(s_{1}-s_{4})\omega_{4}}
d\omega_{1}d\omega_{2}d\omega_{3},
\end{eqnarray*}
into $B_{1}(r_{1},r_{2})$ to give 
\begin{eqnarray*}
&&B_{1}(r_{1},r_{2}) \\
&=& \frac{c_{4}}{(2\pi)^{3}\lambda^{3}}\sum_{k_{1},k_{2}=-a}^{a}g(\omega_{k_{1}})\overline{g(\omega_{k_{2}})}
\int_{\mathbb{R}^{3}}
f_{4}(\omega_{1},\omega_{2},\omega_{3})\int_{[-\lambda/2,\lambda/2]^{4}}
 e^{is_{1}(\omega_{1}+\omega_{2}+\omega_{3}+\omega_{k_{1}})} \\
&&e^{-is_{2}(\omega_{1}+\omega_{k_{1}+r_{1}})}e^{-is_{3}(\omega_{2}+\omega_{k_{2}})} 
e^{is_{4}(-\omega_{3}+\omega_{k_{2}+r_{2}})}ds_{1}ds_{2}ds_{3}ds_{4}d\omega_{1}d\omega_{2}d\omega_{3} \\
 &=&   \frac{c_{4}\lambda}{(2\pi)^{3}}\sum_{k_{1},k_{2}=-a}^{a}g(\omega_{k_{1}})\overline{g(\omega_{k_{2}})}
\int_{\mathbb{R}^{3}}
f_{4}(\omega_{1},\omega_{2},\omega_{3})
\sinc\left(\frac{\lambda(\omega_{1}+\omega_{2}+\omega_{3})}{2}+k_{1}\pi\right)
\sinc\left(\frac{\lambda \omega_{1}}{2} + (k_{1}+r_{1})\pi\right)\\
&&\times \sinc\left(\frac{\lambda \omega_{2}}{2}+k_{2}\pi\right)
\sinc\left(\frac{\lambda \omega_{3}}{2}-(k_{2}+r_{2})\pi\right)
d\omega_{1}d\omega_{2}d\omega_{3}. 
\end{eqnarray*}
Now we make a change of variables and let 
$u_{1}=\frac{\lambda \omega_{1}}{2} + (k_{1}+r_{1})$, 
$u_{2} = \frac{\lambda \omega_{2}}{2}+k_{2}\pi$
and $u_{3}=\frac{\lambda \omega_{3}}{2}-(k_{2}+r_{2})\pi$, this gives 
\begin{eqnarray*}
B_{1}(r_{1},r_{2}) 
 &=&   \frac{c_{4}}{\pi^{3}\lambda^{2}}\sum_{k_{1},k_{2}=-a}^{a}\int_{\mathbb{R}^{3}}g(\omega_{k_{1}})\overline{g(\omega_{k_{2}})}
f_{4}\left(\frac{2u_{1}}{\lambda}-\omega_{k_{1}+r_{1}},\frac{2u_{2}}{\lambda}-\omega_{k_2},\frac{2u_{3}}{\lambda}+\omega_{k_2 + r_2}\right)
 \times \\
&&\times\sinc\left(u_{1}+u_{2}+u_{3}+(r_{2}-r_{1})\pi\right)\sinc(u_{1})
\sinc(u_{2})\sinc(u_{3})du_{1}du_{2}du_{3}.
\end{eqnarray*}
Next we exchange the summand with a double integral and use Lemma \ref{lemma:sum-integral}(iii)
together with Lemma \ref{lemma:1a}, equation (\ref{eq:in}) (in \citeA{p:sub-15})
to obtain 
\begin{eqnarray*}
B_{1}(r_{1},r_{2}) 
 &=&
 \frac{c_{4}}{(2\pi)^{2}\pi^{3}}\int_{2\pi[-a/\lambda,a/\lambda]^{2}}\int_{\mathbb{R}^{3}}g(\omega_{1})\overline{g(\omega_{2})}
f_{4}\left(\frac{2u_{1}}{\lambda}-\omega_{1}-\omega_{r_{1}},
\frac{2u_{2}}{\lambda}-\omega_{2},\frac{2u_{3}}{\lambda}+\omega_{2}+\omega_{r_2}\right)
 \times \\
&&\sinc\left(u_{1}+u_{2}+u_{3}+(r_{2}-r_{1})\pi\right)\sinc(u_{1})
\sinc(u_{2})\sinc(u_{3})du_{1}du_{2}du_{3}d\omega_{1}d\omega_{2} + O\left(\frac{1}{\lambda}\right).
\end{eqnarray*}
By using Lemma \ref{lemma:cum4bound},
we  replace 
$f_{4}(\frac{2u_{1}}{\lambda}-\omega_{1}-\omega_{r_{1}},\frac{2u_{2}}{\lambda}-\omega_{2},\frac{2u_{3}}{\lambda}+\omega_{2}+\omega_{r_2})$ in the integral with 
$f_{4}(-\omega_{1}-\omega_{r_{1}},-\omega_{2},\omega_{2} + \omega_{r_2})$, this gives 
\begin{eqnarray*}
&&B_{1}(r_{1},r_{2})  \\
&=&\frac{c_{4}}{(2\pi)^{2}\pi^{3}}\int_{2\pi[-a/\lambda,a/\lambda]^{2}}g(\omega_{1})\overline{g(\omega_{2})}
 f_{4}(-\omega_{1}-\omega_{r_{1}},-\omega_{2},\omega_{2} +\omega_{r_2})\times \\
&&\int_{\mathbb{R}^{3}}\sinc(u_{1}+u_{2}+u_{3}+(r_{2}-r_{1})\pi)\sinc(u_{1})
\sinc(u_{2})\sinc(u_{3})du_{1}du_{2}du_{3}d\omega_{1}d\omega_{2}
+O\left(\frac{\log^{3}\lambda}{\lambda}\right) \\
&=& \frac{c_{4}I_{r_{1}=r_{2}}}{(2\pi)^{2}}\int_{2\pi[-a/\lambda,a/\lambda]^{2}}g(\omega_{1})\overline{g(\omega_{2})}
 f_{4}(-\omega_{1}-\omega_{r_{1}},-\omega_{2},\omega_{2} +\omega_{r_2})d\omega_{1}d\omega_{2}
+O\left(\frac{\log^{3}\lambda}{\lambda}\right),
\end{eqnarray*}
where the last line follows from the orthogonality relation of the
sinc function in equation (\ref{eq:sinc3}). 

Finally we make one further approximation. In the case $r_{1}=r_{2}$ we replace 
$ f_{4}(-\omega_{1}-\omega_{r_{1}},-\omega_{2},\omega_{2}
+\omega_{r_2})$ with $ f_{4}(-\omega_{1},-\omega_{2},\omega_{2})$
which by using the Lipschitz continuity of $f_{4}$ and  Lemma \ref{lemma:1a}, equation (\ref{eq:in}) gives 
\begin{eqnarray*}
B_{1}(r_{1},r_{1})  
 &=& \frac{c_{4}}{(2\pi)^{2}}\int_{2\pi[-a/\lambda,a/\lambda]^{2}}g(\omega_{1})\overline{g(\omega_{2})}
f_{4}(-\omega_{1},-\omega_{2},\omega_{2})d\omega_{1}d\omega_{2} +
O\left(\frac{\log^{3}\lambda}{\lambda}+\frac{|r_{1}|}{\lambda}\right).
\end{eqnarray*}
Altogether this gives 
\begin{eqnarray*}
&&B_{1}(r_{1},r_{2}) \\
&=& \left\{
\begin{array}{cc}
\frac{1}{(2\pi)^{2}}\int_{2\pi[-a/\lambda,a/\lambda]^{2}}g(\omega_{1})\overline{g(\omega_{2})}
f_{4}(-\omega_{1},-\omega_{2},\omega_{2})d\omega_{1}d\omega_{2} + 
O(\frac{\log^{3}\lambda}{\lambda} + 
\frac{|r_{1}|+|r_{2}|}{\lambda}) & r_{1}=r_{2} \\
O(\frac{\log^{3}\lambda}{\lambda}) & r_{1}\neq r_{2}.
\end{array}
\right..
\end{eqnarray*}

Next we show that $B_{2}(r_{1},r_{2})$ is asymptotically negligible.  
To bound $B_{2}(r_{1},r_{2})$ we decompose $\mathcal{D}_{3}$ into six
sets, the set $\mathcal{D}_{3,1}=\{j_{1},\ldots,j_{4};j_{1}=j_{3}$, 
$j_{2}$ and $j_{4}$ are different)$\}$,
$\mathcal{D}_{3,2}=\{j_{1},\ldots,j_{4};j_{1}=j_{4}$, 
$j_{2}$ and $j_{3}$ are different), 
$\mathcal{D}_{3,3}=\{j_{1},\ldots,j_{4};j_{2}=j_{3}$, 
$j_{1}$ and $j_{4}$ are different), 
$\mathcal{D}_{3,4}=\{j_{1},\ldots,j_{4};j_{2}=j_{4}$, 
$j_{1}$ and $j_{3}$ are different), 
$\mathcal{D}_{2,1}=\{j_{1},\ldots,j_{4};j_{1}=j_{3}$ and
$j_{2}=j_{4}\}$,  
$\mathcal{D}_{2,2}=\{j_{1},\ldots,j_{4};j_{1}=j_{4}$ and
$j_{2}=j_{3}\}$.
Using this decomposition we have $B_{2}(r_{1},r_{2})
=\sum_{j=1}^{4}B_{2,(3,j)}(r_{1},r_{2})+\sum_{j=1}^{2}B_{2,(2,j)}(r_{1},r_{2})$, where  
 \begin{eqnarray*}
 B_{2,(3,1)}(r_{1},r_{2}) 
&=& \frac{|\mathcal{D}_{3,1}|\lambda}{n^{4}}\sum_{k_{1},k_{2}=-a}^{a}
g(\omega_{k_{1}})\overline{g(\omega_{k_{2}})}\times \\
 &&\Ex\bigg[\kappa_{4}(s_{j_2}-s_{j_1},0,s_{j_4}-s_{j_1})e^{is_{j_1}\omega_{k_{1}}}e^{-is_{j_2}\omega_{k_{1}+r_{1}}}
e^{-is_{j_1}\omega_{k_{2}}}e^{is_{j_4}\omega_{k_{2}+r_{2}}} \bigg] \\
\end{eqnarray*}
for $j=2,3,4$, $B_{2,(3,j)}(r_{1},r_{2})$  are defined similarly, 
\begin{eqnarray*}
B_{2,(2,1)}(r_{1},r_{2}) &=& \frac{|\mathcal{D}_{2,1}|\lambda}{n^{4}}\sum_{k_{1},k_{2}=-a}^{a}
g(\omega_{k_{1}})\overline{g(\omega_{k_{2}})}\times \\
 &&\Ex\bigg[\kappa_{4}(s_{j_2}-s_{j_1},0,s_{j_2}-s_{j_1})e^{is_{j_1}\omega_{k_{1}}}e^{-is_{j_2}\omega_{k_{1}+r_{1}}}
e^{-is_{j_1}\omega_{k_{2}}}e^{is_{j_2}\omega_{k_{2}+r_{2}}} \bigg], \\ 
\end{eqnarray*}
$B_{2,(2,2)}(r_{1},r_{2})$ is defined similarly and 
 $|\cdot|$ denotes the cardinality of a set.
By using identical methods to those used to bound $B_{1}(r_{1},r_{2})$ we have 
\begin{eqnarray*}
&& |B_{2,(3,1)}(r_{1},r_{2})| \\
&\leq& \frac{C\lambda}{n(2\pi)^{3}}\sum_{k_{1},k_{2}=-a}^{a}
|g(\omega_{k_{1}})\overline{g(\omega_{k_{2}})}| \int_{\mathbb{R}^{3}}
|f_{4}(\omega_{1},\omega_{2},\omega_{3})|
\left|\sinc\left(\frac{\lambda(\omega_{1}+\omega_{3})}{2}+(k_{2}-k_{1})\pi\right)\right|\times \\
&&\left|\sinc\left(\frac{\lambda \omega_{1}}{2} - (k_{1}+r_{1})\pi\right)\right|
\times \left|\sinc\left(\frac{\lambda \omega_{3}}{2}+(k_{2}+r_{2})\pi\right)\right|
d\omega_{1}d\omega_{2}d\omega_{3} = O\left(\frac{\lambda}{n}\right).
\end{eqnarray*}
Similarly we can show
$|B_{2,(3,j)}(r_{1},r_{2})|=O(\frac{\lambda}{n})$ (for $2\leq j \leq 4$)
and 
\begin{eqnarray*}
&&|B_{2,(2,1)}(r_{1},r_{2})| \\
&\leq&  \frac{\lambda}{(2\pi)^{3}n^{2}}\sum_{k_{1},k_{2}=-a}^{a}
|g(\omega_{k_{1}})\overline{g(\omega_{k_{2}})}|\int_{\mathbb{R}^{3}}
|f_{4}(\omega_{1},\omega_{2},\omega_{3})|\bigg|\sinc\left(\frac{\lambda(\omega_{1}+\omega_{2})}{2}
    + (k_{2}-k_{1})\pi\right)\times \\
&&\sinc\left(\frac{\lambda(\omega_{1}+\omega_{2})}{2} + (k_{2}-k_{1}+r_{2}-r_{1})\pi\right)\bigg|
d\omega_{1}d\omega_{2} d\omega_{3}
=O\left(\frac{a\lambda}{n^{2}}\right). 
\end{eqnarray*}
This immediately gives us
(\ref{eq:nonGaussian1}). To prove (\ref{eq:nonGaussian2}) we use 
identical methods. Thus we obtain the result.  \hfill $\Box$

\vspace{2mm}
Note that the term $B_{2,(2,1)}(r_{1},r_{2})$ in the proof above is important as
it does not seem possible to improve on the bound $O(a\lambda/n^{2})$.

%% file: A_3_nonuniform.tex
\section{Proofs in the case of non-uniform sampling}\label{sec:proof-nonuniform}

In this section we prove the results in Section
\ref{sec:uniform}. Most of the results derived here are 
based on the methodology developed in the uniform
sampling case. 

\vspace{3mm}
{\bf PROOF of Theorem \ref{lemma:nonuniformdft}}
We prove the result for $d=1$.

First we prove (\ref{eq:dft-nonuniform}). By expanding
$\cov\left[J_{n}(\omega_{k_{1}}),J_{n}(\omega_{k_{2}})\right]$ (and
assuming that $\Ex[Z(s)|s]=0$) we have 
\begin{eqnarray*}
&&\cov\left[J_{n}(\omega_{k_{1}}),J_{n}(\omega_{k_{2}})\right] \\
&= &
c_{2}\Ex\left[Z(s_{1})Z(s_{2})\exp(is_{1}\omega_{k_1}-is_{2}\omega_{2})\right]
+ n^{-1}\Ex\left[Z(s)^{2}\exp(is(\omega_{k_1}-\omega_{k_{2}}))\right] \\
&=&
\frac{c_{2}}{\lambda^{}}\int_{[-\lambda/2,\lambda/2]^{2}}c(s_{1}-s_{2})e^{is(\omega_{k_1}-\omega_{k_{2}}}
h\left(\frac{s_1}{\lambda}\right)
h\left(\frac{s_2}{\lambda}\right)ds_1ds_2 
 + \frac{c(0)}{n}\int_{-\lambda/2}^{\lambda/2}h\left(\frac{s}{\lambda}\right)e^{is(\omega_{k_1}-\omega_{k_2})}ds.
\end{eqnarray*}
Replacing $h(s/\lambda)$ with its Fourier representation $h(s/\lambda)
= \sum_{j=-\infty}^{\infty}\gamma_{j}e^{2\pi js/\lambda}$ 
\begin{eqnarray*}
&&\cov\left[J_{n}(\omega_{k_{1}}),J_{n}(\omega_{k_{2}})\right] \\
&=& \frac{c_{2}}{\lambda}\sum_{j_{1},j_{2}=-\infty}^{\infty}\gamma_{j_{1}}\gamma_{j_{2}}
\int_{-\lambda/2}^{\lambda/2}\int_{-\lambda/2}^{\lambda/2}
c(s_{1}-s_{2})e^{is_{1}\omega_{j_{1}}}e^{is_{2}\omega_{j_{2}}}e^{is_{1}\omega_{k_{1}}-is_{2}\omega_{k_{2}}}ds_{1}ds_{2} \\
&& +\frac{ c(0)\gamma_{k_{2}-k_{1}}\lambda}{n}.
\end{eqnarray*}
Changing variables in the integral (using $t=s_{1}-s_2$ and $s=s_2$),
then the limits of the integral gives 
\begin{eqnarray*}
\cov\left[J_{n}(\omega_{k_{1}}),J_{n}(\omega_{k_{2}})\right] &=& 
 \frac{c_{2}}{\lambda}\sum_{j_{1},j_{2}=-\infty}^{\infty}\gamma_{j_{1}}\gamma_{j_{2}}
\left(\int_{-\lambda/2}^{\lambda/2}e^{-is(\omega_{k_{2}}-\omega_{j_{2}}-\omega_{j_{1}}-\omega_{k_{1}})}ds\right)
\left(\int_{-\lambda/2}^{\lambda/2}c(t)e^{it(\omega_{j_{1}}+\omega_{k_{1}})}dt\right)
\\
&& +\frac{ c(0)\gamma_{k_{2}-k_{1}}\lambda}{n} + O\left(\frac{1}{\lambda}\right),
\end{eqnarray*}
where to obtain the remainder $O(\frac{1}{\lambda})$ we use that
$\sum_{j=-\infty}^{\infty}|\gamma_{j}|<\infty$ (note that a similar
change in the limits of an integral is given in the proof of Theorem 2.1, 
\citeA{p:ban-sub-15}).
Next, by using the identity 
\begin{eqnarray*}
\int_{-\lambda/2}^{\lambda/2}e^{-is(\omega_{k_{2}}-\omega_{j_{2}}-\omega_{j_{1}}-\omega_{k_{1}})}ds
= 
\left\{
\begin{array}{cc}
 0 & k_{2}-j_{2}\neq k_{1}+j_{1} \\
\lambda & k_{2}-j_{2}=k_{1}+j_{1}\\
\end{array}
\right.
\end{eqnarray*}
we have 
\begin{eqnarray*}
\cov\left[J_{n}(\omega_{k_{1}}),J_{n}(\omega_{k_{2}})\right] &=& 
 c_{2}\sum_{j=-\infty}^{\infty}\gamma_{j}\gamma_{k_{2}-k_{1}-j}
\int_{-\lambda/2}^{\lambda/2}c(t)e^{it(\omega_{j_{1}}+\omega_{k_{1}})}dtds
 +\frac{ c(0)\gamma_{k_{1}-k_{2}}\lambda}{n} + O(\lambda^{-1}) \\
 &=&
 \sum_{j=-\infty}^{\infty}\gamma_{j}\gamma_{k_{2}-k_{1}-j}f(\omega_{k_{1}+j}) +\frac{ c(0)\gamma_{k_{1}-k_{2}}\lambda}{n}
 + O\left(\frac{1}{\lambda} + \frac{1}{n}\right).
\end{eqnarray*}
Finally, we replace $f(\omega_{k_{1}+j})$ with $f(\omega_{k_{1}})$ and
use the Lipschitz continiuity of $f(\cdot)$ to give 
\begin{eqnarray*}
\left| \sum_{j=-\infty}^{\infty}\gamma_{j}\gamma_{k_{2}-k_{1}-j}\left[f(\omega_{k_{1}})-f(\omega_{k_{1}+j})\right]\right|
 \leq
 \frac{C}{\lambda}\sum_{j=-\infty}^{\infty}|j|\cdot|\gamma_{j}\gamma_{k_{2}-k_{1}-j}|
 = O\left(\frac{1}{\lambda}\right),
\end{eqnarray*}
where the last line follows from $|\gamma_{j}|\leq
C|j|^{-(1+\delta)}I(j\neq 0)$. Altogether, this gives 
\begin{eqnarray*}
\cov\left[J_{n}(\omega_{k_{1}}),J_{n}(\omega_{k_{2}})\right] =f(\omega_{k_{1}})
 \sum_{j=-\infty}^{\infty}\gamma_{j}\gamma_{k_{2}-k_{1}-j} +\frac{ c(0)\gamma_{k_{2}-k_{1}}\lambda}{n}
 + O\left(\frac{1}{\lambda}\right).
\end{eqnarray*}
This completes the proof for $d=1$, the proof for $d>1$ is the
same. 

To prove (\ref{eq:dft-nonuniformOmega}) we use that
\begin{eqnarray*}
&&\cov\left[J_{n}(\omega_{\Omega,k_{1}}),J_{n}(\omega_{\Omega,k_{2}})\right] \\
&=& \frac{c_{2}}{\lambda}\sum_{j_{1},j_{2}=-\infty}^{\infty}\gamma_{j_{1}}\gamma_{j_{2}}
\int_{-\lambda/2}^{\lambda/2}\int_{-\lambda/2}^{\lambda/2}
c(s_{1}-s_{2})e^{is_{1}\omega_{j_{1}}}e^{is_{2}\omega_{j_{2}}}e^{is_{1}\omega_{\Omega,k_{1}}-is_{2}\omega_{\Omega,k_{2}}}ds_{1}ds_{2} \\
&&
+\frac{c(0)}{n}\int_{-\lambda/2}^{\lambda/2}h\left(\frac{s}{\lambda}\right)e^{is(\omega_{\Omega,k_1-k_2})}ds
\\
&=&
\frac{\lambda c_{2}}{(2\pi)^{}}\sum_{j_{1},j_{2}=-\infty}^{\infty}\gamma_{j_{1}}\gamma_{j_{2}}\int_{\mathbb{R}}f(\omega)
\frac{1}{\lambda^{2}}\int_{-\lambda/2}^{\lambda/2}\int_{-\lambda/2}^{\lambda/2}
e^{is_{1}(\omega+\omega_{j_{1}}+\omega_{\Omega,k_1})}e^{is_{2}(\omega_{j_{2}}-\omega-\omega_{\Omega,k_{2}})}
ds_{1}ds_{2} + O\left(\frac{\lambda}{n}\right) \\
&=& \frac{\lambda
c_{2}}{(2\pi)^{}}\sum_{j_{1},j_{2}=-\infty}^{\infty}\gamma_{j_{1}}\gamma_{j_{2}}\int_{\mathbb{R}}f(\omega)
\sinc\left(\frac{\lambda}{2}(\omega+\omega_{j_{1}}+\omega_{\Omega,k_1})\right)
\sinc\left(\frac{\lambda}{2}(\omega-\omega_{j_{2}}+\omega_{\Omega,k_2})\right)d\omega
+ O\left(\frac{\lambda}{n}\right) \\
&=& \frac{1}{\pi}\sum_{j_{1},j_{2}=-\infty}^{\infty}\gamma_{j_{1}}\gamma_{j_{2}}\int_{\mathbb{R}}f\left(\frac{2u}{\lambda}-
\omega_{\Omega,k_{1}}-\omega_{j_{1}}\right)
\sinc\left(u\right)
\sinc\left(u-(j_1 +j_2)\pi+\frac{\lambda}{\Omega}(k_2 - k_1)\pi\right)du
+ O\left(\frac{\lambda}{n}\right) \\
&=&\frac{1}{\pi}\sum_{j_{1},j_{2}=-\infty}^{\infty}\gamma_{j_{1}}\gamma_{j_{2}}f\left(
\omega_{\Omega,k_{1}}+\omega_{j_{1}}\right)\int_{\mathbb{R}}
\sinc\left(u\right)
\sinc\left(u-(j_1 +j_2)\pi+\frac{\lambda}{\Omega}(k_2 - k_1)\pi\right)du
+ O\left(\frac{\lambda}{n}+\frac{\log \lambda}{\lambda}\right) \\
&=&\sum_{j_{1},j_{2}=-\infty}^{\infty}\gamma_{j_{1}}\gamma_{j_{2}}f\left(
\omega_{\Omega,k_{1}}+\omega_{j_{1}}\right)
\sinc\left((j_1 +j_2)\pi+\frac{\lambda}{\Omega}(k_1 - k_2)\pi\right)
+ O\left(\frac{\lambda}{n}+\frac{\log \lambda}{\lambda}\right).
\end{eqnarray*}
Finally, under Assumption \ref{assum:GG}(c) (which gives the Lipschitz
continuity of $f$) we replace $f(\omega_{\Omega,k_1}+\omega_{j_{1}})$
with $f(\omega_{\Omega,k_1})$ to give
\begin{eqnarray*}
\cov\left[J_{n}(\omega_{\Omega,k_{1}}),J_{n}(\omega_{\Omega,k_{2}})\right] =
f\left(\omega_{\Omega,k_{1}}\right)
\sum_{j_{1},j_{2}=-\infty}^{\infty}\gamma_{j_{1}}\gamma_{j_{2}}
\sinc\left((j_1 +j_2)\pi+\frac{\lambda}{\Omega}(k_1 - k_2)\pi\right)
+ O\left(\frac{\lambda}{n}+\frac{\log \lambda}{\lambda}\right),
\end{eqnarray*}
thus giving the required result. 
\hfill $\Box$

\vspace{3mm}

%The proof of Theorem \ref{theorem:nonuniformmean} uses the
%method given in the proof of Theorem \ref{lemma:mean-stat} (see
%Appendix \ref{}). 
%\vspace{3mm}

{\bf PROOF of Theorem \ref{theorem:nonuniformmean}} 
We prove the result for the case $d=1$. Using the same
method used to prove Theorem \ref{lemma:mean-stat} (see the arguments
in Section \ref{sec:3mean}) we obtain 
\begin{eqnarray*}
&&\Ex\left[\widetilde{Q}_{a,\lambda}(g;r)\right] \\
&=& \frac{c_{2}}{2\pi}\sum_{j_{1},j_{2}=-\infty}^{\infty}\gamma_{j_{1}}\gamma_{j_{2}}
\sum_{k=-a}^{a}g(\omega_{k})
\int_{-\infty}^{\infty}f(\omega)\sinc\left(\frac{\lambda \omega}{2}+(k+j_{1})\pi\right)
\sinc\left(\frac{\lambda \omega}{2}+(k+r-j_{2})\pi\right)d\omega. 
\end{eqnarray*}
By the change of variables $y=\frac{\lambda \omega}{2}+(k+j_{1})\pi$
we obtain 
\begin{eqnarray*}
&&\Ex\left[\widetilde{Q}_{a,\lambda}(g;r)\right] \\
&=& \frac{c_{2}}{\lambda\pi}\sum_{j_{1},j_{2}=-\infty}^{\infty}\gamma_{j_{1}}\gamma_{j_{2}}
\sum_{k=-a}^{a}g(\omega_{k})
\int_{-\infty}^{\infty}f\left(\frac{2y}{\lambda}-\omega_{k+j_1}\right)\sinc(y)\sinc(y+(r-j_{1}-j_{2})\pi)dy.
\end{eqnarray*}
Replacing sum with an integral and 
using Lemma \ref{lemma:sum-integral}(ii) gives 
\begin{eqnarray*}
&&\Ex\left[\widetilde{Q}_{a,\lambda}(g;r)\right] \\
&=& \frac{c_{2}}{2\pi^{2}}\sum_{j_{1},j_{2}=-\infty}^{\infty}\gamma_{j_{1}}\gamma_{j_{2}}
\int_{-2\pi a/\lambda}^{2\pi a/\lambda}g(\omega)
f\left(\frac{2y}{\lambda}-\omega_{} - \omega_{j_1}\right)d\omega
\int_{-\infty}^{\infty}\sinc(y)\sinc(y+(r-j_{1}-j_{2})\pi)dy +  O\left(\frac{1}{\lambda}\right).
\end{eqnarray*}
Next, replacing $f(\frac{2y}{\lambda}-\omega)$ with $f(-\omega-\omega_{j_1})$ and
we have
\begin{eqnarray*}
&&\Ex\left[\widetilde{Q}_{a,\lambda}(g;r)\right] \\
&=& \frac{c_{2}}{2\pi^{2}}\sum_{j_{1},j_{2}=-\infty}^{\infty}\gamma_{j_{1}}\gamma_{j_{2}}
\int_{-2\pi a/\lambda}^{2\pi a/\lambda}g(\omega) f\left(\frac{2y}{\lambda}-\omega_{}-\omega_{j_1}\right)d\omega
\int_{-\infty}^{\infty}\sinc(y)\sinc(y+(r-j_{1}-j_{2})\pi)dy +  R_{n} + O\left(\frac{1}{\lambda}\right).
\end{eqnarray*}
where 
\begin{eqnarray*}
R_{n} = \frac{c_{2}}{2\pi^{2}}\sum_{j_{1},j_{2}=-\infty}^{\infty}\gamma_{j_{1}}\gamma_{j_{2}}
\int_{-2\pi a/\lambda}^{2\pi a/\lambda}g(\omega)
\left[f\left(\frac{2y}{\lambda}-\omega_{}-\omega_{j_{1}}\right) - f(-\omega-\omega_{j_{1}})\right]
\int_{-\infty}^{\infty}\sinc(y)\sinc(y+(r-j_{1}-j_{2})\pi)dyd\omega.
\end{eqnarray*}
By using Lemma \ref{lemma:1star} we have 
\begin{eqnarray*}
|R_{n}| &\leq& C\sum_{j_{1},j_{2}=-\infty}^{\infty}|\gamma_{j_{1}}|\cdot|\gamma_{j_{2}}|
\frac{\log\lambda + \log |r-j_{1}-j_{2}|}{\lambda} \leq C\sum_{j_{1},j_{2}=-\infty}^{\infty}|\gamma_{j_{1}}|\cdot|\gamma_{j_{2}}|
\frac{\log\lambda + I(r\neq 0)\log|r|+\log|j_{1}|+\log|j_{2}|}{\lambda},
\end{eqnarray*}
noting that $C$ is a generic constant that changes between
inequalities. 
This gives 
\begin{eqnarray*}
&&\Ex\left[\widetilde{Q}_{a,\lambda}(g;r)\right] \\
&=& \frac{c_{2}}{2\pi^{2}}\sum_{j_{1},j_{2}=-\infty}^{\infty}\gamma_{j_{1}}\gamma_{j_{2}}
\int_{-2\pi a/\lambda}^{2\pi a/\lambda}g(\omega) f(-\omega-\omega_{j_1})d\omega
\int_{-\infty}^{\infty}\sinc(y)\sinc(y+(r-j_{1}-j_{2})\pi)dy + 
 O\left(\frac{\log \lambda + I(r \neq 0)\log |r|}{\lambda}\right).
\end{eqnarray*}
Finally, by the orthogonality of the sinc function at integer shifts (and $f(-\omega)=f(\omega)$) we
have
\begin{eqnarray*}
\Ex\left[\widetilde{Q}_{a,\lambda}(g;r)\right] 
&=& \frac{1}{2\pi}\sum_{j=-\infty}^{\infty}\gamma_{j}\gamma_{r-j}
\int_{-2\pi a/\lambda}^{2\pi a/\lambda}g(\omega) f(\omega+\omega_{j})d\omega + 
 O\left(\frac{\log\lambda + \log |r|}{\lambda}+\frac{1}{n}\right) \\
 &=& \frac{1}{2\pi}\sum_{j=-\infty}^{\infty}\gamma_{j}\gamma_{r-j}
\int_{-2\pi a/\lambda}^{2\pi a/\lambda}g(\omega) f(\omega)d\omega + 
 O\left(\frac{\log\lambda + I(r\neq 0)\log |r|}{\lambda}+\frac{1}{n}\right) 
\end{eqnarray*}
thus we obtain the desired result. \hfill $\Box$
\vspace{3mm}

In order to prove Theorem \ref{theorem:nonuniformvar}(i), we define
the quantities 
$U_{1}(\cdot)$ and $U_{2}(\cdot)$ are defined as 
\begin{eqnarray}
\label{eq:Ur}
U_{1}(\rb_{1},\rb_{2};\ob_{\rb_{1}},\ob_{\rb_{2}}) &=&
U_{1,1}(\rb_{1},\rb_{2};\ob_{\rb_{1}},\ob_{\rb_{2}}) + 
U_{1,2}(\rb_{1},\rb_{2};\ob_{\rb_{1}},\ob_{\rb_{2}}) \nonumber\\
U_{2}(\rb_{1},\rb_{2};\ob_{\rb_{1}},\ob_{\rb_{2}}) &=&
U_{2,1}(\rb_{1},\rb_{2};\ob_{\rb_{1}},\ob_{\rb_{2}}) + 
U_{2,2}(\rb_{1},\rb_{2};\ob_{\rb_{1}},\ob_{\rb_{2}}) 
\end{eqnarray}
with 
\begin{eqnarray*}
%\label{eq:Ur11}
&&U_{1,1}(\rb_{1},\rb_{2};\ob_{\rb_{1}},\ob_{\rb_{2}}) \\
&=& 
\frac{1}{(2\pi)^{d}}\sum_{\jb_{1}+\ldots+\jb_{4}=\rb_{1}-\rb_{2}}\gamma_{\jb_{1}}\gamma_{\jb_{2}}\gamma_{\jb_{3}}\gamma_{\jb_{4}}
\int_{\mathcal{D}_{(\jb_{1}+\jb_{3})}}g(\ob)\overline{g(\ob+\ob_{\jb_{1}+\jb_{3}})} 
f(\ob+\ob_{\jb_{1}})f(\ob+\ob_{\rb_{1}-\jb_{2}})d\ob \\
&&U_{1,2}(\rb_{1},\rb_{2};\ob_{\rb_{1}},\ob_{\rb_{2}}) \\
&=& 
\frac{1}{(2\pi)^{d}}\sum_{\jb_{1}+\jb_{2}+\jb_{3}+\jb_{4} =
  \rb_{1}-\rb_{2}}\gamma_{\jb_{1}}\gamma_{\jb_{2}}\gamma_{\jb_{3}}\gamma_{\jb_{4}}
\int_{\mathcal{D}_{(\rb_{1}-\jb_{2}-\jb_{3})}}g(\ob)\overline{g(-\ob-\ob_{\rb_{2}-\jb_{3}-\jb_{2}})}
f(\ob+\ob_{\jb_{1}})f(\ob+\ob_{\rb_{1}-\jb_{2}})d\ob\\
&&U_{2,1}(\rb_{1},\rb_{2};\ob_{\rb_{1}},\ob_{\rb_{2}}) \\
&=&
\frac{1}{(2\pi)^{d}}\sum_{\jb_{1}+\jb_{2}+\jb_{3}+\jb_{4}=\rb_{1}+\rb_{2}}\gamma_{\jb_{1}}\gamma_{\jb_{2}}\gamma_{\jb_{3}}\gamma_{\jb_{4}}
\int_{\mathcal{D}_{(\jb_{1}+\jb_{3})}}g(\ob)g(\ob_{-\jb_{1}-\jb_{3}}-\ob)  f(\ob+\ob_{\jb_{1}})
f(\ob+\ob_{\rb_{1}-\jb_{2}})d\ob \\
&& U_{2,2}(\rb_{1},\rb_{2};\ob_{\rb_{1}},\ob_{\rb_{2}}) \\
&=&\frac{1}{(2\pi)^{d}}\sum_{\jb_{1}+\jb_{2}+\jb_{3}+\jb_{4}=\rb_{1}+\rb_{2} }
\gamma_{\jb_{1}}\gamma_{\jb_{2}}\gamma_{\jb_{3}}\gamma_{\jb_{4}}
\int_{\mathcal{D}_{(\rb_{1}-\jb_{2}-\jb_{3})}}g(\ob)g(\ob-\ob_{\jb_{2}+\jb_{3}-\rb_{1}})
f(\ob+\ob_{\jb_{1}})f(\ob+\ob_{\rb_{1}-\jb_{2}})d\ob,
\end{eqnarray*}
where the integral is defined as 
\\*
$\int_{\mathcal{D}_{\rb}} = 
\int_{2\pi\max(-a,-a-r_{1})/\lambda}^{2\pi\min(a,a-r_{1})/\lambda}\ldots
\int_{2\pi\max(-a,-a-r_{d})/\lambda}^{2\pi\min(a,a-r_{d})/\lambda}$.

\vspace{3mm}
{\bf PROOF of Theorem \ref{theorem:nonuniformvar}(i)} To prove (i)
we use Theorem \ref{lemma:nonuniformdft}
and  Lemma \ref{lemma:cumulantsA} which immediately gives the
result. 

\vspace{3mm}

{\bf PROOF of Theorem \ref{theorem:nonuniformvar}(ii)} 
We first note that
by using Lemma \ref{lemma:cumulantsA} (generalized to non-uniform sampling),  we can show that 
\begin{eqnarray}
\label{eq:Qgamma}
\lambda\cov\left[\widetilde{Q}_{a,\lambda}(g;r_{1}),\widetilde{Q}_{a,\lambda}(g;r_{2})\right]
&=& A_{1}(r_{1},r_{2}) + A_{2}(r_{1},r_{2}) + O\left(\frac{\lambda}{n}\right) 
\end{eqnarray}
where
\begin{eqnarray*}
A_{1}(r_{1},r_{2})&=& \lambda\sum_{k_{1},k_{2}=-a}^{a}g(\omega_{k_{1}})\overline{g(\omega_{k_{2}})}
\cov\big[Z(s_{1})\exp(is_{1}\omega_{k_{1}}),Z(s_{3})\exp(is_{3}\omega_{k_{2}})\big]
\times \\
&&\cov\big[Z(s_{2})\exp(-is_{2}\omega_{k_{1}+r_{1}}),Z(s_{4})\exp(-is_{4}\omega_{k_{2}+r_{2}})\big]\\
A_{2}(r_{1},r_{2})&=& \lambda\sum_{k_{1},k_{2}=-a}^{a}g(\omega_{k_{1}})\overline{g(\omega_{k_{2}})}
 \cov\big[Z(s_{1})\exp(is_{1}\omega_{k_{1}}),Z(s_{4})\exp(-is_{4}\omega_{k_{2}+r_{2}})\big] \times\\
&&\cov\big[Z(s_{2})\exp(-is_{2}\omega_{k_{1}+r_{1}}),Z(s_{3})\exp(is_{3}\omega_{k_{2}})\big].
\end{eqnarray*}
We first analyze $A_{1}(r_{1},r_{2})$. Conditioning on the locations $s_{1},\ldots,s_{4}$ gives 
\begin{eqnarray*}
&&A_{1}(r_{1},r_{2}) \\
&=&\sum_{j_{1},\ldots,j_{4}=-\infty}^{\infty}\gamma_{j_{1}}\gamma_{j_{2}}\gamma_{j_{3}}\gamma_{j_{4}}
\sum_{k_{1},k_{2}=-a}^{a}g(\omega_{k_{1}})\overline{g(\omega_{k_{2}})}
\frac{1}{\lambda^{3}}\int_{[-\lambda/2,\lambda/2]^{4}}
 c(s_{1}-s_{3})c(s_{2}-s_{4})\\
&& e^{is_{1}\omega_{k_{1}}-is_{3}\omega_{k_{2}}}
e^{-is_{2}\omega_{k_{1}+r_{1}}+is_{4}\omega_{k_{2}+r_{2}}}e^{i(s_{1}\omega_{j_{1}}+s_{2}\omega_{j_{2}}+s_{3}\omega_{j_{3}}+s_{4}\omega_{j_{4}})}
ds_{1}ds_{2}ds_{3}ds_{4}. 
\end{eqnarray*}
By using the spectral representation theorem and integrating out $s_{1},\ldots,s_{4}$
we can write the above as
\begin{eqnarray}
&&A_{1} (r_{1},r_{2}) \nonumber\\
&=& \frac{\lambda}{(2\pi)^{2}}\sum_{j_{1},\ldots,j_{4}=-\infty}^{\infty}\gamma_{j_{1}}\gamma_{j_{2}}\gamma_{j_{3}}\gamma_{j_{4}}
\sum_{k_{1},k_{2}=-a}^{a}g(\omega_{k_{1}})\overline{g(\omega_{k_{2}})} 
\int_{-\infty}^{\infty} \int_{-\infty}^{\infty} f(x)f(y)
\sinc\left(\frac{\lambda x}{2} + (k_{1}+j_{1})\pi\right)\nonumber\\
&&\sinc\left(\frac{\lambda y}{2} - (r_{1}+k_{1}-j_{2})\pi\right)
 \sinc\left(\frac{\lambda x}{2}+(k_{2}-j_{3})\pi\right)
\sinc\left(\frac{\lambda y}{2} - (r_{2}+k_{2}+j_{4})\pi\right)dxdy \nonumber\\
&=& \frac{1}{\pi^{2}\lambda}\sum_{j_{1},\ldots,j_{4}=-\infty}^{\infty}\gamma_{j_{1}}\gamma_{j_{2}}\gamma_{j_{3}}\gamma_{j_{4}}
\sum_{k_{1},k_{2}=-a}^{a}g(\omega_{k_{1}})\overline{g(\omega_{k_{2}})}  \int_{-\infty}^{\infty}
\int_{-\infty}^{\infty} 
f(\frac{2u}{\lambda} - \omega_{k_{1}+j_{1}})
f(\frac{2v}{\lambda} + \omega_{k_{1}+r_{1}-j_{2}})
\nonumber\\
&&\times\sinc(u)
\sinc(u + (k_{2}-k_{1}-j_{1}-j_{3})\pi)
\sinc(v)\sinc(v - (r_{2}+k_{2}+j_{4}-r_{1}-k_{1}+j_{2})\pi)dudv. \qquad\label{eq:A1r1r2gamma}
\end{eqnarray}
By making a change of variables $m=k_{1}-k_{2}$ we have 
\begin{eqnarray*}
&&A_{1} (r_{1},r_{2})\\
& = & \frac{1}{\pi^{2}\lambda}\sum_{j_{1},\ldots,j_{4}=-\infty}^{\infty}\gamma_{j_{1}}\gamma_{j_{2}}\gamma_{j_{3}}\gamma_{j_{4}}
\sum_{m=-\infty}^{\infty}\sum_{k_{1}=-a}^{a}g(\omega_{k_{1}})\overline{g(\omega_{m-k_{1}})}  \int_{-\infty}^{\infty}
\int_{-\infty}^{\infty} 
f(\frac{2u}{\lambda} - \omega_{k_{1}+j_{1}})
f(\frac{2v}{\lambda} + \omega_{k_{1}+r_{1}-j_{2}})
\\
&&\times\sinc(u)
\sinc(u + (m+j_{1}+j_{3})\pi)
\sinc(v)\sinc(v + (m-r_{2}-j_{4}+r_{1}-j_{2})\pi)dudv. 
\end{eqnarray*}
Thus by taking absolutes we have 
\begin{eqnarray*}
&&|A_{1} (r_{1},r_{2})| \\
&\leq &
\frac{1}{\pi^{2}\lambda}\sum_{j_{1},\ldots,j_{4}=-\infty}^{\infty}
|\gamma_{j_{1}}\gamma_{j_{2}}\gamma_{j_{3}}\gamma_{j_{4}}|
\sum_{m=-\infty}^{\infty}\sum_{k_{1}=-a}^{a}
\left|g(\omega_{k_{1}})\overline{g(\omega_{m-k_{1}})}\right|  \int_{-\infty}^{\infty}
\int_{-\infty}^{\infty} 
f(\frac{2u}{\lambda} - \omega_{k_{1}+j_{1}})
f(\frac{2v}{\lambda} + \omega_{k_{1}+r_{1}-j_{2}})
\\
&&\times\left|\sinc(u)
\sinc(u + (m+j_{1}+j_{3})\pi)
\sinc(v)\sinc(v + (m-r_{2}-j_{4}+r_{1}-j_{2})\pi)\right|dudv. 
\end{eqnarray*}
Finally, by following the same series of bounds used to prove
Lemma \ref{lemma:covariance}(iii) we have 
\begin{eqnarray*}
|A_{1}(r_{1},r_{2})|&\leq &\frac{1}{\pi^{2}\lambda}\sum_{j_{1},\ldots,j_{4}=-\infty}^{\infty}
|\gamma_{j_{1}}\gamma_{j_{2}}\gamma_{j_{3}}\gamma_{j_{4}}|
 \sup_{\omega}|g(\omega)|^{2}\|f\|_{2}^{2}\\
&&\times \sum_{m=-\infty}^{\infty}\int_{\mathbb{R}^{2}}
|\sinc(u-(m+j_{1}+j_{3})\pi\sinc(v+(m-r_{2}+r_{1}-j_{4}-j_{2})\pi)
\sinc(u)\sinc(v)|du dv \\
&<& \infty.
\end{eqnarray*}
Similarly we can bound $A_{2}(r_{1},r_{2})$ and 
$\lambda\cov\left[\widetilde{Q}_{a,\lambda}(g;r_{1}),\overline{\widetilde{Q}_{a,\lambda}(g;r_{2})}\right]$,
thus giving the required result. \hfill $\Box$

\vspace{3mm}
{\bf PROOF of Theorem  \ref{theorem:nonuniformvar}(iii)} 
The proof uses the expansion (\ref{eq:Qgamma}). Using this as a basis, 
we will show that 
\begin{eqnarray*}
A_{1}(r_{1},r_{1}) = U_{1}(r_{1},r_{2};\omega_{r_{1}},\omega_{r_{2}})
+ O(\ell_{\lambda,a,n}) \\
A_{2}(r_{1},r_{1}) = U_{2}(r_{1},r_{2};\omega_{r_{1}},\omega_{r_{2}})
+ O(\ell_{\lambda,a,n}).
\end{eqnarray*}
We find an approximation for $A_{1}(r_{1},r_{2})$ starting with the expansion given
in (\ref{eq:A1r1r2gamma}). 
We use the same proof as that used to prove  Lemma
\ref{lemma:var-asymp} to approximate the terms inside the sum
$\sum_{j_{1},\ldots,j_{4}}$. More precisely we let $m=k_{1}-k_{2}$, replace $\omega_{k_{1}}$
with $\omega_{}$ and by using the same methodology given in the proof
of Lemma \ref{lemma:var-asymp}
(and that $\sum_{j}|\gamma_{j}|<\infty$), we have 
\begin{eqnarray*}
&&A_{1} (r_{1},r_{2}) \\
&=& \frac{1}{2\pi^{3}}\sum_{j_{1},\ldots,j_{4}=-\infty}^{\infty}\gamma_{j_{1}}\gamma_{j_{2}}\gamma_{j_{3}}\gamma_{j_{4}}
\sum_{m=-2a}^{2a}\int_{2\pi \max(-a,-a+m)/\lambda}^{2\pi \min(a,a+m)/\lambda} f(-\omega-\omega_{j_{1}})
f(\omega+\omega_{r_{1}-j_{2}}) g(\omega)\overline{g(\omega-\omega_{m})}d\omega  \\
&&\times\int_{\mathbb{R}^{2}}\sinc(u)
\sinc(u + (m+j_{1}+j_{3})\pi)
\sinc(v)\sinc(v - (m-r_{2}-j_{4}+r_{1}-j_{2})\pi)dudv + O(\ell_{\lambda,a,n}). 
\end{eqnarray*}
By orthogonality of the sinc function we see that the above is zero unless
$m=-j_{1}-j_{3}$ {\it and} $m=r_{2}-r_{1}+j_{2}+j_{4}$ (and using that
$f(-\omega-\omega_{j_1})=f(\omega+\omega_{j_1})$), therefore
\begin{eqnarray*}
&&A_{1} (r_{1},r_{2}) \\
&=& \frac{1}{2\pi}\sum_{j_{1}+\ldots+j_{4}=r_{1}-r_{2}}\gamma_{j_{1}}\gamma_{j_{2}}\gamma_{j_{3}}\gamma_{j_{4}}
\int_{2 \pi \max(-a,-a-j_{1}-j_{3})/\lambda}^{2\pi \min(a,a-j_{1}-j_{3})/\lambda}g(\omega)\overline{g(\omega+\omega_{j_{1}+j_{3}})} f(\omega+\omega_{j_{1}})
f(\omega+\omega_{r_{1}-j_{2}})d\omega + O(\ell_{\lambda,a,n}). 
\end{eqnarray*}
This gives us $U_{1,1}(r_{1},r_{2};\omega_{r_{1}},\omega_{r_{2}})$.  Next we consider $A_{2}(r_{1},r_{2})$
\begin{eqnarray*}
A_{2}(r_{1},r_{2})&=& \frac{1}{\lambda^{3}}\sum_{j_{1},j_{2},j_{3},j_{4}\in \mathbb{Z}}\gamma_{j_{1}}\gamma_{j_{2}}\gamma_{j_{3}}\gamma_{j_{4}}\sum_{k_{1},k_{2}=-a}^{a}g(\omega_{k_{1}})\overline{g(\omega_{k_{2}})}
\int_{[-\lambda/2,\lambda/2]^{4}}c(s_{1}-s_{4})c(s_{2}-s_{3})\exp(is_{1}\omega_{k_{1}+j_{1}})\\
&&\exp(is_{4}\omega_{k_{2}+r_{2}+j_{4}})\exp(-is_{2}\omega_{k_{1}+r_{1}-j_{2}})\exp(-is_{3}\omega_{k_{2}-j_{3}})ds_{1}ds_{2}ds_{3}ds_{4}
\\
&=& \lambda\sum_{j_{1},j_{2},j_{3},j_{4}\in
  \mathbb{Z}}\sum_{k_{1},k_{2}=-a}^{a}g(\omega_{k_{1}})\overline{g(\omega_{k_{2}})}\int_{\mathbb{R}^{2}}
f(x)f(y)\sinc\left(\frac{\lambda x}{2} + (k_{1}+j_{1})\pi\right)\times \\
&&\sinc\left(\frac{\lambda x}{2} - (k_{2}+r_{2}+j_{4})\pi\right) 
\sinc\left(\frac{\lambda y}{2} - (k_{1}+r_{1}-j_{2})\pi\right)
\sinc\left(\frac{\lambda y}{2} + (k_{2}-j_{3})\pi\right)dxdy.
\end{eqnarray*}
Making a change of variables $u=\frac{\lambda x}{2} + (k_{1}+j_{1})\pi$
and $v=\frac{\lambda y}{2} - (k_{1}+r_{1}-j_{2})\pi$ we have  
\begin{eqnarray*}
A_{2}(r_{1}.r_{2})&=& \frac{1}{\pi^{2}\lambda}\sum_{j_{1},j_{2},j_{3},j_{4}\in
  \mathbb{Z}}\gamma_{j_{1}}\gamma_{j_{2}}\gamma_{j_{3}}\gamma_{j_{4}}\sum_{k_{1},k_{2}=-a}^{a}g(\omega_{k_{1}})\overline{g(\omega_{k_{2}})}\int_{\mathbb{R}^{2}}
f\left(\frac{2u}{\lambda}-\omega_{k_{1}+j_{1}}\right)f\left(\frac{2v}{\lambda}+\omega_{k_{1}+r_{1}-j_{2}}\right)\times \\
&&\sinc(u)\sinc\left(u-(k_{2}+r_{2}+j_{4}+k_{1}+j_{1})\pi\right)\sinc(v) 
\sinc\left(v+ (k_{2}-j_{3}+k_{1}+r_{1}-j_{2})\pi\right)dudv.
\end{eqnarray*}
Again by using the same proof as that given in  Lemma
\ref{lemma:var-asymp} to approximate the terms inside the sum
$\sum_{j_{1},\ldots,j_{4}}$  (setting $m=k_{1}+k_{2}$ and replacing
$\omega_{k_{1}}$ with $\omega$), we can approximate $A_{2}(r_{1},r_{2})$ with 
\begin{eqnarray*}
&&A_{2}(r_{1},r_{2}) \\
&=&  \frac{1}{2\pi^{3}}\sum_{j_{1},j_{2},j_{3},j_{4}\in
  \mathbb{Z}}\gamma_{j_{1}}\gamma_{j_{2}}\gamma_{j_{3}}\gamma_{j_{4}}\sum_{m=-2a}^{2a}
\int_{2\pi (\max(-a,-a+m)/\lambda}^{2\pi \min(a,a+m)/\lambda}g(\omega)\overline{g(-\omega+\omega_{m})}\int_{\mathbb{R}^{2}}
f(-\omega-\omega_{j_{1}})f(\omega+\omega_{r_{1}-j_{2}})\times \\
&& \sinc(u)\sinc\left(u-(m+r_{2}+j_{4}+j_{1})\pi\right)\sinc(v) 
\sinc\left(v+ (m-j_{3}+r_{1}-j_{2})\pi\right)dudvd\omega+ O(\ell_{\lambda,a,n}).
\end{eqnarray*}
Using the orthogonality of the sinc function, the inner integral is
non-zero when $m-j_{3}+r_{1}-j_{2}=0$ and
$m+r_{2}+j_{4}+j_{1}=0$. Setting $m=-r_{1}+j_{2}+j_{3}$, 
this implies 
\begin{eqnarray*}
&&A_{2}(r_{1},r_{2})\\
&=&  \frac{1}{2\pi}\sum_{j_{1}+j_{2}+j_{3}+j_{4} =
  r_{1}-r_{2}}\gamma_{j_{1}}\gamma_{j_{2}}\gamma_{j_{3}}\gamma_{j_{4}}\int_{2\pi \max(-a,-a-r_{1}+j_{2}+j_{3})/\lambda}^{2\pi \min(a,a-r_{1}+j_{2}+j_{3})/\lambda}
g(\omega)\overline{g(-\omega-\omega_{r_{1}-j_{3}-j_{2}})} \\
&& \times f(\omega+\omega_{j_{1}})f(\omega+\omega_{r_{1}-j_{2}})d\omega+ O(\ell_{\lambda,a,n}),
\end{eqnarray*}
thus giving us $U_{1,2}(\rb_{1},\rb_{2};\ob_{\rb_{1}},\ob_{\rb_{2}})$.

By using Lemma \ref{lemma:cumulantsA}, we can show that 
\begin{eqnarray*}
\lambda\cov\left[\widetilde{Q}_{a,\lambda}(g;r_{1}),\overline{\widetilde{Q}_{a,\lambda}(g;r_{2})}\right]
&=& A_{3}(r_{1},r_{2}) + A_{4}(r_{1},r_{2}) + O\left(\frac{\lambda}{n}\right) 
\end{eqnarray*}
where
\begin{eqnarray*}
A_{3}(r_{1},r_{2})&=& \lambda\sum_{k_{1},k_{2}=-a}^{a}g(\omega_{k_{1}})g(\omega_{k_{2}})
\cov\big[Z(s_{1})\exp(is_{1}\omega_{k_{1}}),Z(s_{3})\exp(-is_{3}\omega_{k_{2}})\big]
\times \\
&&\cov\big[Z(s_{2})\exp(-is_{2}\omega_{k_{1}+r_{1}}),Z(s_{4})\exp(is_{4}\omega_{k_{2}+r_{2}})\big]\\
A_{4}(r_{1},r_{2})&=& \lambda\sum_{k_{1},k_{2}=-a}^{a}g(\omega_{k_{1}})g(\omega_{k_{2}})
 \cov\big[Z(s_{1})\exp(is_{1}\omega_{k_{1}}),Z(s_{4})\exp(is_{4}\omega_{k_{2}+r_{2}})\big] \times\\
&&\cov\big[Z(s_{2})\exp(-is_{2}\omega_{k_{1}+r_{1}}),Z(s_{3})\exp(-is_{3}\omega_{k_{2}})\big].
\end{eqnarray*}
By following the same proof as that used to prove $A_{1}(r_{1},r_{2})$
we have 
\begin{eqnarray*}
&&A_{3} (r_{1},r_{2}) \\
&=& \frac{\lambda}{(2\pi)^{2}}\sum_{j_{1},\ldots,j_{4}=-\infty}^{\infty}\gamma_{j_{1}}\gamma_{j_{2}}\gamma_{j_{3}}\gamma_{j_{4}}
\sum_{k_{1},k_{2}=-a}^{a}g(\omega_{k_{1}})g(\omega_{k_{2}}) 
\int_{-\infty}^{\infty} \int_{-\infty}^{\infty} f(x)f(y)
\sinc\left(\frac{\lambda x}{2} + (k_{1}+j_{1})\pi\right)\\
&& \sinc\left(\frac{\lambda x}{2}-(k_{2}+j_{3})\pi\right)
\sinc\left(\frac{\lambda y}{2} - (r_{1}+k_{1}-j_{2})\pi\right)
\sinc\left(\frac{\lambda y}{2} + (r_{2}+k_{2}-j_{4})\pi\right)dxdy \\
&=& \frac{1}{\pi^{2}\lambda}\sum_{j_{1},\ldots,j_{4}=-\infty}^{\infty}\gamma_{j_{1}}\gamma_{j_{2}}\gamma_{j_{3}}\gamma_{j_{4}}
\sum_{k_{1},k_{2}=-a}^{a}g(\omega_{k_{1}})g(\omega_{k_{2}})  \int_{-\infty}^{\infty}
\int_{-\infty}^{\infty} 
f(\frac{2u}{\lambda} - \omega_{k_{1}+j_{1}})
f(\frac{2v}{\lambda} + \omega_{k_{1}+r_{1}-j_{2}})
\\
&&\times\sinc(u)\sinc(u - (k_{2}+k_{1}+j_{1}+j_{3})\pi)
\sinc(v)\sinc(v + (r_{2}+k_{2}-j_{4}+r_{1}+k_{1}-j_{2})\pi)dudv \\
&=& \frac{1}{\pi^{2}\lambda}\sum_{j_{1},\ldots,j_{4}=-\infty}^{\infty}\gamma_{j_{1}}\gamma_{j_{2}}\gamma_{j_{3}}\gamma_{j_{4}}
\sum_{m=-2a}^{2a}\sum_{k_{1}=\max(-a,-a+m)}^{\min(a,a+m)}g(\omega_{k_{1}})g(\omega_{m-k_{1}})  \int_{-\infty}^{\infty}
\int_{-\infty}^{\infty} 
f(\frac{2u}{\lambda} - \omega_{k_{1}+j_{1}})
f(\frac{2v}{\lambda} + \omega_{k_{1}+r_{1}-j_{2}})
\\
&&\times\sinc(u)\sinc(u - (m+j_{1}+j_{3})\pi)
\sinc(v)\sinc(v + (r_{2}+m-j_{4}+r_{1}-j_{2})\pi)dudv.
\end{eqnarray*}
Again using the method used to bound $A_{1}(r_{1},r_{2})$ gives 
\begin{eqnarray*}
A_{3} (r_{1},r_{2}) 
&=& \frac{1}{2\pi^{3}}\sum_{j_{1}+j_{2}+j_{3}+j_{4}=r_{1}+r_{2}}\gamma_{j_{1}}\gamma_{j_{2}}\gamma_{j_{3}}\gamma_{j_{4}}
\sum_{m}\int_{2\pi \max(-a,-a-j_{1}-j_{3})/\lambda}^{2\pi
  \min(a,a-j_{1}-j_{3})/\lambda}g(\omega)g(\omega_{-j_{1}-j_{3}}-\omega)
\times \\
&& f(\omega+\omega_{j_{1}})f(\omega+\omega_{r_{1}-j_{2}})d\omega +
O(\ell_{\lambda,a,n}) = 
U_{1,2}(\rb_{1},\rb_{2};\ob_{\rb_{1}},\ob_{\rb_{2}})+O(\ell_{\lambda,a,n}).
\end{eqnarray*}
Finally we consider $A_{4}(r_{1},r_{2})$. Using the same expansion as
the above we have  
\begin{eqnarray*}
&&A_{4}(r_{1},r_{2})\\
&=& \frac{1}{\lambda^{3}}\sum_{j_{1},j_{2},j_{3},j_{4}\in \mathbb{Z}}\gamma_{j_{1}}\gamma_{j_{2}}\gamma_{j_{3}}\gamma_{j_{4}}\sum_{k_{1},k_{2}=-a}^{a}g(\omega_{k_{1}})g(\omega_{k_{2}})
\int_{[-\lambda/2,\lambda/2]^{4}}c(s_{1}-s_{4})c(s_{2}-s_{3})\exp(is_{1}\omega_{k_{1}+j_{1}})\\
&&\exp(-is_{4}\omega_{k_{2}+r_{2}-j_{4}})\exp(-is_{2}\omega_{k_{1}+r_{1}-j_{2}})\exp(is_{3}\omega_{k_{2}+j_{3}})ds_{1}ds_{2}ds_{3}ds_{4}
\\
&=& \lambda\sum_{j_{1},j_{2},j_{3},j_{4}\in
  \mathbb{Z}}\sum_{k_{1},k_{2}=-a}^{a}g(\omega_{k_{1}})g(\omega_{k_{2}})\int_{\mathbb{R}^{2}}
f(x)f(y)\sinc\left(\frac{\lambda x}{2} + (k_{1}+j_{1})\pi\right)\times \\
&&\sinc\left(\frac{\lambda x}{2} + (k_{2}+r_{2}-j_{4})\pi\right) 
\sinc\left(\frac{\lambda y}{2} - (k_{1}+r_{1}-j_{2})\pi\right)
\sinc\left(\frac{\lambda y}{2} - (k_{2}+j_{3})\pi\right)dxdy \\
 &=& \frac{1}{\pi^{2}\lambda}\sum_{j_{1},j_{2},j_{3},j_{4}\in
  \mathbb{Z}}\gamma_{j_{1}}\gamma_{j_{2}}\gamma_{j_{3}}\gamma_{j_{4}}\sum_{k_{1},k_{2}=-a}^{a}g(\omega_{k_{1}})g(\omega_{k_{2}})\int_{\mathbb{R}^{2}}
f(\frac{2u}{\lambda}-\omega_{k_{1}+j_{1}})f(\frac{2v}{\lambda}+\omega_{k_{1}+r_{1}-j_{2}})\times \\
&&\sinc(u)\sinc\left(u+(k_{2}+r_{2}-j_{4}-k_{1}-j_{1})\pi\right)\sinc(v) 
\sinc\left(v- (k_{2}+j_{3}-k_{1}-r_{1}+j_{2})\pi\right)dudv \\
&=&  \frac{1}{2\pi^{3}}\sum_{j_{1},j_{2},j_{3},j_{4}\in
  \mathbb{Z}}\gamma_{j_{1}}\gamma_{j_{2}}\gamma_{j_{3}}\gamma_{j_{4}}
\sum_{m=-2a}^{2a}\int_{2\pi\max(-a,-a+m)/\lambda}^{2\pi\min(a,a+m)} g(\omega)g(\omega-\omega_{m})\int_{\mathbb{R}^{2}}
f(-\omega-\omega_{j_{1}})f(\omega+\omega_{r_{1}-j_{2}})\times \\
&& \sinc(u)\sinc\left(u+(r_{2}-j_{4}-m-j_{1})\pi\right)\sinc(v) 
\sinc\left(v- (j_{3}-m-r_{1}+j_{2})\pi\right)du dv d\omega +
O(\ell_{\lambda,a,n}) \\
&=& \frac{1}{2\pi}\sum_{j_{1}+j_{2}+j_{3}+j_{4}=r_{1}+r_{2}
 }\gamma_{j_{1}}\gamma_{j_{2}}\gamma_{j_{3}}\gamma_{j_{4}}
\int_{\max(-a,a-r_{1}+j_{2}+j_{3})/\lambda}^{2\pi \min(a,a-r_{1}+j_{2}+j_{3})/\lambda}
g(\omega)g(\omega-\omega_{j_{2}+j_{3}-r_{2}})\\
&&\times f(\omega+\omega_{j_{1}})f(\omega+\omega_{r_{1}-j_{2}})d\omega +
O(\ell_{\lambda,a,n}).
\end{eqnarray*}
This gives the desired result. \hfill $\Box$

\vspace{3mm}
We now obtain and approximation of $U_{1}$ and $U_{2}$.

\vspace{3mm}
{\bf PROOF of Corollary \ref{cor:C}} By using 
Lipschitz continuity of $g(\cdot)$ and $f(\cdot)$ and $|\gamma_{\jb}| \leq
C\prod_{i=1}^{d}|j_{i}|^{-(1+\delta)}I(j_{i}\neq 0)$ we obtain the
result. \hfill $\Box$

\vspace{3mm}

{\bf PROOF of Theorem \ref{theorem:nonGaussiannonUniform}} 
We prove the result for the case $d=1$ and using
$A_{1}(r_{1},r_{2}),\ldots,A_{4}(r_{1},r_{2})$
defined in proof of Theorem \ref{theorem:nonuniformvar}.
The proof
is identical to the proof of Theorem
\ref{theorem:variance-nongaussian}. Following the same notation in
proof of Theorem \ref{theorem:variance-nongaussian} we have 
\begin{eqnarray*}
\lambda\cov\left[\widetilde{Q}_{a,\lambda}(g;r_{1}),\widetilde{Q}_{a,\lambda}(g;r_{2})\right]
=
A_{1}(r_{1},r_{2}) + A_{2}(r_{1},r_{2}) + B_{1}(r_{1},r_{2}) + B_{2}(r_{1},r_{2})
+ O\left(\frac{\lambda}{n}\right),
\end{eqnarray*}
with $|B_{2}(r_{1},r_{2})| = O((a\lambda)^{2}/n^{2})$ and 
the main term involving the trispectral density is 
\begin{eqnarray*}
&&B_{1}(r_{1},r_{2})  \\
&=& \lambda c_{4}\sum_{k_{1},k_{2}=-a}^{a}g(\omega_{k_{1}})\overline{g(\omega_{k_{2}})}
\Ex\bigg[\kappa_{4}(s_{1}-s_{2},s_{1}-s_{3},s_{1}-s_{4})e^{is_{1}\omega_{k_{1}}}e^{-is_{2}\omega_{k_{1}+r_{1}}}
e^{-is_{3}\omega_{k_{2}}}e^{is_{4}\omega_{k_{2}+r_{2}}} \bigg] 
\\
&=& \frac{c_{4}}{(2\pi)^{3}\lambda^{3}}
\sum_{j_{1},\ldots,j_{4}=-\infty}^{\infty}\gamma_{j_{1}}\gamma_{j_{2}}\gamma_{j_{3}}\gamma_{j_{4}}\sum_{k_{1},k_{2}=-a}^{a}g(\omega_{k_{1}})\overline{g(\omega_{k_{2}})}
\int_{\mathbb{R}^{3}}
f_{4}(\omega_{1},\omega_{2},\omega_{3})\int_{[-\lambda/2,\lambda/2]^{4}}
 e^{is_{1}(\omega_{1}+\omega_{2}+\omega_{3}+\omega_{k_{1}+j_{1}})} \\
&&e^{-is_{2}(\omega_{1}+\omega_{k_{1}+r_{1}-j_{2}})}e^{-is_{3}(\omega_{2}+\omega_{k_{2}-j_{3}})} 
e^{is_{4}(-\omega_{3}+\omega_{k_{2}+r_{2}+j_{4}})}ds_{1}ds_{2}ds_{3}ds_{4}d\omega_{1}d\omega_{2}d\omega_{3} \\
 &=&   \frac{c_{4}\lambda}{(2\pi)^{3}}\sum_{j_{1},\ldots,j_{4}=-\infty}^{\infty}\gamma_{j_{1}}\gamma_{j_{2}}\gamma_{j_{3}}\gamma_{j_{4}}\sum_{k_{1},k_{2}=-a}^{a}g(\omega_{k_{1}})\overline{g(\omega_{k_{2}})}
\int_{\mathbb{R}^{3}}
f_{4}(\omega_{1},\omega_{2},\omega_{3})
\sinc\left(\frac{\lambda(\omega_{1}+\omega_{2}+\omega_{3})}{2}+(k_{1}+j_{1})\pi\right)\\
&&\times \sinc\left(\frac{\lambda \omega_{1}}{2} + (k_{1}+r_{1}-j_{2})\pi\right)
 \sinc\left(\frac{\lambda \omega_{2}}{2}+(k_{2}-j_{3})\pi\right)
\sinc\left(\frac{\lambda \omega_{3}}{2}-(k_{2}+r_{2}+j_{4})\pi\right)
d\omega_{1}d\omega_{2}d\omega_{3}. 
\end{eqnarray*}
Now we make a change of variables and let 
$u_{1}=\frac{\lambda \omega_{1}}{2} + (k_{1}+r_{1}-j_{2})$, 
$u_{2} = \frac{\lambda \omega_{2}}{2}+(k_{2}-j_{3})\pi$
and $u_{3}=\frac{\lambda \omega_{3}}{2}-(k_{2}+r_{2}+j_{4})\pi$, this gives 
\begin{eqnarray*}
&&B_{1}(r_{1},r_{2}) \\
 &=&
 \frac{c_{4}}{\pi^{3}\lambda^{2}}
\sum_{j_{1},\ldots,j_{4}=-\infty}^{\infty}\gamma_{j_{1}}\gamma_{j_{2}}\gamma_{j_{3}}\gamma_{j_{4}}
\sum_{k_{1},k_{2}=-a}^{a}
\int_{\mathbb{R}^{3}}g(\omega_{k_{1}})\overline{g(\omega_{k_{2}})}
f_{4}\left(\frac{2u_{1}}{\lambda}-\omega_{k_{1}+r_{1}-j_{2}},\frac{2u_{2}}{\lambda}-\omega_{k_{2}-j_{3}},
\frac{2u_{3}}{\lambda}+\omega_{k_2 + r_2+j_4}\right)
 \\
&&\times\sinc\left(u_{1}+u_{2}+u_{3}+(r_{2}-r_{1}+j_{1}+j_{2}+j_{3}+j_{4})\pi\right)\sinc(u_{1})
\sinc(u_{2})\sinc(u_{3})du_{1}du_{2}du_{3}.
\end{eqnarray*}
Next we exchange the summand with a double integral and use Lemma \ref{lemma:sum-integral}(iii)
together with Lemma \ref{lemma:1a}, equation (\ref{eq:in})
to obtain 
\begin{eqnarray*}
&&B_{1}(r_{1},r_{2}) \\
 &=&
 \frac{c_{4}}{\pi^{3}(2\pi)^{2}}\sum_{j_{1},\ldots,j_{4}=-\infty}^{\infty}\gamma_{j_{1}}\gamma_{j_{2}}\gamma_{j_{3}}\gamma_{j_{4}}
\int_{2\pi[-a/\lambda,a/\lambda]^{2}}g(\omega_{1})\overline{g(\omega_{2})}\times \\
&&\int_{\mathbb{R}^{3}}f_{4}\left(\frac{2u_{1}}{\lambda}-\omega_{1}-\omega_{r_{1}-j_2},
\frac{2u_{2}}{\lambda}-\omega_{2}-\omega_{j_{3}},\frac{2u_{3}}{\lambda}+\omega_{2}+\omega_{r_2+j_4}\right)\times \\
 &&\sinc\left(u_{1}+u_{2}+u_{3}+(r_{2}-r_{1}+j_{1}+j_{2}+j_{3}+j_{4})\pi\right)
\sinc(u_{1})\sinc(u_{2})\sinc(u_{3})du_{1}du_{2}du_{3}d\omega_{1}d\omega_{2}+ O\left(\frac{1}{\lambda}\right).
\end{eqnarray*}
By using Lemma \ref{lemma:cum4bound},
we  replace 
$f_{4}\left(\frac{2u_{1}}{\lambda}-\omega_{1}-\omega_{r_{1}-j_2},
\frac{2u_{2}}{\lambda}-\omega_{2}-\omega_{j_{3}},\frac{2u_{3}}{\lambda}+\omega_{2}+\omega_{r_2+j_4}\right)$
with $f_{4}\left(-\omega_{1}-\omega_{r_{1}-j_2},
-\omega_{2}-\omega_{j_{3}},\omega_{2}+\omega_{r_2+j_4}\right)$,
 to give
\begin{eqnarray*}
&&B_{1}(r_{1},r_{2}) \\
 &=&
 \frac{c_{4}}{\pi^{3}(2\pi)^{2}}\sum_{j_{1},\ldots,j_{4}=-\infty}^{\infty}\gamma_{j_{1}}\gamma_{j_{2}}\gamma_{j_{3}}\gamma_{j_{4}}
\int_{2\pi[-a/\lambda,a/\lambda]^{2}}g(\omega_{1})\overline{g(\omega_{2})}\times \\
&&\int_{\mathbb{R}^{3}}f_{4}\left(-\omega_{1}-\omega_{r_{1}-j_2},
-\omega_{2}-\omega_{j_{3}},\omega_{2}+\omega_{r_2+j_4}\right)
\sinc\left(u_{1}+u_{2}+u_{3}+(r_{2}-r_{1}+j_{1}+j_{2}+j_{3}+j_{4})\pi\right)\\
&&\times \sinc(u_{1})
\sinc(u_{2})\sinc(u_{3})du_{1}du_{2}du_{3}d\omega_{1}d\omega_{2}+
O\left(\frac{\log^{3}(\lambda)}{\lambda}\right) \\
&=&  \frac{c_{4}}{(2\pi)^{2}}\sum_{j_{1}+j_{2}+j_{3}+j_{4}=r_{2}-r_{1}}\gamma_{j_{1}}\gamma_{j_{2}}\gamma_{j_{3}}\gamma_{j_{4}}
\int_{2\pi[-a/\lambda,a/\lambda]^{2}}g(\omega_{1})\overline{g(\omega_{2})}\times \\
&&f_{4}\left(-\omega_{1}-\omega_{r_{1}-j_2},-\omega_{2}-\omega_{j_{3}},\omega_{2}+\omega_{r_2+j_4}\right)d\omega_{1}d\omega_{2}
+ O\left(\frac{\log^{3}(\lambda)}{\lambda}\right),
\end{eqnarray*}
where the last line follows from (\ref{eq:sinc3}). 

To obtain an expression for
$\lambda\cov\left[\widetilde{Q}_{a,\lambda}(g;r_{1}),\overline{\widetilde{Q}_{a,\lambda}(g;r_{2})}\right]$
we note that 
\begin{eqnarray*}
\lambda\cov\left[\widetilde{Q}_{a,\lambda}(g;r_{1}),\overline{\widetilde{Q}_{a,\lambda}(g;r_{2})}\right]
=
A_{3}(r_{1},r_{2}) + A_{4}(r_{1},r_{2}) + B_{3}(r_{1},r_{2}) + B_{4}(r_{1},r_{2})
+ O\left(\frac{\lambda}{n}\right),
\end{eqnarray*}
just as in the proof of Theorem \ref{theorem:variance-nongaussian} we
can show that $|B_{4}(r_{1},r_{2})| = O((\lambda a)^{d}/n^{2})$ and
the leading term involving the trispectral density is 
\begin{eqnarray*}
&&B_{3}(r_{1},r_{2})  \\
&=& \lambda c_{4}\sum_{k_{1},k_{2}=-a}^{a}g(\omega_{k_{1}})g(\omega_{k_{2}})
\Ex\bigg[\kappa_{4}(s_{1}-s_{2},s_{1}-s_{3},s_{1}-s_{4})e^{is_{1}\omega_{k_{1}}}e^{-is_{2}\omega_{k_{1}+r_{1}}}
e^{is_{3}\omega_{k_{2}}}e^{-is_{4}\omega_{k_{2}+r_{2}}} \bigg] \\
&=& \frac{c_{4}}{(2\pi)^{3}\lambda^{3}}
\sum_{j_{1},\ldots,j_{4}=-\infty}^{\infty}\gamma_{j_{1}}\gamma_{j_{2}}\gamma_{j_{3}}\gamma_{j_{4}}\sum_{k_{1},k_{2}=-a}^{a}g(\omega_{k_{1}})g(\omega_{k_{2}})
\int_{\mathbb{R}^{3}}
f_{4}(\omega_{1},\omega_{2},\omega_{3})\int_{[-\lambda/2,\lambda/2]^{4}}
 e^{is_{1}(\omega_{1}+\omega_{2}+\omega_{3}+\omega_{k_{1}+j_{1}})} \\
&&e^{-is_{2}(\omega_{1}+\omega_{k_{1}+r_{1}-j_{2}})}e^{-is_{3}(\omega_{2}-\omega_{k_{2}+j_{3}})} 
e^{is_{4}(-\omega_{3}-\omega_{k_{2}+r_{2}-j_{4}})}ds_{1}ds_{2}ds_{3}ds_{4}d\omega_{1}d\omega_{2}d\omega_{3} \\
 &=&   \frac{c_{4}\lambda}{(2\pi)^{3}}\sum_{j_{1},\ldots,j_{4}=-\infty}^{\infty}\gamma_{j_{1}}\gamma_{j_{2}}\gamma_{j_{3}}\gamma_{j_{4}}\sum_{k_{1},k_{2}=-a}^{a}g(\omega_{k_{1}})\overline{g(\omega_{k_{2}})}
\int_{\mathbb{R}^{3}}
f_{4}(\omega_{1},\omega_{2},\omega_{3})
\sinc\left(\frac{\lambda(\omega_{1}+\omega_{2}+\omega_{3})}{2}+(k_{1}+j_{1})\pi\right)\\
&&\times \sinc\left(\frac{\lambda \omega_{1}}{2} + (k_{1}+r_{1}-j_{2})\pi\right)
 \sinc\left(\frac{\lambda \omega_{2}}{2}-(k_{2}+j_{3})\pi\right)
\sinc\left(\frac{\lambda \omega_{3}}{2}+(k_{2}+r_{2}-j_{4})\pi\right)
d\omega_{1}d\omega_{2}d\omega_{3}. 
\end{eqnarray*}
We make a change of variables $u_{1}=\frac{\lambda \omega_{1}}{2} + (k_{1}+r_{1}-j_{2})\pi$
$u_{2}=\frac{\lambda \omega_{2}}{2}-(k_{2}+j_{3})\pi$ and
$u_{3}=\frac{\lambda \omega_{3}}{2}+(k_{2}+r_{2}-j_{4})\pi$. 
This gives 
\begin{eqnarray*}
&&B_{3}(r_{1},r_{2}) \\
&=&  \frac{c_{4}}{\lambda^{2}(2\pi)^{3}}\sum_{j_{1},\ldots,j_{4}=-\infty}^{\infty}\gamma_{j_{1}}\gamma_{j_{2}}\gamma_{j_{3}}\gamma_{j_{4}}\sum_{k_{1},k_{2}=-a}^{a}g(\omega_{k_{1}})g(\omega_{k_{2}})
\int_{\mathbb{R}^{3}}
f_{4}(\frac{2u_{1}}{\lambda}-\omega_{k_{1}+r_{1}-j_{2}},\frac{2u_{2}}{\lambda}+\omega_{k_{2}+j_{3}},\frac{2u_{3}}{\lambda}-
\omega_{k_{2}-r_{2}+j_{4}})\\
&&\times \sinc\left(u_{1}+u_{2}+u_{3}+(j_{1}+j_{2}+j_{3}+j_{4}-r_{1}-r_{2})\pi\right)
 \sinc(u_{1})\sinc(u_{2})\sinc(u_{3})du_{1}du_{2}du_{3} \\
&=& \frac{1}{(2\pi)^{2}}\sum_{j_{1}+j_{2}+j_{3}+j_{4}=r_{1}+r_{2}}
\gamma_{j_{1}}\gamma_{j_{2}}\gamma_{j_{3}}\gamma_{j_{4}}\int_{2\pi[-a/\lambda,a/\lambda]^{2}}g(\omega_{1})g(\omega_{2})
f_{4}(-\omega_{1}-\omega_{r_{1}-j_{2}},\omega_{2}+\omega_{j_{3}},-\omega_{2}
-\omega_{-r_{2}+j_{4}})d\omega_{1}d\omega_{2}.
\end{eqnarray*}

Finally, by replacing
$f_{4}\left(-\omega_{1}-\omega_{r_{1}-j_2},-\omega_{2}-\omega_{j_{3}},\omega_{2}+\omega_{r_2+j_4}\right)$
with
$f_{4}\left(-\omega_{1},-\omega_{2},\omega_{2}\right)$
in $B_{1}(r_{1},r_{2})$ and $f_{4}(-\omega_{1}-\omega_{r_{1}-j_{2}},\omega_{2}+\omega_{j_{3}},-\omega_{2}
-\omega_{-r_{2}+j_{4}})$ with
$f_{4}(-\omega_{1},\omega_{2},-\omega_{2})$ and using the pointwise Lipschitz
continuity of $f_{4}$ and that $|\gamma_{j}|\leq CI(j\neq
0)|j|^{-(1+\delta)}$ we obtain 
$B_{1}(r_{1},r_{2}) =
D_{1}+O\left(\frac{\log^{3}(\lambda)}{\lambda}+\frac{1}{\lambda}\right)$
and  $B_{3}(r_{1},r_{2}) =
D_{3}+O\left(\frac{\log^{3}(\lambda)}{\lambda}+\frac{1}{\lambda}\right)$.
Thus giving the required result. \hfill $\Box$

%% file: B_1_spatial.tex
\section{Technical Lemmas}\label{appendix:tech-proofs}

We first prove Lemma \ref{lemma:orthogonal}, then state four lemmas, which form an important component in
the proofs of this paper. Through out this section we use $C$ to
denote a finite generic constant. It is worth mentioning that many of
these results build on the work of  T. Kawata (see
\citeA{p:kaw-59}).

{\bf PROOF  of Lemma \ref{lemma:orthogonal}} We first prove (\ref{eq:sincA}). 
By using partial fractions and the definition of the sinc function we
have  
\begin{eqnarray*}
\int_{-\infty}^{\infty} \sinc(u) \sinc(u+x)du &=&  
\frac{1}{x}\int_{-\infty}^{\infty}\sin(u)\sin(u+x)\bigg(\frac{1}{u} -
\frac{1}{u+x} \bigg)du \\
&=&\frac{1}{x}\int_{-\infty}^{\infty}\frac{\sin(u)\sin(u+x)}{u}du - 
\frac{1}{x}\int_{-\infty}^{\infty}\frac{\sin(u)\sin(u+x)}{u+x}du.  
\end{eqnarray*}
For the second integral we make a change of variables $u^{\prime} =
u+x$, this gives 
\begin{eqnarray*}
\int_{-\infty}^{\infty} \sinc(u) \sinc(u+x)du 
&=&\frac{1}{x}\int_{-\infty}^{\infty}\frac{\sin(u)\sin(u+x)}{u}du - 
\frac{1}{x}\int_{-\infty}^{\infty}\frac{\sin(u^{\prime})\sin(u^{\prime}-x)}{u^{\prime}}du^{\prime}
\\
&= &  \frac{1}{x}\int_{-\infty}^{\infty}\frac{\sin(u)}{u}\bigg(\sin(u+x) -
\sin(u-x)\bigg)du \\
&=&
\frac{2\sin(x)}{x}\int_{-\infty}^{\infty}\frac{\cos(u)\sin(u)}{u}du =
\frac{\pi\sin(x)}{x}. 
\end{eqnarray*}
To prove (\ref{eq:sincB}), it is clear that for $x=s\pi$ (with $s\in
\mathbb{Z}\{0\}$) 
$\frac{\pi\sin(s\pi)}{s\pi}=0$, which gives the result. \hfill $\Box$

\vspace{3mm}

The following result is used to obtain bounds for the variance and
higher order cumulants. 
\begin{lemma}\label{lemma:1a}
 Define the function $\ell_{p}(x) = C/e$ for
  $|x|\leq e$ and $\ell_{p}(x)=C\log^{p}|x|/|x|$ for $|x|\geq e$.
\begin{itemize}
\item[(i)] We have 
\begin{eqnarray}
\label{eq:i1}
  \int_{-\infty}^{\infty}\frac{|\sin^{}(x)\sin(x+y)|}{|x(x+y)|}dx
  \leq 
\left\{
\begin{array}{cc}
C\frac{\log|y|}{|y|} & |y|\geq e \\
C & |y|< e  \\
\end{array}
\right. = \ell_{1}(y),
\end{eqnarray} 
\begin{eqnarray}
\label{eq:i2}
\int_{-\infty}^{\infty}|\sinc(x)|\ell_{p}(x+y)dx \leq \ell_{p+1}(y)
\end{eqnarray}
and 
 \begin{eqnarray}
\label{eq:in}
\int_{\mathbb{R}^{p}}\left|\sinc\left(\sum_{j=1}^{p}x_{j}\right)\prod_{j=1}^{p}\sinc(x_{j})\right|dx_{1}\ldots
dx_{p} \leq C, 
\end{eqnarray}
\item[(ii)] 
\begin{eqnarray*}
\sum_{m=-a}^{a}\int_{-\infty}^{\infty}\left|\frac{\sin^{2}(x)}{x(x+m\pi)}\right|dx \leq
C\log^{2}a
\end{eqnarray*}
\item[(iii)]
\begin{eqnarray}
\label{eq:sinx}
\int_{-\infty}^{\infty}\int_{-\infty}^{\infty}\sum_{m=-a}^{a}\frac{\sin^{2}(x)}{|x(x+m\pi)|}
\frac{\sin^{2}(y)}{|y(y+m\pi)|}dxdy \leq C, 
\end{eqnarray}
\item[(iv)] 
\begin{eqnarray*}
&& \sum_{m_{1},\ldots,m_{q-1}=-a}^{a}\int_{\mathbb{R}^{q}}\bigg|\prod_{j=1}^{q-1}\sinc(x_{j})\sinc(x_{j}+m_{j}\pi)
\times \\
&&\sinc(x_{q})\sinc(x_{q}+\pi\sum_{j=1}^{q-1}m_{j})\bigg|\prod_{j=1}^{q}dx_{j}
\leq  C\log^{2(q-2)}(a), 
\end{eqnarray*}
\end{itemize}
where $C$ is a finite generic constant which is independent of  $a$. 
\end{lemma}
PROOF.  We first prove (i), equation (\ref{eq:i1}). It is clear that for 
$|y|\leq e$ that $\int_{-\infty}^{\infty}\frac{|\sin(x)\sin(x+y)|}{|x(x+y)|}dx\leq C$. Therefore
we now consider the case $|y|>e$, without loss of generality we prove the result for $y>e$.
Partitioning the integral we have 
\begin{eqnarray*}
&&\int_{-\infty}^{\infty}\frac{|\sin(x)\sin(x+y)|}{|x(x+y)|}dx = I_{1}+I_{2}+I_{3}+I_{4}+I_{5}, 
\end{eqnarray*}
where
\begin{eqnarray*}
I_{1} &=& \int_{0}^{y}\frac{|\sin(x)\sin(x+y)|}{|x(x+y)|}dx \qquad
I_{2} = \int_{-y}^{0}\frac{|\sin(x)\sin(x+y)|}{|x(x+y)|}dx \\
I_{3} &=& \int^{-y}_{-2y}\frac{|\sin(x)\sin(x+y)|}{|x(x+y)|}dx \qquad
I_{4} = \int^{-2y}_{-\infty}\frac{|\sin(x)\sin(x+y)|}{|x(x+y)|}dx \\
I_{5} &=& \int_{y}^{\infty}\frac{|\sin(x)\sin(x+y)|}{|x(x+y)|}dx.
\end{eqnarray*}
To bound $I_{1}$ we note that for $y>1$ and $x>0$ that $|\sin(x+y)/(x+y)|\leq 1/y$, thus 
\begin{eqnarray*}
I_{1} &=& \frac{1}{y}\int_{0}^{y}\frac{|\sin(x)|}{|x|}dx \leq C\frac{\log y}{y}.
\end{eqnarray*}
To bound $I_{2}$, we further partition the integral 
\begin{eqnarray*}
I_{2} &=&  \int_{-y}^{-y/2}\frac{|\sin(x)\sin(x+y)|}{|x(x+y)|}dx + \int_{-y/2}^{0}\frac{|\sin(x)\sin(x+y)|}{|x(x+y)|}dx \\
 &\leq& \frac{2}{y}\int_{-y}^{-y/2}\frac{|\sin(x+y)|}{|(x+y)|}dx + \frac{2}{y}\int_{-y/2}^{0}\frac{|\sin(x)|}{|x|}dx 
  \leq C\frac{\log y}{y}. 
\end{eqnarray*}
To bound $I_{3}$, we use the bound 
\begin{eqnarray*}
I_{3} \leq \frac{1}{y}\int_{-2y}^{-y}\frac{|\sin(x+y)|}{|(x+y)|}dx \leq C\frac{\log y}{y}.
\end{eqnarray*}
To bound $I_{4}$ we use that for $y>0$, $\int_{y}^{\infty}x^{-2}dx \leq C|y|^{-1}$, thus 
\begin{eqnarray*}
I_{4} \leq  \int^{-y}_{-\infty}\frac{1}{x^{2}}dx \leq C|y|^{-1}
\end{eqnarray*}
and using a similar argument we have $I_{5}\leq C|y|^{-1}$. Altogether, this gives (\ref{eq:i1}).

We now prove (\ref{eq:i2}). It is clear that for 
$|y|\leq e$ that $\int_{-\infty}^{\infty}|\sinc(x)|\ell_{p}(x+y)|dx\leq C$. Therefore
we now consider the case $|y|>e$, without loss of generality we prove the result for $y>e$.
As in (\ref{eq:i1}) we partition the integral 
\begin{eqnarray*}
\int_{-\infty}^{\infty}|\sinc(x)\ell_{p}(x+y)|dx = II_{1}+\ldots+II_{5}, 
\end{eqnarray*}
where $II_{1},\ldots,II_{5}$ are defined in the same way as $I_{1},\ldots,I_{5}$ just with 
$|\sinc(x)\ell_{p}(x+y)|$ replacing $\frac{|\sin(x)\sin(x+y)|}{|x(x+y)|}$. 
To bound $II_{1}$ we note that 
\begin{eqnarray*}
II_{1} &=& \int_{0}^{y}|\sinc(x)\ell_{p}(x+y)|dx \leq 
\frac{\log^{p}(y)}{y}\int_{0}^{y}|\sinc(x)|dx \leq \frac{\log^{p+1}(y)}{y},
\end{eqnarray*}
we use similar method to show $II_{2}\leq C\frac{\log^{p+1}(y)}{y}$ and $II_{3}\leq C\frac{\log^{p+1}(y)}{y}$. Finally to 
bound $II_{4}$ and $II_{5}$ we note that by using a change of variables
$ x=yz$, we have 
\begin{eqnarray*}
II_{5} &=& \int_{y}^{\infty}\frac{|\sin(x)|\log^{p}(x+y)}{x(x+y)}dx \leq\int_{y}^{\infty}\frac{\log^{p}(x+y)}{x(x+y)}dx 
  = \frac{1}{y}\int_{1}^{\infty}\frac{[\log(y)+\log(z+1)]^{p}}{z(z+1)}dz \leq C\frac{\log^{p}(y)}{y}. 
\end{eqnarray*}
Similarly we can show that $II_{4}\leq C\frac{\log^{p}(y)}{y}$. Altogether, this gives the result. 

To prove (\ref{eq:in}) we recursively apply (\ref{eq:i2}) to give 
\begin{eqnarray*}
\int_{\mathbb{R}^{p}}|\sinc(x_{1}+\ldots+x_{p})|\prod_{j=1}^{p}|\sinc(x_{j})|dx_{1}\ldots dx_{p} 
&\leq&
\int_{\mathbb{R}^{p-1}}|\ell_{1}(x_{1}+\ldots+x_{p-1})\prod_{j=1}^{p-1}\sinc(x_{j})|
dx_{1}\ldots dx_{p-1} \\
&\leq& \int_{\mathbb{R}}|\ell_{p-1}(x_{1})\sinc(x_{1})|dx_{1} = O(1), 
\end{eqnarray*}
thus we have the required the result. 

To bound (ii), without loss of generality we derive a bound over $\sum_{m=1}^{a}$, the bounds for 
$\sum_{-a}^{m}$ is identical. Using (\ref{eq:i1}) we have 
\begin{eqnarray*}
&&\sum_{m=1}^{a}\int_{-\infty}^{\infty}\frac{\sin^{2}(x)}{|x(x+m\pi)|}dx = 
 \sum_{m=1}^{a}\int_{-\infty}^{\infty}\frac{|\sin(x)\sin(x+m\pi)|}{|x(x+m\pi)|}dx \\
&\leq& \sum_{m=1}^{a} \ell_{1}(m\pi)  = C\sum_{m=1}^{a} \frac{\log(m\pi)}{m\pi} \\
&\leq& C\log(a\pi)\sum_{m=1}^{a} \frac{1}{m\pi} = C\log(a\pi)\log(a) \leq C\log^{2}a.
\end{eqnarray*}
Thus we have shown (ii). 

To prove (iii) we use (\ref{eq:i1}) to give 
\begin{eqnarray*}
&&\sum_{m=-a}^{a} \left(\int_{-\infty}^{\infty}\frac{\sin^{2}(x)}{|x(x+m\pi)|}dx\right)
\left(\int_{-\infty}^{\infty}\frac{\sin^{2}(y)}{|y(y+m\pi)|}dy\right) \\
&\leq& C\sum_{m=-a}^{a}\big(\frac{\log m}{m}\big)^{2}  \leq C\sum_{m=-\infty}^{\infty}\big(\frac{\log m}{m}\big)^{2} \leq C.  
\end{eqnarray*}

To prove (iv) we apply (\ref{eq:i1}) to each of the integrals this gives 
\begin{eqnarray*}
&&\sum_{m_{1},\ldots,m_{q-1}=-a}^{a}\int_{\mathbb{R}^{q}}\bigg|\prod_{j=1}^{q-1}\sinc(x_{j})\sinc(x_{j}+m_{j}\pi)\sinc(x_{q})\sinc(x_{q}+\pi\sum_{j=1}^{q-1}m_{j})\bigg|\prod_{j=1}^{q}dx_{j} \\
&\leq& 
\sum_{m_{1},\ldots,m_{q-1}=-a}^{a}\ell_{1}(m_{1}\pi)\ldots\ell_{1}(m_{q-1}\pi)\ell_{1}(m_{q-1}\pi)\ell_{1}(\pi\sum_{j=1}^{q-1}m_{j}) 
 \leq  C\log^{2(q-2)}a, 
\end{eqnarray*} 
thus we obtain the desired result. \hfill $\Box$.

\vspace{3mm}

\vspace{3mm}
The proofs of Theorem
\ref{lemma:mean-stat}, Theorems \ref{theorem:asymptotic} (ii), Theorem \ref{lemma:nonuniformdft} 
involve integrals of $\sinc(u)\sinc(u+m\pi)$, where $|m|\leq
a+|r_{1}|+|r_{2}|$ and that $a\rightarrow \infty$ as
$\lambda\rightarrow \infty$. In the following lemma we obtain bounds
for these integrals. We note that all these results involve an inner integral
difference of the form
\begin{eqnarray*}
\int_{-b}^{b}g(\omega)\bigg[
h\left(\omega+\frac{2u}{\lambda}\right) - 
h(\omega)\bigg]d\omega du\bigg|.
\end{eqnarray*}
By using the mean value theorem heuristically it is clear that
$O(\lambda^{-1})$ should come out of the integral however the $2u$ in
the numerator of $h\left(\omega+\frac{2u}{\lambda}\right)$ makes the
analysis quite involved. 

\begin{lemma}
\label{lemma:1star}
Suppose $h$ is a function which  is absolutely integrable 
and $|h^{\prime}(\omega)|\leq
\beta_{}(\omega)$ (where $\beta_{}$ is a montonically decreasing
function that is absolutely integrable),
$m\in \mathbb{R}$ and  $g(\omega)$ is a bounded function. 
$b$ can take any value, and the bounds given below are independent of $b$.
Then we have 
\begin{eqnarray}
\label{eq:lemma1star1}
&&\bigg| \int_{-\infty}^{\infty}\sinc(u)\sinc(u+m\pi)
\int_{-b}^{b}g(\omega)\bigg[
h\left(\omega+\frac{2u}{\lambda}\right) - 
h(\omega)\bigg]d\omega du\bigg| \nonumber\\  
&\leq& C\frac{\log(\lambda)+I(|m|> e)\log(|m|)}{\lambda},
\end{eqnarray}
where $C$ is a finite constant independent of $m$ and 
$b$. If $g(\omega)$ is a bounded function with a bounded first derivative, then we have 
\begin{eqnarray}
\label{eq:lemma1star1b}
&&\bigg| \int_{-\infty}^{\infty}\sinc(u)\sinc(u+m\pi)
\int_{-b}^{b}h(\omega)\bigg(
g(\omega+\frac{2u}{\lambda}) - 
g(\omega)\bigg)d\omega du\bigg|  \nonumber\\
&\leq& C\frac{\log(\lambda)+I(|m|> e)\log(|m|)}{\lambda}.
\end{eqnarray}
In the case of double integrals, we assume 
 $h(\cdot,\cdot)$ is such that $\int_{\mathbb{R}^{2}}|h(\omega_{1},\omega_{2})|d\omega_{1}d\omega_{2}<\infty$,
and $|\frac{\partial h(\omega_{1},\omega_{2})}{\partial \omega_{1}}|\leq
\beta_{}(\omega_{1})\beta_{}(\omega_{2})$,  $|\frac{\partial h(\omega_{1},\omega_{2})}{\partial \omega_{2}}|\leq
\beta_{}(\omega_{1})\beta_{}(\omega_{2})$ and $|\frac{\partial^{2}
  h(\omega_{1},\omega_{2})}{\partial \omega_{1}\partial \omega_{2}}|\leq
\beta_{}(\omega_{1})\beta_{}(\omega_{2})$
(where $\beta_{}$ is a montonically decreasing
function that is absolutely integrable)
and  $g(\cdot,\cdot)$ is a bounded function. Then if  $m_{1}\neq 0$ and
$m_{2}\neq 0$ and $m_{1},m_{2}\in \mathbb{Z}$  we have 
\begin{eqnarray}
\label{eq:lemma1star1c}
&&\bigg| \int_{-\infty}^{\infty}\int_{-\infty}^{\infty}\sinc(u_{1})\sinc(u_{1}+m_{1}\pi) \sinc(u_{2})\sinc(u_{2}+m_{2}\pi)
\int_{-b_{}}^{b_{}}\int_{-b}^{b}g(\omega_{1},\omega_{2})\times  \nonumber\\
&& \bigg[h\left(\omega_{1}+\frac{2u_{1}}{\lambda},\omega_{2}+\frac{2u_{2}}{\lambda}\right) -
h(\omega_{1},\omega_{2})\bigg]d\omega_{1}d\omega_{2} du_{1}du_{2}\bigg| 
\leq C\frac{\prod_{i=1}^{2}[\log(\lambda)+\log(|m_{i}|)]}{\lambda^{2}}.\qquad
\end{eqnarray}
and
\begin{eqnarray}
\label{eq:lemma1star1d}
&&\bigg| \int_{-\infty}^{\infty}\int_{-\infty}^{\infty}\sinc(u_{1})\sinc(u_{1}+m_{1}\pi) \sinc(u_{2})\sinc(u_{2}+m_{2}\pi)
\frac{1}{\lambda^{2}}\sum_{-b\lambda}^{b\lambda}\sum_{-b\lambda}^{b\lambda}g(\omega_{k_{1}},\omega_{k_{2}})\times  \nonumber\\
&& \bigg[h\left(\omega_{k_{1}}+\frac{2u_{1}}{\lambda},\omega_{k_2}+\frac{2u_{2}}{\lambda}\right) -
h(\omega_{k_1},\omega_{k_2})\bigg]du_{1}du_{2}\bigg| 
\leq C\frac{\prod_{i=1}^{2}[\log(\lambda)+\log(|m_{i}|)]}{\lambda^{2}}.\qquad
\end{eqnarray}
where $a\rightarrow \infty$ as $\lambda \rightarrow \infty$.
\end{lemma}
PROOF. To simplify the notation in the proof, we'll prove (\ref{eq:lemma1star1})
for  $m>0$ (the proof for $m\leq 0$ is identical). 

The proof is based on considering the cases that $|u|\leq \lambda$  and
$|u|>\lambda$ separately. For $|u|\leq\lambda$ we apply the mean value
theorem to the difference $h(\omega+\frac{2u}{\lambda}) - 
h(\omega)$ and for $|u|>\lambda$ we exploit that the integral 
$\int_{u>|\lambda|}|\sinc(u)\sinc(u+m\pi)|du$ decays as
$\lambda\rightarrow \infty$. We now make these argument precise. 
We start by partitioning the integral 
\begin{eqnarray}
\int_{-\infty}^{\infty}\sinc(u)\sinc(u+m\pi)
\int_{-b}^{b}g(\omega)\bigg(
h(\omega+\frac{2u}{\lambda}) - 
h(\omega)\bigg)d\omega du = I_{1} + I_{2},\label{eq:mmdecomp}
\end{eqnarray}
where 
\begin{eqnarray*}
I_{1} &=& \int_{|u|>\lambda}\sinc(u)\sinc(u+m\pi)
\int_{-b}^{b}g(\omega)\bigg(
h(\omega+\frac{2u}{\lambda}) - 
h(\omega)\bigg)d\omega du \\
I_{2} &=& \int_{-\lambda}^{\lambda}\sinc(u)\sinc(u+m\pi)
\int_{-b}^{b}g(\omega)\bigg(
h(\omega+\frac{2u}{\lambda}) - 
h(\omega)\bigg)d\omega du.
\end{eqnarray*}
We further partition the integral $I_{1} = I_{11} + I_{12}+I_{13}$, where
\begin{eqnarray*}
I_{11} & =&  \int_{\lambda}^{\infty}\sinc(u)\sinc(u+m\pi)
\int_{-b}^{b}g(\omega)\bigg(
h(\omega+\frac{2u}{\lambda}) - 
h(\omega)\bigg)d\omega du \\
I_{12} &=& \int_{-\lambda-m\pi}^{-\lambda}\sinc(u)\sinc(u+m\pi)
\int_{-b}^{b}g(\omega)\bigg(
h(\omega+\frac{2u}{\lambda}) - 
h(\omega)\bigg)d\omega du \\
I_{13} &=& \int_{-\infty}^{-\lambda-m\pi}\sinc(u)\sinc(u+m\pi)
\int_{-b}^{b}g(\omega)\bigg(
h(\omega+\frac{2u}{\lambda}) - 
h(\omega)\bigg)d\omega du
\end{eqnarray*}
and partition $I_{2}=I_{21}+I_{22}+I_{23}$, where  
\begin{eqnarray*}
I_{21} &=& \int_{-\lambda}^{\lambda}\sinc(u)\sinc(u+m\pi)
\int_{-\min(4|u|/\lambda,b)}^{\min(4|u|/\lambda,b)}g(\omega)\bigg(
h(\omega+\frac{2u}{\lambda}) - 
h(\omega)\bigg)d\omega du \\
I_{22} &=& \int_{-\lambda}^{\lambda}\sinc(u)\sinc(u+m\pi)
\int^{-\min(4|u|/\lambda,b)}_{-b}g(\omega)\bigg(
h(\omega+\frac{2u}{\lambda}) - 
h(\omega)\bigg)d\omega du \\
I_{23} &=& \int_{-\lambda}^{\lambda}\sinc(u)\sinc(u+m\pi)
\int_{\min(4|u|/\lambda,b)}^{b}g(\omega)\bigg(
h(\omega+\frac{2u}{\lambda}) - 
h(\omega)\bigg)d\omega du.
\end{eqnarray*}
We start by bounding $I_{1}$. Taking absolutes of $I_{11}$, and using
that $h(\omega)$ is absolutely integrable we have  
\begin{eqnarray*}
|I_{11}| &\leq&  2\Gamma\int_{\lambda}^{\infty}\frac{\sin^{2}(u)}{u(u+m\pi)}du,
\end{eqnarray*}
where $\Gamma = \sup_{u}|g(u)|\int_{0}^{\infty}|h(u)|du$. Since $m>0$, it is
straightforward to show that
$\int_{\lambda}^{\infty}\frac{\sin^{2}(u)}{u(u+m\pi)}du\leq
C\lambda^{-1}$, where $C$ is some finite constant. This implies 
$|I_{11}|\leq 2C\Gamma\lambda^{-1}$. Similarly it can be shown that 
$|I_{13}|\leq 2C\Gamma\lambda^{-1}$. To bound $I_{12}$ we note that 
\begin{eqnarray*}
|I_{12}| &\leq&
\frac{2\Gamma}{\lambda}\int_{-\lambda-m\pi}^{-\lambda}\frac{\sin^{2}(u)}{|u+m\pi|} du
 \leq
\frac{2\Gamma}{\lambda}\times\left\{
\begin{array}{cc}
\log\left(\frac{\lambda}{\lambda - m\pi}\right) & m\pi <\lambda \\
\log\lambda + \log(m\pi -\lambda) & m\pi >\lambda \\
\end{array}
\right.. 
\end{eqnarray*}
Thus, we have $|I_{12}| \leq 2\Gamma\lambda^{-1}\big[\log
\lambda + \log m \big]$.
Altogether, the bounds for $I_{11},I_{12},I_{13}$ give 
\begin{eqnarray*}
|I_{1}|\leq \frac{C(\log \lambda + \log m)}{\lambda}.
\end{eqnarray*}
To bound $I_{2}$ we apply the mean value theorem to $h(\omega+\frac{2u}{\lambda}) - 
h(\omega)=
\frac{2u}{\lambda}h^{\prime}(\omega+\zeta(\omega,u)\frac{2u}{\lambda})$,
where $0\leq \zeta(\omega,u)|\leq 1$. Substituting this into $I_{23}$ gives 
\begin{eqnarray*}
|I_{23}| &\leq & \frac{2}{\lambda}\int_{-\lambda}^{\lambda}\frac{\sin^{2}(u)}{|u+m\pi|}
\int_{\min(4|u|/\lambda,b)}^{b}\left|h^{\prime}(\omega+\zeta(\omega,u)\frac{2u}{\lambda})\right| 
d\omega du.
\end{eqnarray*}
Since the limits of the inner integral are greater than $4u/\lambda$,
and the derivative is bounded by $\beta(\omega)$, 
this means $|h^{\prime}(\omega+\zeta(\omega,u)\frac{2u}{\lambda})|
\leq \max[\beta(\omega),\beta(\omega+\frac{2u}{\lambda})]=
\beta(\omega)$. Altogether,
  this gives 
\begin{eqnarray*}
|I_{23}|&\leq&\frac{2}{\lambda}
\left( \int_{\min(4|u|/\lambda,b)}^{b}\beta(\omega) d\omega\right)
\int_{-\lambda}^{\lambda}\frac{\sin^{2}(u)}{|u+m\pi|}du 
 \leq \frac{2\Gamma\log (\lambda+m\pi)}{\lambda}.
\end{eqnarray*}
Using the same method we obtain $I_{22} \leq \frac{2\Gamma \log
  (\lambda+m\pi)}{\lambda}$. 
Finally, to bound $I_{21}$, we cannot bound 
$h^{\prime}(\omega+\zeta(\omega,u)\frac{2u}{\lambda})$
by a monotonic function since $\omega$ and 
$\omega+\frac{2u}{\lambda}$ can have different signs. 
Therefore we simply bound
$h^{\prime}(\omega+\zeta(\omega,u)\frac{2u}{\lambda})$ with a
constant, this gives 
\begin{eqnarray*}
|I_{21}|\leq
\frac{8C}{\lambda^{2}}\int_{-\lambda}^{\lambda}\frac{|u|\sin^{2}(u)}{|u+m\pi|}du
\leq \frac{16C}{\lambda}.
\end{eqnarray*}
Altogether, the bounds for $I_{21},I_{22},I_{23}$ give 
\begin{eqnarray*}
|I_{2}|\leq C\Gamma\frac{\log \lambda + \log (m)+\log(\lambda+m\pi)}{\lambda}.
\end{eqnarray*}
Finally, we recall that if $\lambda >2$ and $m\pi>2$, then
$(\lambda+m\pi)< \lambda m\pi$, thus $\log(\lambda + m\pi)\leq \log \lambda + \log(m\pi)$,
Therefore, we obtain (\ref{eq:lemma1star1}). 
The proof of (\ref{eq:lemma1star1b}) is similar, but avoids some of the awkward details 
that are required to prove (\ref{eq:lemma1star1}).

Now we prove (\ref{eq:lemma1star1c}).  We note that both $m_{1}\neq 0$ and
$m_{2}\neq 0$ \emph{and} $m_{1},m_{2}\in \mathbb{Z}$
(if one of these values were zero or non-integer valued only the slower bound given
in (\ref{eq:lemma1star1}) holds). Without loss of generality we will prove the result
for $m_{1}>0$ and $m_{2}>0$. We first note since $m_{1}\neq 0$ and
$m_{2}\neq 0$, by
orthogonality of the sinc function at integer shifts (see Lemma \ref{lemma:orthogonal})
we have 
\begin{eqnarray}
&&\int_{-\infty}^{\infty}\int_{-\infty}^{\infty}\sinc(u_{1})\sinc(u_{1}+m_{1}\pi) \sinc(u_{2})\sinc(u_{2}+m_{2}\pi)
\int_{-b_{}}^{b_{}}\int_{-b}^{b}g(\omega_{1},\omega_{2})\times  \nonumber\\
&& \bigg[h\left(\omega_{1}+\frac{2u_{1}}{\lambda},\omega_{2}+\frac{2u_{2}}{\lambda}\right) -
h(\omega_{1},\omega_{2})\bigg]d\omega_{1}d\omega_{2} du_{1}du_{2} \nonumber\\
&=& \int_{-\infty}^{\infty}\sinc(u_{2})\sinc(u_{2}+m_{2}\pi) \int_{-\infty}^{\infty}\sinc(u_{1})\sinc(u_{1}+m_{1}\pi)
\int_{-b_{}}^{b_{}}\int_{-b}^{b}g(\omega_{1},\omega_{2})\nonumber\\
&&\times h\left(\omega_{1}+\frac{2u_{1}}{\lambda},\omega_{2}+\frac{2u_{2}}{\lambda}\right) 
d\omega_{1}du_{1}d\omega_{2}du_{2}, \label{eq:51a}
\end{eqnarray}
since in the last line $h(\omega_{1},\omega_{2})$ comes outside the
integral over $u_{1}$ and $u_{2}$. 
To further the proof, we note that for $m_{2}\neq 0$ we have
\begin{eqnarray}
%&&\int_{-\infty}^{\infty}\sinc(u_{2})\sinc(u_{2}+m_{2}\pi)
%\int_{|u_{1}|\leq \lambda}^{}\sinc(u_{1})\sinc(u_{1}+m_{1}\pi)
%\int_{-b_{}}^{b_{}}\int_{-b}^{b}g(\omega_{1},\omega_{2})\nonumber\\
%&&\times h\left(\omega_{1},\omega_{2}+\frac{2u_{2}}{\lambda}\right)d\omega_{1}du_{1}d\omega_{2}du_{2}
%= 0 \label{eq:zero1}\\
&&\int_{-\infty}^{\infty}\sinc(u_{2})\sinc(u_{2}+m_{2}\pi)
\int_{|u_{1}|> \lambda}^{}\sinc(u_{1})\sinc(u_{1}+m_{1}\pi)
\int_{-b_{}}^{b_{}}\int_{-b}^{b}g(\omega_{1},\omega_{2})\nonumber\\
&&\times h\left(\omega_{1}+\frac{2u_{1}}{\lambda},\omega_{2}\right)d\omega_{1}du_{1}d\omega_{2}du_{2}
= 0. \label{eq:zero1}
\end{eqnarray}
We use this zero-equality at relevant parts in the proof. 

Now we subtract (\ref{eq:zero1}) from (\ref{eq:51a}). 
We use the same decomposition and notation used in (\ref{eq:mmdecomp})
to decompose the integral  over $u_{1}$ into 
$\int_{|u_{1}|>\lambda}+\int_{|u_{1}|\leq \lambda}$,  to give 
\begin{eqnarray*}
&=& \int_{-\infty}^{\infty}\sinc(u_{2})\sinc(u_{2}+m_{2}\pi) \bigg(\int_{-\infty}^{\infty}\sinc(u_{1})\sinc(u_{1}+m_{1}\pi)
\int_{-b_{}}^{b_{}}\int_{-b}^{b}g(\omega_{1},\omega_{2})\times  \nonumber\\
&& \bigg[h\left(\omega_{1}+\frac{2u_{1}}{\lambda},\omega_{2}+\frac{2u_{2}}{\lambda}\right) -
h\left(\omega_{1} + \frac{2u_{1}}{\lambda},\omega_{2}\right)\bigg]d\omega_{1}du_{1}\bigg)d\omega_{2}du_{2}
\\
&=& \int_{-\infty}^{\infty}\sinc(u_{2})\sinc(u_{2}+m_{2}\pi) \int_{-b_{}}^{b_{}}\bigg(\int_{-\infty}^{\infty}\sinc(u_{1})\sinc(u_{1}+m_{1}\pi)
\int_{-b}^{b}g(\omega_{1},\omega_{2})\times  \nonumber\\
&& \bigg[h\left(\omega_{1}+\frac{2u_{1}}{\lambda},\omega_{2}+\frac{2u_{2}}{\lambda}\right) -
h\left(\omega_{1}+\frac{2u_{1}}{\lambda},\omega_{2}\right)\bigg]d\omega_{1}du_{1}\bigg)d\omega_{2}du_{2}
\\
&=&
\int_{-\infty}^{\infty}\sinc(u_{2})\sinc(u_{2}+m_{2}\pi)\int_{-b}^{b}\left[I_{1}(u_{2},\omega_{2})
  + I_{2}(u_{2},\omega_{2}) 
  \right]du_{2}d\omega_{2} = J_{1} + J_{2} ,
\end{eqnarray*}
where 
\begin{eqnarray*}
I_{1}(u_{2},\omega_{2})
&=&\int_{|u_{1}|>\lambda}\sinc(u_{1})\sinc(u_{1}+m_{1}\pi)\times \\
&&\int_{-b_{}}^{b_{}}g(\omega_{1},\omega_{2})
\left[h\left(\omega_{1}+\frac{2u_{1}}{\lambda},\omega_{2}+\frac{2u_{2}}{\lambda}\right)
- h\left(\omega_{1}+\frac{2u_{1}}{\lambda},\omega_{2}\right)\right]
d\omega_{1}du_{1} \\
I_{2}(u_{2},\omega_{2}) 
&=&\int_{|u_{1}|\leq\lambda}\sinc(u_{1})\sinc(u_{1}+m_{1}\pi)\int_{-b_{}}^{b_{}}g(\omega_{1},\omega_{2}) \frac{\partial h}{\partial 
\omega_{2}}\left(\omega_{1} +
\frac{2u_{1}}{\lambda},\omega_{2}+\zeta_{1}(\omega_{2},u_{2})
\frac{2u_{2}}{\lambda}\right)d\omega_{1}du_{1}
\end{eqnarray*}
and $|\zeta_{1}(\omega_{1},u_{1})|\leq 1$. Note that the expression
for $I_{2}(u_{2},\omega_{2})$ applies the mean value theorem to the
difference $h\left(\omega_{1}+\frac{2u_{1}}{\lambda},\omega_{2}+\frac{2u_{2}}{\lambda}\right)
- h\left(\omega_{1}+\frac{2u_{1}}{\lambda},\omega_{2}\right)=\frac{\partial h}{\partial 
\omega_{2}}\left(\omega_{1} +
\frac{2u_{1}}{\lambda},\omega_{2}+\zeta_{1}(\omega_{2},u_{2})
\frac{2u_{2}}{\lambda}\right)$; there 
is a slight abuse of  notation since the derivative is applied to
$h$ before evaluating it at $(\omega_{1} +
\frac{2u_{1}}{\lambda},\omega_{2}+\zeta_{1}(\omega_{2},u_{2})
\frac{2u_{2}}{\lambda})$. 
Next to bound $J_{1}$ we
decompose  the outer integral over $u_{2}$ into $\int_{|u_{2}|>\lambda}$ and
$\int_{|u_{2}|\leq \lambda}$ to give 
$J_{1} = J_{11} + J_{12}$, where 
\begin{eqnarray*}
J_{11} &=&
\int_{|u_{1}|> \lambda}\int_{|u_{2}|> \lambda}\sinc(u_{1})\sinc(u_{1}+m_{1}\pi)\sinc(u_{2})\sinc(u_{2}+m_{2}\pi)\times \\
&&\int_{-b_{}}^{b_{}}\int_{-b_{}}^{b_{}}g(\omega_{1},\omega_{2})
\left[h\left(\omega_{1}+\frac{2u_{1}}{\lambda},\omega_{2}+\frac{2u_{2}}{\lambda}\right)
- h\left(\omega_{1}+\frac{2u_{1}}{\lambda},\omega_{2}\right)\right]
d\omega_{1}du_{1} 
du_{2}d\omega_{2}\\
\\
J_{12} &=&
\int_{|u_{1}|>  \lambda}\int_{-b}^{b}\sinc(u_{1})\sinc(u_{1}+m_{1}\pi)\times\\
&&
\int_{|u_{2}|\leq\lambda}\sinc(u_{2})\sinc(u_{2}+m_{2}\pi)
\int_{-b_{}}^{b_{}}g(\omega_{1},\omega_{2})\frac{\partial h}{\partial 
\omega_{2}}\left(\omega_{1} + \frac{2u_{1}}{\lambda},\omega_{2} + \zeta_{2}(\omega_{2},u_{2})\frac{2u_{2}}{\lambda}\right)
d\omega_{1}du_{1}
du_{2}d\omega_{2},
\end{eqnarray*}
and $|\zeta_{2}(\omega_{2},u_{2})|<1$
 (applying the mean
value theorem to $h\left(\omega_{1}+\frac{2u_{1}}{\lambda},\omega_{2}+\frac{2u_{2}}{\lambda}\right)
- h\left(\omega_{1}+\frac{2u_{1}}{\lambda},\omega_{2}\right)$).
By using the same methods used to bound $I_{1}$ and $I_{2}$ in
(\ref{eq:mmdecomp}) we can show that 
\begin{eqnarray*}
|J_{11}|\textrm{ and }|J_{12}| = O\left(\frac{(\log(\lambda)+\log|m_{1}|)(\log(\lambda)+\log|m_{2}|)}{\lambda^{2}}\right).
\end{eqnarray*}
%$|J_{12}| = O([\log(\lambda)+\log|m_{2}|]/\lambda^{2})$  and  $|J_{21}| =O(\lambda^{-2})$. 
To bound
$J_{2}$ we again use that $m_{2}\neq 0$ and orthogonality of the sinc
function to subtract the  term $\frac{\partial }{\partial 
\omega_{2}}h\left(\omega_{1},\omega_{2} +\zeta_{1}(\omega_{2},u_{2})\frac{2u_{2}}{\lambda}\right)$
(whose total contribution is zero, see (\ref{eq:zero1})) from $J_{2}$ to give 
\begin{eqnarray*}
J_{2} &=& \int_{-\infty}^{\infty}\int_{-b_{}}^{b_{}}\sinc(u_{2})\sinc(u_{2}+m_{2}\pi)
\int_{|u_{1}|\leq\lambda}\sinc(u_{1})\sinc(u_{1}+m_{1}\pi)\times \\
&&\int_{-b_{}}^{b_{}}g(\omega_{1},\omega_{2}) \left[\frac{\partial h}{\partial 
\omega_{2}}\left(\omega_{1} +\frac{2u_{1}}{\lambda},\omega_{2}+
\zeta_{1}(\omega_{2},u_{2})\frac{2u_{2}}{\lambda}\right) - 
\frac{\partial h}{\partial 
\omega_{2}}\left(\omega_{1}, \omega_{2}+\zeta_{1}(\omega_{2},u_{2})\frac{2u_{2}}{\lambda}\right) 
\right]
d\omega_{1}d\omega_{2}du_{1}du_{2}.
\end{eqnarray*}
By decomposing the outer integral of $J_{2}$ over $u_{2}$ into $\int_{|u_{2}|>\lambda}$ and
$\int_{|u_{2}|\leq \lambda}$ we have 
$J_{2}= J_{21} + J_{22}$,
where 
\begin{eqnarray*}
J_{21} &=& \int_{|u_{2}|>\lambda}\int_{-b_{}}^{b_{}}\sinc(u_{2})\sinc(u_{2}+m_{2}\pi)
\int_{|u_{1}|\leq\lambda}\sinc(u_{1})\sinc(u_{1}+m_{1}\pi)\times \\
&&\int_{-b_{}}^{b_{}}g(\omega_{1},\omega_{2}) \left[\frac{\partial h}{\partial 
\omega_{2}}\left(\omega_{1}+\frac{2u_{1}}{\lambda},\omega_{2}+
\zeta_{1}(\omega_{2},u_{2})\frac{2u_{2}}{\lambda}\right) - 
\frac{\partial h}{\partial 
\omega_{2}}\left(\omega_{1},\omega_{2}+\zeta_{1}(\omega_{2},u_{2})\frac{2u_{2}}{\lambda}\right)
\right]
d\omega_{1}d\omega_{2}du_{1}du_{2}
\end{eqnarray*}
and 
\begin{eqnarray*}
J_{22} &=& \int_{|u_{2}|\leq \lambda}\int_{-b_{}}^{b_{}}\sinc(u_{2})\sinc(u_{2}+m_{2}\pi)
\int_{|u_{1}|\leq\lambda}\sinc(u_{1})\sinc(u_{1}+m_{1}\pi)\times \\
&&\int_{-b_{}}^{b_{}}g(\omega_{1},\omega_{2}) \frac{\partial^{2} h}{\partial 
\omega_{1}\partial \omega_{2}}\left(\omega_{1} +
\zeta_{3}(\omega_{1},u_{1})\frac{2u_{1}}{\lambda},\omega_{2}+
\zeta_{1}(\omega_{2},u_{2})\frac{2u_{2}}{\lambda}\right) 
d\omega_{1}d\omega_{2}du_{1}du_{2}
\end{eqnarray*}
with $|\zeta_{3}(\omega_{2},u_{2})|<1$. Note that the expression for
$J_{22}$ is obtained by applying the mean value theorem to 
$\frac{\partial h}{\partial 
\omega_{2}}\left(\omega_{1} +
\frac{2u_{1}}{\lambda},\omega_{2}+\zeta_{1}(\omega_{2},u_{2})
\frac{2u_{2}}{\lambda}\right) - 
\frac{\partial h}{\partial 
\omega_{2}}\left(\omega_{1} ,\omega_{2}+\zeta_{1}(\omega_{2},u_{2})\frac{2u_{2}}{\lambda}\right)$.

Again by using the methods used to bound $I_{1}$ and $I_{2}$ in 
(\ref{eq:mmdecomp}) we can show that $|J_{21}| =
O(\prod_{i=1}^{2}[\log(\lambda)+\log|m_{i}|]/\lambda^{2})$ and 
$|J_{22}| =
O(\prod_{i=1}^{2}[\log(\lambda)+\log|m_{i}|]/\lambda^{2})$. Altogether
this proves (\ref{eq:lemma1star1c}).

The proof of (\ref{eq:lemma1star1d}) follows exactly the same method
used to prove (\ref{eq:lemma1star1c}) the only difference is that the
summand rather than the integral makes the notation more
cumbersome. For this reason we omit the details. 
\hfill $\Box$

\vspace{3mm}

The following result is used in to obtain expression for the fourth
order cumulant term (in the case that the spatial random field is not
Gaussian). It is used in the proofs of Theorems
\ref{theorem:nonGaussiannonUniform} and \ref{theorem:variance-nongaussian}. 

\begin{lemma}\label{lemma:cum4bound}
Suppose $h$ is a function which is absolutely integrable 
and $|h^{\prime}(\omega)|\leq
\beta(\omega)$ (where $\beta$ is a monotonically decreasing
function that is absolutely integrable), $m\in \mathbb{Z}$ and  $g(\omega)$ is a bounded
function. Then we have 
\begin{eqnarray}
\label{eq:cum4bound1}
&&\int_{\mathbb{R}^{3}}\sinc(u_{1}+u_{2}+u_{3}+m\pi)\sinc(u_{1})\sinc(u_{2})\sinc(u_{3})
\int_{-a/\lambda}^{a/\lambda}
g(\omega)\left[h\left(\frac{2u_{1}}{\lambda}-\omega\right) - h(-\omega)
\right]d\omega du_{1}du_{2}du_{3} \nonumber\\
&=& O\left(\frac{[\log(\lambda)+\log|m|]^{3}}{\lambda}\right). 
\end{eqnarray}
\begin{eqnarray}
\label{eq:cum4bound2}
&&\int_{\mathbb{R}^{2}}\sinc(u_{1}+u_{2}+m\pi)\sinc(u_{1})\sinc(u_{2})
\int_{-a/\lambda}^{a/\lambda}
g(\omega)\left[h\left(\frac{2u_{1}}{\lambda}-\omega\right) - h(-\omega)
\right]d\omega du_{1}du_{2} \nonumber\\
&=& O\left(\frac{[\log(\lambda)+\log |m|]^{3}}{\lambda}\right). 
\end{eqnarray}
\end{lemma}
PROOF. The proof of (\ref{eq:cum4bound1}) is very similar to the proof of Lemma
\ref{lemma:1star}. 
%Since $m<<\lambda$ to simplify notation we prove the result for $m=0$. 
We start by partitioning the integral over $u_{1}$ into
$\int_{|u_{1}|\leq \lambda}$ and $\int_{|u_{1}|>\lambda}$ 
\begin{eqnarray*}
&&\int_{\mathbb{R}^{3}}\sinc(u_{1}+u_{2}+u_{3}+m\pi)\sinc(u_{1})\sinc(u_{2}) \sinc(u_{3})\times\\
&&\int_{-a/\lambda}^{a/\lambda}
g(\omega)\left(h(\frac{2u_{1}}{\lambda}-\omega) - h(-\omega)
\right)d\omega du_{1}du_{2}du_{3}  = I_{1} + I_{2}, 
\end{eqnarray*}
where 
\begin{eqnarray*}
&&I_{1}= \\
&&\int_{|u_{1}|>\lambda}\int_{\mathbb{R}^{3}}\sinc(u_{1}+u_{2}+u_{3}+m\pi)\sinc(u_{1})\sinc(u_{2}) \sinc(u_{3})\\
&& \times
\int_{-a/\lambda}^{a/\lambda}
g(\omega)\left(h(\frac{2u_{1}}{\lambda}-\omega) - h(-\omega)
\right)d\omega du_{1}du_{2}du_{3}  \\
&&I_{2}=\\
&&\int_{|u_{1}|\leq\lambda}\int_{\mathbb{R}^{3}}\sinc(u_{1}+u_{2}+u_{3}+m\pi)\sinc(u_{1})\sinc(u_{2}) \sinc(u_{3})\\
&&\times\int_{-a/\lambda}^{a/\lambda}
g(\omega)\left(h(\frac{2u_{1}}{\lambda}-\omega) - h(-\omega)
\right)d\omega du_{1}du_{2}du_{3}. 
\end{eqnarray*}
Taking absolutes of $I_{1}$ and using Lemma \ref{lemma:1a}, equations (\ref{eq:i1}) and (\ref{eq:i2})
we have
\begin{eqnarray*}
|I_{1}| \leq
4\Gamma\int_{|u_{1}|>\lambda}|\sinc(u_{1}+m\pi)|\ell_{2}(u_{1})du_{1} \leq  
C\int_{|u|>\lambda}\frac{\log^{2}(u)}{|u|}\times
\frac{|\sinc(u+m\pi)|}{|u+m\pi|}du,
\end{eqnarray*}
where $\Gamma = \sup_{\omega}|g(\omega)|\int_{0}^{\infty}|h(\omega)|d\omega$ and $C$ is a
finite constant (which has absorbed $\Gamma$). Decomposing the above
integral we have 
\begin{eqnarray*}
|I_{1}| \leq C\int_{|u|>\lambda}\frac{\log^{2}(u)}{|u|}\times
\frac{|\sinc(u+m\pi)|}{|u+m\pi|}du  = I_{11} + I_{12},
\end{eqnarray*}
where 
\begin{eqnarray*}
I_{11} &=& \int_{\lambda <|u| \leq \lambda(1+|m|)} \frac{\log^{2}(u)}{|u|}\times
\frac{|\sinc(u_{1}+m\pi)|}{|u_{1}+m\pi|}du \\
I_{12} &=& \int_{|u|> \lambda(1+|m|)} \frac{\log^{2}(u)}{|u|}\times
\frac{|\sinc(u_{1}+m\pi)|}{|u_{1}+m\pi|}du.  
\end{eqnarray*}
We first bound $I_{11}$ 
\begin{eqnarray*}
I_{11} &\leq& \frac{\log^{2}[\lambda(1+|m|\pi)]}{\lambda}\int_{\lambda <|u| \leq \lambda(1+|m|)} 
\frac{|\sinc(u_{1}+m\pi)|}{|u_{1}+m\pi|}du \\
&\leq&
\frac{2C\log^{2}[\lambda(1+|m|\pi)]}{\lambda}\times \log(\lambda+m\pi)
= C\frac{(\log|m| + \log \lambda)^{3}}{\lambda}
\end{eqnarray*}
To bound $I_{12}$ we make a change of variables $u=\lambda z$, the above becomes
\begin{eqnarray*}
I_{12}\leq
\frac{C}{\lambda}\int_{|z|>1+|m|\lambda}\frac{[\log\lambda +
  \log z]^{2}}{z(z+\frac{m}{\lambda})}dz
 = O\left(\frac{\log^{2}(\lambda)}{\lambda}\right). 
\end{eqnarray*}
Altogether the bounds for $I_{11}$ and $I_{12}$ give 
$|I_{1}| \leq C(\log|m| + \log \lambda)^{3}/\lambda$. 

To bound $I_{2}$, just as in Lemma \ref{lemma:1star}, equation (\ref{eq:lemma1star1}), we decompose it into
three parts $I_{2} = I_{21}+I_{22}+I_{23}$, where using Lemma \ref{lemma:1a}, equations (\ref{eq:i1}) and (\ref{eq:i2})
we have the bounds
\begin{eqnarray*}
&&|I_{21}|\leq \int_{|u_{}|\leq\lambda}|\sinc(u_{} + m\pi)|\ell_{2}(u_{})
\int_{-\min(a,4|u|)/\lambda}^{\min(a,4u)/\lambda}
|g(\omega)|\left|h\left(\frac{2u_{}}{\lambda}-\omega\right) - h(-\omega)
\right|d\omega du \\
&&|I_{22}|\leq\int_{|u_{}|\leq\lambda}|\sinc(u_{} + m\pi)|\ell_{2}(u_{})
\int^{-\min(a,4|u|)/\lambda}_{-a/\lambda}
|g(\omega)|\left|h\left(\frac{2u_{}}{\lambda}-\omega\right) - h(-\omega)
\right|d\omega du \\
&&|I_{23}|\leq\int_{|u_{}|\leq\lambda}|\sinc(u_{}+m\pi)|\ell_{2}(u_{})
\int_{\min(a,4|u|)/\lambda}^{a/\lambda}
|g(\omega)|\left|h\left(\frac{2u_{}}{\lambda}-\omega\right) - h(-\omega)
\right|d\omega du. 
\end{eqnarray*}
Using the same method used to bound $|I_{21}|,|I_{22}|, |I_{23}|$ in Lemma \ref{lemma:1star}, 
we have $|I_{21}|,|I_{22}|,|I_{23}|\leq C[\log(\lambda)+\log(|m|)]^{3}/\lambda$.
Having bounded all partitions of the integral, we have the result. 

The proof of (\ref{eq:cum4bound2}) is identical and we omit the
details. 
\hfill $\Box$

\subsection{Lemmas required to prove Lemma
  \ref{lemma:var-asymp} and 
Theorem \ref{theorem:asymptotic}}\label{appendix:lemma-asymptotic}

In this section we give the proofs of the three results used in Lemma
\ref{lemma:var-asymp} (which in turn proves Theorem
\ref{theorem:asymptotic}). 
\begin{lemma}\label{lemma:A8}
Suppose Assumptions \ref{assum:G}(ii) and  \ref{assum:GG}(b,c) 
holds. Then for $r_{1},r_{2}\in \mathbb{Z}$ we have 
\begin{eqnarray*}
|A_{1}(r_{1},r_{2})-B_{1}(r_{1}-r_{2};r_{1})| = O\left(\frac{\log^{2}(a)}{\lambda}\right)
\end{eqnarray*}
\end{lemma}
PROOF. 
To obtain a bound for the difference we use Lemma
\ref{lemma:sum-integral}(ii) to give 
\begin{eqnarray*}
&& |A_{1}(r_{1},r_{2})-B_{1}(r_{1}-r_{2};r_{1})|\\
&=& \int_{-\infty}^{\infty}\int_{-\infty}^{\infty}\sum_{m=-2a}^{2a}
|\sinc(u)\sinc(u-m\pi)\sinc(v)\sinc(v+(m+r_{1}-r_{2})\pi)| \\
 && \underbrace{\bigg|H_{m,\lambda}(\frac{2u}{\lambda},\frac{2v}{\lambda};r_{1})  - 
H_{m}(\frac{2u}{\lambda},\frac{2v}{\lambda};r_{1})\bigg|}_{\leq C/\lambda}
dudv \\
&\leq& \frac{C}{\lambda}\sum_{m=-2a}^{2a}\int_{-\infty}^{\infty}
|\sinc(u)\sinc(u-m\pi)|du \underbrace{\int_{-\infty}^{\infty}
|\sinc(v)\sinc(v+(m+r_{1}-r_{2})\pi)|dv}_{<\infty} \\
&\leq & \frac{C}{\lambda}\sum_{m=-2a}^{2a}\int_{-\infty}^{\infty}
|\sinc(u)\sinc(u-m\pi)|du \textrm{ (from Lemma }\ref{lemma:1a}(ii))\\
&=& O(\frac{\log^{2}a}{\lambda}),
\end{eqnarray*}
thus giving the desired result. \hfill $\Box$

\vspace{3mm}

\begin{lemma}\label{lemma:A9}
Suppose Assumptions \ref{assum:G}(ii) and  \ref{assum:GG}(b,c) 
holds (with $0\leq |r_{1}|, |r_{2}|< C|a|$). Then we have 
\begin{eqnarray*}
|B_{1}(s;r)-C_{1}(s;r)| = 
O\bigg(\log^{2}(a)\bigg[
\frac{\log a+\log \lambda}{\lambda}\bigg]\bigg). 
\end{eqnarray*}
Note, if we relax the assumption on $r_{1},r_{2}$ to $r_{1},r_{2}\in
\mathbb{Z}$ then the above bound requires the additional term 
$\log^{2}(a)[I(r_{1}\neq 0)\log|r_{1}| + I(r_{2}\leq 0)\log|r_{2}|]/\lambda$.
\end{lemma}
PROOF. Taking differences, it is 
easily seen that
\begin{eqnarray*}
&&B_{1}(s,r) -C_{1}(s,r) \\
 &=&\int_{\mathbb{R}^{2}}\sum_{m=-2a}^{2a}\sinc(u)\sinc(u-m\pi)
\sinc(v)\sinc(v+(m+s)\pi)\\
&&\times\frac{1}{(2\pi)}\int_{2\pi\max(-a,-a+m)/\lambda}^{2\pi\min(a,a+m)/\lambda}
g(\omega)\overline{g(\omega+\omega_{m})}
\left[ f(\omega-\frac{2u}{\lambda})f(\omega+\frac{2v}{\lambda}+ \omega_{r}) - 
f(\omega)f(\omega + \omega_{r})\right]d\omega dv du\\
 &=& I_{1} + I_{2} 
\end{eqnarray*}
where 
\begin{eqnarray*}
I_{1} &=& \int_{\mathbb{R}^{2}}\sum_{m=-2a}^{2a}\sinc(u)\sinc(u-m\pi)
\sinc(v)\sinc(v+(m+s)\pi) \\
&&\times \frac{1}{(2\pi)}\int_{2\pi\max(-a,-a+m)/\lambda}^{2\pi\min(a,a+m)/\lambda}
g(\omega)
\overline{g(\omega+\omega_{m})}f(\omega+\frac{2v}{\lambda}+ \omega_{r})
\left[ f(\omega-\frac{2u}{\lambda}) -f(\omega)\right]
d\omega dv du \\
  &=& \int_{\mathbb{R}}\sum_{m=-2a}^{2a}
\sinc(v)\sinc(v+(m+s)\pi)D_{m}(v)dv 
\end{eqnarray*}
and 
\begin{eqnarray*}
I_{2} &=& \int_{\mathbb{R}^{2}}\sum_{m=-2a}^{2a}\sinc(u)\sinc(u-m\pi)
\sinc(v)\sinc(v+(m+s)\pi)\\
&& \times \frac{1}{(2\pi)}\int_{2\pi\max(-a,-a+m)/\lambda}^{2\pi\min(a,a+m)/\lambda}
g(\omega)
\overline{g(\omega+\omega_{m})}f(\omega)
\left[ f(\omega+\frac{2v}{\lambda}+ \omega_{r})  -f(\omega+\omega_{r})\right]
d\omega dv du \\
 &=& \sum_{m=-2a}^{2a}d_{m}\int_{\mathbb{R}}\sinc(u)\sinc(u-m\pi)du
\end{eqnarray*}
with 
\begin{eqnarray*}
&&D_{m}(v) = \\
&&\int_{\mathbb{R}^{}}\sinc(u)\sinc(u-m\pi)\frac{1}{(2\pi)}\int_{2\pi\max(-a,-a+m)/\lambda}^{2\pi\min(a,a+m)/\lambda}
g(\omega)
\overline{g(\omega+\omega_{m})}f(\omega+\frac{2v}{\lambda}+ \omega_{r})\\
&&\times\left[ f(\omega-\frac{2u}{\lambda}) -f(\omega)\right]
d\omega du 
\end{eqnarray*}
and 
\begin{eqnarray*}
&& d_{m} =  \\
&&\int_{\mathbb{R}^{}}
\sinc(v)\sinc(v+(m+s)\pi)\frac{1}{(2\pi)}\int_{2\pi\max(-a,-a+m)/\lambda}^{2\pi\min(a,a+m)/\lambda}
g(\omega)
\overline{g(\omega+\omega_{m})}f(\omega)\\
&&\times\left[ f(\omega+\frac{2v}{\lambda}+ \omega_{r})  -f(\omega+\omega_{r})\right]
d\omega dv.
\end{eqnarray*}
Since the functions $f(\cdot)$ and $g(\cdot)$ satisfy the conditions stated in  
Lemma \ref{lemma:1star}, the lemma can be used to show that
\begin{eqnarray*}
\max_{|m|\leq a}\sup_{v}|D_{m}(v)|\leq C\bigg( \frac{\log\lambda+\log a}{\lambda}\bigg)  
\end{eqnarray*}
and 
\begin{eqnarray*}
\max_{|m|\leq a}|d_{m}|\leq C\bigg(\frac{\log\lambda+\log a}{\lambda}\bigg).
\end{eqnarray*}
Substituting these bounds into $I_{1}$ and $I_{2}$ give 
\begin{eqnarray*}
|I_{1}|&\leq& C\bigg(\frac{\log\lambda+\log a}{\lambda}\bigg)
\sum_{m=-2a}^{2a}\int_{-\infty}^{\infty}
|\sinc(v)\sinc(v+(m+s)\pi)|dv \\
|I_{2}|&\leq & C\bigg(\frac{\log\lambda+\log a}{\lambda}\bigg)
\sum_{m=-2a}^{2a}\int_{-\infty}^{\infty}|\sinc(u)\sinc(u-m\pi)|du.
\end{eqnarray*}
Therefore, by using Lemma \ref{lemma:1a}(ii) we have 
\begin{eqnarray*}
|I_{1}|\textrm{ and }|I_{2}|&=& O\bigg(\log^{2}(a)
\frac{\log a+\log \lambda}{\lambda}\bigg).
\end{eqnarray*}
Since $|B_{1}(s;r)-C_{1}(s;r)|\leq |I_{1}|+|I_{2}|$ this gives the desired result. 
\hfill $\Box$

%% file: B_2_cumulant.tex
\section{Approximations to the covariance and cumulants of 
$\widetilde{Q}_{a,\lambda}(g;\rb)$}\label{sec:cumulants}

In this section, our objective is to obtain bounds for 
$\cum_{q}\big(\widetilde{Q}_{a,\lambda}(g;\rb_{1}),\ldots,
\widetilde{Q}_{a,\lambda}(g;\rb_{q})\big)$, 
these results will be used to prove the asymptotic expression for the
variance of $\widetilde{Q}_{a,\lambda}(g;\rb)$ (given in Section \ref{sec:variance}) and asymptotic normality of
$\widetilde{Q}_{a,\lambda}(g;\rb)$. \citeA{p:fox-taq-87},
\citeA{p:dah-89}, \citeA{p:gir-sur-90} (see also \citeA{b:taq-11}) have
developed techniques for dealing with the cumulants of sums of
periodograms of Gaussian (discrete time) time series, and one would
have expected that these results could be used here. However, 
in our setting there are a few differences that we now describe
(i) despite the spatial random being Gaussian the locations
are randomly sampled, thus the composite process $Z(\ub)$ is not
Gaussian (we can only exploit the Gaussianity when we condition on
the location) (ii) the random field is defined over $\mathbb{R}^{d}$
(not $\mathbb{Z}^{d}$) (iii) the number of terms in the sums 
$\widetilde{Q}_{a,\lambda}(\cdot)$ is not necessarily the sample
size. Unfortunately, these differences make it difficult to apply the
above mentioned results to our setting. Therefore, in this section we
consider cumulant based results for spatial data observed at irregular
locations. In order to reduce cumbersome notation we focus on the case
that the locations are from a uniform distribution. 

As a simple motivation we first consider
$\var[\widetilde{Q}_{a,\lambda}(1,0)]$.
 By using indecomposable partitions we have 
\begin{eqnarray}
&&\var[\widetilde{Q}_{a,\lambda}(1,0)] \nonumber\\
&=&
\frac{1}{n^{4}}\sum_{\stackrel{j_{1},j_{2},j_{3},j_{4}=1}{j_{1}\neq
  j_{2},j_{3}\neq j_{4}}}^{n}\sum_{k_{1},k_{2}=-a}^{a}
\cov\left[Z(s_{j_{1}})Z(s_{j_{2}})\exp(i\omega_{k_{1}}(s_{j_{1}}-s_{j_{2}})),
Z(s_{j_{3}})Z(s_{j_{4}})\exp(i\omega_{k_{2}}(s_{j_{3}}-s_{j_{4}}))
\right]\nonumber\\
&=&
\frac{1}{n^{4}}\sum_{\stackrel{j_{1},j_{2},j_{3},j_{4}=1}{r_{1}\neq
    r_{2},r_{3}\neq r_{4}}}^{n}\sum_{k_{1},k_{2}=-a}^{a}\bigg(
\cum\big[Z(s_{j_1})e^{is_{j_1}\omega_{k_{1}}},Z(s_{j_3})e^{-is_{j_3}\omega_{k_{2}}}\big]
\cum\big[Z(s_{j_2})e^{-is_{j_2}\omega_{k_{1}}},Z(s_{j_4})e^{is_{j_4}\omega_{k_{2}}}
\big] \nonumber\\
&& +\cum\big[Z(s_{j_1})e^{is_{j_1}\omega_{k_{1}}},Z(s_{j_4})e^{is_{j_4}\omega_{k_{2}}}\big]
\cum\big[Z(s_{j_2})e^{-is_{j_2}\omega_{k_{1}}},Z(s_{j_3})e^{-is_{j_3}\omega_{k_{2}}}\big]\nonumber\\
&&  +\cum\big[Z(s_{j_1})e^{is_{j_1}\omega_{k_{1}}},
Z(s_{j_2})e^{-is_{j_2}\omega_{k_{1}}},Z(s_{j_3})e^{-is_{j_3}\omega_{k_{2}}},Z(s_{j_4})e^{is_{j_4}\omega_{k_{2}}}\big]\bigg).
\label{eq:Qalambda}
\end{eqnarray}
In order to evaluate the covariances in the above we condition
on the locations $\{s_{j}\}$. To evaluate the fourth order cumulant of the
above we appeal to a generalisation of the conditional variance
method. This expansion was first derived in \citeA{p:bri-69}, and in
the general setting it is stated as 
\begin{eqnarray}
\label{eq:conditional-cumulants}
\cum(Y_{1},Y_{2},\ldots,Y_{q}) = \sum_{\pi}\cum\left[\cum(Y_{\pi_{1}}|s_{1},\ldots,s_{q}),\ldots,\cum(Y_{\pi_{b}}|s_{1},\ldots,s_{q})\right],
\end{eqnarray}
where the sum is over all partitions $\pi$ of $\{1,\ldots,q\}$ and 
$\{\pi_{1},\ldots,\pi_{b}\}$ are all the blocks in the partition
$\pi$. We use (\ref{eq:conditional-cumulants}) to evaluate 
$\cum[Z(s_{j_{1}})e^{is_{j_{1}}\omega_{k_{1}}},\ldots,Z(s_{j_{q}})e^{is_{j_{q}}\omega_{k_{q}}}]$,
where $Y_{i} = Z(s_{j_{i}})e^{is_{j_{i}}\omega_{k_{i}}}$ and we condition on
the locations $\{s_{j}\}$.  Using this decomposition we observe that
because the spatial process is Gaussian, 
$\cum[Z(s_{j_{1}})e^{is_{j_{1}}\omega_{k_{1}}},\ldots,Z(s_{j_{q}})e^{is_{j_{q}}\omega_{k_{q}}}]$ can only be composed of
cumulants of covariances conditioned on the locations. Moreover, if $s_{1},\ldots,s_{q}$
are independent then by using the same reasoning we see that 
$\cum[Z(s_{1})e^{is_{1}\omega_{k_{1}}},\ldots,Z(s_{q})e^{is_{q}\omega_{k_{q}}}]=0$.
Therefore,  $\cum\big[Z(s_{j_1})e^{is_{j_1}\omega_{k_{1}}},
Z(s_{j_2})e^{-is_{j_2}\omega_{k_{1}}},Z(s_{j_3})e^{is_{j_3}\omega_{k_{2}}},Z(s_{j_4})e^{-is_{j_4}(\omega_{k_{2}}}\big]$
will only be non-zero if some elements of $s_{j_{1}},s_{j_{2}},s_{j_{3}},s_{j_{4}}$
are dependent. Using these rules we have 
\begin{eqnarray}
&&\var[\widetilde{Q}_{a,\lambda}(1,0)] \nonumber\\
&=& \frac{1}{n^{4}}\sum_{j_{1},j_{2},j_{3},j_{4}\in \mathcal{D}_{4}}\sum_{k_{1},k_{2}=-a}^{a}\bigg(
\cum\big[Z(s_{j_1})e^{is_{j_1}\omega_{k_{1}}},Z(s_{j_3})e^{-is_{j_3}\omega_{k_{2}}}\big]
\cum\big[Z(s_{j_2})e^{-is_{j_2}\omega_{k_{1}}},Z(s_{j_4})e^{is_{j_4}\omega_{k_{2}}}
\big] \nonumber\\
&& +\cum\big[Z(s_{j_1})e^{is_{j_1}\omega_{k_{1}}},Z(s_{j_4})e^{is_{j_4}\omega_{k_{2}}}\big]
\cum\big[Z(s_{j_2})e^{-is_{j_{2}}\omega_{k_{1}}},Z(s_{j_3})e^{-is_{j_3}\omega_{k_{2}}}\big]\bigg)\nonumber\\
&&  + \frac{1}{n^{4}}\sum_{j_{1},j_{2},j_{3},j_{4}\in \mathcal{D}_{3}}\sum_{k_{1},k_{2}=-a}^{a}\bigg(
\cum\big[Z(s_{j_1})e^{is_{j_1}\omega_{k_{1}}},Z(s_{j_3})e^{-is_{j_3}\omega_{k_{2}}}\big]
\cum\big[Z(s_{j_2})e^{-is_{j_2}\omega_{k_{1}}},Z(s_{j_4})e^{is_{j_4}\omega_{k_{2}}}
\big] \nonumber\\
&& +\cum\big[Z(s_{j_1})e^{is_{j_1}\omega_{k_{1}}},Z(s_{j_4})e^{is_{j_4}\omega_{k_{2}}}\big]
\cum\big[Z(s_{j_2})e^{-is_{j_{2}}\omega_{k_{1}}},Z(s_{j_3})e^{-is_{j_3}\omega_{k_{2}}}\big]\bigg)\nonumber\\
&& +\frac{1}{n^{4}}\sum_{j_{1},j_{2},j_{3},j_{4}\in \mathcal{D}_{3}}\sum_{k_{1},k_{2}=-a}^{a}\cum\big[Z(s_{j_1})e^{is_{j_1}\omega_{k_{1}}},
Z(s_{j_2})e^{-is_{j_2}\omega_{k_{1}}},Z(s_{j_3})e^{is_{j_3}\omega_{k_{2}}},
Z(s_{j_4})e^{-is_{j_4}\omega_{k_{2}}}\big], \label{eq:varQ}
\end{eqnarray}
where $\mathcal{D}_{4} = \{j_{1},\ldots,j_{4}=\textrm{all $j$s are
  different}\}$, $\mathcal{D}_{3}=\{j_{1},\ldots,j_{4};\textrm{two $j$s
are the same but }j_{1}\neq  j_{2} \textrm{ and }j_{3}\neq j_{4}\}$ (noting that by definition of $\widetilde{Q}_{a,\lambda}(1,0)$ more
than two elements in $\{j_{1},\ldots,j_{4}\}$ cannot be the same).
We observe that $|\mathcal{D}_{4}| = O(n^{4})$ and
$|\mathcal{D}_{3}|=O(n^{3})$,  where $|\cdot|$ denotes the cardinality
of a set. We will show that the second and third terms are
asymptotically negligible with respect to the first term. To show this
we require the following lemma.

%Before we elaborate on further rules we prove the following result,
%which is used in the proof of Lemma \ref{lemma:covariance}. 
\begin{lemma}\label{lemma:cumulantsA}
Suppose Assumptions \ref{assum:S}, \ref{assum:uniform}, \ref{assum:G}
and \ref{assum:GG}(b,c) hold (note we only use Assumption
\ref{assum:GG}(c) to get `neater expressions in the proofs' it is not needed to obtain
the same order). Then we have 
%\begin{eqnarray}
%\label{eq:cum1}
%\sum_{k_{1},k_{2}=-n}^{n}\cum\bigg(c(s_{1}-s_{3})\exp(is_{1}\omega_{k_{1}}- is_{3}\omega_{k_{3}},
%c(s_{2}-s_{1})\exp(-is_{2}\omega_{k_{1}}+ is_{1}\omega_{k_{3}})\bigg) = O(1)
%\end{eqnarray}
\begin{eqnarray}
\label{eq:cum2}
\sup_{a}\sum_{k_{1},k_{2}=-a}^{a}\cum[Z(s_{1})e^{is_{1}\omega_{k_{1}}},Z(s_{2})e^{-is_{2}\omega_{k_{1}}},
Z(s_{3})e^{-is_{3}\omega_{k_{2}}},Z(s_{1})e^{is_{1}\omega_{k_{2}}}] = O(1)
\end{eqnarray}
\begin{eqnarray}
\label{eq:cum3}
\sup_{a}\sum_{k_{1},k_{2}=-a}^{a}\cum\left[Z(s_{1})e^{is_{1}\omega_{k_{1}}},Z(s_{1})e^{is_{1}\omega_{k_{2}}}\right]
\cum\left[Z(s_{2})e^{-is_{2}\omega_{k_{1}}},Z(s_{3})e^{-is_{3}\omega_{k_{2}}}\right]= O(1).
\end{eqnarray}
%\begin{eqnarray}
%\label{eq:cum3}
%\cum(Z(s_{1})e^{is_{1}\omega_{k_{1}}},Z(s_{2})e^{-is_{2}\omega_{k_{1}}},
%Z(s_{3})e^{-is_{3}\omega_{k_{2}}},Z(s_{4})e^{is_{4}\omega_{k_{2}}}) = 0
%\end{eqnarray}
%\begin{eqnarray}
%\label{eq:cum4}
%\cum(Z(s)\exp(is\omega_{r_{1}}),Z(s),Z(s)\exp(-is\omega_{r_{2}})) = 
%\left\{
%\begin{array}{cc}
%0 & r_{1}\neq r_{2} \\
%2c(0)^{2} & \textrm{otherwise}
%\end{array}
%\right.
%\end{eqnarray}
\end{lemma}
PROOF. 
To show (\ref{eq:cum2}) we use conditional cumulants (see
(\ref{eq:conditional-cumulants})). By using the conditional cumulant expansion and 
Gaussianity of $Z(s)$ conditioned on the location we have 
\begin{eqnarray*}
&&\sum_{k_{1},k_{2}=-a}^{a}\cum[Z(s_{1})e^{is_{1}\omega_{k_{1}}},Z(s_{2})e^{-is_{2}\omega_{k_{1}}},
Z(s_{3})e^{-is_{3}\omega_{k_{2}}},Z(s_{1})e^{is_{1}\omega_{k_{2}}}] \\
&=& \sum_{k_{1},k_{2}=-a}^{a}\cum[c(s_{1}-s_{2})e^{is_{1}\omega_{k_{1}}-is_{2}\omega_{k_{1}}},
c(s_{1}-s_{3})e^{-is_{3}\omega_{k_{2}}+is_{1}\omega_{k_{2}}}] + \\
&& \sum_{k_{1},k_{2}=-a}^{a}\cum[c(s_{1}-s_{3})e^{is_{1}\omega_{k_{1}}-is_{3}\omega_{k_{2}}},
c(s_{1}-s_{2})e^{-is_{2}\omega_{k_{1}}+is_{1}\omega_{k_{2}}}] +\\
 &&  \underbrace{\sum_{k_{1},k_{2}=-a}^{a}\cum[c(0)e^{is_{1}(\omega_{k_{1}}+\omega_{k_{2}})},
c(s_{2}-s_{3})e^{-is_{2}\omega_{k_{1}}-is_{3}\omega_{k_{2}}}]}_{=0
\textrm{ (since $s_{1}$ is independent of $s_{2}$ and $s_{3}$)}}  = I_{1} + I_{2}.
\end{eqnarray*}
Writing $I_{1}$ in terms of expectations and using the
spectral representation of the covariance we have 
\begin{eqnarray*}
I_{1}&=&\sum_{k_{1},k_{2}=-a}^{a}\bigg( 
\Ex\big[c(s_{1}-s_{2})c(s_{1}-s_{3})
e^{is_{1}(\omega_{k_{1}}+\omega_{k_{2}})}
e^{-is_{3}\omega_{k_{2}}}e^{-is_{2}\omega_{k_{1}}}\big]  - 
\Ex\big[c(s_{1}-s_{2})e^{is_{1}\omega_{k_{1}}- is_{2}\omega_{k_{1}}}\big] \\
&& \times\Ex\big[c(s_{1}-s_{3})e^{is_{1}\omega_{k_{2}}- is_{3}\omega_{k_{2}}}\big]\bigg) \\
&=&\frac{1}{(2\pi)^{2}}\sum_{k_{1},k_{2}=-a}^{a}\int\int f(x)f(y)\sinc\bigg(\frac{\lambda}{2}(x+y)+(k_{1}+k_{2})\pi\bigg)
\sinc\bigg(\frac{\lambda}{2}x+k_{1}\pi\bigg)\sinc\bigg(\frac{\lambda}{2}y+k_{2}\pi\bigg)dxdy -\\
&& \frac{1}{(2\pi)^{2}}\sum_{k_{1},k_{2}=-a}^{a}\int\int f(x)f(y)\sinc\bigg(\frac{\lambda}{2}x+k_{2}\pi\bigg)
\sinc\bigg(\frac{\lambda}{2}x+k_{2}\pi\bigg)
\sinc\bigg(\frac{\lambda}{2}y+k_{1}\pi\bigg)
\sinc\bigg(\frac{\lambda}{2}y+k_{1}\pi\bigg)dxdy \\
&=& E_{1} - E_{2}.
\end{eqnarray*}
To bound $E_{1}$ we make a change of variables $u=\frac{\lambda
  x}{2}+k_{1}\pi$, $v=\frac{\lambda y}{2}+k_{2}\pi$, 
and replace sum with integral (and use Lemma \ref{lemma:sum-integral}) to give 
\begin{eqnarray*}
E_{1} &=& \frac{1}{\pi^{2}\lambda^{2}}\int \int \sum_{k_{1},k_{2}=-a}^{a}
f(\frac{2u}{\lambda} - \omega_{k_{1}}) 
f(\frac{2v}{\lambda} - \omega_{k_{2}})\sinc(u+v)\sinc(u)\sinc(v)dudv
\\
&=& \frac{1}{(2\pi)^{2}\pi^{2}}\int \int\sinc(u+v)\sinc(u)\sinc(v)
\bigg( \int_{-2\pi a/\lambda}^{2\pi a/\lambda}\int_{-2\pi
  a/\lambda}^{2\pi a/\lambda}
f(\frac{2u}{\lambda} - \omega_{1}) 
f(\frac{2v}{\lambda} - \omega_{2})d\omega_{1}d\omega_{2}
\bigg)dudv + O(\frac{1}{\lambda}).
\end{eqnarray*}
Let $G(\frac{2u}{\lambda}) = \frac{1}{2\pi}\int_{-2\pi a/\lambda}^{2\pi a/\lambda}
f(\frac{2u}{\lambda} - \omega) d\omega$, then substituting this into
the above and using equation (\ref{eq:in}) in 
Lemma \ref{lemma:1a} we have
\begin{eqnarray*}
E_{1}= \frac{1}{\pi^{2}}\int
\int\sinc(u+v)\sinc(u)\sinc(v)G\left(\frac{2u}{\lambda}\right)G\left(\frac{2v}{\lambda}\right)dudv
+ O(\frac{1}{\lambda}) = O(1).
\end{eqnarray*}
To bound $E_{2}$ we use a similar technique and Lemma \ref{lemma:1a}(iii) to give
$E_{2}=O(1)$. Altogether, this gives $I_{1}=O(1)$. The same proof can be used to 
show that $I_{2} = O(1)$. Altogether this gives (\ref{eq:cum2}).

To bound (\ref{eq:cum3}), we observe that if $k_{1}\neq -k_{2}$, then  
$\cum\big[Z(s_{1})e^{is_{1}\omega_{k_{1}}},Z(s_{1})e^{is_{1}\omega_{k_{2}}}\big]
=\Ex[c(0)e^{is(\omega_{k_{1}}+\omega_{k_{2}})}]=0$ otherwise
$\cum\big[Z(s_{1})e^{is_{1}\omega_{k}},Z(s_{1})e^{-is_{1}\omega_{k}}\big]=c(0)$. Using
this,  (\ref{eq:cum3}) can be reduced to   
\begin{eqnarray*}
&&\sum_{k_{1},k_{2}=-a}^{a}\cum\left[Z(s_{1})e^{is_{1}\omega_{k_{1}}},Z(s_{1})e^{is_{1}\omega_{k_{2}}}\right]
\cum\left[Z(s_{2})e^{-is_{2}\omega_{k_{1}}},Z(s_{3})e^{-is_{3}\omega_{k_{2}}}\right] \\
&=& c(0)\sum_{k=-a}^{a}
\cum\left[Z(s_{2})e^{is_{2}\omega_{k_{}}},Z(s_{3})e^{-is_{3}\omega_{k_{}}}\right] \\
&=& \frac{c(0)}{2\pi}\int_{-\infty}^{\infty}\sum_{k=-a}^{a}f(x)
\sinc\left(\frac{\lambda x }{2}+k\pi\right)\sinc\left(\frac{\lambda
    x}{2} + k\pi\right)dx \qquad \left(\textrm{change variables using }
\omega = \frac{\lambda x}{2}+k\pi\right) \\
 &=& c(0)\int_{-\infty}^{\infty}\frac{1}{\pi\lambda}\sum_{k=-a}^{a}
f\left(\frac{2\omega}{\lambda}-\omega_{k}\right)
\sinc^{2}(\omega)d\omega \\
 &=&  \frac{c(0)}{2\pi^{2}}\int_{-\infty}^{\infty}
\sinc^{2}(\omega)\bigg(\int_{-2\pi a/\lambda}^{2\pi a/\lambda}
f\left(\frac{2 \omega}{\lambda}-x\right)dx\bigg)d\omega + O\left(\frac{1}{\lambda}\right) = O(1),
\end{eqnarray*}
thus proving (\ref{eq:cum3}). \hfill $\Box$
%To prove (\ref{eq:cum3}) once again we apply the conditional cumulant
%to $\cum(Z(s_{1})e^{is_{1}\omega_{k_{1}}},Z(s_{2})e^{-is_{2}\omega_{k_{1}}},
%Z(s_{3})e^{-is_{3}\omega_{k_{2}}},Z(s_{4})e^{is_{4}\omega_{k_{2}}})$. However,
%since $s_{1},\ldots,s_{4}$ are independent the covariances involving 
%$s_{1},\ldots,s_{4}$ are zero, thus giving (\ref{eq:cum3}).

%To prove (\ref{eq:cum4}) we use the conditional cumulant expansion 
%\begin{eqnarray*}
%&&\cum(Z(s)e^{is\omega_{r_{1}}},Z(s),
%Z(s)e^{-is\omega_{r_{2}}},Z(s)) \\
%&=& 
%2\cum(c(0)e^{is\omega_{r_{1}}},c(0)e^{-is\omega_{r_{2}}}) +
%\cum(c(0)e^{is\omega_{r_{1}-r_{2}}},c(0)).
%\end{eqnarray*}
%It is immediately clear that the second term is zero. To analyse the
%first term of the above we write 
%$\cum(c(0)e^{is\omega_{r_{1}}},c(0)e^{-is\omega_{r_{2}}})=
%\Ex(c(0)^{2}e^{is\omega_{r_{1}-r_{2}}})-
%\Ex(c(0)e^{is\omega_{r_{1}}})\Ex(c(0)e^{-is\omega_{r_{2}}})$. It is
%clear that when $r_{1}\neq r_{2}$ then this term is zero. On the other
%hand for $r_{1}=r_{2}$ (and they are not zero) the above expansion gives 
%$\cum(Z(s)e^{is\omega_{r_{1}}},Z(s),
%Z(s)e^{-is\omega_{r_{1}}},Z(s)) = 2c(0)^{2}$. 
\vspace{3mm}

We now derive an expression for $\var[\widetilde{Q}_{a,\lambda}(1,0)]$, by using
Lemma \ref{lemma:cumulantsA} we have 
\begin{eqnarray}
&&\var[\widetilde{Q}_{a,\lambda}(1,0)] \nonumber\\
&=& \frac{1}{n^{4}}\sum_{j_{1},j_{2},j_{3},j_{4}\in \mathcal{D}_{4}}\sum_{k_{1},k_{2}=-a}^{a}\bigg(
\cum\big[Z(s_{j_1})e^{is_{j_1}\omega_{k_{1}}},Z(s_{j_3})e^{-is_{j_3}\omega_{k_{2}}}\big]
\cum\big[Z(s_{j_2})e^{-is_{j_2}\omega_{k_{1}}},Z(s_{j_4})e^{is_{j_4}\omega_{k_{2}}}
\big] \nonumber\\
&& +\cum\big[Z(s_{j_1})e^{is_{j_1}\omega_{k_{1}}},Z(s_{j_4})e^{is_{j_4}\omega_{k_{2}}}\big]
\cum\big[Z(s_{j_2})e^{-is_{j_2}\omega_{k_{1}}},Z(s_{j_3})e^{-is_{j_3}\omega_{k_{2}}}\big]\bigg) +O\left(\frac{1}{n}\right). 
\label{eq:Q10}
\end{eqnarray}
In  Lemma \ref{lemma:covariance} we have shown that the covariance terms above are of order $O(\lambda^{-1})$, 
thus dominating the fourth order cumulant terms which is of order $O(n^{-1})$ (so long as $\lambda<<n$). 

\begin{lemma}\label{lemma:basic}
Suppose that $\{Z(\ub);\ub\in \mathbb{R}^{d}\}$ is a Gaussian random process and 
$\{\ub_{j}\}$ are iid random variables. Then the following hold:
\begin{itemize}
\item[(i)]  As we mentioned above, by using 
that the spatial process is Gaussian and  (\ref{eq:conditional-cumulants}),
\\*
$\cum[Z(\ub_{j_{1}})e^{i\ub_{j_{1}}\ob_{\kb_{1}}},\ldots,Z(\ub_{j_{q}})e^{i\ub_{j_{q}}\ob_{\kb_{q}}}]$
can be written as the sum of products of 
cumulants of the spatial covariance conditioned on location. Therefore,  it is easily seen that the
odd order cumulant
\\*
  $\cum_{2q+1}[Z(\ub_{j_1})\exp(i\ub_{j_1}\ob_{\kb_{1}}),\ldots,Z(\ub_{j_{2q+1}}\exp(i\ub_{j_{2q+1}}
\ob_{\kb_{2q+1}}))]=0$ for all $q$ and regardless of $\{\ub_{j}\}$ being
dependent or not. 
\item[(ii)] By the same reasoning given above, 
we observe that if more than $(q+1)$ locations $\{s_{j};j=1,\ldots,2q\}$ are
  independent, then 
\\*
$\cum_{2q}[Z(\ub_{j_1})\exp(i\ub_{j_1}\ob_{\kb_{1}}),\ldots,Z(\ub_{j_{2q}})\exp(i\ub_{j_{2q}}
\ob_{\kb_{2q}}))]=0$. 
%\item[(iii)] As illustrated in the bounds for $E_{1}$ and $E_{2}$ in
 % the proof of (\ref{eq:cum2}), the sum of product of cumulant 
%terms is dominated by the corresponding product of expectations. The
%proof follows from expanding the cumulant in terms of
%expectations. 
\end{itemize}
\end{lemma}

\begin{lemma}\label{lemma:cumulant3}
Suppose Assumptions \ref{assum:S}, \ref{assum:uniform} \ref{assum:G} and
\ref{assum:GG}(b) are satisfied, and $d=1$. 
Then we have 
\begin{eqnarray}
\label{eq:cum3a1}
\cum_{3}[\widetilde{Q}_{a,\lambda}(g,r)]
&=&  O\bigg(\frac{\log^{2}(a)}{\lambda^{2}}\bigg)
%  + 
% \frac{\log^{2}(n)}{\lambda n } + \frac{1}{n^{2}}\bigg).
\end{eqnarray}
with $\lambda^{d}/(\log^{2}(a)n)\rightarrow 0$ as
$\lambda \rightarrow
\infty$, $n\rightarrow \infty$ and $a\rightarrow \infty$.
\end{lemma}
PROOF. We prove the result for $\cum_{3}[\widetilde{Q}_{a,\lambda}(1,0)]$, noting that the proof is identical
for general $g$ and $r$. 
We first expand $\cum_{3}[\tilde{Q}_{a,\lambda}(1,0)]$ using
indecomposable partitions. Using  Lemma \ref{lemma:basic}(i) we note that 
the third order cumulant is zero, therefore
\begin{eqnarray*}
&&\cum_{3}[\tilde{Q}_{a,\lambda}(1,0)] \\
&=&\frac{1}{n^{6}}\sum_{\underline{j}\in \mathcal{D}}^{}\sum_{k_{1},k_{2},k_{3}=-a}^{a} 
\cum\big[Z(s_{j_{1}})Z(s_{j_{2}})\exp(i\omega_{k_{1}}(s_{j_{1}}-s_{j_{2}})),\\
&& Z(s_{j_{3}})Z(s_{j_{4}}) 
\exp(i\omega_{k_{2}}(s_{j_{3}}-s_{j_{4}})),Z(s_{j_{5}})Z(s_{j_{6}})\exp(i\omega_{k_{3}}(s_{j_{5}}-s_{j_{6}}))\big]\\
&=& \frac{1}{n^{6}}\sum_{\underline{j}\in
  \mathcal{D}}^{}\sum_{\pi_{(2,2,2),} \in
  \mathcal{P}_{2,2,2}}A_{2,2,2}^{\underline{j}}(\pi_{(2,2,2)})  +
\frac{1}{n^{6}}\sum_{\underline{j}\in
  \mathcal{D}}^{}\sum_{\pi_{4,2} \in
  \mathcal{P}_{4,2}}A_{4,2}^{\underline{j}}(\pi_{(4,2)}) +  \frac{1}{n^{6}}\sum_{\underline{j}\in
  \mathcal{D}}^{}  A_{6}^{\underline{j}} \\ 
&=& B_{2,2,2} + B_{4,2} + B_{6}, 
\end{eqnarray*}
where $\mathcal{D} = \{j_{1},\ldots,j_{6}\in \{1,\ldots,n\} \textrm{
  but } j_{1}\neq j_{2},j_{3}\neq j_{4},j_{5}\neq j_{6}\}$,
$A_{2,2,2}^{\underline{j}}$   consists of only the product  of cumulants of order two and 
$\mathcal{P}_{2,2,2}$ is the set of all cumulants of order two from
the set of indecomposable partitions
of $\{(1,2),(3,4),(5,6)\}$, $A_{4,2}^{\underline{j}}$ consists 
of only the product of 4th and 2nd order cumulants and
$\mathcal{P}_{4,2}$ is the set of all 4th
order and 2nd order cumulant  indecomposable partitions of
$\{(1,2),(3,4),(5,6)\}$, finally $A_{6}^{\underline{j}}$ is the 6th order cumulant. 
Examples of $A$'s are given below 
\begin{eqnarray}
A_{2,2,2}^{\underline{j}}(\pi_{(2,2,2),1}) 
&=&\sum_{k_{1},k_{2},k_{3}=-a}^{a} 
 \cum\left[Z(s_{j_1})e^{is_{j_1}\omega_{k_{1}}},Z(s_{j_3})e^{is_{j_3}\omega_{k_{2}}}\right]
\cum\left[Z(s_{j_2})e^{-is_{j_2}\omega_{k_{1}}},Z(s_{j_5})e^{is_{j_5}\omega_{k_{3}}}\right]\nonumber\\
&&\times\cum\left[Z(s_{j_4})e^{-is_{j_4}\omega_{k_{2}}},Z(s_{j_6})e^{-is_{j_6}\omega_{k_{3}}}\right] \nonumber\\
&=& \sum_{k_{1},k_{2},k_{3}=-a}^{a}
\Ex[c(s_{j_1}-s_{j_3})e^{i(s_{j_1}\omega_{k_{1}}+s_{j_3}\omega_{k_{2}})}]
\Ex[c(s_{j_2}-s_{j_5})e^{i(-s_{j_2}\omega_{k_{1}}+s_{j_5}\omega_{k_{3}})}]\times \nonumber\\
&&\Ex[c(s_{j_4}-s_{j_6}) e^{-i(s_{j_4}\omega_{k_{2}}+s_{j_6}\omega_{k_{3}})}], \label{eq:partA1}
\end{eqnarray}
\begin{eqnarray}
A_{4,2}^{\underline{j}}(\pi_{(4,2),1}) &=&\sum_{k_{1},k_{2},k_{3}=-a}^{a} 
\cum[Z(s_{j_4})\exp(-is_{j_4}\omega_{k_{2}}),Z(s_{j_6})\exp(-is_{j_6}\omega_{k_{3}})]\nonumber\\
&&\cum[Z(s_{j_1}) \exp(is_{j_1}\omega_{k_{1}}),Z(s_{j_2})\exp(-is_{j_2}\omega_{k_{1}}),
Z(s_{j_3})\exp(is_{j_3}\omega_{k_{2}}),
Z(s_{j_5})\exp(is_{j_5}\omega_{k_{3}})], \nonumber\\
&& \label{eq:partA2}\\
A_{6}^{\underline{j}}   
&=& 
\sum_{k_{1},k_{2},k_{3}=-a}^{a}\cum[Z(s_{j_1})
\exp(is_{j_1}\omega_{k_{1}}),Z(s_{j_2})\exp(-is_{j_2}\omega_{k_{1}}),
Z(s_{j_3})\exp(is_{j_3}\omega_{k_{2}}),
\nonumber\\
&&Z(s_{j_4})\exp(-is_{j_4}\omega_{k_{2}}),Z(s_{j_5})\exp(is_{j_5}\omega_{k_{3}}),
Z(s_{j_6})\exp(-is_{j_6}\omega_{k_3})], \label{eq:partA3}
\end{eqnarray}
where $\underline{j} = (j_{1},\ldots,j_{6})$. 

\subsubsection*{Bound for $B_{222}$}

We will show that $B_{222}$ is 
the leading term in $\cum_{3}(\widetilde{Q}_{a,\lambda}(g;0))$. The set
$\mathcal{D}$ is split into four sets, 
$\mathcal{D}_{6}$ where all the elements of $\underline{j}$ are different,
and for $3\leq i \leq 5$,
$\mathcal{D}_{i}$ where $i$ elements in $\underline{j}$ are different, such that 
\begin{eqnarray*}
B_{2,2,2}=\frac{1}{n^{6}}\sum_{i=0}^{3}\sum_{\underline{j}\in
  \mathcal{D}_{6-i}}^{}\sum_{\pi_{(2,2,2)} \in
  \mathcal{P}_{(2,2,2)}}A_{2,2,2}^{\underline{j}}(\pi_{(2,2,2)}).
\end{eqnarray*}
We start by bounding the partition given in (\ref{eq:partA1}), we later explain
how the same bounds can be obtained for other  
indecomposable partitions in $\mathcal{P}_{2,2,2}$. By using the
spectral representation of the covariance and that $|\mathcal{D}_{6}|=O(n^{6})$, 
it is straightforward to show that 
\begin{eqnarray}
&&\frac{1}{n^{6}}\sum_{\underline{j}\in
  \mathcal{D}_{6}}^{}A_{2,2,2}^{\underline{j}}(\pi_{(2,2,2),1}) \nonumber\\
&=&\sum_{k_{1},k_{2},k_{3}}\frac{c_{6}}{(2\pi)^{3}\lambda^{6}}\int_{\mathbb{R}^{3}}
\int_{[-\lambda/2,\lambda/2]^{6}}f(x)f(y)f(z) 
\exp(i\omega_{k_{1}}(s_{1}-s_{2}))\times \nonumber\\
&& \exp(i\omega_{k_{2}}(s_{3}-s_{4}))
\exp(i\omega_{k_{3}}(s_{5}-s_{6}))\nonumber\\
&&\exp(ix(s_{1}-s_{3}))\exp(iy(s_{4}-s_{6}))
\exp(iz(s_{2}-s_{5})) \prod_{j=1}^{3}ds_{2j-1}ds_{2j}dx dy dz \nonumber\\
&=& \frac{c_{6}}{(2\pi)^{3}}\sum_{k_{1},k_{2},k_{3}}\int_{\mathbb{R}^{3}}
f(x)f(y)f(z)
\sinc\left(\frac{\lambda x}{2}+k_{1}\pi\right)
\sinc\left(\frac{\lambda z}{2}-k_{1}\pi\right)\times \nonumber\\
&&\sinc\left(\frac{\lambda x}{2}-k_{2}\pi\right) 
\sinc\left(\frac{\lambda y}{2}-k_{2}\pi\right)
\sinc\left(\frac{\lambda z}{2}-k_{3}\pi\right)
\sinc\left(\frac{\lambda y}{2}+k_{3}\pi\right)dxdydz, \label{eq:Aexpan}
\end{eqnarray}
where $c_{6}=n(n-1)\ldots (n-5)/n^{6}$. 
By changing variables $x=\frac{\lambda x}{2}+k_{1}\pi$,
$y=\frac{\lambda y}{2}-k_{2}\pi$ and $z=\frac{\lambda z}{2}-k_{3}\pi$
we have 
\begin{eqnarray}
&&\frac{1}{n^{6}}\sum_{\underline{j}\in
  \mathcal{D}_{6}}^{}A_{2,2,2}^{\underline{j}}(\pi_{(2,2,2),1}) \nonumber\\
&=& \frac{c_{6}}{(2\pi)^{3}}\sum_{k_{1},k_{2},k_{3}}\frac{1}{\lambda^{3}}\int_{\mathbb{R}^{3}}
f(\frac{2x}{\lambda}-\omega_{k_{1}})
f(\frac{2y}{\lambda}+\omega_{k_{2}})f(\frac{2z}{\lambda}+\omega_{k_{3}})
\sinc(x)\sinc(z)\sinc(y) \nonumber\\
&&\sinc(x-(k_{2}+k_{1})\pi) \sinc(y+(k_{3}+k_{2})\pi)
\sinc(z-(k_{1}-k_{3})\pi)
dxdydz. \label{eq:D61}
\end{eqnarray}
In order to understand how this case can generalise to other partitions in
$\mathcal{P}_{1}$, 
we represent the $k$s inside the sinc function using the the linear equations 
\begin{eqnarray}
\label{eq:rank}
\left(
\begin{array}{ccc}
-1 & -1 & 0 \\
-1 & 0 & 1 \\
 0 & 1 & 1 \\ 
\end{array}
\right)
\left(
\begin{array}{c}
k_{1} \\
k_{2} \\
k_{3}\\
\end{array}
\right), 
\end{eqnarray}
where we observe that the above is a rank two matrix. Based on this we make the following 
change of variables $k_{1}=k_{1}$, 
$m_{1}=k_{2}+k_{1}$ and $m_{2} = k_{1}-k_{3}$, and rewrite the sum as 
\begin{eqnarray}
&&\frac{1}{n^{6}}\sum_{\underline{j}\in
  \mathcal{D}_{6}}^{}A_{2,2,2}^{\underline{j}}(\pi_{(2,2,2),1})  \nonumber\\
&=& \frac{c_{6}}{(2\pi)^{3}}\sum_{k_{1},m_{1},m_{2}}\frac{1}{\lambda^{3}}\int_{\mathbb{R}^{3}}
f\left(\frac{2x}{\lambda}-\omega_{k_{1}}\right)
f\left(\frac{2y}{\lambda}+\omega_{m_{1}-k_{1}}\right)
f\left(\frac{2z}{\lambda}+\omega_{k_{1}-m_{2}}\right)
\sinc(x)\sinc(z)\sinc(y) \nonumber\\
&&\sinc(x-m_{1}\pi) 
\sinc(y-(m_{2}-m_{1})\pi)
\sinc(z-m_{2}\pi)dxdydz \nonumber\\
&=& \frac{c_{6}}{\lambda^{2}(2\pi)^{3}}\sum_{m_{1},m_{2}}\int_{\mathbb{R}^{3}}
\underbrace{\sinc(x)\sinc(x-m_{1}\pi)  
\sinc(y)\sinc(y+(m_{1}-m_{2})\pi)\sinc(z)\sinc(z-m_{2}\pi)}_{\textrm{contains
two independent $m$ terms}}\nonumber\\
&& \times\frac{1}{\lambda}\sum_{k_{1}}f\left(\frac{2 x}{\lambda}-\omega_{k_{1}}\right)
f\left(\frac{2 y}{\lambda}+\omega_{m_{1}-k_{1}}\right)f\left(\frac{2 z}{\lambda}-\omega_{k_{1}-m_{2}}\right)
dxdydz. \label{eq:D6}
\end{eqnarray}
Finally, we apply Lemma \ref{lemma:1a}(iv)  to give
\begin{eqnarray}
\label{eq:A1}
\frac{1}{n^{6}}\sum_{\underline{j}\in
  \mathcal{D}_{6}}^{}A_{2,2,2}^{\underline{j}}(\pi_{(2,2,2),1})= O\bigg(\frac{\log^{2}a}{\lambda^{2}}\bigg).
\end{eqnarray}
% To generalize to higher dimensions we note that
% $f(\omega)\rightarrow f(\omega_{1},\omega_{2})$ and the 6 sinc
% functions and 3 summands become 12 sinc functions and 6 summands. 
%  
The above only gives the bound for one partition of $\mathcal{P}_{1}$, but we now 
show that the same bound applies to all the other partitions. Looking back at 
(\ref{eq:D61}) and comparing with (\ref{eq:D6}), the reason that only one of the three
$\lambda$s in the denominator of (\ref{eq:D61}) gets `swallowed' is
because the matrix in (\ref{eq:rank}) has rank two. Therefore,  
there are two independent $m$s in the sinc function of (\ref{eq:D61}),
thus by  applying Lemma \ref{lemma:1a}(iv) the sum
$\sum_{m_{1},m_{2}}$ only grows with rate $O(\log^{2}a)$. 
Moreover, it can be shown that all
indecomposable partitions of $\mathcal{P}_{2,2,2}$ 
correspond to rank two matrices (for a proof see equation (A.13) in \citeA{p:che-deo-00}). 
Thus all indecomposable partitions in $\mathcal{P}_{2,2,2}$
will have the same order, which altogether gives 
 \begin{eqnarray*}
\frac{1}{n^{6}}\sum_{\underline{j}\in
  \mathcal{D}_{6}}^{}\sum_{\pi_{(2,2,2)} \in
  \mathcal{P}_{1}}A_{2,2,2}^{\underline{j}}(\pi_{(2,2,2)})
= O\bigg(\frac{\log^{2}a}{\lambda^{2}}\bigg).
\end{eqnarray*}
Now we consider the case that $\underline{j}\in \mathcal{D}_{5}$. In this case, there
are two `typical' cases $\underline{j} =
(j_{1},j_{2},j_{3},j_{4},j_{1},j_{6})$, which gives 
\begin{eqnarray*}
&&A_{2,2,2}^{\underline{j}}(\pi_{(2,2,2),1}) =\\ 
&&\cum\left[Z(s_{j_1})e^{is_{j_1}\omega_{k_{1}}},Z(s_{j_3})e^{is_{j_3}\omega_{k_{2}}}\right]
\cum\left[Z(s_{j_2})e^{-is_{j_2}\omega_{k_{1}}},Z(s_{j_1})e^{is_{j_1}\omega_{k_{3}}}\right]
\cum\left[Z(s_{j_4})e^{-is_{j_4}\omega_{k_{2}}},Z(s_{j_6})e^{-is_{j_6}\omega_{k_{3}}}\right] 
\end{eqnarray*}
and $\underline{j} = (j_{1},j_{2},j_{1},j_{4},j_{5},j_{6})$,  which
gives 
\begin{eqnarray*}
&&A_{2,2,2}^{\underline{j}}(\pi_{(2,2,2),1}) = \\
&&\cum\left[Z(s_{j_1})e^{is_{j_1}\omega_{k_{1}}},Z(s_{j_1})e^{is_{j_1}\omega_{k_{2}}}\right]
\cum\left[Z(s_{j_2})e^{-is_{j_2}\omega_{k_{1}}},Z(s_{j_5})e^{is_{j_5}\omega_{k_{3}}}\right]
\cum\left[Z(s_{j_4})e^{-is_{j_4}\omega_{k_{2}}},Z(s_{j_6})e^{-is_{j_6}\omega_{k_{3}}}\right].
\end{eqnarray*}
Using the same method used to bound (\ref{eq:A1}), when $\underline{j} =
(j_{1},j_{2},j_{3},j_{4},j_{1},j_{6})$ 
$A_{2,2,2}^{\underline{j}}(\pi_{(2,2,2),1}) =
O(\frac{\log^{2}a}{\lambda^{2}})$. 
However, when $\underline{j} =
(j_{1},j_{2},j_{1},j_{4},j_{5},j_{6})$ we use the
same proof used to prove (\ref{eq:cum3}) to give 
$A_{2,2,2}^{\underline{j}}(\pi_{(2,2,2),1}) =
O(\frac{1}{\lambda})$. As we get similar expansions for all 
$\underline{j}\in \mathcal{D}_{5}$ and 
$|\mathcal{D}_{5}|=O(n^{5})$ we have 
\begin{eqnarray*}
\frac{1}{n^{6}}\sum_{\underline{j}\in
  \mathcal{D}_{5}}^{}\sum_{\pi_{(2,2,2)} \in
  \mathcal{P}_{2,2,2}}A_{2,2,2}^{\underline{j}}(\pi_{(2,2,2)})
= O\bigg(\frac{1}{\lambda n} + \frac{\log^{2}(a)}{n\lambda^{2}}\bigg).
\end{eqnarray*}
Similarly we can show that 
\begin{eqnarray*}
\frac{1}{n^{6}}\sum_{\underline{j}\in
  \mathcal{D}_{4}}^{}\sum_{\pi_{(2,2,2)} \in
  \mathcal{P}_{2,2,2}}A_{2,2,2}^{\underline{j}}(\pi_{(2,2,2)})
= O\bigg(\frac{1}{\lambda n^{2}} + \frac{\log^{2}(a)}{n^{2}\lambda^{2}}\bigg).
\end{eqnarray*}
and 
\begin{eqnarray*}
\frac{1}{n^{6}}\sum_{\underline{j}\in
  \mathcal{D}_{3}}^{}\sum_{\pi_{(2,2,2),} \in
  \mathcal{P}_{2,2,2}}A_{2,2,2}^{\underline{j}}(\pi_{(2,2,2)})
= O\bigg(\frac{1}{n^{3}} + \frac{\log^{2}(a)}{n^{3}\lambda^{2}}\bigg).
\end{eqnarray*}
Therefore, if $n>>\lambda/\log^{2}(a)$ we have 
$B_{2,2,2}
= O(\frac{\log^{2}(a)}{\lambda^{2}})$.

\subsubsection*{Bound for $B_{4,2}$}

To bound $B_{4,2}$ we consider the `typical' partition given in (\ref{eq:partA2}). 
Since $A_{4,2}^{\underline{j}}(\pi_{(4,2),1})$ involves fourth order
cumulants by Lemma \ref{lemma:basic}(ii) 
it will be zero in the case that the $\underline{j}$ are all different. 
Therefore, only a maximum of five terms
in $\underline{j}$ can be different, which gives 
\begin{eqnarray*}
\frac{1}{n^{6}}\sum_{\underline{j}\in
  \mathcal{D}}^{}A_{4,2}^{\underline{j}}(\pi_{(4,2),1}) = 
\frac{1}{n^{6}}\sum_{i=1}^{3}\sum_{\underline{j}\in
  \mathcal{D}_{6-i}}^{}A_{4,2}^{\underline{j}}(\pi_{(4,2),1}).
\end{eqnarray*}
We will show that for $\underline{j}\in \mathcal{D}_{5}$,
$A_{4,2}^{\underline{j}}(\pi_{(4,2),1})$ 
will not be as small as $O(\log^{2}(a)/\lambda^{2})$, however, this
will be compensated by $|\mathcal{D}_{5}|=O(n^{5})$ (noting that $|\mathcal{D}_{6}|=O(n^{6})$). Let 
$\underline{j}=(j_{1},j_{2},j_{3},j_{4},j_{1},j_{6})$, then   
expanding the fourth order cumulant in
$A_{4,2}^{\underline{j}}(\pi_{(4,2),1})$ and using conditional cumulants (see
(\ref{eq:conditional-cumulants})) we have 
\begin{eqnarray}
&&   A_{4,2}^{\underline{j}}(\pi_{(4,2),1}) \nonumber\\
&=& \sum_{k_{1},k_{2},k_{3}=-a}^{a} 
\cum[Z(s_{j_4})\exp(-is_{j_4}\omega_{k_{2}}),Z(s_{j_6})\exp(-is_{j_6}\omega_{k_{3}})]\nonumber\\
&&\cum[Z(s_{j_1}) \exp(is_{j_1}\omega_{k_{1}}),Z(s_{j_2})\exp(-is_{j_2}\omega_{k_{1}}),
Z(s_{j_3})\exp(is_{j_3}\omega_{k_{2}}),
Z(s_{j_5})\exp(is_{j_5}\omega_{k_{3}})] \nonumber\\
&=& \sum_{k_{1},k_{2},k_{3}=-a}^{a}
\bigg\{ \cum\big[c(s_{1}-s_{2})e^{i\omega_{k_{1}}(s_{1}-s_{2})},c(s_{3}-s_{1})
e^{is_{3}\omega_{k_{2}}+is_{1}\omega_{k_{3}}}\big]\nonumber\\
&& + \cum\big[c(s_{1}-s_{3})e^{i(s_{1}\omega_{k_{1}}+s_{3}\omega_{k_{2}})},
c(s_{1}-s_{2})e^{is_{1}\omega_{k_{3}}-is_{2}\omega_{k_{1}}}\big]  \nonumber\\
&& +\cum\big[c(0)e^{is_{1}\omega_{k_{1}}+is_{1}\omega_{k_{3}}},
c(s_{2}-s_{3})e^{-is_{2}\omega_{k_{1}}+is_{3}\omega_{k_{2}}}\big]
\bigg\} 
\Ex\big[c(s_{4}-s_{6})e^{-is_{4}\omega_{k_{2}}-is_{6}\omega_{k_{3}}}\big] \nonumber\\
&=& A_{4,2}^{\underline{j}}(\pi_{(4,2),1},\Omega_{1}) + A_{4,2}^{\underline{j}}(\pi_{(4,2),1},\Omega_{2}) + A_{4,2}^{\underline{j}}(\pi_{(4,2),1},
\Omega_{3}), \label{eq:a42pi21}
\end{eqnarray}
where we use the notation $\Omega$ to denote the partition of the fourth order 
cumulant into it's conditional cumulants. To bound each term we expand the 
covariances as expectations, this gives 
\begin{eqnarray*}
&&A_{4,2}^{\underline{j}}(\pi_{(4,2),1},\Omega_{1}) \nonumber\\
&=&\sum_{k_{1},k_{2},k_{3}=-a}^{a}\bigg\{
\Ex\big[c(s_{1}-s_{2})c(s_{3}-s_{1})e^{i\omega_{k_{1}}(s_{1}-s_{2})}
e^{is_{3}\omega_{k_{2}}+is_{1}\omega_{k_{3}}}\big]\Ex\big[c(s_{4}-s_{6})e^{-is_{4}\omega_{k_{2}}-is_{6}\omega_{k_{3}}}\big]  - \nonumber\\ 
&& \Ex\big[c(s_{1}-s_{2})e^{i\omega_{k_{1}}(s_{1}-s_{2})}\big]\Ex\big[c(s_{3}-s_{1})
e^{is_{3}\omega_{k_{2}}+is_{1}\omega_{k_{3}}}\big]\Ex\big[c(s_{4}-s_{6})e^{-is_{4}\omega_{k_{2}}-is_{6}\omega_{k_{3}}}\big] \bigg\}\nonumber\\
&=& A_{4,2}^{\underline{j}}(\pi_{(4,2),1},\Omega_{1},\Pi_{1}) + A_{4,2}^{\underline{j}}(\pi_{(4,2),1},\Omega_{1},\Pi_{2}),\label{eq:a42pi21omega}
\end{eqnarray*}
where we use the notation $\Pi$ to denote the expansion of the
cumulants of the spatial covariances expanded into expectations. 
To bound $A_{4,2}^{\underline{j}}(\pi_{(4,2),1},\Omega_{1},\Pi_{1}) $ we use the spectral representation theorem to give 
\begin{eqnarray*}
&&A_{4,2}^{\underline{j}}(\pi_{(4,2),1},\Omega_{1},\Pi_{1}) = \\
&&\frac{1}{(2\pi)^{3}}
  \int_{\mathbb{R}^{3}} \sum_{k_{1},k_{2},k_{3}=-a}^{a} f(x)f(y)f(z)
\sinc\bigg(\frac{\lambda(x-y)}{2}+(k_{1}+k_{3})\pi\bigg)\sinc\bigg(\frac{\lambda
x}{2}+k_{1}\pi\bigg)\times\\
&&\sinc\bigg(\frac{\lambda
y}{2}+k_{2}\pi\bigg) \sinc\left(\frac{\lambda z}{2}-k_{2}\pi\right)\sinc\left(\frac{\lambda
z}{2}+k_{3}\pi\right)dxdydz. 
\end{eqnarray*}
By changing variables %and applying  Lemma \ref{lemma:1a} we have  
\begin{eqnarray*}
A_{4,2}^{\underline{j}}(\pi_{(4,2),1},\Omega_{1},\Pi_{1}) 
&=&\frac{1}{\pi^{3}\lambda^{3}}
\sum_{k_{1},k_{2},k_{3}=-a}^{a}
\int_{\mathbb{R}^{3}}
f\left(\frac{2u}{\lambda}-\omega_{k_{1}}\right)
f\left(\frac{2v}{\lambda}-\omega_{k_{2}}\right)
f\left(\frac{2w}{\lambda}-\omega_{k_{3}}\right)\times\\
&&\sinc(u-v+(k_{2}+k_{3})\pi)\sinc(u)\sinc(v)\sinc(w)\sinc(w-(k_{2}+k_{3})\pi)dudvdw. 
\end{eqnarray*}
Just as in the bound for $B_{2,2,2}$, 
we represent the $k$s inside the sinc function as a set of linear equations
\begin{eqnarray}\label{eq:rank1}
\left(
\begin{array}{ccc}
0 & 1 & 1 \\
0 & 0 & 0\\
0 & 1 & 1 \\
\end{array}
\right)
\left(
\begin{array}{c}
k_{1} \\
k_{2} \\
k_{3} \\
\end{array}
\right),
\end{eqnarray}
observing the matrix has rank one. We make a change of variables $m=k_{2}+k_{3}$,
$k_{1}=k_{1}$ and $k_{2}=k_{2}$ to give 
\begin{eqnarray*}
&&A_{4,2}^{\underline{j}}(\pi_{(4,2),1},\Omega_{1},\Pi_{1}) \\ 
&=&\frac{1}{\pi^{3}\lambda}
\int_{\mathbb{R}^{3}}\frac{1}{\lambda^{2}}\sum_{k_{1},k_{2}=-a}^{a}
f\left(\frac{2u}{\lambda}-\omega_{k_{1}}\right)
f\left(\frac{2v}{\lambda}-\omega_{k_{2}}\right)
f\left(\frac{2w}{\lambda}-\omega_{m_{1}-k_{2}}\right)\times\\
&&\sum_{m_{}}\sinc(u-v+m_{}\pi)\sinc(u)\sinc(v)\sinc(w)\sinc(w-m_{}\pi)dudvdw
\\
 &=& \frac{2^{3}}{\lambda}
\int_{\mathbb{R}^{3}}
\sum_{m_{}}G_{\lambda,m}\left(\frac{2u}{\lambda},\frac{2v}{\lambda},\frac{2w}{\lambda}\right)
\sinc(u-v+m_{}\pi)\sinc(u)\sinc(v)\sinc(w)\sinc(w-m_{}\pi)dudvdw,
\end{eqnarray*}
where $G_{\lambda,m}(\frac{2u}{\lambda},\frac{2v}{\lambda},\frac{2w}{\lambda}) = \frac{1}{\lambda^{2}}\sum_{k_{1},k_{2}=-a}^{a}
f(\frac{2u}{\lambda}-\omega_{k_{1}})
f(\frac{2v}{\lambda}-\omega_{k_{2}})
f(\frac{2w}{\lambda}-\omega_{m_{}-k_{2}})$. Taking absolutes gives 
\begin{eqnarray*}
|A_{4,2}^{\underline{j}}(\pi_{(4,2),1},\Omega_{1},\Pi_{1}) |\leq \frac{C}{\lambda}
\int_{\mathbb{R}^{3}}^{}
\sum_{m_{}}|\sinc(u-v+m_{}\pi)\sinc(w-m_{}\pi)\sinc(u)\sinc(v)\sinc(w)|dudvdw
\end{eqnarray*}
Since the above contains $m$ in the sinc function we use Lemma
\ref{lemma:1a}(i), equations (\ref{eq:i1}) and (\ref{eq:i2}), to show 
\begin{eqnarray*}
|A_{4,2}^{\underline{j}}(\pi_{(4,2),1},\Omega_{1},\Pi_{1}) |&\leq &\frac{C}{\lambda}\int_{\mathbb{R}^{3}}
\sum_{m_{}}|\sinc(u-v+m_{}\pi)\sinc(w-m_{}\pi)\sinc(u)\sinc(v)\sinc(w)|dudvdw \\
&\leq& 
\frac{C}{\lambda}\sum_{m}\ell_{1}(m\pi)\ell_{2}(m\pi) = O(\lambda^{-1}),
\end{eqnarray*}
where the functions $\ell_{1}(\cdot)$ and $\ell_{2}(\cdot)$ are defined in  Lemma
\ref{lemma:1a}(i).
Thus $|A_{4,2}^{\underline{j}}(\pi_{(4,2),1},\Omega_{1},\Pi_{1}) | =
O(\frac{1}{\lambda})$. We use the same method used to 
bound  (\ref{eq:A1}) to show that $A_{4,2}^{\underline{j}}(\pi_{(4,2),1},\Omega_{1},\Pi_{2}) 
= O(\frac{\log^{2}(a)}{\lambda^{2}})$ and
$|A_{4,2}^{\underline{j}}(\pi_{(4,2),1},\Omega_{2})| =
O(\frac{1}{\lambda})$. Furthermore, it is straightforward to see that by the
independence of $s_{1}$ and $s_{2}$ and $s_{3}$ that
$A_{4,2}^{\underline{j}}(\pi_{(4,2),1},\Omega_{3})=0$ (recalling that
$A_{4,2}^{\underline{j}}(\pi_{(4,2),1},\Omega_{3})$ is defined in equation (\ref{eq:a42pi21})).
Thus altogether we have for $\underline{j}=(j_{1},j_{2},j_{3},j_{4},j_{1},j_{6})$ and the partition $\pi_{(4,2),1}$,
that $|A_{4,2}^{\underline{j}}(\pi_{(4,2),1})| = O(\frac{1}{\lambda})$. However, it is important to note
for all other $\underline{j}\in \mathcal{D}_{5}$ and partitions in $\mathcal{P}_{4,2}$ the
same method will lead to a similar decomposition given in (\ref{eq:a42pi21}) and the rank one
matrix given in (\ref{eq:rank1}). The rank one matrix means one `free' $m$ in the sinc functions
and this $|A_{4,2}^{\underline{j}}(\pi_{4,2})| = O(\frac{1}{\lambda})$ for all 
$\underline{j}\in \mathcal{D}_{5}$ and $\pi_{4,2}\in \mathcal{P}_{4,2}$. Thus, since  
 $|\mathcal{D}_{5}|=O(n^{5})$ we have 
\begin{eqnarray*}
\frac{1}{n^{6}}\sum_{\underline{j}\in
  \mathcal{D}_{5}} \sum_{\pi_{4,2}\in \mathcal{P}_{4,2}}
  A_{4,2}^{\underline{j}}(\pi_{4,2}) =
  O\left(\frac{1}{\lambda n}+\frac{\log^{2}a}{\lambda^{2}n}\right) = O\left(\frac{\log^{2}a}{\lambda^{2}}\right) 
\end{eqnarray*}
if $n>>\lambda/\log^{2}(a)$
(ie. $\frac{\lambda}{n\log^{2}a}\rightarrow 0$). For $\underline{j}\in \mathcal{D}_{4}$ and 
$\underline{j}\in \mathcal{D}_{3}$ we use the same argument, noting that the number of 
free $m$'s in the sinc functions goes down but to compensate, $|\mathcal{D}_{4}|=O(n^{4})$
and $|\mathcal{D}_{3}|=O(n^{3})$. Therefore,  if
$n>>\log^{2}(a)/\lambda$, then
\begin{eqnarray*}
B_{4,2} = \frac{1}{n^{6}}\sum_{i=1}^{3}\sum_{\underline{j}\in
  \mathcal{D}_{6-i}}^{}\sum_{\pi_{4,2}\in \mathcal{P}_{4,2}}A_{4,2}^{\underline{j}}(\pi_{4,2}) = 
 O\left(\frac{\log^{2}a}{\lambda^{2}}\right).
\end{eqnarray*}

\subsubsection*{Bound for $B_{6}$}

Finally, we bound $B_{6}$.
By using Lemma \ref{lemma:basic}(ii) we 
observe that $A_{6}^{\underline{j}}(k_{1},k_{2},k_{3})=0$ if more than
four elements of $\underline{j}$ are different. Thus 
\begin{eqnarray*}
B_{6} = \frac{1}{n^{6}}\sum_{i=2}^{3}\sum_{\underline{j}\in
  \mathcal{D}_{6-i}}A_{6}^{\underline{j}}.
\end{eqnarray*}
We start by considering the case
that $\underline{j}=(j_{1},j_{2},j_{1},j_{4},j_{1},j_{6})$ (three
elements in $\underline{j}$ are the same), then by using 
conditional cumulants we have 
\begin{eqnarray*}
&&  A_{6}^{\underline{j}} \\
&=& \sum_{k_{1},k_{2},k_{3}=-a}^{a}\cum[Z(s_{1})e^{is_{1}\omega_{k_{1}}},
Z(s_{2})e^{-is_{2}\omega_{k_{1}}},Z(s_{1})e^{is_{1}\omega_{k_{2}}},Z(s_{4})e^{-is_{4}\omega_{k_{2}}},
Z(s_{1})e^{is_{1}\omega_{k_{3}}},Z(s_{6})e^{-is_{6}\omega_{k_{3}}}] \\
&=& \sum_{k_{1},k_{2},k_{3}=-a}^{a}\sum_{\Omega \in
  \mathcal{R}_{}}A_{6}^{\underline{j}}(\Omega),
\end{eqnarray*}
where $\mathcal{\mathcal{R}}_{}$ is the set of all pairwise partitions of 
$\{1,2,1,4,1,6\}$, for example
\begin{eqnarray*}
A_{6}^{\underline{j}}(\Omega_{1}) =\cum\big[
c(s_{1}-s_{2})e^{i\omega_{k_{1}}(s_{1}-s_{2})},c(s_{1}-s_{4})e^{i\omega_{k_{2}}(s_{1}-s_{4})},c(s_{1}-s_{6})
e^{i\omega_{k_{3}}(s_{1}-s_{6})}
\big]. 
\end{eqnarray*}
We will first bound the above and then explain how this generalises to
the other $\Omega \in \mathcal{R}_{}$ and $\underline{j}\in \mathcal{D}_{4}$. 
Expanding the above third order cumulant in terms of expectations gives 
\begin{eqnarray*}
&&A_{6}^{\underline{j}}(\Omega_{1}) \\
&=& \sum_{k_{1},k_{2},k_{3}=-a}^{a}\bigg\{ 
\Ex\big[
c(s_{1}-s_{2})e^{i\omega_{k_{1}}(s_{1}-s_{2})}c(s_{1}-s_{4})e^{i\omega_{k_{2}}(s_{1}-s_{4})}c(s_{1}-s_{6})
e^{i\omega_{k_{3}}(s_{1}-s_{6})}\big] - \\
 && \Ex\big[
c(s_{1}-s_{2})e^{i\omega_{k_{1}}(s_{1}-s_{2})}\big]
\Ex\big[c(s_{1}-s_{4})e^{i\omega_{k_{2}}(s_{1}-s_{4})}c(s_{1}-s_{6})
e^{i\omega_{k_{3}}(s_{1}-s_{6})}\big] - \\
 && \Ex\big[
c(s_{1}-s_{2})e^{i\omega_{k_{1}}(s_{1}-s_{2})}c(s_{1}-s_{4})e^{i\omega_{k_{2}}(s_{1}-s_{4})}\big)
\Ex\big[c(s_{1}-s_{6})
e^{i\omega_{k_{3}}(s_{1}-s_{6})}\big] - \\
&& \Ex\big[
c(s_{1}-s_{2})e^{i\omega_{k_{1}}(s_{1}-s_{2})}c(s_{1}-s_{6})
e^{i\omega_{k_{3}}(s_{1}-s_{6})}\big]\Ex\big[
c(s_{1}-s_{4})e^{i\omega_{k_{2}}(s_{1}-s_{4})}\big]+ \\
&& 2\Ex\big[
c(s_{1}-s_{2})e^{i\omega_{k_{1}}(s_{1}-s_{2})}\big]\Ex\big[c(s_{1}-s_{4})e^{i\omega_{k_{2}}(s_{1}-s_{4})}\big] 
\Ex\big[c(s_{1}-s_{6})e^{i\omega_{k_{3}}(s_{1}-s_{6})}\big] \bigg\} \\
 &=& \sum_{\ell=1}^{5}A_{6}^{\underline{j}}(\Omega_{1},\Pi_{\ell}).
\end{eqnarray*}
We observe that for $2\leq \ell \leq 5$, 
$A_{6}^{\underline{j}}(\Omega_{1},\Pi_{\ell})$ resembles 
$A_{4,2}^{\underline{j}}(\pi_{(4,2),1},\Omega_{1},\Pi_{1})$ defined in
(\ref{eq:a42pi21omega}), thus the same proof used to bound the terms
in (\ref{eq:a42pi21omega}) can be use to show that for $2\leq \ell
\leq 5$ 
 $A_{6}^{\underline{j}}(\Omega_{1},\Pi_{\ell}) =
 O(\frac{1}{\lambda})$. However, the first term $A_{6}^{\underline{j}}(\Omega_{1},\Pi_{1})$
 involves just one expectation, and is not included in the previous
 cases.  By using the spectral representation theorem we have
\begin{eqnarray*}
&&A_{6}^{\underline{j}}(\Omega_{1},\Pi_{1})\\
&=& \sum_{k_{1},k_{2},k_{3}=-a}^{a}\Ex\left[
c(s_{1}-s_{2})e^{i\omega_{k_{1}}(s_{1}-s_{2})}c(s_{1}-s_{4})e^{i\omega_{k_{2}}(s_{1}-s_{4})}c(s_{1}-s_{6})
e^{i\omega_{k_{3}}(s_{1}-s_{6})}\right] \\
&=& \frac{1}{(2\pi)^{3}}\sum_{k_{1},k_{2},k_{3}=-a}^{a}\int\int\int
f(x)f(y)f(z)\sinc\left(\frac{\lambda(x+y+z)}{2}+(k_{1}+k_{2}+k_{3})\pi\right)
\sinc\left(\frac{\lambda x}{2}+k_{1}\pi\right) \\
&& \sinc\left(\frac{\lambda y}{2}+k_{2}\pi\right) \sinc\left(\frac{\lambda
  z}{2}+k_{3}\pi\right)dxdydz \\
&=&
\frac{2^{3}}{(2\pi)^{3}\lambda^{3}}\sum_{k_{1},k_{2},k_{3}=-a}^{a}\int\int\int
f(\frac{2u}{\lambda}-\omega_{k_{1}})
f(\frac{2v}{\lambda}-\omega_{k_{2}})
f(\frac{2w}{\lambda}-\omega_{k_{3}})\times \\
&&\sinc(u+v+w)\sinc(u)\sinc(v)\sinc(w)dudvdw.
\end{eqnarray*}
It is obvious that the $k$s within the sinc function correspond to
a rank zero matrix, and thus $A_{6}(\Omega_{1},\Pi_{1})=O(1)$. Therefore,
  $A_{6}^{\underline{j}}(\Omega_{1}) = O(1)$. A similar bound holds
  for all $\underline{j}\in \mathcal{D}_{4}$, this we have  
\begin{eqnarray*}
\frac{1}{n^{6}}\sum_{\underline{j}\in
  \mathcal{D}_{4}}^{}  A_{6}^{\underline{j}}(\Omega_{1}) = O\left(\frac{1}{n^{2}}\right),
\end{eqnarray*}
since $|\mathcal{D}_{4}|=O(n^{4})$. Indeed, the same argument applies to the other partitions  
$\Omega$ and $\underline{j}\in \mathcal{D}_{3}$, thus altogether we
have 
\begin{eqnarray*}
\mathcal{B}_{6}=\frac{1}{n^{6}}\sum_{\Omega\in \mathcal{R}}\sum_{i=2}^{3}\sum_{\underline{j}\in
  \mathcal{D}_{6-i}}^{}   A_{6}^{\underline{j}}(\Omega) = O\left(\frac{1}{n^{2}}\right).
\end{eqnarray*}

Altogether, using the  bounds derived for $\mathcal{B}_{2,2,2}$, $\mathcal{B}_{4,2}$ and $\mathcal{B}_{6}$  
we have 
\begin{eqnarray*}
\cum_{3}(\widetilde{Q}_{a,\lambda}(1,0)) =
O\bigg(\frac{\log^{2}(a)}{\lambda^{2}} + \frac{1}{n\lambda} + \frac{\log^{2}(a)}{\lambda^{2}n} +\frac{1}{n^{2}}\bigg)
=O\left(\frac{\log^{2}(a)}{\lambda^{2}}\right),
\end{eqnarray*}
where the last bound is due to the conditions on $a,n$ and
$\lambda$. This gives the result. 
\hfill $\Box$

\vspace{3mm}
We now generalize the above results to higher order cumulants. 
\begin{lemma}\label{lemma:cumulantq}
Suppose Assumptions \ref{assum:S}, \ref{assum:uniform},
\ref{assum:G} and \ref{assum:GG}(b) are satisfied, and $d=1$. 
Then for $q\geq 3$ we have 
\begin{eqnarray}
\label{eq:cumqa1}
\cum_{q}[\widetilde{Q}_{a,\lambda}(g,r)]
&=&  O\bigg(\frac{\log^{2(q-2)}(a)}{\lambda^{q-1}}\bigg),
%  + 
% \frac{\log^{2}(n)}{\lambda n } + \frac{1}{n^{2}}\bigg).
\end{eqnarray}
\begin{eqnarray}
\label{eq:cumqa2}
\cum_{q}[\widetilde{Q}_{a,\lambda}(g,r_{1}),\ldots,\widetilde{Q}_{a,\lambda}(g,r_{q})]
&=&  O\bigg(\frac{\log^{2(q-2)}(a)}{\lambda^{q-1}}\bigg)
%  + 
% \frac{\log^{2}(n)}{\lambda n } + \frac{1}{n^{2}}\bigg).
\end{eqnarray}
and in the case $d>1$, we have 
\begin{eqnarray}
\label{eq:cumqa2d}
\cum_{q}[\widetilde{Q}_{a,\lambda}(g,\rb_{1}),\ldots,\widetilde{Q}_{a,\lambda}(g,\rb_{q})]=
 O\bigg(\frac{\log^{2d(q-2)}(a)}{\lambda^{d(q-1)}}\bigg)
\end{eqnarray}
with $\lambda^{d}/(\log^{2}(a)n)\rightarrow 0$ as
$\lambda \rightarrow
\infty$, $n\rightarrow \infty$ and $a\rightarrow \infty$.
\end{lemma}
PROOF. The proof  essentially follows the  same method used to bound the second and third
cumulants. We first prove (\ref{eq:cumqa1}). 
To simplify the notation we prove the result for $g=1$ and $r=0$, noting that the proof
in the general case is identical. Expanding out
$\cum_{q}[\widetilde{Q}_{a,\lambda}(1,0)]$ and using indecomposable
partitions gives 
\begin{eqnarray*}
&&\cum_{q}[\widetilde{Q}_{a,\lambda}(1,0)]\\
&=&
\frac{1}{n^{2q}}\sum_{j_{1},\ldots,j_{2q}\in \mathcal{D}}\sum_{k_{1},\ldots,k_{q}=-a}^{a}
\cum\left[Z(s_{j_1})Z(s_{j_2})e^{i\omega_{k_{1}}(s_{j_1}-s_{j_2})},\ldots,
Z(s_{j_{2q-1}})Z(s_{j_{2q}})e^{i\omega_{k_{q}}(s_{j_{2q-1}}-s_{j_{2q}})}
\right]\nonumber\\
&=& \frac{1}{n^{2q}}\sum_{\underline{j}\in \mathcal{D}}
\sum_{\underline{b}\in \mathcal{B}_{q}}\sum_{\pi_{\underline{b}}\in
  \mathcal{P}_{2\underline{b}}}A_{2\underline{b}}^{\underline{j}}(\pi_{2\underline{b}}) 
%=  \frac{1}{n^{2q}}\sum_{\underline{j}\in \mathcal{D}}
%\sum_{\underline{b}\in \mathcal{B}_{q}}\sum_{\pi_{\underline{b}}\in
%  \mathcal{P}_{\underline{b}}}A_{2\underline{b}}^{\underline{j}}(\pi_{\underline{b}}), 
\end{eqnarray*}
where $\mathcal{B}_{q}$
corresponds to the set of integer partitions of $q$ (more precisely,
each partition is a sequence of positive integers which 
sums to $q$). The notation used here is simply a generalization of the
notation used in Lemma \ref{lemma:cumulant3}. 
Let
$\underline{b}=(b_{1},\ldots,b_{m})$ (noting $\sum_{j}b_{j}=q$) denote one of these partitions,
then $\mathcal{P}_{2\underline{b}}$ is the set of all 
indecomposable partitions of $\{(1,2),(3,4),\ldots,(2q-1,2q)\}$ where
the size of each partition is $2b_{1},2b_{2},\ldots,2b_{m}$. For example, if $q=3$, then one
example of an element of $\mathcal{B}_{3}$ is 
$\underline{b}=(1,1,1)$ and  $\mathcal{P}_{2\underline{b}}=\mathcal{P}_{(2,2,2)}$
corresponds to all pairwise indecomposable partitions of 
$\{(1,2),(3,4),(5,6)\}$. Finally,
$A_{2\underline{b}}^{\underline{j}}(\pi_{2\underline{b}})$ corresponds
to the product of one indecomposable partition of the 
 cumulant $\cum\big[Z(s_{j_1})Z(s_{j_2})e^{i\omega_{k_{1}}(s_{j_1}-s_{j_2})},\ldots,
Z(s_{j_{2q-1}})Z(s_{j_{2q}})e^{i\omega_{k_{q}}(s_{j_{2q-1}}-s_{j_{2q}})}
\big]$, where the cumulants are of order 
$2b_{1},2b_{2},\ldots,2b_{m}$  (examples, in the case $q=3$ are given
in equation (\ref{eq:partA1})-(\ref{eq:partA3})). Let 
\begin{eqnarray*}
B_{2\underline{b}}=\frac{1}{n^{2q}}\sum_{\underline{j}\in
  \mathcal{D}}\sum_{\pi_{2\underline{b}}\in
  \mathcal{P}_{2\underline{b}}}A_{2\underline{b}}^{\underline{j}}(\pi_{2\underline{b}}), 
\end{eqnarray*}
therefore $\cum_{q}[\widetilde{Q}_{a,\lambda}(1,0)]=\sum_{\underline{b}\in
  \mathcal{B}_{q}}B_{2\underline{b}}$. 
%We will derive bounds for two the
%extreme cases 
%\begin{eqnarray*}
%&&B_{1,1,\ldots,1} \\
%&=& \frac{1}{n^{2q}}\sum_{\underline{j}\in
%  \mathcal{D}}\sum_{k_{1},\ldots,k_{q}=-a}^{a}
%\bigg(\cov(Z(s_{j_{1}})e^{is_{j_{1}}\omega_{k_{1}}},
%Z(s_{j_{2q}})e^{-is_{j_{2q}}\omega_{k_{q}}})
%\prod_{c=1}^{q-1}\cov(Z(s_{j_{2c}})e^{-is_{j_{2c}}\omega_{k_{c}}},Z(s_{j_{2c+1}})e^{is_{j_{2c+1}}\omega_{k_{c+1}}}) \\
%&&+\ldots+\cov(Z(s_{j_{2}})e^{-is_{j_{2}}\omega_{k_{1}}},
%Z(s_{j_{2q-1}})e^{is_{j_{2q}}\omega_{k_{q}}})\prod_{c=1}^{q-1}\cov(Z(s_{j_{2c-1}})e^{is_{j_{2c-1}}\omega_{k_{c-1}}},
%Z(s_{j_{2c}})e^{is_{j_{2c}}\omega_{k_{c}}}) 
%\bigg)
%\end{eqnarray*}
%and
%\begin{eqnarray*}
%B_{2q} &=& \frac{1}{n^{2q}}\sum_{\underline{j}\in
%  \mathcal{D}}\sum_{k_{1},\ldots,k_{q}=-a}^{a}
%\cum_{2q}(Z(s_{j_{1}})e^{is_{j_{1}}\omega_{k_{1}}},\ldots,
%Z(s_{j_{2q}})e^{-is_{j_{2q}}\omega_{k_{q}}}), 
%\end{eqnarray*}
%and show that the other cases lie `between' these two cases. 

Just as in the proof of Lemma \ref{lemma:cumulant3}, we will show that under
the condition $n>>\lambda/\log^{2}(a)$, the pairwise decomposition   
$B_{2,\ldots,2}$ is the denominating term. We start with a `typical'
decomposition
$\pi_{(2,\ldots,2),1}\in  \mathcal{P}_{2,\ldots,2}$, 
\begin{eqnarray*}
A_{2,\ldots,2}^{\underline{j}}(\pi_{(2,\ldots,2),1}) = \sum_{k_{1},\ldots,k_{q}=-a}^{a}
\cum[Z(s_{1})e^{is_{1}\omega_{k_{1}}},
Z(s_{2q})e^{-is_{2q}\omega_{k_{q}}}]
\prod_{c=1}^{q-1}\cum[Z(s_{2c})e^{-is_{2c}\omega_{k_{c}}},Z(s_{2c+1})e^{is_{2c+1}\omega_{k_{c+1}}}]. 
\end{eqnarray*}
and 
\begin{eqnarray*}
\frac{1}{n^{2q}}\sum_{\underline{j}\in
  \mathcal{D}}A_{2,\ldots,2}^{\underline{j}}(\pi_{(2,\ldots,2),1}) = \frac{1}{n^{2q}}\sum_{i=0}^{q-1}\sum_{\underline{j}\in
  \mathcal{D}_{2q-i}}A_{2,\ldots,2}^{\underline{j}}(\pi_{(2,\ldots,2),1}), 
\end{eqnarray*}
where $\mathcal{D}_{2q}$ denotes the set where all elements of $\underline{j}$ are different
and $\mathcal{D}_{2q-i}$ denotes the set that $(2q-i)$ elements in $\underline{j}$
are different. We first consider the case that $\underline{j}=(1,2,\ldots,2q)\in \mathcal{D}_{2q}$. 
Using identical arguments to those used for
$\var[\widetilde{Q}_{a,\lambda}(1,0)]$ and
 $\cum_{3}[\widetilde{Q}_{a,\lambda}(1,0)]$ we can show that 
\begin{eqnarray}
%&&\sum_{k_{1},\ldots,k_{q}=-a}^{a}
%\cov(Z(s_{1})e^{is_{1}\omega_{k_{1}}},
%Z(s_{2q})e^{-is_{2q}\omega_{k_{q}}})
%\prod_{c=1}^{q-1}\cov(Z(s_{2c})e^{-is_{2c}\omega_{k_{c}}},Z(s_{2c+1})e^{is_{2c+1}\omega_{k_{c+1}}}) \\
A_{2,\ldots,2}^{\underline{j}}(\pi_{(2,\ldots,2),1}) &=& \sum_{k_{1},\ldots,k_{q}=-a}^{a}\int_{\mathbb{R}^{q}}
f(x_{q})\sinc\left(\frac{2x_{q}}{\lambda}+k_{1}\pi\right)\sinc\left(\frac{2x_{q}}{\lambda}+k_{q}\pi\right)\times\nonumber\\
&&\prod_{c=1}^{q-1}f(x_{c})\sinc\left(\frac{2x_{c}}{\lambda}-k_{c}\pi\right)
\sinc\left(\frac{2x_{c}}{\lambda}-k_{c+1}\right)\prod_{c=1}^{q}
dx_{c}. \label{eq:A222}
\end{eqnarray}
By a change of variables we get 
\begin{eqnarray*}
A_{2,\ldots,2}^{\underline{j}}(\pi_{(2,\ldots,2),1}) &=&\frac{1}{\lambda^{q}}\sum_{k_{1},\ldots,k_{q}=-a}^{a}
\int_{\mathbb{R}^{q}}\prod_{c=1}^{q-1}f\left(\frac{\lambda
  u_{c}}{2}+\omega_{c}\right)\sinc(u_{c})\sinc(u_{c}-(k_{c+1}-k_{c})\pi) \\
&&\qquad\times f\left(\frac{\lambda
  u_{q}}{2}-\omega_{1}\right)\sinc(u_{q})\sinc(u_{q}+(k_{q}-k_{1})\pi)
%\sinc(\frac{2x_{c}}{\lambda}-k_{c}\pi)
%\sinc(\frac{2x_{c}}{\lambda}+k_{c+1})
%\sinc(\frac{2x_{q}}{\lambda}+k_{1}\pi)\sinc(\frac{2x_{q}}{\lambda}+k_{p}\pi)
\prod_{c=1}^{q}dx_{c}.
\end{eqnarray*}
As in the proof of the third order cumulant we can rewrite the $k$s in the above as a
matrix equation
\begin{eqnarray*}
\left(
\begin{array}{cccccc}
-1    & 1   &  0      & \ldots &  0 & 0\\
0         & -1   & 1  & \ldots &  0 & 0\\
\vdots & \vdots & \vdots & \vdots & 0 & 0\\
-1         & 0        &  \ldots & \ldots & 0 & 1 \\
\end{array}
\right)
\left(
\begin{array}{c}
k_{1} \\
k_{2} \\
\vdots \\
k_{q} \\
\end{array}
\right),
\end{eqnarray*}
noting that that above is a $(q-1)$-rank matrix. Therefore applying
the same arguments that were used in the proof of  $\cum_{3}[\widetilde{Q}_{a,\lambda}(1,0)]$ and also Lemma \ref{lemma:1a}(iii) we can show that  
$A^{\underline{j}}_{2,\ldots,2}(\pi_{(2,\ldots,2),1}) 
=O(\frac{\log^{2(q-2)}(a)}{\lambda^{q-1}})$. Thus for
$\underline{j}\in \mathcal{D}_{2q}$ we have 
$\frac{1}{n^{2q}}\sum_{\underline{j}\in \mathcal{D}_{2q}}A^{\underline{j}}_{2,\ldots,2}(\pi_{(2,\ldots,2),1}) 
=O(\frac{\log^{2(q-2)}(a)}{\lambda^{q-1}})$.

In the case that $\underline{j}\in \mathcal{D}_{2q-1}$ 
($(2q-1)$-terms in $\underline{j}$ are different)  by using the 
same arguments as those used to bound $A_{2,2,2}$ (in the proof of $\cum[\widetilde{Q}_{a,\lambda}(1,0)]$) we have 
$A^{\underline{j}}_{2,\ldots,2}(\pi_{(2,\ldots,2),1}) = O(\frac{\log^{2(q-3)}(a)}{\lambda^{q-2}})$, similarly if 
$(2q-2)$-terms in $\underline{j}$ are different, then 
$A^{\underline{j}}_{2,\ldots,2}(\pi_{(2,\ldots,2),1}) = O(\frac{\log^{2(q-4)}(a)}{\lambda^{q-3}})$ and so forth. Therefore, since
$|\mathcal{D}_{2q-i}|=O(n^{2q-i})$ we have 
\begin{eqnarray*}
\frac{1}{n^{2q}}\sum_{\underline{j}\in
  \mathcal{D}}|A_{2,\ldots,2}^{\underline{j}}(\pi_{(2,\ldots,2),1})| \leq 
C\sum_{i=0}^{q}\frac{\log^{2(q-2-i)}(a)}{\lambda^{q-1-i}n^{i}}.
\end{eqnarray*}
Now by using that 
$n>>\lambda/\log^{2}(a)$ we have
\begin{eqnarray*}
\frac{1}{n^{2q}}\sum_{\underline{j}\in
  \mathcal{D}}|A_{2,\ldots,2}^{\underline{j}}(\pi_{(2,\ldots,2),1})| =O
\left(\frac{\log^{2(q-2)}(a)}{\lambda^{q-1}}\right).
\end{eqnarray*}
The same argument holds
for every other second order cumulant indecomposable partition,  because the
corresponding matrix will always have rank $(q-1)$ in the case that 
$\underline{j}\in \mathcal{D}_{2q}$ or for $\underline{j}\in
\mathcal{D}_{2q-i}$ and the dependent $s_{j}$'s 
lie in different cumulants (see \citeA{p:che-deo-00}), thus 
$B_{2,\ldots,2}=O(\frac{\log^{2(q-1)}(a)}{\lambda^{q-1}})$. 

Now, we bound the other extreme $B_{2q}$. Using the conditional
cumulant expansion (\ref{eq:conditional-cumulants}) and 
noting that $\cum_{2q}[Z(s_{j_{1}})e^{is_{j_{1}}\omega_{k_{1}}},\ldots,
Z(s_{j_{2q}})e^{-is_{j_{2q}}\omega_{k_{q}}}]$ is non-zero, only when
at  most $(q+1)$ elements of $\underline{j}$ are different we have 
\begin{eqnarray*}
B_{2q} &=& \frac{1}{n^{2q}}\sum_{\underline{j}\in
  \mathcal{D}}\sum_{k_{1},\ldots,k_{q}=-a}^{a}
\cum_{2q}[Z(s_{j_{1}})e^{is_{j_{1}}\omega_{k_{1}}},\ldots,
Z(s_{j_{2q}})e^{-is_{j_{2q}}\omega_{k_{q}}}]   \\
 &=& \frac{1}{n^{2q}}\sum_{\underline{j}\in
  \mathcal{D}_{q+1}}\sum_{\Omega \in
  \mathcal{R}_{2q}}A^{\underline{j}}_{2q}(\Omega).
%\cum_{2q}(Z(s_{j_{1}})e^{is_{j_{1}}\omega_{k_{1}}},\ldots,
%Z(s_{j_{2q}})e^{-is_{j_{2q}}\omega_{k_{q}}})
\end{eqnarray*}
where $\mathcal{R}_{2q}$ is the set of all pairwise partitions of
$\{1,2,1,3,\ldots,1,q\}$. We consider a `typical' partition 
$\Omega_{1}\in \mathcal{R}_{2q}$
\begin{eqnarray}
\label{eq:Ajj}
A^{\underline{j}}_{2q}(\Omega_{1})=\sum_{k_{1},\ldots,k_{q}=-a}^{a}
\cum\big[c(s_{1}-s_{2})e^{i(s_{1}-s_{2})\omega_{k_{1}}},\ldots,c(s_{1}-s_{q+1})e^{i(s_{1}-s_{q+1})\omega_{k_{q}}}\big].
\end{eqnarray}
By expanding the above the cumulant  as the sum of the product of
expectations we have 
 \begin{eqnarray*}
A^{\underline{j}}_{2q}(\Omega_{1}) =\frac{1}{n^{2q}}\sum_{\underline{j}\in
  \mathcal{D}_{q+1}}\sum_{\Omega \in \mathcal{R}_{2q}}\sum_{\Pi \in \mathcal{S}_{q}}A_{2q}^{\underline{j}}(\Omega_{1},\Pi),
%\cum_{2q}(Z(s_{j_{1}})e^{is_{j_{1}}\omega_{k_{1}}},\ldots,
%Z(s_{j_{2q}})e^{-is_{j_{2q}}\omega_{k_{q}}})
\end{eqnarray*}
where $\mathcal{S}_{q}$ is the set of all partitions of
$\{1,\ldots,q\}$. As we have seen in both the 
$\var[\widetilde{Q}_{a,\lambda}(1,0)]$ and  $\cum_{3}[\widetilde{Q}_{a,\lambda}(1,0)]$
calculations, the leading term in the cumulant expansion is
the expectation over all the covariance terms. 
The same result holds true for higher order cumulants, the expectation
over all the covariances in that cumulant is the leading term because
the it gives the linear equation of the $k$s in the sinc function with
the lowest order rank (we recall the lower the rank the less `free' $\lambda$s). 
Based on this we will only derive bounds for the expectation over all
the covariances. Let $\Pi_{1}\in
\mathcal{S}_{q}$, where 
\begin{eqnarray*}
A^{\underline{j}}_{2q}(\Omega_{1},\Pi_{1})=\sum_{k_{1},\ldots,k_{q}=-a}^{a}
\Ex\big[\prod_{c=1}^{q}c(s_{1}-s_{c+1})e^{i(s_{1}-s_{c+1})\omega_{k_{c}}}\big].
\end{eqnarray*}
Representing the
above expectation as an integral and using the spectral representation
theorem and a change of variables gives
\begin{eqnarray*}
A^{\underline{j}}_{2q}(\Omega_{1},\Pi_{1})
&=& \sum_{k_{1},\ldots,k_{q}=-a}^{a}
\Ex\big[\prod_{c=1}^{q}c(s_{1}-s_{c+1})e^{i(s_{1}-s_{c+1})\omega_{k_{c}}}\big] \\
&=& \frac{1}{(2\pi)^{q}}\sum_{k_{1},k_{2},k_{3}=-a}^{a}\int_{\mathbb{R}^{q}}
\sinc\bigg(\frac{\lambda(\sum_{c=1}^{q}x_{c})}{2}+\pi\sum_{c=1}^{q}k_{c}\bigg)
\prod_{c=1}^{q}f(x_{c})\sinc(x_{c}+k_{c}\pi)\prod_{c=1}^{q}dx_{c}
 \\
&=&
\frac{2^{q}}{(2\pi)^{q}\lambda^{q}}\sum_{k_{1},\ldots,k_{q}=-a}^{a}\int_{\mathbb{R}^{q}}
\sinc(\sum_{c=1}^{q}u_{c})
\prod_{c=1}^{q}f\left(\frac{2u_{c}}{\lambda}-\omega_{k_{c}}\right)\sinc(u_{c})du_{c}
=O(1),
\end{eqnarray*}
where the last line follows from Lemma \ref{lemma:1a}, equation
(\ref{eq:in}). Therefore, $A^{\underline{j}}_{2q}(\Omega_{1}) =
O(1)$. By using the same method on every partition $\Omega \in
\mathcal{R}_{q+1}$ and $\underline{j}\in\mathcal{D}_{q+1}$ and 
$|\mathcal{D}_{q+1}|=O(n^{q+1})$, we have 
\begin{eqnarray*}
B_{2q}=\frac{1}{n^{2q}}\sum_{\underline{j}\in
  \mathcal{D}_{q+1}}\sum_{\Omega \in \mathcal{R}_{2q}}\sum_{\Pi \in \mathcal{S}_{2q}}A^{\underline{j}}(\Omega_{1},\Pi_{1}) = 
O(\frac{1}{n^{q-1}}).
\end{eqnarray*}

Finally, we briefly discuss the terms $B_{2\underline{b}}$ which lie between the two extremes $B_{2,\ldots,2}$ and 
$B_{2q}$. Since $B_{2\underline{b}}$ is the product of
$2b_{1},\ldots,2b_{m}$ cumulants, 
by Lemma \ref{lemma:basic}(ii)
at most 
$\sum_{j=1}^{m}(b_{j}+1) = q+m$ elements of $\underline{j}$ can be
different. Thus 
\begin{eqnarray*}
B_{2\underline{b}}=\frac{1}{n^{2q}}\sum_{i=q}^{q+m}\sum_{\underline{j}\in
  \mathcal{D}_{i}}\sum_{\pi_{\underline{2b}}\in
  \mathcal{P}_{\underline{2b}}}A_{2\underline{b}}^{\underline{j}}(\pi_{2\underline{b}}).
\end{eqnarray*} 
%where $\mathcal{D}_{q+m}$ is the set where at most $\sum_{j=1}^{m}(b_{j}+1)$ elements of 
%$\underline{j}$ are different.  
By expanding the cumulants in terms of the cumulants of covariances
conditioned on the location 
(which is due to Gaussianity of the random field, see for example, (\ref{eq:Ajj})) we have 
\begin{eqnarray*}
B_{2\underline{b}}=\frac{1}{n^{2q}}\sum_{i=q}^{q+m}\sum_{\underline{j}\in
  \mathcal{D}_{i}}\sum_{\pi_{2\underline{b}}\in
  \mathcal{P}_{2\underline{b}}}\sum_{\Omega \in \mathcal{R}_{2\underline{b}}}
A_{2\underline{b}}^{\underline{j}}(\pi_{2\underline{b}},\Omega),
\end{eqnarray*}
where $\mathcal{R}_{2\underline{b}}$ is the set of all paired partitions of $\{(1,\ldots,2b_{1}),
\ldots,(2b_{m-1}+1,\ldots,2b_{m})\}$. The leading terms are the highest order expectations. 
This term leads to a matrix equation for the $k$'s within the sinc
functions, where the rank of the corresponding matrix is at least $(m-1)$ (we do not give a formal proof of this). 
Therefore, $B_{2\underline{b}} = O(\frac{\log^{2(m-1)}(a)}{n^{q-m}\lambda^{m-1}}) = O(\frac{\log^{2(q-2)}(a)}{\lambda^{q-1}})$
(since $n>>\lambda/\log^{2}(a)$). This concludes the proof of (\ref{eq:cumqa1}). 

The proof of (\ref{eq:cumqa2}) is identical and we omit the
details. 

To prove the result for $d>1$, (\ref{eq:cumqa2d}) we use the same 
method, the main difference is that the spectral density
function in (\ref{eq:A222}) is a multivariate function of dimension
$d$, there are $2dp$ sinc functions and the integral is over
$\mathbb{R}^{dp}$, however the analysis is identical. \hfill $\Box$

\vspace{3mm}
{\bf PROOF of Theorem \ref{theorem:CLT}}
By using the well known identities
\begin{eqnarray}
\label{eq:cov-standard}
\cov(\Re A, \Re B) &=&  \frac{1}{2}\big( \Re \cov(A,B) + \Re
\cov(A,\bar{B})\big) \nonumber\\
\cov(\Im A, \Im B) &=&  \frac{1}{2}\big( \Re \cov(A,B) -
\Re\cov(A,\bar{B})\big), \nonumber\\ 
\cov(\Re A, \Im B) &=&  \frac{-1}{2}\big( \Im \cov(A,B) - \Im
\cov(A,\bar{B})\big),
\end{eqnarray}
and equation (\ref{eq:asymptoticUniform}), we immediately obtain
\begin{eqnarray*}
\lambda^{d}\var[\Re\widetilde{Q}_{a,\lambda}(g;0)] =
\frac{1}{2}\left(C_{1}(0) + \Re C_{2}(0) \right) + O(\ell_{\lambda,a,n}),
\end{eqnarray*} 
\begin{eqnarray*}
\lambda^{d}\cov\left[\Re\widetilde{Q}_{a,\lambda}(g;\rb_{1}),\Re\widetilde{Q}_{a,\lambda}(g;\rb_{2})\right]=
\left\{
\begin{array}{cc}
\frac{\Re}{2}C_{1}(\ob_{\rb}) + O(\ell_{\lambda,a,n})
& \rb_{1}=\rb_{2}(=\rb) \\
\frac{\Re}{2}C_{2}(\ob_{\rb}) + O(\ell_{\lambda,a,n}) &
\rb_{1} = -\rb_{2}(=\rb) \\
O(\ell_{\lambda,a,n}) & \textrm{ otherwise }
\end{array}
\right.
\end{eqnarray*}
\begin{eqnarray*}
\lambda^{d}\cov\left[\Im\widetilde{Q}_{a,\lambda}(g;\rb_{1}),\Im\widetilde{Q}_{a,\lambda}(g;\rb_{2})\right]=
\left\{
\begin{array}{cc}
\frac{\Re}{2}C_{1}(\ob_{\rb}) + O(\ell_{\lambda,a,n})
& \rb_{1}=\rb_{2}(=\rb) \\
\frac{-\Re}{2}C_{2}(\ob_{\rb}) + O(\ell_{\lambda,a,n})) &
\rb_{1} = -\rb_{2}(=\rb) \\
O(\ell_{\lambda,a,n}) & \textrm{ otherwise }
\end{array}
\right.
\end{eqnarray*}
and 
\begin{eqnarray*}
\lambda^{d}\cov\left[\Re\widetilde{Q}_{a,\lambda}(g;\rb_{1}),\Im\widetilde{Q}_{a,\lambda}(g;\rb_{2})\right]=
\left\{
\begin{array}{cc}
O(\ell_{\lambda,a,n}) & \rb_{1}\neq -\rb_{2} \\
\frac{\Im}{2}C_{2}(\ob_{\rb}) + O(\ell_{\lambda,a,n}) &
\rb_{1} = -\rb_{2}(=\rb) \\
\end{array}
\right.
\end{eqnarray*}
Similar expressions for the covariances of  $\Re
Q_{a,\lambda}(g;\rb_{})$ and $\Im
Q_{a,\lambda}(g;\rb_{})$ can also be derived. 

Finally, asymptotic normality of
$\Re\widetilde{Q}_{a,\lambda}(g;\rb_{})$ and
$\Im\widetilde{Q}_{a,\lambda}(g;\rb_{})$ follows from Lemma \ref{lemma:cumulantq}. Thus giving (\ref{eq:cltQtilde}). 
%Similar
%results as those given in Theorem
%\ref{theorem:cumulants} can also be derived for the cumulants of
%$Q_{a,\lambda}(g;\rb)$ (we do not give the details here, but they
%follow the same argument used for
%$\widetilde{Q}_{a,\lambda}(g;\rb)$), thus giving (\ref{eq:cltQ}).
\hfill $\Box$

\vspace{3mm}

{\bf PROOF of Theorem \ref{theorem:cumulantsnonuniform}}  The proof is
  identical to the proof Lemma \ref{lemma:cumulantq}, we omit
  the details. 
\hfill $\Box$

\vspace{3mm}
{\bf PROOF of Theorem \ref{theorem:CLTnu}} The proof is similar to the
proof of Theorem \ref{theorem:CLT}, we omit the details.

%% file: B_3.tex
\section{Additional proofs}\label{sec:b3}

In this section we prove  the remaining results required  in this paper. 

For example,  Theorems \ref{lemma:mean-stat}, \ref{theorem:nonuniformmean},
\ref{theorem:nonuniformvar}(i,iii),
\ref{theorem:nonGaussiannonUniform} and 
\ref{theorem:variance-nongaussian}, 
involve replacing
sums with integrals. In the case that the frequency domain is
increasing stronger assumptions are required than in the case of fixed
frequency domain. We state the required result in the following
lemma. 
\begin{lemma}\label{lemma:sum-integral}
Let us suppose the function $g_{1},g_{2}$ are 
bounded ($\sup_{\ob\in \mathbb{R}^{d}}|g_{1}(\ob)|<\infty$ and $\sup_{\ob\in \mathbb{R}^{d}}|g_{2}(\ob)|<\infty$) and for all
$1\leq j\leq d$, $\sup_{\ob\in \mathbb{R}^{d}}|\frac{\partial
  g_{1}(\ob)}{\partial \omega_{j}}|<\infty$ and $\sup_{\ob\in \mathbb{R}^{d}}|\frac{\partial
  g_{2}(\ob)}{\partial \omega_{j}}|<\infty$. 
\begin{itemize}
\item[(i)] Suppose $a/\lambda = C$ (where $C$ is a fixed finite
  constant) and $h$ is a bounded function whose first partial
  derivative $1\leq j\leq d$, $\sup_{\ob\in \mathbb{R}^{d}}|\frac{\partial
  h(\ob)}{\partial \omega_{j}}|<\infty$. Then we have 
\begin{eqnarray*}
\left|\frac{1}{\lambda^{d}}\sum_{\kb=-a}^{a}g_{1}(\ob_{\kb})h(\ob_{\kb}) -
\frac{1}{(2\pi)^{d}}\int_{2\pi[-C,C]^{d}}g_{1}(\ob)h(\ob)\ob
\right| \leq K\lambda^{-1},
\end{eqnarray*}
where $K$ is a finite constant independent of $\lambda$.
\item[(ii)] Suppose $a/\lambda \rightarrow \infty$ as $\lambda
  \rightarrow \infty$. Furthermore, 
$h(\ob)\leq \beta_{\delta}(\ob)$ and  for all $1 \leq j\leq d$ 
the partial derivatives satisfy 
$|\frac{\partial h(\ob)}{\partial \omega_{j}}|\leq 
\beta_{\delta}(\ob)$. Then uniformly over $a$ we have
\begin{eqnarray*}
&&\bigg|\frac{1}{\lambda^{d}}
\sum_{k=-a}^{a}g_{1}(\ob_{\kb})
h\left(\ob_{k}\right)
-   \frac{1}{(2\pi)^{d}}\int_{2\pi[-a/\lambda,a/\lambda]^{d}}g_{1}(\ob_{})
h\left(\ob\right)d\ob\bigg|
\leq K\lambda^{-1}
\end{eqnarray*}
\item[(iii)] Suppose $a/\lambda \rightarrow \infty$ as $\lambda
  \rightarrow \infty$. Furthermore, 
$f_{4}(\ob_{1},\ob_{2},\ob_{3})\leq \beta_{\delta}(\ob_{1})\beta_{\delta}(\ob_{2})\beta_{\delta}(\ob_{3})$ and  for all $1 \leq j\leq 3d$ 
the partial derivatives satisfy 
$|\frac{\partial f_{4}(\ob_{1},\ob_{2},\ob_{3})}{\partial \omega_{j}}|\leq 
\beta_{\delta}(\ob)$. 
\begin{eqnarray*}
&&\bigg|\frac{1}{\lambda^{2d}}\sum_{\kb_{1},\kb_{2}=-a}^{a}g_{1}(\ob_{\kb_{1}})g_{2}(\ob_{\kb_{2}})
f_{4}\left(\ob_{\kb_{1}+\rb_{1}},\ob_{\kb_2},\ob_{\kb_2+ \rb_2}\right)-\\
&&\qquad\frac{1}{(2\pi)^{2d}}\int_{2\pi[-a/\lambda,a/\lambda]^{d}}\int_{2\pi[-a/\lambda,a/\lambda]^{d}}g_{1}(\ob_{1})g_{2}(\ob_{2})
f_{4}\left(\ob_{1}+\ob_{\rb_{1}},\ob_{2},\ob_{2}+
    \omega_{\rb_2}\right)\bigg| \leq K\lambda^{-1}.
\end{eqnarray*}
\end{itemize}
\end{lemma}
PROOF. We first prove the result in the univariate case. 
We expand the difference between sum and integral
\begin{eqnarray*}
&&\frac{1}{\lambda}\sum_{k=-a}^{a}g_{1}(\omega_{k})h(\omega_{k}) -
\frac{1}{2\pi}\int_{-2\pi a/\lambda}^{2\pi a/\lambda}g_{1}(\omega)h(\omega)d\omega \\
&=&\frac{1}{\lambda}\sum_{k=-a}^{a}g_{1}(\omega_{k})h(\omega_{k}) -
\frac{1}{2\pi}\sum_{k=-a}^{a-1}\int_{0}^{2\pi/\lambda}g_{1}(\omega_{k}+\omega)h(\omega_{k}+\omega)d\omega.
\end{eqnarray*}
By applying the mean value theorem for integrals to the integral above
we have 
\begin{eqnarray*}
&=&\frac{1}{\lambda}\sum_{k=-a}^{a}g_{1}(\omega_{k})h(\omega_{k}) -
\frac{1}{2\pi}\sum_{k=-a}^{a-1}\int_{0}^{2\pi/\lambda}g_{1}(\omega_{k}+\omega)h(\omega_{k}+\omega)d\omega\\
&=&\frac{1}{\lambda}\sum_{k=-a}^{a-1}\left[g_{1}(\omega_{k})h(\omega_{k}) -
g_{1}(\omega_{k}+\overline{\omega_{k}})h(\omega_{k}+\overline{\omega_{k}})\right]
+ \frac{1}{\lambda}g_{1}(\omega_{a})h(\omega_{a})
\end{eqnarray*}
where $\overline{\omega_{k}}\in [0,\frac{2\pi}{\lambda}]$.
Next, by applying the mean value theorem to the difference above we have
\begin{eqnarray}
\label{eq:hgtilde}
\left|\frac{1}{\lambda}\sum_{k=-a}^{a}\left[g_{1}(\omega_{k})h(\omega_{k}) -
g_{1}(\omega_{k}+\overline{\omega_{k}})h(\omega_{k}+\overline{\omega_{k}})\right]\right| 
\leq \frac{1}{\lambda^{2}}\sum_{k=-a}^{a}\left[
  g^{\prime}_{1}(\widetilde{\omega}_{k})h(\widetilde{\omega}_{k})+
g_{1}(\widetilde{\omega}_{k})
h^{\prime}(\widetilde{\omega}_{k})\right]
\end{eqnarray}
where $\widetilde{\omega}_{k}\in [\omega_{k},\omega_{k}+\omega]$ (note this is analogous to the
expression given in \citeA{b:bri-81}, Exercise 1.7.14). 
%We note that
%for a given $k$ we have
%\begin{eqnarray}
%\label{eq:hgtilde1}
%&&\left|\frac{1}{2\pi}\int_{0}^{2\pi/\lambda}
%\left(g^{\prime}_{1}(\widetilde{\omega}_{k,\omega})h(\widetilde{\omega}_{k,\omega})+g_{1}(\widetilde{\omega}_{k,\omega})
%h^{\prime}(\widetilde{\omega}_{k,\omega})\right)d\omega\right| \nonumber\\
%&\leq& \frac{1}{\lambda}\sup_{\omega_{k}\leq
 % \omega \leq
 % \omega_{k+1}}(|g^{\prime}_{1}(\omega)h(\omega)|+|g_{1}(\omega)h^{\prime}(\omega)|),
%\end{eqnarray}
%which will be useful in the proofs below. 

Under the condition that $a =
C\lambda$, $g_{1}(\omega)$ and $h(\omega)$ are bounded and using (\ref{eq:hgtilde})  it is
clear that 
\begin{eqnarray*}
&&\left|\frac{1}{\lambda}\sum_{k=-a}^{a}g_{1}(\omega_{k})h(\omega_{k}) -
\frac{1}{2\pi}\int_{-2\pi a/\lambda}^{2\pi
  a/\lambda}g_{1}(\omega)h(\omega)d\omega \right|\\
&\leq& \frac{\sup_{\omega}|h^{\prime}(\omega)g_{1}(\omega)|+\sup_{\omega}|h(\omega)g_{1}^{\prime}(\omega)|}{\lambda^{2}}\sum_{k=-a}^{a}1
= C\lambda^{-1}.
\end{eqnarray*}
For $d=1$, this proves (i). 

In the case that $a/\lambda \rightarrow \infty$ as $\lambda
\rightarrow \infty$, we use that 
$h$ and $h^{\prime}$ are dominated by a
montonic function and that $g_{1}$ is bounded. Thus by using
(\ref{eq:hgtilde})  we have 
\begin{eqnarray*}
&&\left|\frac{1}{\lambda}\sum_{k=-a}^{a}g_{1}(\omega_{k})h(\omega_{k}) -
\frac{1}{2\pi}\int_{-2\pi a/\lambda}^{2\pi
  a/\lambda}g_{1}(\omega)h(\omega)d\omega\right| 
\leq  \frac{1}{\lambda^{2}}\sum_{k=-a}^{a}\sup_{\omega_{k}\leq
  \omega \leq
  \omega_{k+1}}(|g^{\prime}_{1}(\omega)h(\omega)|+|g_{1}(\omega)h^{\prime}(\omega)|)
\\
&\leq&\frac{2(\sup_{\omega}|g_{1}(\omega)|+\sup_{\omega}|g_{1}^{\prime}(\omega)|)}{\lambda^{2}}\sum_{k=0}^{a}
\beta_{\delta}(\omega_{k}) \leq
\frac{C}{\lambda}\int_{0}^{\infty}\beta_{\delta}(\omega)d\omega = O(\lambda^{-1}).
\end{eqnarray*}
%where $C$ is a finite constant which does not depend on $\lambda$ or $a$.
For $d=1$,  this proves (ii).

To prove the result for $d=2$ we take
differences, that is 
\begin{eqnarray*}
&&\frac{1}{\lambda^{2}}\sum_{k_{1}=-a}^{a}\sum_{k_{2}=-a}^{a}g(\omega_{k_1},\omega_{k_2})h(\omega_{k_1},\omega_{k_{2}})
-
\frac{1}{(2\pi)^{2}}\int_{-a/\lambda}^{a/\lambda}\int_{-a/\lambda}^{a/\lambda}g(\omega_{1},\omega_{2})
h(\omega_{1},\omega_{2})d\omega_{1}d\omega_{2}
\\
&=& \frac{1}{\lambda}\sum_{k_{1}=-a}^{a}\left(\frac{1}{\lambda}
\sum_{k_{2}=-a}^{a}g(\omega_{k_1},\omega_{k_2})h(\omega_{k_1},\omega_{k_{2}})d\omega_{2}
-
\frac{1}{2\pi}\int_{-a/\lambda}^{a/\lambda}g(\omega_{k_1},\omega_{2})h(\omega_{k_1},\omega_{2})\right)\\
&& + \frac{1}{2\pi}\int_{-a/\lambda}^{a/\lambda}\left(
  \frac{1}{\lambda}\sum_{k_{1}=-a}^{a}
g(\omega_{k_1},\omega_{2})h(\omega_{k_1},\omega_{2}) - \frac{1}{(2\pi)^{2}}
\int_{-a/\lambda}^{a/\lambda}g(\omega_{1},\omega_{2})h(\omega_{1},\omega_{2})d\omega_{2}\right)d\omega_{1}.
\end{eqnarray*}
For each of the terms above we apply the method described for the case
$d=1$; for the first term we take the partial derivative over
$\omega_{2}$ and the for the second term we take the partial
derivative over $\omega_{1}$. This method can easily be generalized to
the case $d>2$.  The proof of (iii) is identical to the proof of (ii). 

We mention that the assumptions on the derivatives (used replace sum
with integral) can  be
relaxed to that of bounded variation of the function. However, since
we require the bounded derivatives to decay at certain rates (to prove
other results) we do not relax the assumption here. \hfill $\Box$

\vspace{3mm}

We now prove the claims at the start of Appendix \ref{sec:variance}. 

\begin{lemma}\label{lemma:covariance}
Suppose Assumptions \ref{assum:S}, \ref{assum:uniform} and  
\begin{itemize}
\item[(i)] Assumptions \ref{assum:G}(i) and \ref{assum:GG}(a,c) hold. Then we have 
\begin{eqnarray*}
&&\lambda^{d}\cov\left[\widetilde{Q}_{a,\lambda}(g;\rb_{1}),\widetilde{Q}_{a,\lambda}(g;\rb_{2})\right]=    
\left\{
\begin{array}{cl}
C_{1}(\ob_{\rb})
+ O(\frac{1}{\lambda^{}}+\frac{\lambda^{d}}{n}) & \rb_{1} = \rb_{2} (=\rb) \\
O(\frac{1}{\lambda^{}}+\frac{\lambda^{d}}{n}) & \rb_{1} \neq \rb_{2} 
\end{array}
\right.
\end{eqnarray*}
and 
\begin{eqnarray*}
&& \lambda^{d}\cov\left[\widetilde{Q}_{a,\lambda}(g;\rb_{1}),\overline{\widetilde{Q}_{a,\lambda}(g;\rb_{2})}\right]=  
\left\{
\begin{array}{cl}
C_{2}(\ob_{\rb})
+ O(\frac{1}{\lambda^{}}+\frac{\lambda^{d}}{n}) & \rb_{1} = -\rb_{2} (= \rb) \\
O(\frac{1}{\lambda^{}}+\frac{\lambda^{d}}{n}) & \rb_{1} \neq -\rb_{2} 
\end{array}
\right.
\end{eqnarray*}
where
\begin{eqnarray}
\label{eq:Cr}
C_{1}(\ob_{\rb})  = C_{1,1}(\ob_{\rb}) + C_{1,2}(\ob_{\rb})
\textrm{ and } C_{2}(\ob_{\rb})  = C_{2,1}(\ob_{\rb}) + C_{2,2}(\ob_{\rb}),
\end{eqnarray}
with 
\begin{eqnarray}
\label{eq:Cr11}
C_{1,1}(\ob_{\rb}) 
&=&\frac{1}{(2\pi)^{d}}\int_{2\pi[-a/\lambda,a/\lambda]^{d}}f(\ob)f(\ob+\ob_{\rb}) 
|g(\ob)|^{2} d\ob, \nonumber\\
C_{1,2}(\ob_{\rb}) 
&=&\frac{1}{(2\pi)^{d}} \int_{\mathcal{D}_{\rb}}f(\ob)f(\ob+\ob_{\rb}) 
g(\ob)\overline{g(-\ob-\ob_{\rb})} d\ob, \nonumber\\
C_{2,1}(\ob_{\rb}) 
&=&\frac{1}{(2\pi)^{d}}\int_{2\pi[-a/\lambda,a/\lambda]^{d}}f(\ob)f(\ob+\ob_{\rb}) 
g(\ob)g(-\ob) d\ob,\nonumber\\ 
 C_{2,2}(\ob_{\rb}) 
&=&\frac{1}{(2\pi)^{d}}\int_{\mathcal{D}_{\rb}}f(\ob)f(\ob+\ob_{\rb}) 
g(\ob)g(\ob+\ob_{\rb}) d\ob, 
\end{eqnarray}
where the integral is defined as 
\\*
$\int_{\mathcal{D}_{\rb}} = 
\int_{2\pi\max(-a,-a-r_{1})/\lambda}^{2\pi\min(a,a-r_{1})/\lambda}\ldots
\int_{2\pi\max(-a,-a-r_{d})/\lambda}^{2\pi\min(a,a-r_{d})/\lambda}$
(note that $C_{1,1}(\ob_{\rb}) $ and $C_{1,2}(\ob_{\rb})$ are real).
\item[(ii)] Assumptions \ref{assum:G}(ii) and \ref{assum:GG}(b) hold. Then 
\begin{eqnarray*}
\lambda^{d}\cov\left[\widetilde{Q}_{a,\lambda}(g;\rb_{1}),\widetilde{Q}_{a,\lambda}(g;\rb_{2})\right] =
A_{1}(\rb_{1},\rb_{2})+A_{2}(\rb_{1},\rb_{2}) + O(\frac{\lambda^{d}}{n}),
\end{eqnarray*}
\begin{eqnarray*}
\lambda^{d}\cov\left[\widetilde{Q}_{a,\lambda}(g;\rb_{1}),\overline{\widetilde{Q}_{a,\lambda}(g;\rb_{2})}\right]
= 
A_{3}(\rb_{1},\rb_{2})+A_{4}(\rb_{1},\rb_{2}) + O(\frac{\lambda^{d}}{n}),
\end{eqnarray*}
where 
\begin{eqnarray*}
A_{1}(\rb_{1},\rb_{2}) &=& \frac{1}{\pi^{2d}\lambda^{d}}
\sum_{\mb=-2a}^{2a}\sum_{\kb=\max(-a,-a+\mb)}^{\min(a,a+\mb)} 
\int_{\mathbb{R}^{2d}} f(\frac{2\ubb}{\lambda} -\ob_{\kb_{}})
f(\frac{2\vbb}{\lambda} + \ob_{\kb_{}}+\ob_{\rb_{1}}) \\
&& g(\ob_{\kb})\overline{g(\ob_{\kb}-\ob_{\mb})}
\Sinc(\ubb-\mb \pi)\Sinc(\vbb+(\mb+\rb_{1}-\rb_{2})\pi)
\Sinc(\ubb)\Sinc(\vbb)d\ubb d\vbb
\end{eqnarray*}
\begin{eqnarray*}
A_{2}(\rb_{1},\rb_{2}) &=& \frac{1}{\pi^{2d}\lambda^{d}}
\sum_{\mb=-2a}^{2a}\sum_{\kb=\max(-a,-a+\mb)}^{\min(a,a+\mb)} 
\int_{\mathbb{R}^{2d}} f(\frac{2\ubb}{\lambda} -\ob_{\kb_{}})
f(\frac{2\vbb}{\lambda} + \ob_{\kb_{}}+\ob_{\rb_{1}}) \\
&& g(\ob_{\kb})\overline{g(\ob_{\mb}-\ob_{\kb})}
\Sinc(\ubb-(\mb+\rb_{2}) \pi)\Sinc(\vbb+(\mb+\rb_{1})\pi)
\Sinc(\ubb)\Sinc(\vbb)d\ubb d \vbb
\end{eqnarray*}
\begin{eqnarray*}
A_{3}(\rb_{1},\rb_{2}) &=& \frac{1}{\pi^{2d}\lambda^{d}}
\sum_{\mb=-2a}^{2a}\sum_{\kb=\max(-a,-a+\mb)}^{\min(a,a+\mb)} 
\int_{\mathbb{R}^{2d}} f(\frac{2\ubb}{\lambda} -\ob_{\kb_{}})
f(\frac{2\vbb}{\lambda} + \ob_{\kb_{}}+\ob_{\rb_{1}}) \\
&& g(\ob_{\kb})g(\ob_{\mb}-\ob_{\kb})
\Sinc(\ubb+\mb \pi)\Sinc(\vbb+(\mb+\rb_{2}+\rb_{1})\pi)
\Sinc(\ubb)\Sinc(\vbb)d\ubb d\vbb
\end{eqnarray*}
\begin{eqnarray*}
A_{4}(\rb_{1},\rb_{2}) &=& \frac{1}{\pi^{2d}\lambda^{d}}
\sum_{\mb=-2a}^{2a}\sum_{\kb=\max(-a,-a+\mb)}^{\min(a,a+\mb)} 
\int_{\mathbb{R}^{2d}} f(\frac{2\ubb}{\lambda} -\ob_{\kb_{}})
f(\frac{2\vbb}{\lambda} + \ob_{\kb_{}}+\ob_{\rb_{1}}) \\
&& g(\ob_{\kb})g(\ob_{\kb}-\ob_{\mb})
\Sinc(\ubb-(\mb-\rb_{2}) \pi)\Sinc(\vbb+(\mb+\rb_{1})\pi)
\Sinc(\ubb)\Sinc(\vbb)d\ubb d\vbb
\end{eqnarray*}
where $\kb=\max(-a,-a-\mb)=\{k_{1}=\max(-a,-a-m_{1}),\ldots,k_{d}=\max(-a,-a-m_{d})\}$, 
$\kb=\min(-a+\mb)=\{k_{1}=\max(-a,-a+m_{1}),\ldots,k_{d}=\min(-a+m_{d})\}$.
% $C_{4} = n(n-1)(n-2)(n-3)/n^{4}$.  
\item[(iii)]  Assumptions \ref{assum:G}(ii) and \ref{assum:GG}(b) hold. Then we have 
\begin{eqnarray*}
\lambda^{d}\sup_{a}\left|\cov\left[\widetilde{Q}_{a,\lambda}(g;\rb_{}),\widetilde{Q}_{a,\lambda}(g;\rb_{})\right]\right|
<\infty \quad\textrm{and}\quad
\lambda^{d}\sup_{a}\left|\cov\left[\widetilde{Q}_{a,\lambda}(g;\rb_{}),
\overline{\widetilde{Q}_{a,\lambda}(g;-\rb_{})}\right]\right|<\infty,
\end{eqnarray*}
if $\lambda^{d}/n\rightarrow c$ ($0\leq c<\infty$) as $\lambda\rightarrow \infty$ and
$n\rightarrow \infty$. 
\end{itemize}
\end{lemma}
{\bf PROOF} We prove the result in the case $d=1$ (the proof for $d>1$ is
identical).
We first prove (i). By using indecomposable partitions,  Theorem \ref{lemma:nonuniformdft} and 
Lemma \ref{lemma:cumulantsA} and noting that the fourth
order cumulant is of order $O(1/n)$, it is
straightforward to show that 
\begin{eqnarray*}
&&\lambda\cov\left[\widetilde{Q}_{a,\lambda}(g;r_{1}),\widetilde{Q}_{a,\lambda}(g;r_{2})\right]  \\
&=&  
\frac{1}{\lambda}\sum_{k_{1},k_{2}=-a}^{a}g(\omega_{k_{1}})\overline{g(\omega_{k_{2}})}\bigg[
\cov\bigg(J_{n}(\omega_{k_{1}})\overline{J_{n}(\omega_{k_{1}+r_{1}})}-\frac{\lambda^{d}}{n}\sum_{j=1}^{n}Z(s_{j})^{2}e^{-is_{j}\omega_{r}}, \\
&& J_{n}(\omega_{k_{2}})\overline{J_{n}(\omega_{k_{2}+r_{2}})}-\frac{\lambda^{d}}{n}\sum_{j=1}^{n}Z(s_{j})^{2}e^{-is_{j}\omega_{r}}\bigg) \bigg]\\
&=& 
\left\{
\begin{array}{cc}
\frac{1}{\lambda}\sum_{k=-a}^{a}
f(\omega_{k})f(\omega_{k}+\omega_{r})g(\omega_{k})\overline{g(\omega_{k})} & \\
 +
 \frac{1}{\lambda}\sum_{k=\max(-a,-a-r)}^{\min(a,a-r)}f(\omega_{k})f(\omega_{k}+\omega_{r})g(\omega_{k})\overline{g(-\omega_{k+r})}+ 
O(\frac{a}{\lambda^{2}}+\frac{\lambda}{n}) & r_{1}=r_{2} \\
O(\frac{a}{\lambda^{2}}+\frac{\lambda}{n}) & r_{1} \neq r_{2} 
\end{array}
\right.
\end{eqnarray*} 
Since $a=C\lambda$ we have 
$O(\frac{a}{\lambda^{2}}) = O(\frac{1}{\lambda})$ and 
by replacing sum with integral (using Lemma \ref{lemma:sum-integral}(i)) we have 
\begin{eqnarray*}
&&\lambda\cov\left[\widetilde{Q}_{a,\lambda}(g;r_{1}),\widetilde{Q}_{a,\lambda}(g;r_{2})\right]  =
\left\{
\begin{array}{cc}
C_{1}(\omega_{r})
+ O(\frac{1}{\lambda^{}}+\frac{\lambda}{n}) 
& r_{1}=r_{2}(=r) \\
O(\frac{1}{\lambda^{}}+\frac{\lambda}{n}) & r_{1} \neq r_{2} 
\end{array}
\right.
\end{eqnarray*} 
where 
\begin{eqnarray*}
&&C_{1}(\omega_{r}) = C_{11}(\omega_{r})+C_{12}(\omega_{r}) + O\left(\frac{1}{\lambda}\right), 
\end{eqnarray*}
with 
\begin{eqnarray*}
C_{11}(\omega_{r})&=&\frac{1}{2\pi}\int_{-2\pi a/\lambda}^{2\pi a/\lambda}
f(\omega)f(\omega+\omega_{r})g(\omega)\overline{g(\omega)}d\omega \\
C_{12}(\omega_{r})&=& \frac{1}{2\pi}\int_{2\pi\max(-a,-a-r)/\lambda}^{2\pi\min(a,a-r)/\lambda}f(\omega)f(\omega+\omega_{r})g(\omega)\overline{g(-\omega-\omega_{r})}
d\omega.
\end{eqnarray*}

We now show that $C_{1}(\omega_{r}) = C_{11}(\omega_{r}) +
C_{12}(\omega_{r})$ is real. It is clear
that $C_{11}(\omega_{r})$ is real. Thus we focus on $C_{12}(\omega_{r})$.
To do this, we write $g(\omega) = g_{1}(\omega)+ig_{2}(\omega)$,
therefore 
\begin{eqnarray*}
\Im g(\omega)\overline{g(-\omega-\omega_{r})} =
%[g_{1}(\omega)g_{1}(-\omega-\omega_{r})+g_{2}(\omega)g_{2}(-\omega-\omega_{r})]+
[g_{2}(\omega)g_{1}(-\omega-\omega_{r})-g_{1}(\omega)g_{2}(-\omega-\omega_{r})].
\end{eqnarray*} 
Substituting the above into $\Im C_{12}(\omega_{r})$ gives us $\Im C_{12}(r)=  [C_{121}(r)+C_{122}(r)]$ where
\begin{eqnarray*}
C_{121}(r)&=&\frac{1}{2\pi}\int_{2\pi\max(-a,-a-r)/\lambda}^{2\pi\min(a,a-r)/\lambda}f(\omega)f(\omega+\omega_{r})g_{2}(\omega)
g_{1}(-\omega-\omega_{r})d\omega \\
C_{122}(r)&=& - \frac{1}{2\pi}
\int_{2\pi\max(-a,-a-r)/\lambda}^{2\pi\min(a,a-r)/\lambda}f(\omega)f(\omega+\omega_{r})g_{1}(\omega)g_{2}(-\omega-\omega_{r})d\omega
\end{eqnarray*} 
We will show that $C_{122}(r)=-C_{121}(r)$. Focusing on 
$C_{122}$ and making the change of variables $u=-\omega-\omega_{r}$ gives us 
\begin{eqnarray*}
C_{122} = \frac{1}{2\pi}\int^{2\pi\max(-a,-a-r)/\lambda}_{2\pi\min(a,a-r)/\lambda}f(u+\omega_{r})f(-u)g_{1}(-u-\omega_{r})g_{2}(u)du,
\end{eqnarray*}
noting that the spectral density function is symmetric with $f(-u)=f(u)$, and that 
$\int^{2\pi\max(-a,-a-r)/\lambda}_{2\pi\min(a,a-r)/\lambda} = 
-\int_{2\pi\max(-a,-a-r)/\lambda}^{2\pi\min(a,a-r)/\lambda}$. Thus we have $\Im
C_{12}(r)=0$, which shows that
$C_{1}(\omega_{r})$ is real. 
The proof of
$\lambda^{d}\cov[\widetilde{Q}_{a,\lambda}(g;\rb_{1}),\overline{\widetilde{Q}_{a,\lambda}(g;\rb_{2})}]$
is the same. Thus, we have proven (i).

To prove (ii) we first expand
$\cov[\widetilde{Q}_{a,\lambda}(g;r_{1}),\widetilde{Q}_{a,\lambda}(g;r_{2})]$
to give 
\begin{eqnarray*}
&&\lambda\cov[\widetilde{Q}_{a,\lambda}(g;r_{1}),\widetilde{Q}_{a,\lambda}(g;r_{2})]
\nonumber\\
&=& \lambda\sum_{k_{1},k_{2}=-a}^{a}g(\omega_{k_{1}})\overline{g(\omega_{k_{2}})}
\sum_{j_{1},\ldots,j_{4}=1}^{n}\bigg( 
\Ex\big[c(s_{j_{1}}-s_{j_{3}})e^{is_{j_1}\omega_{k_{1}}-is_{j_3}\omega_{k_{2}}}\big] 
\Ex\big[c(s_{j_{2}}-s_{j_{4}})e^{-is_{j_2}\omega_{k_{1}+r_{1}}+is_{j_4}\omega_{k_{2}+r_{2}}}\big] +\\ 
&&\Ex\big[c(s_{j_{1}}-s_{j_{4}})e^{is_{j_1}\omega_{k_{1}}+is_{j_4}\omega_{k_{2}+r_{2}}}\big]
\Ex\big[c(s_{j_{2}}-s_{j_{3}})e^{-is_{j_2}\omega_{k_{1}+r_{1}}-is_{j_3}\omega_{k_{2}}}\big] + \\
&& \cum\big[Z(s_{j_{1}})e^{i\omega_{k_{1}}s_{j_{1}}},
Z(s_{j_{2}})e^{-is_{j_{2}}(\omega_{k_{1}}+\omega_{r_{2}})},
Z(s_{j_{3}})e^{-i\omega_{k_{2}}s_{j_{3}}},Z(s_{j_{4}})e^{is_{j_{4}}(\omega_{k_{2}}+\omega_{r_{2}})}
\big]\bigg).
\end{eqnarray*}
By applying Lemma \ref{lemma:cumulantsA}, (\ref{eq:cum3}) in \citeA{p:sub-14}, to the
fourth order cumulant we can 
show that 
\begin{eqnarray}
\label{eq:ahatcov}
\lambda\cov\left[\widetilde{Q}_{a,\lambda}(r_{1}),\widetilde{Q}_{a,\lambda}(g;r_{2})\right]
&=& A_{1}(r_{1},r_{2}) + A_{2}(r_{1},r_{2}) + O\left(\frac{\lambda}{n}\right) 
\end{eqnarray}
where
\begin{eqnarray*}
A_{1}(r_{1},r_{2})&=& \lambda\sum_{k_{1},k_{2}=-a}^{a}g(\omega_{k_{1}})\overline{g(\omega_{k_{2}})}
\Ex\big[c(s_{1}-s_{3})\exp(is_{1}\omega_{k_{1}}-is_{3}\omega_{k_{2}})\big]
\times \\
&&\Ex\big[c(s_{2}-s_{4})\exp(-is_{2}\omega_{k_{1}+r_{1}}+is_{4}\omega_{k_{2}+r_{2}})\big]\\
A_{2}(r_{1},r_{2})&=& \lambda\sum_{k_{1},k_{2}=-a}^{a}g(\omega_{k_{1}})\overline{g(\omega_{k_{2}})}
 \Ex\big[c(s_{1}-s_{4})\exp(is_{1}\omega_{k_{1}}+is_{4}\omega_{k_{2}+r_{2}})\big] \times\\
&&\Ex\big[c(s_{2}-s_{3})\exp(-is_{2}\omega_{k_{1}+r_{1}}-is_{3}\omega_{k_{2}})\big].
\end{eqnarray*}
Note that the $O(\lambda/n)$ term include the error $n^{-1}[A_{1}(r_{1},r_{2})]+A_{2}(r_{1},r_{2})]$
(we show below that $A_{1}(r_{1},r_{2})$ and $A_{2}(r_{1},r_{2})$ are both bounded over
$\lambda$ and $a$).
To write $A_{1}(r_{1},r_{2})$ in  the form stated in the lemma we
integrate over $s_{1},s_{2},s_{3}$ and $s_{4}$ to give 
\begin{eqnarray*}
&&A_{1}(r_{1},r_{2}) \\
&=&\lambda
\sum_{k_{1},k_{2}=-a}^{a}g(\omega_{k_{1}})\overline{g(\omega_{k_{2}})}
\frac{1}{\lambda^{4}}\int_{[-\lambda/2,\lambda/2]^{4}}
 c(s_{1}-s_{3})c(s_{2}-s_{4})e^{is_{1}\omega_{k_{1}}-is_{3}\omega_{k_{2}}}
e^{-is_{2}\omega_{k_{1}+r_{1}}+is_{4}\omega_{k_{2}+r_{2}}}ds_{1}\ldots ds_{4}. 
%\times 
% \nonumber\\
%&&\cov\big[Z(s_{2})\exp(-is_{2}\omega_{k_{1}+r_{1}}),Z(s_{4})\exp(-is_{4}\omega_{k_{2}+r_{2}})\big]
\end{eqnarray*}
By using the spectral representation theorem and integrating out $s_{1},\ldots,s_{4}$
we can write the above as
\begin{eqnarray*}
&&A_{1} (r_{1},r_{2}) \\
&=& \frac{\lambda}{(2\pi)^{2}}
\sum_{k_{1},k_{2}=-a}^{a}g(\omega_{k_{1}})\overline{g(\omega_{k_{2}})} 
\int_{-\infty}^{\infty} \int_{-\infty}^{\infty} f(x)f(y)
\sinc\left(\frac{\lambda x}{2} + k_{1}\pi\right)\\
&&\sinc\left(\frac{\lambda y}{2} - (r_{1}+k_{1})\pi\right)
 \sinc\left(\frac{\lambda x}{2}+k_{2}\pi\right)
\sinc\left(\frac{\lambda y}{2} - (r_{2}+k_{2})\pi\right)dxdy \\
&=& \frac{1}{\pi^{2}\lambda}
\sum_{k_{1},k_{2}=-a}^{a} g(\omega_{k_{1}})\overline{g(\omega_{k_{2}})} \int_{-\infty}^{\infty}
\int_{-\infty}^{\infty} 
f(\frac{2u}{\lambda} - \omega_{k_{1}})
f(\frac{2v}{\lambda} + \omega_{k_{1}}+\omega_{r_{1}})\\
&&\times\sinc(u)
\sinc(u + (k_{2}-k_{1})\pi)
\sinc(v)\sinc(v + (k_{1}-k_{2}+r_{1}-r_{2})\pi)dudv, 
\end{eqnarray*}
where the second equality is due to the change of variables 
$u=\frac{\lambda x}{2} + k_{1}\pi$ and $v=\frac{\lambda y}{2} -
(r_{1}+k_{1})\pi$. Finally, by making a change of variables $k =
k_{1}$ and 
$m=k_{1}-k_{2}$ ($k_{1}=-k_{2}+m$) we obtain the expression for 
$A_{1}(r_{1},r_{2})$ given in Lemma \ref{lemma:covariance}.

A similar method can be used to obtain the expression for
$A_{2}(r_{1},r_{2})$
\begin{eqnarray*}
&&A_{2} (r_{1},r_{2}) \\
&=& \frac{\lambda}{(2\pi)^{2}}
\sum_{k_{1},k_{2}=-a}^{a}g(\omega_{k_{1}})\overline{g(\omega_{k_{2}})} 
\int_{-\infty}^{\infty} \int_{-\infty}^{\infty} f(x)f(y)
\sinc\left(\frac{\lambda x}{2} + k_{1}\pi\right)\\
&&\sinc\left(\frac{\lambda y}{2} - (r_{1}+k_{1})\pi\right)
 \sinc\left(\frac{\lambda y}{2}+k_{2}\pi\right)
\sinc\left(\frac{\lambda x}{2} - (r_{2}+k_{2})\pi\right)dxdy \\
&=& \frac{1}{\pi^{2}\lambda}
\sum_{k_{1},k_{2}=-a}^{a} g(\omega_{k_{1}})\overline{g(\omega_{k_{2}})} \int_{-\infty}^{\infty}
\int_{-\infty}^{\infty} 
f(\frac{2u}{\lambda} - \omega_{k_{1}})
f(\frac{2v}{\lambda} + \omega_{k_{1}}+\omega_{r_{1}})\\
&&\times\sinc(u)
\sinc(u - (k_{2}+r_{2}+k_{1})\pi)
\sinc(v)\sinc(v + (k_{2}+k_{1}+r_{1})\pi)dudv. 
\end{eqnarray*}
By making the change of variables $k=k_{1}$ and $m=k_{1}+k_{2}$
($k_{1} = -k_{2}+m$) we
obtain 
the stated expression for $A_{2}(r_{1},r_{2})$.

Finally following the same steps as those above we obtain 
\begin{eqnarray*}
\lambda\cov(\widetilde{Q}_{a,\lambda}(g;r_{1}),\overline{\widetilde{Q}_{a,\lambda}(g;r_{2})})
= A_{3}(r_{1},r_{2}) + A_{4}(r_{1},r_{2}) + O(\frac{\lambda}{n}), 
\end{eqnarray*}
where 
\begin{eqnarray*}
A_{3}(r_{1},r_{2})&=& \frac{1}{\lambda^{3}}
\sum_{k_{1},k_{2}=-a}^{a}g(\omega_{k_{1}})g(\omega_{k_{2}})
\int_{[-\lambda/2,\lambda/2]^{4}}e^{i(s_{1}\omega_{k_{1}}+s_{3}\omega_{k_{2}})}
e^{-is_{2}\omega_{k_{1}+r_{1}}-is_{4}\omega_{k_{2}+r_{2}}}
\\
 &&\times c(s_{1}-s_{3})c(s_{2}-s_{4})ds_{1}ds_{2}ds_{3}ds_{4} \\
A_{4}(r_{1},r_{2}) &=& \frac{1}{\lambda^{3}}
\sum_{k_{1},k_{2}=-a}^{a}g(\omega_{k_{1}})g(\omega_{k_{2}})
\int_{[-\lambda/2,\lambda/2]^{4}}e^{i(s_{1}\omega_{k_{1}}+s_{3}\omega_{k_{2}})}
e^{-is_{2}\omega_{k_{1}+r_{1}}-is_{4}\omega_{k_{2}+r_{2}}}
\\
 && \times c(s_{1}-s_{4})c(s_{2}-s_{3})ds_{1}ds_{2}ds_{3}ds_{4}. 
\end{eqnarray*}
Again by replacing the covariances in  $A_{3}(r_{1},r_{2})$ and
$A_{4}(r_{1},r_{4})$ with their spectral representation gives (ii)  for $d=1$. The result for
$d>1$ is  identical.

It is clear that (iii) is true under Assumption
\ref{assum:G}(a,b). To prove (iii) under Assumption \ref{assum:GG}(a,b)
we will show 
that for $1\leq j \leq 4$, $\sup_{a}|A_{j}(\rb_{1},\rb_{2})|<\infty$.
To do this, we first note that 
by the Cauchy Schwarz inequality we have 
\begin{eqnarray*}
&& \sup_{a,\lambda}
\frac{1}{\pi^{2d}}\left|\frac{1}{\lambda^{d}}\sum_{\kb=\max(-a,-a+\mb)}^{\min(-a+\mb)}
g(\omega_{\kb})\overline{g(\omega_{\kb}-\ob_{m})}  
f(\frac{2\ubb}{\lambda} -\ob_{\kb_{}})
f(\frac{2\vbb}{\lambda} + \ob_{\kb_{}}+\ob_{\rb_{1}})\right| \\
&\leq& C\sup_{\ob}|g(\ob)|^{2}\|f\|_{2}^{2},
\end{eqnarray*}
where $\|f\|_{2}$ is the $L_{2}$ norm of the spectral density
function and $C$ is a finite constant. Thus by taking absolutes of $A_{1}(\rb_{1},\rb_{2})$ we have 
\begin{eqnarray*}
&&|A_{1}(\rb_{1},\rb_{2})| \leq \\ 
&&C\sup_{\ob}|g(\ob)|^{2}\|f\|_{2}^{2}\sum_{\mb=-\infty}^{\infty}\int_{\mathbb{R}^{2d}}
|\Sinc(\ubb-\mb \pi)\Sinc(\vbb+(\mb+\rb_{1}-\rb_{2})\pi)
\Sinc(\ubb)\Sinc(\vbb)|d\ubb d\vbb. 
\end{eqnarray*}
Finally, by using Lemma \ref{lemma:1a}(iii), \citeA{p:sub-14}, we have that 
$\sup_{a}|A_{1}(\rb_{},\rb_{})|<\infty$. By using the same method we can show that 
$\sup_{a}|A_{2}(\rb_{},\rb_{})|,A_{3}(\rb,-\rb),\sup_{a}|A_{4}(\rb_{},\rb_{})|<\infty$. This
completes the proof. \hfill
\hfill $\Box$

\vspace{3mm}

{\bf PROOF of Theorem \ref{theorem:var-est}} 
Making the classical variance-bias decomposition we have 
\begin{eqnarray*}
\Ex\left(\widetilde{V} -
    \lambda^{d}\var[\widetilde{Q}_{a,\lambda}(g;0)]\right)^{2} =
  \var[\widetilde{V} ] + \left(\Ex[\widetilde{V}] - \lambda^{d}\var[\widetilde{Q}_{a,\lambda}(g;0)]\right)^{2}.
\end{eqnarray*}
We first analysis the bias term, in particular $\Ex[\widetilde{V}]$. We
note that by using the expectation and variance result in Theorem
\ref{lemma:mean-stat} and equation
(\ref{eq:asymptoticUniform}), respectively, we have 
 \begin{eqnarray*}
\Ex[\widetilde{V}] &=& 
\frac{\lambda^{d}}{|\mathcal{S}|}\sum_{\rb\in
  \mathcal{S}}^{}\var[\widetilde{Q}_{a,\lambda}(g;\rb)]+\frac{\lambda^{d}}{|\mathcal{S}|}\sum_{\rb\in
  \mathcal{S}}^{}
\underbrace{\big|\Ex[\widetilde{Q}_{a,\lambda}(g;\rb)]\big|^{2}}_{=O(\lambda^{-2d}\prod_{j=1}^{d}(\log\lambda
  + \log|\rb_{j}|)^{2})}\\
 &=& \frac{1}{|\mathcal{S}|}\sum_{\rb \in \mathcal{S}}C_{1}(\ob_{\rb})
 + O\left(\ell_{\lambda,a,n} + \frac{[\log\lambda + \log
   M]^{d}}{\lambda^{d}}\right) \\ 
 &=& C_{1} + O\left(\ell_{a,\lambda,n} + \frac{[\log\lambda + \log
   M]^{d}}{\lambda^{d}} + \frac{|M|}{\lambda}\right).
\end{eqnarray*}
Next we consider $\var[\widetilde{V} ]$, by using the classical
cumulant decomposition we have 
\begin{eqnarray*}
\var[\widetilde{V}] &=& \frac{\lambda^{2d}}{|\mathcal{S}|^{2}}\sum_{\rb_{1},\rb_{2}\in \mathcal{S}}\bigg(
\left|\cov\left[\widetilde{Q}_{a,\lambda}(g;\rb_{1}),\widetilde{Q}_{a,\lambda}(g;\rb_{2})\right]\right|^{2}
+\left|\cov\left[\widetilde{Q}_{a,\lambda}(g;\rb_{1}),\overline{\widetilde{Q}_{a,\lambda}(g;\rb_{2})}\right]\right|^{2} \\
&& + \frac{\lambda^{2d}}{|\mathcal{S}|^{2}}\sum_{\rb_{1},\rb_{2}\in
  \mathcal{S}}\cum\left(\widetilde{Q}_{a,\lambda}(g;\rb_{1}),\overline{\widetilde{Q}_{a,\lambda}(g;\rb_{1})},
\overline{\widetilde{Q}_{a,\lambda}(g;\rb_{2})},\widetilde{Q}_{a,\lambda}(g;\rb_{2})\right)
\\
&&+ O\left( \frac{\lambda^{2d}}{|\mathcal{S}|^{2}}\sum_{\rb_{1},\rb_{2}\in
  \mathcal{S}}\left|\cum\left(\widetilde{Q}_{a,\lambda}(g;\rb_{1}),\overline{\widetilde{Q}_{a,\lambda}(g;\rb_{1})},
\overline{\widetilde{Q}_{a,\lambda}(g;\rb_{2})}\right)\Ex\left(\widetilde{Q}_{a,\lambda}(g;\rb_{2})\right)\right|\right).
\end{eqnarray*}
By substituting the variance/covariance results for
$\widetilde{Q}_{a,\lambda}(\cdot)$ based on uniformly sampled
locations in (\ref{eq:asymptoticUniform}) and 
the cumulant bounds in Lemma \ref{lemma:cumulantq} into the above we have 
\begin{eqnarray*}
\var[\widetilde{V}] &=& \frac{\lambda^{2d}}{|\mathcal{S}|^{2}}\sum_{\rb\in\mathcal{S}}
\left|\var\left[\widetilde{Q}_{a,\lambda}(g;\rb_{})\right]\right|^{2}
+ O\left(\ell_{\lambda,a,n}+\frac{\log^{4d}(a)}{\lambda^{d}}\right) = 
O\left(\frac{1}{|\mathcal{S}|} + \ell_{\lambda,a,n}+\frac{\log^{4d}(a)}{\lambda^{d}}\right).
\end{eqnarray*}
Thus altogether we have the result. \hfill $\Box$

%% file: B_5_grids.tex
\section{Sampling properties of $\widetilde{Q}_{a,\Omega,\lambda}(g;\rb)$}\label{sec:fine-grid}

In this section we expand on and prove the results stated in Section
\ref{sec:summary}. We mainly assume that the spatial process is observed
at $\{\sbb_{j}\}_{j=1}^{d}$ where $\sbb_{j}$ are iid uniformly
distributed random variable defined on
$[-\lambda/2,\lambda/2]^{d}$. We consider the estimator 
\begin{eqnarray}
\label{eq:Qomega}
\widetilde{Q}_{a,\Omega,\lambda}(g;\rb) = \frac{1}{\Omega^{d}}\sum_{k_{1},\ldots,k_{d}=-a}^{a}
g(\ob_{\Omega,\kb})J_{n}(\ob_{\Omega,\kb})\overline{J_{n}(\ob_{\Omega,\kb+\rb})} 
- \frac{\lambda^{d}}{\Omega^{d}n}\sum_{\kb=-a}^{a}g(\ob_{\Omega,\kb})
\frac{1}{n}\sum_{j=1}^{n}Z(\sbb_{j})^{2}\exp(-i\sbb^{\prime}_{j}\ob_{\Omega,\rb}),
\end{eqnarray}
and $\ob_{\Omega,\kb} = 2\pi \kb/\Omega$. In the previous sections we
considered the case $\Omega = \lambda$ here we consider the case that
$\Omega \neq \lambda$, noting that when $\Omega>\lambda$ the frequency
finer grid is finer and $\Omega<\lambda$ corresponds to a coarser frequency grid.  

First  we consider the mean of 
$\widetilde{Q}_{a,\Omega,\lambda}(g;\rb)$. We note that this result
holds under Assumption
\ref{assum:nonuniform} (and is not just under the uniform distribution).
\begin{lemma}\label{lemma:meangrid}
Suppose Assumptions \ref{assum:S}(i),  Assumptions  \ref{assum:G}(i)
or (ii) and \ref{assum:GG}(b,c) hold. 
\begin{itemize}
\item[(i)] If  Assumption \ref{assum:uniform} holds then
\begin{eqnarray*}
\Ex\left[\widetilde{Q}_{a,\Omega,\lambda}(g;\rb)\right] 
 &=&  \frac{c_{2}}{\pi^{d}(2\pi)^{d}}\int_{\mathbb{R}^{d}}
\sinc({\boldsymbol y})\sinc\left({\boldsymbol y}+\frac{\lambda\pi \rb}{\Omega}\right)d\ob
\int_{[-a/\Omega,a/\Omega]^{d}}g(\ob)f(\ob)d\ob
dy \nonumber\\
&& + 
O\left(\frac{\log \lambda + I(|\lambda r/\Omega|>e)
\log\left(|\frac{r\lambda}{\Omega}|\right)}{\lambda}+\frac{1}{\Omega}\right).
\end{eqnarray*}
\item[(ii)] If  Assumption \ref{assum:nonuniform} holds then
\begin{eqnarray*}
\Ex\left[\widetilde{Q}_{a,\Omega,\lambda}(g;\rb)\right] 
&=& \frac{c_{2}}{\pi^{d}(2\pi)^{d}}\sum_{\jb_{1},\jb_{2}=-\infty}^{\infty}\gamma_{\jb_{1}}\gamma_{\jb_{2}}
\int_{-[2\pi a/\Omega,2\pi a/\Omega]^{d}}g(\ob)f\left(-\ob_{}-\ob_{j_1}\right)d\ob\times \\
&&\int_{\mathbb{R}^{d}}\sinc({\boldsymbol y})
\sinc\left({\boldsymbol y}+\frac{\pi\lambda \rb}{\Omega}-(\jb_{1}+\jb_{2})\pi\right)d{\boldsymbol y} \\
&& +   O\left(\frac{1}{\Omega} +  
\frac{\log\lambda+\sum_{j=1}^{d}I(\frac{\lambda |r_{j}|}{\Omega} >e)\log\frac{\lambda  |r_{j}|}{\Omega}}{\lambda}\right).
\end{eqnarray*}
\item[(iii)] If  Assumption \ref{assum:nonuniform} holds then
\begin{eqnarray*}
\Ex\left[\widetilde{Q}_{a,\Omega,\lambda}(g;0)\right] 
 &=&  \frac{1}{(2\pi)^{d}}\sum_{\jb}|\gamma_{\jb}|^{2}\int_{\mathbb{R}^{d}}
\int_{[-2\pi a/\Omega,2\pi a/\Omega]^{d}}g(\ob)f(\ob)d\ob
dy  + O\left(\frac{\log \lambda}{\lambda}+\frac{1}{\Omega}+\frac{1}{n}\right).
\end{eqnarray*}
\end{itemize}
\end{lemma}
PROOF. We first prove (i). To simplify notation, we let $d=1$.

Taking expectations in (\ref{eq:Qomega}) (for $d=1$) we have 
\begin{eqnarray*}
\Ex\left[\widetilde{Q}_{a,\Omega,\lambda}(g;r)\right] &=& \frac{c_{2}\lambda}{2\pi\Omega}\sum_{\kb=-a}^{a}g(\omega_{\Omega,k})\int_{\mathbb{R}^{}}f(\omega)
\sinc\left(\frac{\lambda}{2}\left[\omega+\omega_{\Omega,k}\right]\right)\Sinc\left(\frac{\lambda}{2}\left[\omega+\omega_{\Omega,k+r}\right]
\right)d\omega \\
 &=&  \frac{c_{2}}{\pi^{}}\int_{\mathbb{R}^{}}
\sinc(y)\sinc\left(y+\frac{\lambda\pi r}{\Omega}\right)\left[\frac{1}{\Omega}
\sum_{k=-a}^{a}g(\omega_{\Omega,k})f(\frac{2y}{\lambda}-\omega_{\Omega,k})\right]
dy.
\end{eqnarray*}
Replacing summand with integral and 
$f(\frac{2y}{\lambda}-\omega_{\Omega,k})$ with $f(-\omega_{\Omega,k})$
(and using Lemma \ref{lemma:1star}, equation (\ref{eq:lemma1star1}))
we have 
\begin{eqnarray*}
\Ex\left[\widetilde{Q}_{a,\Omega,\lambda}(g;r)\right] &=&
\frac{2c_{2}}{\pi^{}}\int_{\mathbb{R}^{}}
\sinc(y)\sinc\left(y+\frac{\lambda\pi r}{\Omega}\right)
\int_{-2\pi a/\Omega}^{2\pi a/\Omega}g(\omega)f(\frac{2y}{\lambda}-\omega)d\omega
dy + O\left(\frac{1}{\Omega}\right).
 \\
 &=&  \frac{c_{2}}{\pi^{d}}\int_{\mathbb{R}^{d}}
\sinc(y)\sinc\left(y+\frac{\lambda\pi r}{\Omega}\right)
\int_{-2\pi a/\Omega}^{2\pi a/\Omega}g(\omega)f(\omega)d\omega
dy \\
&& + O\left(\frac{\log \lambda + I(|\lambda r/\Omega|>e)\log(|\lambda r/\Omega|)}{\lambda}+\frac{1}{\Omega}\right).
\end{eqnarray*}
Thus proving (i).

We now prove (ii). The proof is very similar to the proof of Theorem
\ref{theorem:nonuniformmean}. Taking expectation
\begin{eqnarray*}
&&\Ex\left[\widetilde{Q}_{a,\Omega,\lambda}(g;r)\right] \\
&=& \frac{c_{2}\lambda }{2\pi \Omega}\sum_{j_{1},j_{2}=-\infty}^{\infty}\gamma_{j_{1}}\gamma_{j_{2}}
\sum_{k=-a}^{a}g(\omega_{\Omega,k})
\int_{-\infty}^{\infty}f(\omega)\sinc\left(\frac{\lambda}{2}\left(\omega+\omega_{\Omega,k}+\omega_{\lambda,j_{1}}\right)\right)
\sinc\left(\frac{\lambda}{2}\left(\omega+\omega_{\Omega,k}-\omega_{\lambda,j_{1}}\right)\right)
d\omega. 
\end{eqnarray*}
By the change of variables $y=\frac{\lambda}{2}(\omega+\omega_{\Omega,k}+\omega_{\lambda,j_{1}})$
we obtain 
\begin{eqnarray*}
&&\Ex\left[\widetilde{Q}_{a,\Omega,\lambda}(g;r)\right] \\
&=& \frac{c_{2}}{\pi \Omega}\sum_{j_{1},j_{2}=-\infty}^{\infty}\gamma_{j_{1}}\gamma_{j_{2}}
\sum_{k=-a}^{a}g(\omega_{\Omega,k})
\int_{-\infty}^{\infty}f\left(\frac{2y}{\lambda}-\omega_{\lambda,j_1}-\omega_{\Omega,k}\right)\sinc(y)
\sinc\left(y+\frac{\lambda \pi r}{\Omega} - (j_{1}+j_{2})\pi\right)dy.
\end{eqnarray*}
Replacing sum with an integral and 
using Lemma \ref{lemma:sum-integral}(ii) gives 
\begin{eqnarray*}
&&\Ex\left[\widetilde{Q}_{a,\Omega,\lambda}(g;r)\right] \\
&=& \frac{c_{2}}{\pi}\sum_{j_{1},j_{2}=-\infty}^{\infty}\gamma_{j_{1}}\gamma_{j_{2}}
\int_{-2\pi a/\Omega}^{2\pi a/\Omega}g(\omega)
\int_{-\infty}^{\infty}f\left(\frac{2y}{\lambda}-\omega_{\lambda,j_1}-\omega\right)\sinc(y)
\sinc\left(y+\frac{\lambda \pi r}{\Omega} - (j_{1}+j_{2})\pi\right)dy+
O\left(\frac{1}{\Omega}\right).
\end{eqnarray*}
Next, replacing $f(\frac{2y}{\lambda}-\omega_{\lambda,j_1}-\omega)$
with $f(-\omega-\omega_{j_1})$ gives 
\begin{eqnarray*}
&&\Ex\left[\widetilde{Q}_{a,\Omega,\lambda}(g;r)\right] \\
&=& \frac{c_{2}}{2\pi^{2}}\sum_{j_{1},j_{2}=-\infty}^{\infty}\gamma_{j_{1}}\gamma_{j_{2}}
\int_{-2\pi a/\Omega}^{2\pi a/\Omega}g(\omega) f\left(\frac{2y}{\lambda}-\omega_{}-\omega_{j_1}\right)d\omega
\int_{-\infty}^{\infty}\sinc(y)\sinc\left(y+\frac{\pi\lambda r}{\Omega}
-(j_{1}+j_{2})\pi\right)dy \\
&& +  R_{n} + O\left(\frac{1}{\Omega}\right).
\end{eqnarray*}
where 
\begin{eqnarray*}
R_{n} &=& \frac{c_{2}}{2\pi^{2}}\sum_{j_{1},j_{2}=-\infty}^{\infty}\gamma_{j_{1}}\gamma_{j_{2}}
\int_{-2\pi a/\Omega}^{2\pi a/\Omega}g(\omega)
\left[f\left(\frac{2y}{\lambda}-\omega_{}-\omega_{\lambda,j_{1}}\right)
  - 
f(-\omega-\omega_{\lambda,j_{1}})\right] \\
&&\times \int_{-\infty}^{\infty}\sinc(y)\sinc\left(y+\frac{\pi\lambda r}{\Omega}
-(j_{1}+j_{2})\pi\right)dyd\omega.
\end{eqnarray*}
By using Lemma \ref{lemma:1star} 
\begin{eqnarray*}
|R_{n}| &\leq& C\sum_{j_{1},j_{2}=-\infty}^{\infty}|\gamma_{j_{1}}|\cdot|\gamma_{j_{2}}|
\frac{\log\lambda + \log |\frac{\lambda}{\Omega}r -
  (j_{1}+j_{2})|}{\lambda} \\
&\leq& C\sum_{j_{1},j_{2}=-\infty}^{\infty}|\gamma_{j_{1}}|\cdot|\gamma_{j_{2}}|
\frac{\log\lambda + I(|\lambda r/\Omega| >e)\log|\lambda
  r/\Omega|+\log|j_{1}|+\log|j_{2}|}{\lambda} \\
&=& 
O\left( \frac{\log \lambda +   I(|\lambda r/\Omega| >e)\log|\lambda
  r/\Omega|}{\lambda}\right).
\end{eqnarray*}
This gives
\begin{eqnarray*}
&&\Ex\left[\widetilde{Q}_{a,\Omega,\lambda}(g;r)\right] \\
&=& \frac{c_{2}}{2\pi^{2}}\sum_{j_{1},j_{2}=-\infty}^{\infty}\gamma_{j_{1}}\gamma_{j_{2}}
\int_{-2\pi a/\Omega}^{2\pi a/\Omega}g(\omega) f\left(\frac{2y}{\lambda}-\omega_{}-\omega_{j_1}\right)d\omega
\int_{-\infty}^{\infty}\sinc(y)\sinc\left(y+\frac{\pi\lambda r}{\Omega}
-(j_{1}+j_{2})\pi\right)dy \\
&& +   O\left(\frac{1}{\Omega} +  \frac{\log \lambda +   I(|\lambda r/\Omega| >e)\log|\lambda
  r/\Omega|}{\lambda}\right).
\end{eqnarray*}
This proves (ii). To prove (iii) we note that when $r=0$, the above
reduces to 
\begin{eqnarray*}
&&\Ex\left[\widetilde{Q}_{a,\Omega,\lambda}(g;r)\right] \\
&=& \frac{c_{2}}{2\pi^{2}}\sum_{j_{1},j_{2}=-\infty}^{\infty}\gamma_{j_{1}}\gamma_{j_{2}}
\int_{-2\pi a/\Omega}^{2\pi a/\Omega}g(\omega) f\left(\frac{2y}{\lambda}-\omega_{}-\omega_{j_1}\right)d\omega
\int_{-\infty}^{\infty}\sinc(y)\sinc\left(y
-(j_{1}+j_{2})\pi\right)dy \\
&& +   O\left(\frac{1}{\Omega} +  \frac{\log \lambda +   I(|\lambda r/\Omega| >e)\log|\lambda
  r/\Omega|}{\lambda}\right).
\end{eqnarray*}
By using the orthogonality of the sinc function at integer shifts (and $f(-\omega)=f(\omega)$) we
have
\begin{eqnarray*}
\Ex\left[\widetilde{Q}_{a,\Omega,\lambda}(g;r)\right] 
&=& \frac{1}{2\pi}\sum_{j=-\infty}^{\infty}|\gamma_{j}|^{2}
\int_{-2\pi a/\Omega}^{2\pi a/\Omega}g(\omega) f(\omega+\omega_{j})d\omega + 
 O\left(\frac{\log\lambda}{\lambda}+ \frac{1}{\Omega}+\frac{1}{n}\right) \\
 &=& \frac{1}{2\pi}\sum_{j=-\infty}^{\infty}|\gamma_{j}|^{2}
\int_{-2\pi a/\Omega}^{2\pi a/\Omega}g(\omega) f(\omega)d\omega + 
 O\left(\frac{\log\lambda}{\lambda}+\frac{1}{\Omega}+\frac{1}{n}\right) 
\end{eqnarray*}
thus we obtain (iii).  \hfill $\Box$

\vspace{3mm}

We now obtain the variance of
$\widetilde{Q}_{a,\Omega,\lambda}(g;\rb_{})$. The rate of convergence
is based on whether the ratio $\lambda/\Omega < 1$ or $\lambda/\Omega
\geq 1$. The reason behind this difference is due to the following result. 

\begin{lemma}\label{lemma:triangle}
Let $T(\cdot)$ denote the triangle kernel where $T(x) = (1-|x|)$ for
$|x|\leq 1$ else $T(x)=0$ and suppose $\alpha >0$. Then we have 
\begin{eqnarray}
\label{eq:TaSum}
\sum_{k=-a}^{a} \sinc^{2}\left(\alpha k \pi\right) =
\frac{1}{\alpha}\int_{-\alpha}^{\alpha}
T\left(\frac{x}{\alpha}\right)\frac{\sin(2\pi (a+\frac{1}{2})x)}{\sin(\pi x)}dx,
\end{eqnarray}
\begin{eqnarray}
\label{eq:TaSum1}
\sum_{k\in \mathbb{Z}} \sinc^{2}\left(\alpha k \pi\right) =
\frac{1}{\alpha}\sum_{|k| \leq \lfloor \alpha \rfloor} T\left(\frac{k}{\alpha}\right)
\end{eqnarray}
and 
\begin{eqnarray}
\label{eq:TaSum2}
\sum_{k\in \mathbb{Z}} \sinc^{2}\left(\alpha k \pi\right) \rightarrow 1
\end{eqnarray}
as $\alpha \rightarrow \infty$.
\end{lemma}
PROOF. We recall that $\sinc(\alpha x/2) = \frac{1}{\alpha}
\int_{-\alpha/2}^{\alpha/2}e^{ix\omega}dx$ (Fourier transform of the
rectangle kernel $\alpha^{-1}I_{-\alpha/2,\alpha/2}(x)$). Thus $\sinc^{2}(\alpha x/2)$ 
is the convolution of two rectangle kernels (which is the triangle
kernel defined on $[-\alpha,\alpha]$) and   
\begin{eqnarray*}
\sinc^{2}\left(\frac{\alpha}{2} \omega\right) =
\frac{1}{\alpha}\int_{-\alpha}^{\alpha}
T\left(\frac{x}{\alpha}\right)\exp(i\omega x)dx,
\end{eqnarray*}
thus with $\omega = k\pi $
\begin{eqnarray*}
\sinc^{2}\left(\alpha k \pi \right) =
\frac{1}{\alpha}\int_{-\alpha}^{\alpha}
T\left(\frac{x}{\alpha}\right)\exp(i2\pi k x)dx.
\end{eqnarray*}
Substituting this into the sum gives 
\begin{eqnarray*}
\sum_{k=-a}^{a} \sinc^{2}\left(\alpha k \pi\right) =
\frac{1}{\alpha}\int_{-\alpha}^{\alpha}
T\left(\frac{x}{\alpha}\right)\left[\sum_{k=-a}^{a}\exp(i2\pi k x)\right]dx =
\frac{1}{\alpha}\int_{-\alpha}^{\alpha}
T\left(\frac{x}{\alpha}\right)\frac{\sin(2\pi (a+\frac{1}{2})x)}{\sin(\pi x)}dx,
\end{eqnarray*}
thus proving (\ref{eq:TaSum}). Next letting the limit in the sum
$a\rightarrow \infty$ we have that the Dirichlet kernel
 limits to the generalized function
\begin{eqnarray*}
\frac{\sin(2\pi (a+\frac{1}{2})x)}{\sin(\pi x)}\rightarrow \sum_{m\in
  \mathbb{Z}}\delta_{m}(x) \qquad \textrm{ as }a\rightarrow \infty
\end{eqnarray*}
where $\delta_{m}(x)$ is the dirac delta function which is zero
everywhere but $m$. Substituting this into the the above integral gives
(\ref{eq:TaSum1}).  

To prove (\ref{eq:TaSum2}) we note that using (\ref{eq:TaSum1}) we
have  
\begin{eqnarray*}
\sum_{k\in \mathbb{Z}} \sinc^{2}\left(\alpha k \pi\right) &=&
\frac{1}{\alpha}\sum_{|k| \leq \lfloor \alpha \rfloor}
T\left(\frac{k}{\alpha}\right)=
\frac{1}{\alpha}\left(1 +2\sum_{|k| \leq \lfloor \alpha \rfloor} 
\left[1-\frac{|k|}{\alpha}\right] 
\right) \\
&=& \frac{1}{\alpha}\left(1+\frac{2}{\alpha}\sum_{|k| \leq \lfloor
    \alpha \rfloor} \left(\alpha - |k|\right)\right) \\
&=& \frac{1}{\alpha}\left(1+\frac{2}{\alpha}
\lfloor \alpha\rfloor\left(\alpha - 
\lfloor \alpha\rfloor \right) + \frac{1}{\alpha}\lfloor \alpha \rfloor
(\lfloor \alpha \rfloor +1)
\right) = 1 + O\left(\frac{1}{\alpha}\right),
\end{eqnarray*}
where $\lfloor x \rfloor$ denotes the smallest integer less than or
equal to $x$. Thus we see that $\sum_{k\in \mathbb{Z}}
\sinc^{2}\left(\alpha k \pi\right)$ is uniformly bounded for all
$\alpha \geq 1$ and  as
$\alpha \rightarrow \infty$ $\sum_{k\in \mathbb{Z}}
\sinc^{2}\left(\alpha k \pi\right) \rightarrow 1$, thus proving
(\ref{eq:TaSum2}).  
\hfill $\Box$

\vspace{3mm}
We now summarize the pertinent points of the above lemma.
If $0<\alpha\leq 1$ we have 
\begin{eqnarray*}
\alpha\sum_{k\in \mathbb{Z}} \sinc^{2}\left(\alpha k \pi\right) = T(0)
=1 \Rightarrow \sum_{k\in \mathbb{Z}} \sinc^{2}\left(\alpha k \pi\right) =\frac{1}{\alpha}.
\end{eqnarray*}
On the other hand, if $\alpha > 1$ then
$\sum_{k\in \mathbb{Z}} \sinc^{2}\left(\alpha k \pi\right)$ is
uniformly bounded for all $\alpha$. We show in the lemmas below this
simple result determines the optimal choice of frequency grid.

In the following lemma we obtain the first approximation under the
assumption of Gaussianity of the spatial process. 

\begin{lemma}\label{lemma:covarianceFine}
Suppose Assumptions \ref{assum:S}, \ref{assum:uniform} and  
Assumptions \ref{assum:G}(ii) and \ref{assum:GG}(b) hold. 
Then
\begin{itemize}
\item[(i)] If $\lambda<\Omega$ (finer frequency grid)
\begin{eqnarray*}
\lambda^{d}\cov\left[\widetilde{Q}_{a,\Omega,\lambda}(g;\rb_{1}),
\widetilde{Q}_{a,\Omega,\lambda}(g;\rb_{2})\right] =
\frac{\lambda^{d}}{\Omega^{2d}}\left[A_{1}(\rb_{1},\rb_{2})+A_{2}(\rb_{1},\rb_{2})\right] +F_{1,\mathrm{fine}},
\end{eqnarray*}
\begin{eqnarray*}
\lambda^{d}\cov\left[\widetilde{Q}_{a,\Omega,\lambda}(g;\rb_{1}),
\overline{\widetilde{Q}_{a,\Omega,\lambda}(g;\rb_{2})}\right]
= 
\frac{\lambda^{d}}{\Omega^{2d}}\left[A_{3}(\rb_{1},\rb_{2})+A_{4}(\rb_{1},\rb_{2})\right] + F_{2,\mathrm{fine}},
\end{eqnarray*}
with 
\begin{eqnarray*}
F_{1,\mathrm{fine}} \textrm{ and } F_{2,\mathrm{fine}} = 
O\left(\frac{\lambda^{d}}{n}+
\frac{\lambda^{d} \log^{3}(a)}{\Omega^{d} n }\right) 
\end{eqnarray*}
\item[(ii)] If $\lambda>\Omega$ (coarser frequency grid)
\begin{eqnarray*}
\Omega^{d}\cov\left[\widetilde{Q}_{a,\Omega,\lambda}(g;\rb_{1}),
\widetilde{Q}_{a,\Omega,\lambda}(g;\rb_{2})\right] =
\frac{1}{\Omega^{d}}\left[A_{1}(\rb_{1},\rb_{2})+A_{2}(\rb_{1},\rb_{2})\right] +F_{1,\mathrm{coarse}},
\end{eqnarray*}
\begin{eqnarray*}
\Omega^{d}\cov\left[\widetilde{Q}_{a,\Omega,\lambda}(g;\rb_{1}),
\overline{\widetilde{Q}_{a,\Omega,\lambda}(g;\rb_{2})}\right]
= 
\frac{1}{\Omega^{d}}\left[A_{3}(\rb_{1},\rb_{2})+A_{4}(\rb_{1},\rb_{2})\right] + F_{2,\mathrm{coarse}},
\end{eqnarray*}
with 
\begin{eqnarray*}
F_{1,\mathrm{coarse}} \textrm{ and } F_{2,\mathrm{coarse}} = 
O\left(\frac{\lambda^{d}}{n} +\frac{
  \log^{}(a)\log^{2}(\frac{\lambda}{\Omega}a)}{n}I_{\frac{\lambda}{\Omega}\notin
\mathbb{Z}}\right).
\end{eqnarray*}
\end{itemize}
where 
\begin{eqnarray*}
&&A_{1}(\rb_{1},\rb_{2}) = \frac{1}{\pi^{2d}}
\sum_{\mb=-2a}^{2a}\sum_{\kb=\max(-a,-a+\mb)}^{\min(a,a+\mb)} 
\int_{\mathbb{R}^{2d}} f(\frac{2\ubb}{\lambda} -\ob_{\Omega,\kb_{}})
f(\frac{2\vbb}{\lambda} + \ob_{\Omega,\kb_{}}+\ob_{\Omega,\rb_{1}}) \\
&& g(\ob_{\Omega,\kb})\overline{g(\ob_{\Omega,\kb}-\ob_{\Omega,\mb})}
\Sinc\left(\ubb-\frac{\lambda}{\Omega}\mb \pi\right)
\Sinc\left(\vbb+\frac{\lambda}{\Omega}(\mb+\rb_{1}-\rb_{2})\pi\right)
\Sinc(\ubb)\Sinc(\vbb)d\ubb d\vbb
\end{eqnarray*}
\begin{eqnarray*}
&&A_{2}(\rb_{1},\rb_{2}) = \frac{1}{\pi^{2d}}
\sum_{\mb=-2a}^{2a}\sum_{\kb=\max(-a,-a+\mb)}^{\min(a,a+\mb)} 
\int_{\mathbb{R}^{2d}} f(\frac{2\ubb}{\lambda} -\ob_{\Omega,\kb_{}})
f(\frac{2\vbb}{\lambda} + \ob_{\Omega,\kb_{}}+\ob_{\Omega,\rb_{1}}) \\
&& g(\ob_{\Omega,\kb})\overline{g(\ob_{\Omega,\mb-\kb})}
\Sinc\left(\ubb-\frac{\lambda}{\Omega}(\mb+\rb_{2}) \pi\right)\Sinc\left(\vbb+\frac{\lambda}{\Omega}(\mb+\rb_{1})\pi\right)
\Sinc(\ubb)\Sinc(\vbb)d\ubb d \vbb
\end{eqnarray*}
\begin{eqnarray*}
&&A_{3}(\rb_{1},\rb_{2}) = \frac{1}{\pi^{2d}}
\sum_{\mb=-2a}^{2a}\sum_{\kb=\max(-a,-a+\mb)}^{\min(a,a+\mb)} 
\int_{\mathbb{R}^{2d}} f(\frac{2\ubb}{\lambda} -\ob_{\Omega,\kb_{}})
f(\frac{2\vbb}{\lambda} + \ob_{\Omega,\kb}+\ob_{\Omega,\rb_{1}}) \\
&& g(\ob_{\Omega,\kb})g(\ob_{\Omega,\mb-\kb})
\Sinc\left(\ubb+\frac{\lambda}{\Omega}\mb \pi\right)\Sinc\left(\vbb+\frac{\lambda}{\Omega}(\mb+\rb_{2}+\rb_{1})\pi\right)
\Sinc(\ubb)\Sinc(\vbb)d\ubb d\vbb
\end{eqnarray*}
\begin{eqnarray*}
&&A_{4}(\rb_{1},\rb_{2}) = \frac{1}{\pi^{2d}}
\sum_{\mb=-2a}^{2a}\sum_{\kb=\max(-a,-a+\mb)}^{\min(a,a+\mb)} 
\int_{\mathbb{R}^{2d}} f(\frac{2\ubb}{\lambda} -\ob_{\Omega,\kb_{}})
f(\frac{2\vbb}{\lambda} + \ob_{\Omega,\kb_{}}+\ob_{\Omega,\rb_{1}}) \\
&& g(\ob_{\Omega,\kb})g(\ob_{\Omega,\kb-\mb})
\Sinc\left(\ubb-\frac{\lambda}{\Omega}(\mb-\rb_{2}) \pi\right)\Sinc\left(\vbb+\frac{\lambda}{\Omega}(\mb+\rb_{1})\pi\right)
\Sinc(\ubb)\Sinc(\vbb)d\ubb d\vbb
\end{eqnarray*}
where $\kb=\max(-a,-a-\mb)=\{k_{1}=\max(-a,-a-m_{1}),\ldots,k_{d}=\max(-a,-a-m_{d})\}$, 
$\kb=\min(-a+\mb)=\{k_{1}=\max(-a,-a+m_{1}),\ldots,k_{d}=\min(-a+m_{d})\}$.
% $C_{4} = n(n-1)(n-2)(n-3)/n^{4}$
\end{lemma}
PROOF. To simplify notation, we prove the result for $d=1$.
% We start with the case $\lambda/\Omega <1$ (using a ``fine'' frequency grid).

Using a similar expansion to that in (\ref{eq:varQ}), for  
the case $\lambda/\Omega <1$ (``fine'' frequency grid) we have 
\begin{eqnarray}
&&\lambda\cov\left[\widetilde{Q}_{a,\Omega,\lambda}(g;r_{1}),\widetilde{Q}_{a,\Omega,\lambda}(g;r_{2})
\right] \nonumber\\
&=&
\frac{\lambda^{3}}{\Omega^{2}n^{4}}\sum_{j_{1},j_{2},j_{3},j_{4}\in
  \mathcal{D}_{4}}\sum_{k_{1},k_{2}=-a}^{a}
g(\omega_{\Omega,k_{1}})\overline{g(\omega_{\Omega,k_{2}})}
\bigg(\nonumber\\
&&\cum\big[Z(s_{j_1})e^{is_{j_1}\omega_{\Omega,k_{1}}},
Z(s_{j_3})e^{-is_{j_3}\omega_{\Omega,k_{2}}}\big] 
 \cum\big[ Z(s_{j_{2}}) e^{-is_{j_2}\omega_{\Omega,k_{1}+r_{1}}},
Z(s_{j_4})e^{is_{j_4}\omega_{\Omega,k_{2}+r_{2}}}\big] \nonumber\\
&& +\cum\big[Z(s_{j_1})e^{is_{j_1}\omega_{\Omega,k_{1}}},Z(s_{j_4})e^{is_{j_4}\omega_{\Omega,k_{2}+r_2}}\big]
\cum\big[Z(s_{j_2})e^{-is_{j_{2}}\omega_{\Omega,k_{1}+r_2}},Z(s_{j_3})e^{-is_{j_3}\omega_{\Omega,k_{2}}}\big]\bigg)
+ F_{1,\textrm{fine}} \qquad \label{eq:sumMain}
%\\
%&=& A_{1}(r_{1},r_{2}) + A_{2}(r_{1},r_{2}) + F_{1}
\end{eqnarray}
whereas for the case $\lambda/\Omega >1$ (``coarse'' frequency grid)
\begin{eqnarray}
&&\Omega\cov\left[\widetilde{Q}_{a,\Omega,\lambda}(g;r_{1}),\widetilde{Q}_{a,\Omega,\lambda}(g;r_{2})
\right] \nonumber\\
&=&
\frac{\lambda^{2}}{\Omega^{}n^{4}}\sum_{j_{1},j_{2},j_{3},j_{4}\in
  \mathcal{D}_{4}}\sum_{k_{1},k_{2}=-a}^{a}
g(\omega_{\Omega,k_{1}})\overline{g(\omega_{\Omega,k_{2}})}
\bigg( \nonumber\\
&&\cum\big[Z(s_{j_1})e^{is_{j_1}\omega_{\Omega,k_{1}}},
Z(s_{j_3})e^{-is_{j_3}\omega_{\Omega,k_{2}}}\big] 
  \cum\big[ Z(s_{j_{2}}) e^{-is_{j_2}\omega_{\Omega,k_{1}+r_{1}}},
Z(s_{j_4})e^{is_{j_4}\omega_{\Omega,k_{2}+r_{2}}}\big] \nonumber\\
&& +\cum\big[Z(s_{j_1})e^{is_{j_1}\omega_{\Omega,k_{1}}},Z(s_{j_4})e^{is_{j_4}\omega_{\Omega,k_{2}+r_2}}\big]
\cum\big[Z(s_{j_2})e^{-is_{j_{2}}\omega_{\Omega,k_{1}+r_2}},Z(s_{j_3})e^{-is_{j_3}\omega_{\Omega,k_{2}}}\big]\bigg)
+ F_{1,\textrm{coarse}} \qquad \label{eq:sumMain2}
%\\
%&=& A_{1}(r_{1},r_{2}) + A_{2}(r_{1},r_{2}) + F_{1}
\end{eqnarray}
where 
\begin{eqnarray*}
F_{1,\textrm{fine}}&=& \lambda F_{1}\quad  \textrm{ and }\quad 
F_{1,\textrm{coarse}} = \Omega F_{1}
\end{eqnarray*}
with 
\begin{eqnarray*}
F_{1} &=& \frac{\lambda^{2}}{\Omega^{2}n^{4}}\sum_{j_{1},j_{2},j_{3},j_{4}\in
  \mathcal{D}_{3}}\sum_{k_{1},k_{2}=-a}^{a}
g(\omega_{\Omega,k_{1}})\overline{g(\omega_{\Omega,k_{2}})} \\
&& \times \bigg(
\cum\big[Z(s_{j_1})e^{is_{j_1}\omega_{\Omega,k_{1}}},Z(s_{j_3})e^{-is_{j_3}\omega_{\Omega,k_{2}}}\big]
\cum\big[Z(s_{j_2})e^{-is_{j_2}\omega_{\Omega,k_{1}+r_1}},Z(s_{j_4})e^{is_{j_4}\omega_{\Omega,k_{2}+r_2}}
\big] \nonumber\\
&& +\cum\big[Z(s_{j_1})e^{is_{j_1}\omega_{\Omega,k_{1}}},Z(s_{j_4})e^{is_{j_4}\omega_{\Omega,k_{2}+r_2}}\big]
\cum\big[Z(s_{j_2})e^{-is_{j_{2}}\omega_{\Omega,k_{1}+r_1}},
Z(s_{j_3})e^{-is_{j_3}\omega_{\Omega,k_{2}}}\big]\nonumber\\
&& +\cum\big[Z(s_{j_1})e^{is_{j_1}\omega_{\Omega,k_{1}}},
Z(s_{j_2})e^{-is_{j_2}\omega_{\Omega,k_{1}+r_1}},Z(s_{j_3})e^{-is_{j_3}\omega_{\Omega,k_{2}}},
Z(s_{j_4})e^{is_{j_4}\omega_{\Omega,k_{2}+r_2}}\big]\bigg).
\end{eqnarray*} 
By using the same techniques used in the proof of Lemma
\ref{lemma:covariance} we can show that the first sum in
(\ref{eq:sumMain}) and (\ref{eq:sumMain2}) is equal to the $A_{1}(r_{1},r_{2})$ and
$A_{2}(r_{1},r_{2})$ given at the start of the lemma.
We now  bound the remainders $F_{1,\textrm{fine}}=\lambda F_{1}$ and
$F_{1,\textrm{coarse}} = \Omega F_{1}$.

To bound $F_{1}$ 
we note that it is mainly comprised of terms which 
take the form
\begin{eqnarray*}
D_{1}&=&\sum_{k_{1},k_{2}=-a}^{a}
g(\omega_{\Omega,k_{1}})\overline{g(\omega_{\Omega,k_{2}})} 
\cum[Z(s_{1})e^{is_{1}\omega_{\Omega,k_{1}}},Z(s_{2})e^{-is_{2}\omega_{\Omega,k_{1}+r_{1}}},
Z(s_{3})e^{-is_{3}\omega_{\Omega,k_{2}}},Z(s_{1})e^{is_{1}\omega_{\Omega,k_{2}+r_{2}}}] 
\end{eqnarray*}
and
\begin{eqnarray*}
D_{2}&=&\sum_{k_{1},k_{2}=-a}^{a}g(\omega_{\Omega,k_{1}})\overline{g(\omega_{\Omega,k_{2}})} 
\cum\left[Z(s_{1})e^{is_{1}\omega_{\Omega,k_{1}}},Z(s_{1})e^{-is_{1}\omega_{\Omega,k_{2}}}\right]
\cum\left[Z(s_{2})e^{-is_{2}\omega_{\Omega,k_{1}+r_{1}}},
Z(s_{3})e^{is_{3}\omega_{\Omega,k_{2}+r_{2}}}\right],
\end{eqnarray*}
see Lemma \ref{lemma:cumulantsA}.
We start by bounding $D_{1}$.  Expanding the cumulant (using
(\ref{eq:conditional-cumulants})) gives 
\begin{eqnarray*}
D_{1}
&=&\sum_{k_{1},k_{2}=-a}^{a}g(\omega_{\Omega,k_{1}})\overline{g(\omega_{\Omega,k_{2}})} 
\cum[c(s_{1}-s_{2})e^{is_{1}\omega_{\Omega,k_{1}}-is_{2}\omega_{\Omega,k_{1}+r_{1}}},
c(s_{1}-s_{3})e^{-is_{3}\omega_{\Omega,k_{2}}+is_{1}\omega_{\Omega,k_{2}+r_{2}}}] + \\
&& \sum_{k_{1},k_{2}=-a}^{a}g(\omega_{\Omega,k_{1}})\overline{g(\omega_{\Omega,k_{2}})} 
\cum[c(s_{1}-s_{3})e^{is_{1}\omega_{\Omega,k_{1}}-is_{3}\omega_{\Omega,k_{2}}},
c(s_{1}-s_{2})e^{-is_{2}\omega_{\Omega,k_{1}+r_{1}}+is_{1}\omega_{\Omega,k_{2}+r_{2}}}]  +\\
&& \sum_{k_{1},k_{2}=-a}^{a}g(\omega_{\Omega,k_{1}})\overline{g(\omega_{\Omega,k_{2}})} 
\cum[c(0)e^{is_{1}(\omega_{\Omega,k_{1}}+\omega_{\Omega,k_{2}})},
c(s_{2}-s_{3})e^{-is_{2}\omega_{\Omega,k_{1}}-is_{3}\omega_{\Omega,k_{2}}}]
= I_{1} + I_{2} + I_{3}.
\end{eqnarray*}
Expanding $I_{1}$ just as was done in the proof of Lemma
\ref{lemma:cumulantsA} gives $I_{1} = E_{1} - E_{2}$ where 
\begin{eqnarray*}
E_{1}&=&\frac{1}{(2\pi)^{2}}\sum_{k_{1},k_{2}=-a}^{a}
g(\omega_{\Omega,k_{1}})\overline{g(\omega_{\Omega,k_{2}})} 
\int\int
f(x)f(y)\sinc\bigg(\frac{\lambda}{2}\left[
(x+y)+\omega_{\Omega,k_{1}+k_{2}+r_{2}}\right]\bigg) \\
&&\times \sinc\bigg(\frac{\lambda}{2}\left[x+\omega_{\Omega,k_{1}+r_{2}}\right]\bigg)\sinc\bigg(\frac{\lambda}{2}\left[
y+\omega_{\Omega,k_{2}}\right]\bigg)dxdy 
%-\\
%&& \frac{1}{(2\pi)^{2}}\sum_{k_{1},k_{2}=-a}^{a}\int\int f(x)f(y)\sinc\bigg(\frac{\lambda}{2}x+k_{2}\pi\bigg)
%\sinc\bigg(\frac{\lambda}{2}x+k_{2}\pi\bigg)
%\sinc\bigg(\frac{\lambda}{2}y+k_{1}\pi\bigg)
%\sinc\bigg(\frac{\lambda}{2}y+k_{1}\pi\bigg)dxdy \\
%&=& E_{1} - E_{2}.
\end{eqnarray*}
and $E_{2}$ is defined similarly. Changing variables 
$u = \frac{\lambda}{2}\left[x+\omega_{\Omega,k_{1}+r_{2}}\right]$ and 
$v=\frac{\lambda}{2}\left[
y+\omega_{\Omega,k_{2}}\right]$ gives 
\begin{eqnarray*}
E_{1} &=& \frac{1}{\pi^{2}\lambda^{2}}
\sum_{k_{1},k_{2}=-a}^{a}g(\omega_{\Omega,k_{1}})\overline{g(\omega_{\Omega,k_{2}})} 
\int\int f\left(\frac{2u}{\lambda} -
  \omega_{\Omega,k_{1}+r_{1}} \right) \\
&&f\left(\frac{2v}{\lambda} -
  \omega_{\Omega,k_{2}} \right)\sinc\left(u+v+\frac{\lambda}{\Omega}(r_{2}-r_{2})\pi \right)\sinc(u)\sinc(v)dudv.
\end{eqnarray*}
Since $\frac{1}{\Omega^{2}}\sum_{k_{1},k_{2}} |g(\omega_{\Omega,k_{1}})\overline{g(\omega_{\Omega,k_{2}})} |
f(\frac{2u}{\lambda} -  \omega_{\Omega,k_{1}+r_{1}}) 
f(\frac{2v}{\lambda} - \omega_{\Omega,k_{2}})<\infty$ we have 
 $E_{1} =
O(\frac{\Omega^{2}}{\lambda^{2}})$. Similarly we can show that $E_{2}
= O(\frac{\Omega^{2}}{\lambda^{2}})$, thus $I_{1} = O(\frac{\Omega^{2}}{\lambda^{2}})$.  
The same proof can
be used to show that $I_{2} = O(\frac{\Omega^{2}}{\lambda^{2}})$. Further, $I_{3}=0$, since $s_{1},s_{2}$
and $s_{3}$ are independent. 
Thus 
\begin{eqnarray}
|D_{1}|= O\left(\frac{\Omega^{2}}{\lambda^{2}}\right) \label{eq:cumulantOmega1}
\end{eqnarray}
Next we consider $D_{2}$
\begin{eqnarray*}
&&D_{2}\\
&=&\sum_{k_{1},k_{2}=-a}^{a}g(\omega_{\Omega,k_{1}})\overline{g(\omega_{\Omega,k_{2}})} 
\cum\left[Z(s_{1})e^{is_{1}\omega_{\Omega,k_{1}}},Z(s_{1})e^{-is_{1}\omega_{\Omega,k_{2}}}\right]
\cum\left[Z(s_{2})e^{-is_{2}\omega_{\Omega,k_{1}+r_{1}}},Z(s_{3})e^{is_{3}\omega_{\Omega,k_{2}+r_{2}}}\right] \\
&=& c(0) \sum_{k_{1},k_{2}=-a}^{a}g(\omega_{\Omega,k_{1}})\overline{g(\omega_{\Omega,k_{2}})} \Ex\left[e^{is_{1}\omega_{\Omega,k_{1}-k_{2}}}\right]
\cum\left[Z(s_{2})e^{-is_{2}\omega_{\Omega,k_{1}+r_{1}}},Z(s_{3})e^{is_{3}\omega_{\Omega,k_{2}+r_{2}}}\right]
\\
&=&
\sum_{k_{1},k_{2}=-a}^{a}g(\omega_{\Omega,k_{1}})\overline{g(\omega_{\Omega,k_{2}})} \sinc\left(\frac{\lambda}{\Omega}(k_{1}-k_{2})\right)\int
  f(x)\sinc\left(\frac{\lambda}{2}(x+\omega_{\Omega,k_{1}+r_{1}})
\right)\sinc\left(\frac{\lambda}{2}(x+\omega_{\Omega,k_{2}+r_{2}}) \right)dx.
\end{eqnarray*}
If the ratio $\lambda/\Omega\in \mathbb{Z}$, then
$\sinc\left(\frac{\lambda}{\Omega}(k_{1}-k_{2})\pi\right)=0$ if
$k_{1}\neq k_{2}$, this reduces the above double summand to a single
summand and 
\begin{eqnarray*}
D_{2}
&=& \frac{c(0)\Omega}{\lambda}
\sum_{k=-a}^{a}|g(\omega_{\Omega,k_{}})|^{2}\int
  f(x)\sinc\left(\frac{\lambda}{2}(x+\omega_{\Omega,k_{}+r_{1}})
\right)\sinc\left(\frac{\lambda}{2}(x+\omega_{\Omega,k_{}+r_{2}}
\right)dx \\
&=& \frac{c(0)\Omega}{\lambda^{2}}
\sum_{k=-a}^{a}|g(\omega_{\Omega,k_{}})|^{2}\int
  f\left(\frac{2u}{\lambda} - \omega_{\Omega,k+r_{1}} \right)\sinc\left(u
\right)\sinc\left(u+\frac{\lambda}{\Omega}(r_2-r_2)\pi\right)du =
O\left(\frac{\Omega^{2}}{\lambda^{2}} \right).
\end{eqnarray*}
In the general case that $\lambda/\Omega\notin \mathbb{Z}$, $D_{2}$
remains a double sum.  In this case we make a change
of variables $u = \frac{\lambda}{2}(x+\omega_{\Omega,k_{1}+r_{1}})$ 
\begin{eqnarray*}
D_{2}
&=&\frac{c(0)}{\lambda}
\sum_{k_{1},k_{2}=-a}^{a}g(\omega_{\Omega,k_{1}})\overline{g(\omega_{\Omega,k_{2}})} 
\sinc\left(\frac{\lambda}{\Omega}(k_{1}-k_{2})\right)\\
&&\times \int
  f\left(\frac{2u}{\lambda}-\omega_{\Omega,k_{1}+r_{1}}\right)\sinc\left(u
\right)\sinc\left(u+\frac{\lambda}{\Omega}(k_{2}-k_{1}+r_{2}-r_{1})\pi
\right)du \\
&=& \frac{c(0)}{\lambda}\sum_{m=-2a}^{2a}
\sum_{k_{1}=\max(-a,-a+m)}^{\min(a,a+m)}g(\omega_{\Omega,k_{1}})\overline{g(\omega_{\Omega,m-k_{1}})} 
\sinc\left(\frac{\lambda}{\Omega}m\right)\\
&&\times \int
  f\left(\frac{2u}{\lambda}-\omega_{\Omega,k_{1}+r_{1}}\right)\sinc\left(u+\frac{\lambda}{\Omega}(m+r_{2}-r_{1})\pi
\right)\sinc\left(u
\right)du \\
&=& \frac{\Omega c(0)}{\lambda}\sum_{m=-2a}^{2a}\sinc\left(\frac{\lambda}{\Omega}(k_{1}-k_{2})\right)
\int_{\mathbb{R}}
\sinc\left(\frac{\lambda}{\Omega}m\right)\sinc\left(u+\frac{\lambda}{\Omega}(m+r_{2}-r_{1})\pi\right)\\
&&\times 
\frac{1}{\Omega}\sum_{k_{1}=\max(-a,-a+m)}^{\min(a,a+m)}g(\omega_{\Omega,k_{1}})\overline{g(\omega_{\Omega,m-k_{1}})} 
  f\left(\frac{2u}{\lambda}-\omega_{\Omega,k_{1}+r_{1}}\right)du 
\end{eqnarray*}
Replacing the sum $\frac{1}{\Omega}\sum_{k_1}$ in $D_{2}$
with an integral we have
\begin{eqnarray*}
D_{2}^{(1)} &=& \frac{c(0)\Omega}{\lambda 2\pi}\sum_{m=-2a}^{2a}
\sinc\left(\frac{\lambda}{\Omega}m\right)\int_{\mathbb{R}}
\int_{\max(-a,-a+m)/\Omega}^{\min(a,a+m)/\Omega}
g(\omega)\overline{g(\omega-\omega_{\Omega,k_{1}})} 
  f\left(\frac{2u}{\lambda}-\omega+\omega_{\Omega,r_{1}}\right)\\
&& \sinc\left(u
\right)\sinc\left(u+\frac{\lambda}{\Omega}(m+r_{2}-r_{1})\pi
\right)d\omega du.
\end{eqnarray*}
This gives the error
\begin{eqnarray*}
\left|D_{2} -  D_{2}^{(1)}\right| &\leq& 
\frac{C}{\lambda 2\pi}\sum_{m=-2a}^{2a}
\left|\sinc\left(\frac{\lambda}{\Omega}m\right)\right|\ell_{1}\left(\frac{\lambda}{\Omega}(m+r_2-r_1)\right)
%&=&O\left(\frac{\Omega}{\lambda^{2}}\log(a)\left[I_{\lambda/\Omega <1}\log(a) +
%    I_{\lambda/\Omega >1}\log (\frac{\lambda}{\Omega}a)\right]\right)
%\quad (\textrm{using }(\ref{eq:bounds1088}))
\end{eqnarray*}
In order to bound the above we use
\begin{eqnarray}
\sum_{m=-a}^{a}|\sinc(\alpha m \pi)\ell_{1}(\alpha m \pi)| &=& 
\sum_{|m|< e/\alpha}^{}|\sinc(\alpha m \pi)\ell_{1}(\alpha m \pi)|
+ \sum_{e/\alpha\leq |m|\leq a}^{}|\sinc(\alpha m \pi)\ell_{1}(\alpha
m \pi)|\nonumber\\
&\leq& C\sum_{|m|< e/\alpha}^{}|\frac{|\sin(\alpha m \pi)|}{|\alpha m \pi|}
+ C\log(\alpha a \pi)\sum_{e/\alpha\leq |m|\leq a}^{}\frac{1}{|\alpha
  \pi m|}\nonumber\\
&\leq& 
C\frac{[I_{\alpha \leq 1}\log a + I_{\alpha >1}\log(\alpha a)]\log a}{\alpha},\label{eq:bounds1088}
\end{eqnarray}
where $C$ is a finite constant and $\ell_{1}(m) = C\log(|m|)/m$ if
$|m|\geq e$ else $\ell_{1}(m)=C$. Thus 
\begin{eqnarray*}
\left|D_{2} -  D_{2}^{(1)}\right| &=& 
%\frac{C}{\lambda 2\pi}\sum_{m=-2a}^{2a}
%\left|\sinc\left(\frac{\lambda}{\Omega}m\right)\right|\ell_{1}(\frac{\lambda}{\Omega}(m+r_2-r_1))
O\left(\frac{\Omega}{\lambda^{2}}\log(a)\left[I_{\lambda/\Omega <1}\log(a) +
    I_{\lambda/\Omega >1}\log (\frac{\lambda}{\Omega}a)\right]\right).
\end{eqnarray*}
Replacing  $f\left(\frac{2u}{\lambda}-\omega+\omega_{r_{1}}\right)$
with $f\left(-\omega+\omega_{r_{1}}\right)$  in 
$D_{2}^{(1)}$ and using Lemma
\ref{lemma:1star}, equation (\ref{eq:lemma1star1}) gives 
\begin{eqnarray*}
D_{2} 
&=& D_{2}^{(2)} +
O\left(\frac{\Omega}{\lambda^{2}}\log(a)\left[I_{\lambda/\Omega \leq 1}\log(a) +
    I_{\lambda/\Omega >1}\log
    (\frac{\lambda}{\Omega}a)\right]^{2}\right) \\
\end{eqnarray*}
where 
\begin{eqnarray*}
&&D_{2}^{(2)} \\
&=& \frac{\Omega}{\lambda 2\pi}\sum_{m=-2a}^{2a}
\sinc\left(\frac{\lambda}{\Omega}m\right) \sinc\left(\frac{\lambda}{\Omega}(m+r_{2}-r_{1})\pi
\right)
\int_{\max(-a,-a+m)/\Omega}^{\min(a,a+m)/\Omega}
g(\omega)\overline{g(\omega-\omega_{\Omega,k_{1}})} 
  f\left(-\omega+\omega_{r_{1}}\right)d\omega.
\end{eqnarray*}
By using Lemma \ref{lemma:triangle} we have 
\begin{eqnarray*}
D_{2}^{(2)} = 
\left\{
\begin{array}{cc}
O\left(\frac{\Omega^{2}}{\lambda^{2}}\right) & \frac{\lambda}{\Omega}
< 1 \\
O\left( \frac{\Omega}{\lambda}\right) & \frac{\lambda}{\Omega} \geq 1\\
\end{array}
\right.
\end{eqnarray*}
Thus altogether the bound for $D_{2}$ is 
\begin{eqnarray}
D_{2} = 
\left\{
\begin{array}{cc}
O\left(\frac{\Omega^{2}}{\lambda^{2}}+\frac{\Omega \log^{3}(a)}{\lambda^{2}}\right) & \frac{\lambda}{\Omega}
< 1 \\
O\left( \frac{\Omega}{\lambda} +
\frac{\Omega
  \log^{}(a)\log^{2}(\frac{\lambda}{\Omega}a)}{\lambda^{2}}I_{\frac{\lambda}{\Omega}\notin
\mathbb{Z}}  \right) & \frac{\lambda}{\Omega} \geq 1\\
\end{array}
\right.. \label{eq:cumulantOmega2}
\end{eqnarray}
Since $F_{1} = O(\frac{\lambda^{2}}{\Omega^{2}n}[D_{1}+D_{2}])$ using
(\ref{eq:cumulantOmega1}) and (\ref{eq:cumulantOmega2})
gives 
\begin{eqnarray*}
F_{1} = 
\frac{\lambda^{2}}{\Omega^{2}n}\left\{
\begin{array}{cc}
O\left(\frac{\Omega^{2}}{\lambda^{2}}+\frac{\Omega \log^{3}(a)}{\lambda^{2}}\right) & \frac{\lambda}{\Omega}
< 1 \\
O\left( \frac{\Omega}{\lambda} +
\frac{\Omega
  \log^{}(a)\log^{2}(\frac{\lambda}{\Omega}a)}{\lambda^{2}}I_{\frac{\lambda}{\Omega}\notin
\mathbb{Z}}  \right) & \frac{\lambda}{\Omega} \geq 1\\
\end{array}
\right..
\end{eqnarray*}
Finally, since $F_{1,\mathrm{fine}} = \lambda F_{1}$ and 
$F_{1,\mathrm{coarse}} = \Omega F_{1}$ we have 
\begin{eqnarray*}
F_{1,\mathrm{fine}} = O\left(\frac{\lambda}{n}+
\frac{\lambda \log^{3}(a)}{\Omega n }\right) \quad \textrm{and}\quad
F_{1,\mathrm{coarse}} = O\left(\frac{\lambda}{n} +\frac{
  \log^{}(a)\log^{2}(\frac{\lambda}{\Omega}a)}{n}I_{\frac{\lambda}{\Omega}\notin
\mathbb{Z}}\right).
\end{eqnarray*}
A similar expression holds for 
$\lambda\cov\left[\widetilde{Q}_{a,\Omega,\lambda}(g;r_{1}),
\overline{\widetilde{Q}_{a,\Omega,\lambda}(g;r_{2})}\right]$. 
Thus we obtain the required result. 
\hfill $\Box$

\vspace{3mm}
Now we obtain approximations to
$A_{1}(\rb_1,\rb_2),\ldots,A_{4}(\rb_1,\rb_2)$ by separating the
$\sinc$ function from the spectral density. 
Let 
\begin{eqnarray*}
&& C_{11}\left(\frac{a}{\Omega},\ob_{\Omega,\rb} \right) = 
\frac{1}{(2\pi)^{d}}\int_{[-2\pi a/\Omega,2\pi a/\Omega]^{d}} 
f(\ob_{})f(\ob_{}+\ob_{\Omega,\rb_{}})
|g(\ob)|^{2}d\ob \\ 
&&C_{12}\left(\frac{a}{\Omega},\ob_{\Omega,\rb} \right) = 
\frac{1}{2\pi}\int_{[-2\pi a/\Omega,2\pi a/\Omega]^{d}} 
f(\ob_{})f(\ob_{}+\ob_{\Omega,\rb_{1}})
g(\ob)\overline{g(-\ob)}d\ob 
\end{eqnarray*}
\begin{eqnarray*}
&&C_{21}\left(\frac{a}{\Omega},\ob_{\Omega,\rb} \right) = 
\frac{1}{(2\pi)^{d}}\int_{[-2\pi a/\Omega,2\pi a/\Omega]^{d}} 
f(\ob)f(\ob_{}+\ob_{\Omega,\rb_{}})
g(\ob)^{2}d\ob \\
&&C_{22}\left(\frac{a}{\Omega},\ob_{\Omega,\rb} \right) = 
\frac{1}{(2\pi)^{d}}\int_{[-2\pi a/\Omega, 2\pi a/\Omega]^{d}} 
f(\ob_{})f(\ob_{}+\ob_{\Omega,\rb_{}})
g(\ob)g(-\ob)d\ob
\end{eqnarray*}
and define
\begin{eqnarray}
\label{eq:Awidetilde1}
\widetilde{A}_{1}(\rb_1,\rb_2) &=& C_{11}\left(\frac{a}{\Omega},\ob_{\Omega,\rb} \right)
\sum_{\mb=-2a}^{2a}\Sinc\left( \frac{\lambda}{\Omega}\mb
  \pi\right)\Sinc\left(
  \frac{\lambda}{\Omega}(\mb+\rb_{1}-\rb_{2})\pi\right) 
\end{eqnarray}
\begin{eqnarray}
\label{eq:Awidetilde2}
\widetilde{A}_{2}(\rb_1,\rb_2) &=& 
C_{12}\left(\frac{a}{\Omega},\ob_{\Omega,\rb} \right)
\sum_{\mb=-2a}^{2a}\Sinc\left(\frac{\lambda}{\Omega}(\mb+\rb_{1})\pi \right)
\Sinc\left(\frac{\lambda}{\Omega}(\mb+\rb_{2})\pi \right) 
\end{eqnarray}
\begin{eqnarray}
\label{eq:Awidetilde3}
\widetilde{A}_{3}(\rb_{1},\rb_{2}) 
&=& 
C_{21}\left(\frac{a}{\Omega},\ob_{\Omega,\rb} \right)
\sum_{\mb=-2a}^{2a}\Sinc\left( \frac{\lambda}{\Omega}\mb
  \pi\right)\Sinc\left( \frac{\lambda}{\Omega}(\mb+\rb_{1}+\rb_{2})\pi\right)
\end{eqnarray}
and 
\begin{eqnarray}
\label{eq:Awidetilde4}
\widetilde{A}_{4}(\rb_1,\rb_2) &=&
C_{22}\left(\frac{a}{\Omega},\ob_{\Omega,\rb} \right)
\sum_{\mb=-2a}^{2a}\Sinc\left(\frac{\lambda}{\Omega}(\mb+\rb_{1})\pi \right)
\Sinc\left(\frac{\lambda}{\Omega}(\mb-\rb_{2})\pi \right).
\end{eqnarray}

\begin{lemma}\label{lemma:fine-grid-gaussian}
Suppose Assumptions \ref{assum:S}, \ref{assum:uniform} and  
Assumptions \ref{assum:G}(ii) and \ref{assum:GG}(b) hold. 
Let
$\widetilde{A}_{1}(\rb_1,\rb_2),\ldots,\widetilde{A}_{4}(\rb_1,\rb_2)$
be defined as in (\ref{eq:Awidetilde1}), (\ref{eq:Awidetilde2}),
(\ref{eq:Awidetilde3}) and (\ref{eq:Awidetilde4}). 
Then for $|r_1|,|r_2|\leq a$ and $1\leq j \leq 4$
\begin{itemize}
\item[(i)] if $\lambda/\Omega<1$ then 
\begin{eqnarray*}
\frac{\lambda^{d}}{\Omega^{2d}}A_{j}(\rb_{1},\rb_{2}) &=& 
\frac{\lambda^{d}}{\Omega^{d}}\widetilde{A}_{j}(\rb_1,\rb_2) + 
O\left(\frac{\log^{2}(a)[\log a +\log \lambda]}{\lambda} \right) 
\end{eqnarray*}
\item[(ii)] if $\lambda/\Omega \geq 1$ then
\begin{eqnarray*}
\frac{1}{\Omega^{d}}A_{j}(\rb_{1},\rb_{2}) &=& 
\widetilde{A}_{j}(\rb_1,\rb_2) + 
O\left(
\frac{(\log^{2}a)(\log \lambda + \log \frac{\lambda a}{\Omega})}{\Omega} \right) 
\end{eqnarray*}
\end{itemize}
\end{lemma}
PROOF.  We prove the result for $j=1$ and $d=1$. 
We first define a sequence of approximations. Let 
\begin{eqnarray*}
A_{1}^{(1)}(r_{1},r_{2}) &=& \frac{1}{\pi^{2}}
\sum_{m=-2a}^{2a}\frac{1}{2\pi}\int_{2\pi\max(-a,-a+m)/\Omega}^{2\pi\min(a,a+m)/\Omega} 
\int_{\mathbb{R}^{2}} f(\frac{2u}{\lambda} -\omega_{})
f(\frac{2v}{\lambda} + \omega_{}+\omega_{\Omega,r_{1}})\times \\
&& g(\omega)\overline{g(\omega-\omega_{\Omega,m})}
\sinc\left(u-\frac{\lambda}{\Omega}m \pi\right)\sinc\left(v+\frac{\lambda}{\Omega}(m+r_{1}-r_{2})\pi\right)
\sinc(u)\sinc(v)du dv d\omega \\
A_{1}^{(2)}(r_{1},r_{2}) &=& \frac{1}{\pi^{2}}
\sum_{m=-2a}^{2a}\frac{1}{2\pi}\int_{2\pi\max(-a,-a+m)/\Omega}^{2\pi\min(a,a+m)/\Omega} 
f(\omega_{})
f(\omega_{}+\omega_{\Omega,r_{1}})
g(\omega)\overline{g(\omega-\omega_{\Omega,m})} d\omega \\
&&\times \int_{\mathbb{R}^{2}} 
\sinc\left(u-\frac{\lambda}{\Omega}m \pi\right)
\sinc\left(v+\frac{\lambda}{\Omega}(m+r_{1}-r_{2})\pi\right)
\sinc(u)\sinc(v)du dv\\
 &=& \frac{1}{\pi^{2}}
\sum_{m=-2a}^{2a}\sinc\left( \frac{\lambda}{\Omega}m
  \pi\right)\sinc\left( \frac{\lambda}{\Omega}(m+r_{1}-r_{2})\pi\right)\\
&&\times \frac{1}{2\pi}\int_{2\pi\max(-a,-a+m)/\Omega}^{2\pi\min(a,a+m)/\Omega} 
f(\omega_{})
f(\omega_{}+\omega_{\Omega,r_{1}})
g(\omega)\overline{g(\omega-\omega_{\Omega,m})} d\omega. 
\end{eqnarray*}
We first prove part (i) ($\lambda/\Omega<1$).
By using Lemma \ref{lemma:sum-integral} (to replace summand by
integral) and Lemma \ref{lemma:1a}, equation (\ref{eq:i1}), we have 
\begin{eqnarray*}
&& \frac{\lambda}{\Omega}\left|\frac{1}{\Omega} A_{1}(r_{1},r_{2})  - A_{1}^{(1)}(r_{1},r_{2})  \right| \\
&\leq& \frac{C\lambda^{}}{\Omega^{2}}
\sum_{m=-2a}^{2a}\int_{\mathbb{R}^{2}}\left|
\sinc\left(u-\frac{\lambda}{\Omega}m \pi\right)\sinc\left(v+\frac{\lambda}{\Omega}(m+r_{1}-r_{2})\pi\right)
\sinc(u)\sinc(v)\right|du dv \\
&\leq& \frac{C\lambda^{}}{\Omega^{2}}
\sum_{m=-2a}^{2a}\ell_{1}\left(\frac{\lambda}{\Omega}m\pi
\right)\ell_{1}\left( \frac{\lambda}{\Omega}(m+r_{1}-r_{2})\pi\right)
= O\left( \frac{\log^{2}(a)}{\lambda} \right).
\end{eqnarray*}
We replace
$f(\frac{2u}{\lambda}-\omega)f(\frac{2v}{\lambda}+\omega+\omega_{\Omega,r_1})$
with $f(-\omega)f(\omega+\omega_{\Omega,r_1})$. 
By using Lemma \ref{lemma:1star}, equation (\ref{eq:lemma1star1})
(twice) we have 
\begin{eqnarray*}
&&\frac{\lambda}{\Omega}\left|A_{1}^{(1)}(r_{1},r_{2})  - A_{1}^{(2)}(r_{1},r_{2})\right| \\
&\leq& \frac{C\lambda^{}}{\Omega^{}}
\sum_{m=-2a}^{2a}\left[\ell_{1}\left(\frac{\lambda m
      \pi}{\Omega}\right) \frac{\log \lambda +
    I(\frac{\lambda}{\Omega}|m+r_{1}-r_{2}|>e)\log \frac{\lambda}{\Omega}|m+r_{1}-r_{2}| }{\lambda}
\right]+\\
&& \frac{C\lambda^{}}{\Omega^{}}
\sum_{m=-2a}^{2a}\left[\ell_{1}\left(\frac{\lambda (m+r_{1}-r_{2})
      \pi}{\Omega}\right) \frac{\log \lambda +
    I(\frac{\lambda}{\Omega}|m|>e)\log \frac{\lambda}{\Omega}|m| }{\lambda}
\right] = O\left(\frac{\log^{2}(a)[\log \lambda + \log^{}(a)]}{\lambda}\right).
\end{eqnarray*}
The final approximation is based on 
replacing $\overline{g(\omega-\omega_{\Omega,m})}$ with 
$\overline{g(\omega)}$ and the limits of the integral,
$\int_{\max(-a,-a+m)/\Omega}^{\min(a,a+m)/\Omega} $ with $\int_{-a/\Omega}^{a/\Omega}$. To do this we note 
\begin{eqnarray*}
&&\frac{\lambda}{\Omega}\left( A_{1}^{(2)}(r_{1},r_{2}) -
  \widetilde{A}_{1}(r_{1},r_{2}) \right) = 
\frac{\lambda}{\Omega^{}}
\sum_{m=-2a}^{2a}\sinc\left( \frac{\lambda}{\Omega}m
  \pi\right)\sinc\left(
  \frac{\lambda}{\Omega}(m+r_{1}-r_{2})\pi\right)\times \\
&& \frac{1}{2\pi}\left( \int_{\max(-a,-a+m)/\Omega}^{\min(a,a+m)/\Omega} 
f(\omega_{})^{2}g(\omega)\overline{g(\omega-\omega_{\Omega,m})}
d\omega
- 
\int_{-a/\Omega}^{a/\Omega} 
f(\omega_{})^{2}|g(\omega)|^{2}
\right)d\omega.
\end{eqnarray*}
%Under Assumption \ref{assum:G} $g$ is Lipschitz continuous, thus by
%the 
Next we bound (we assume $m>0$, though the same proof holds for $m<0$)
\begin{eqnarray*}
&&\left| \int_{2\pi\max(-a,-a+m)/\Omega}^{2\pi\min(a,a+m)/\Omega} 
f(\omega_{})f(\omega+\omega_{\Omega,r})g(\omega)\overline{g(\omega-\omega_{\Omega,m})}
d\omega
- 
\int_{-a/\Omega}^{a/\Omega} 
f(\omega_{})f(\omega+\omega_{\Omega,r})|g(\omega)|^{2}
\right| \nonumber\\
&\leq& \left| \int_{2\pi\max(-a,-a+m)/\Omega}^{2\pi\min(a,a+m)/\Omega} 
f(\omega_{})f(\omega+\omega_{\Omega,r})g(\omega)
\left[\overline{g(\omega-\omega_{\Omega,m})}-
\overline{g(\omega)}\right]d\omega\right| +\nonumber\\
&&\left| \int_{-2\pi a/\Omega}^{-2\pi(a-m)/\Omega} 
f(\omega_{})f(\omega+\omega_{\Omega,r})|g(\omega)|^{2}d\omega
\right|.
\end{eqnarray*}
By using the Lipschitz continuity of $g$ and the mean value
theorem twice we have 
\begin{eqnarray}
\label{eq:sum-int-proof}
&&\left| \int_{2\pi\max(-a,-a+m)/\Omega}^{2\pi\min(a,a+m)/\Omega} 
f(\omega_{})f(\omega+\omega_{\Omega,r})g(\omega)\overline{g(\omega-\omega_{\Omega,m})}
d\omega
- 
\int_{-a/\Omega}^{a/\Omega} 
f(\omega_{})f(\omega+\omega_{\Omega,r})|g(\omega)|^{2}
\right| \leq C\frac{|m|}{\Omega}\qquad 
\end{eqnarray}
where $C$ is a finite constant which only depends on $f$ and $g$. This
bound gives
\begin{eqnarray*}
\frac{\lambda}{\Omega}\left|A_{1}^{(2)}(r_1,r_2) - \widetilde{A}_{1}(r_1,r_2)\right| &\leq&  
\frac{C}{\lambda^{}}
\sum_{m=-2a}^{2a}\left|\sin\left( \frac{\lambda}{\Omega} m\pi\right)
\sin\left(
  \frac{\lambda}{\Omega}(m+r_{1}-r_{2})\pi\right)\right|\ell_{0}(|m+r_1-r_2|) \\
&=& 
O\left(\frac{\log (a+|r_1|+|r_2|)}{\lambda}\right).
\end{eqnarray*}
Thus altogether the three bounds above give 
\begin{eqnarray*}
\frac{\lambda}{\Omega^{2}}A_{j}(r_{1},r_{2}) &=& 
\frac{\lambda}{\Omega}\widetilde{A}_{j}(r_1,r_2) + 
O\left(\frac{\log^{2}(a)[\log a +\log \lambda]}{\lambda} \right) 
\end{eqnarray*}
and same sequence of bounds apply to $A_{j}(r_1,r_2)$ for $2\leq j \leq
4$, thus proving (i).

Now we prove (ii) ($\lambda/\Omega >1$).
From Lemma \ref{lemma:triangle} we observe that the rate of
growth of $A_{j}(r_1,r_2)$ and $\widetilde{A}_{j}(r_1,r_2)$ are different when 
$\lambda/\Omega <1$ and $\lambda/\Omega >1$, thus requiring different
standardisations.  
By using Lemma \ref{lemma:sum-integral} (to replace summand by
integral) and Lemma \ref{lemma:1a}, equation (\ref{eq:i1}), we have 
\begin{eqnarray*}
&& \left|\frac{1}{\Omega} A_{1}(r_{1},r_{2})  - A_{1}^{(1)}(r_{1},r_{2})  \right| \\
&\leq& \frac{C}{\Omega}
\sum_{m=-2a}^{2a}\int_{\mathbb{R}^{2}}\left|
\sinc\left(u-\frac{\lambda}{\Omega}m \pi\right)\sinc\left(v+\frac{\lambda}{\Omega}(m+r_{1}-r_{2})\pi\right)
\sinc(u)\sinc(v)\right|du dv \\
&\leq& \frac{C}{\Omega^{}}
\sum_{m=-2a}^{2a}\ell_{1}\left(\frac{\lambda}{\Omega}m\pi
\right)\ell_{1}\left( \frac{\lambda}{\Omega}(m+r_{1}-r_{2})\pi\right)
= O\left(\frac{\Omega\log^{2}(a)}{\lambda^{2}} \right).
\end{eqnarray*}
By using Lemma \ref{lemma:1star}, equation (\ref{eq:lemma1star1})
(twice) we have 
\begin{eqnarray*}
&&\left|A_{1}^{(1)}(r_{1},r_{2})  - A_{1}^{(2)}(r_{1},r_{2})\right| \\
&\leq& 
C\sum_{m=-2a}^{2a}\left[\ell_{1}\left(\frac{\lambda m
      \pi}{\Omega}\right) \frac{\log \lambda +
    I(\frac{\lambda}{\Omega}|m+r_{1}-r_{2}|>e)\log \frac{\lambda}{\Omega}|m+r_{1}-r_{2}| }{\lambda}
\right]+\\
&& C\sum_{m=-2a}^{2a}\left[\ell_{1}\left(\frac{\lambda (m+r_{1}-r_{2})
      \pi}{\Omega}\right) \frac{\log \lambda +
    I(\frac{\lambda}{\Omega}|m|>e)\log \frac{\lambda}{\Omega}|m| }{\lambda}
\right] = O\left(\frac{\log^{2}(a)[\log \lambda + \log^{}(\frac{\lambda}{\Omega}a)]^{}}{\Omega}\right).
\end{eqnarray*}
Replacing $\overline{g(\omega-\omega_{\Omega,m})}$ with 
$\overline{g(\omega)}$ and limits of the integral,
$\int_{\max(-a,-a+m)/\Omega}^{\min(a,a+m)/\Omega} $ with $\int_{-a/\Omega}^{a/\Omega}$
\begin{eqnarray*}
&&A_{1}^{(2)}(r_{1},r_{2}) -
  \widetilde{A}_{1}(r_{1},r_{2}) = 
\sum_{m=-2a}^{2a}\sinc\left( \frac{\lambda}{\Omega}m
  \pi\right)\sinc\left(
  \frac{\lambda}{\Omega}(m+r_{1}-r_{2})\pi\right)\times \\
&& \frac{1}{2\pi}\left( \int_{2\pi\max(-a,-a+m)/\Omega}^{2\pi\min(a,a+m)/\Omega} 
f(\omega_{}) f(\omega+\omega_{\Omega,r}) g(\omega)\overline{g(\omega-\omega_{\Omega,m})}
d\omega
- 
\int_{-a/\Omega}^{a/\Omega} 
f(\omega_{})f(\omega+\omega_{\Omega,r})|g(\omega)|^{2}
\right)d\omega,
\end{eqnarray*}
where $\widetilde{A}_{1}(r_{1},r_{2})$ is defined in (\ref{eq:Awidetilde1}).
Finally, using (\ref{eq:sum-int-proof}) we have 
\begin{eqnarray*}
\left|A_{1}^{(2)}(r_1,r_2) - \widetilde{A}_{1}(r_1,r_2)\right| &\leq&  
\frac{\Omega}{\lambda^{2}}
\sum_{m=-2a}^{2a}\left|\sin\left( \frac{\lambda}{\Omega} m\pi\right)
\sin\left(
  \frac{\lambda}{\Omega}(m+r_{1}-r_{2})\pi\right)\right|\ell_{0}(m+r_1-r_2) \\
&=& 
O\left(\frac{\Omega\log (\frac{\lambda}{\Omega}a)}{\lambda^{2}}\right).
\end{eqnarray*}
Thus the three bounds together give 
\begin{eqnarray*}
\frac{1}{\Omega}A_{1}(r_1,r_2) = \widetilde{A}_{1}(r_1,r_2) + O\left( 
\frac{(\log^{2}a)[\log \lambda+\log(a\frac{\lambda}{\Omega})]}{\Omega}
\right).
\end{eqnarray*}
The same sequence bounds apply to $A_{j}(r_1,r_2)$ for $2\leq j \leq
4$. Thus proving (ii). \hfill $\Box$
\vspace{3mm}

We now summarize the key points in the above result. Using the above result we have
\begin{eqnarray*}
&&\var\left[\widetilde{Q}_{a,\Omega,\lambda}(g;\rb)\right] \\
 &=&  C_{11}\left(\frac{a}{\lambda},\ob_{\Omega,\rb
}\right)
\frac{1}{\Omega^{d}}\sum_{\mb=-2a}^{2a}\sinc^{2}\left( \frac{\lambda}{\Omega}\mb
  \pi\right) +C_{12}\left(\frac{a}{\lambda},\ob_{\Omega,\rb
}\right)
\frac{1}{\Omega^{d}}\sum_{m=-2a}^{2a}\sinc^{2}\left( \frac{\lambda}{\Omega}(\mb+\rb)
  \pi\right) \\
&& + 
O\left(\widetilde{\ell}_{a,\Omega,\lambda}\left[\frac{1}{\lambda^{d}}I_{\frac{\lambda}{\Omega}<1}
  + \frac{1}{\Omega^{d}}I_{\frac{\lambda}{\Omega}\geq 1}\right]\right).
\end{eqnarray*}
where 
\begin{eqnarray}
\label{eq:widetildeell}
\widetilde{\ell}_{a,\Omega,\lambda} = 
\left\{
\begin{array}{cc}
\frac{\lambda^{d}}{n} + \frac{\lambda^{d}\log^{3}(a)}{\Omega^{d} n} + \frac{\log^{2}(a)[\log
  \lambda + \log a]}{\lambda^{}} & \frac{\lambda}{\Omega} <1\\
\frac{\lambda^{d}}{ n} + \frac{\log
  a\log^{2}(\frac{\lambda}{\Omega}a)}{n}I_{\frac{\lambda}{\Omega}\notin
\mathbb{Z}} + \frac{\log^{2}(a)[\log
\lambda+\log(\frac{\lambda}{\Omega}a)]}{\Omega^{}} &
\frac{\lambda}{\Omega}\geq 1 \\
\end{array}
\right..
\end{eqnarray}
This result shows that the rate of convergence of the variance is 
determined by the term 
\begin{eqnarray*}
\frac{1}{\Omega}\sum_{m=-2a}^{2a}\sinc^{2}\left( \frac{\lambda}{\Omega}m
  \pi\right).
\end{eqnarray*}
Using Lemma \ref{lemma:triangle} we see that the rate that the above
term converges to zero depends on whether $\Omega <\lambda$ or $\Omega \geq \lambda$.
In particular setting $a=\infty$
\begin{eqnarray*}
\frac{1}{\Omega}\sum_{m=-\infty}^{\infty}\sinc^{2}\left( \frac{\lambda}{\Omega}m
  \pi\right) = 
\left\{
\begin{array}{cc}
\frac{1}{\lambda} & \frac{\lambda}{\Omega} <1 \\
\frac{1}{\Omega} & \frac{\lambda}{\Omega} \in \mathbb{Z} \\
O\left(\frac{1}{\Omega}\right) & \frac{\lambda}{\Omega} > 1 \textrm{
  and }\frac{\lambda}{\Omega} \notin \mathbb{Z}
\end{array}
\right.
\end{eqnarray*}
and 
$\sum_{m=-\infty}^{\infty}\sinc^{2}\left( \frac{\lambda}{\Omega}m
  \pi\right)\rightarrow 1$ as $\frac{\lambda}{\Omega}\rightarrow \infty$ . 
Therefore for the
fine frequency grid with $\lambda < \Omega$ we have 
\begin{eqnarray*}
&&\lambda^{d}\var[\widetilde{Q}_{a,\Omega,\lambda}(g;\rb)] \\
&=&  C_{11}\left(\frac{a}{\lambda},\ob_{\Omega,\rb}\right)
\frac{\lambda^{d}}{\Omega^{d}}\sum_{\mb=-2a}^{2a}
\Sinc^{2}\left( \frac{\lambda}{\Omega}\mb\pi\right) +C_{12}\left(\frac{a}{\lambda},\ob_{\Omega,\rb}\right)
\frac{\lambda^{d}}{\Omega^{d}}\sum_{m=-2a}^{2a}\Sinc^{2}\left( \frac{\lambda}{\Omega}(\mb+\rb)
  \pi\right) 
 + 
O\left(\widetilde{\ell}_{a,\Omega,\lambda}\right).
\end{eqnarray*}
On the other hand if a coarse frequency grid is used with $\lambda
\geq \Omega$ then the rate of convergence is worse with  
\begin{eqnarray*}
&&\Omega^{d}\var[\widetilde{Q}_{a,\Omega,\lambda}(g;\rb)] \\
&=&  C_{11}\left(\frac{a}{\lambda},\ob_{\Omega,\rb}\right)
\sum_{\mb=-2a}^{2a}
\Sinc^{2}\left( \frac{\lambda}{\Omega}\mb\pi\right) +C_{12}\left(\frac{a}{\lambda},\ob_{\Omega,\rb}\right)
\sum_{m=-2a}^{2a}\Sinc^{2}\left( \frac{\lambda}{\Omega}(\mb+\rb)
  \pi\right) 
 + 
O\left(\widetilde{\ell}_{a,\Omega,\lambda}\right).
\end{eqnarray*}
The above result implies that 
\begin{eqnarray*}
\var[\widetilde{Q}_{a,\Omega,\lambda}(g;\rb)] = 
\left\{
\begin{array}{cc}
O(\lambda^{-d}) & \frac{\lambda}{\Omega} <1 \\
O(\Omega^{-d}) & \frac{\lambda}{\Omega} \geq 1\\
\end{array}
\right.
\end{eqnarray*}

The above results assume the spatial process is Gaussian. We now relax
the assumption of Gaussianity. 

\begin{theorem}
Let us suppose that  $\{Z(\ub);\ub \in\mathbb{R}^{d}\}$ is a fourth order stationary spatial 
random field that satisfies Assumption \ref{assum:S}(i),  
\ref{assum:uniform}, \ref{assum:G}, \ref{assum:nonGaussian} and \ref{assum:GG}(a,c) or \ref{assum:GG}(b,c)  are satisfied. 
%Suppose Assumptions \ref{assum:S}, \ref{assum:uniform} and  
%Assumptions \ref{assum:G}(ii) and \ref{assum:GG}(b) hold. 
Then for $|\rb_1|,|\rb_2|\leq a$ 
\begin{itemize}
\item[(i)] If $\lambda< \Omega$ 
\begin{eqnarray*}
\lambda^{d}\cov\left[\widetilde{Q}_{a,\Omega,\lambda}(g;\rb_{1}),
\widetilde{Q}_{a,\Omega,\lambda}(g;\rb_{2})\right] 
= \frac{\lambda^{d}}{\Omega^{d}}\left(A_{1}(\rb_{1},\rb_{2}) + 
A_{2}(\rb_{1},\rb_{2})\right) + B_{1}(\rb_{1},\rb_{2}) +
+ O\left(\widetilde{\ell}^{(2)}_{a,\Omega,\lambda}
\right),
\end{eqnarray*}
\begin{eqnarray*}
\lambda^{d}\cov\left[\widetilde{Q}_{a,\Omega,\lambda}(g;\rb_{1}),
\overline{\widetilde{Q}_{a,\Omega,\lambda}(g;\rb_{2})}\right] 
= \frac{\lambda^{d}}{\Omega^{d}}\left(A_{3}(\rb_{1},\rb_{2}) +
  A_{4}(\rb_{1},\rb_{2})\right) + 
B_{3}(\rb_{1},\rb_{2}) 
+ O\left(\widetilde{\ell}^{(2)}_{a,\Omega,\lambda}
\right),
\end{eqnarray*}
\item[(ii)]  If $\lambda \geq \Omega$ 
\begin{eqnarray*}
\Omega^{d}\cov\left[\widetilde{Q}_{a,\Omega,\lambda}(g;\rb_{1}),
\widetilde{Q}_{a,\Omega,\lambda}(g;\rb_{2})\right] 
= \frac{1}{\Omega^{d}}\left(A_{1}(\rb_{1},\rb_{2}) + A_{2}(\rb_{1},\rb_{2})\right) + \frac{\Omega^{d}}{\lambda^{d}}
B_{1}(\rb_{1},\rb_{2}) + 
+ O\left(\widetilde{\ell}^{(2)}_{a,\Omega,\lambda}
\right),
\end{eqnarray*}
\begin{eqnarray*}
\Omega^{d}\cov\left[\widetilde{Q}_{a,\Omega,\lambda}(g;r_{1}),
\overline{\widetilde{Q}_{a,\Omega,\lambda}(g;r_{2})}\right] 
= \frac{1}{\Omega^{d}}\left( A_{3}(\rb_{1},\rb_{2}) + A_{4}(\rb_{1},\rb_{2})\right) +
\frac{\Omega^{d}}{\lambda^{d}}B_{3}(\rb_{1},\rb_{2}) 
+ O\left(\widetilde{\ell}^{(2)}_{a,\Omega,\lambda}
\right),
\end{eqnarray*}
\end{itemize}
where 
\begin{eqnarray}
\label{eq:B1r1r2}
B_{1}(\rb_{1},\rb_{2}) 
 &=&   \frac{c_{4}}{\pi^{2d+1}\Omega^{2d}}\sum_{\kb_{1},\kb_{2}=-a}^{a}\int_{\mathbb{R}^{3d}}g(\ob_{\kb_{1}})\overline{g(\ob_{\kb_{2}})}
f_{4}\left(\frac{2\ubb_{1}}{\lambda}-\ob_{\Omega,\kb_{1}+\rb_{1}},\frac{2\ubb_{2}}{\lambda}-\ob_{\Omega,\kb_2},\frac{2\ubb_{3}}{\lambda}+
\ob_{\Omega,\kb_2 + \rb_2}\right)
 \times \nonumber\\
&&\times\Sinc\left(\ubb_{1}+\ubb_{2}+\ubb_{3}+\frac{\lambda}{\Omega}(\rb_{2}-\rb_{1})\pi\right)\Sinc(\ubb_{1})
\Sinc(\ubb_{2})\Sinc(\ubb_{3})d\ubb_{1}d\ubb_{2}d\ubb_{3},
\end{eqnarray}
\begin{eqnarray*}
%\label{eq:B1r1r2}
B_{3}(\rb_{1},\rb_{2}) 
 &=&   \frac{c_{4}}{\pi^{2d+1}\Omega^{2d}}\sum_{\kb_{1},\kb_{2}=-a}^{a}\int_{\mathbb{R}^{3d}}g(\ob_{\kb_{1}})g(\ob_{\kb_{2}})
f_{4}\left(\frac{2\ubb_{1}}{\lambda}-\ob_{\Omega,\kb_{1}+\rb_{1}},\frac{2\ubb_{2}}{\lambda}+\ob_{\Omega,\kb_2},\frac{2\ubb_{3}}{\lambda}-
\ob_{\Omega,\kb_2 + \rb_2}\right)
 \times \nonumber\\
&&\times\Sinc\left(\ubb_{1}+\ubb_{2}+\ubb_{3}-\frac{\lambda}{\Omega}(\rb_{2}+\rb_{1})\pi\right)\Sinc(\ubb_{1})
\Sinc(\ubb_{2})\Sinc(\ubb_{3})d\ubb_{1}d\ubb_{2}d\ubb_{3}
\end{eqnarray*}
and 
\begin{eqnarray}
\label{eq:widetildeell2}
\widetilde{\ell}_{a,\Omega,\lambda}^{(2)} = 
\left\{
\begin{array}{cc}
\frac{\lambda^{d}}{n} + \frac{\lambda^{d}\log^{3}(a)}{\Omega^{d} n} + \frac{\log^{2}(a)[\log
  \lambda + \log a]}{\lambda^{}} +\frac{\lambda^{2d}a^{2}}{n^{2}\Omega^{d}} & \frac{\lambda}{\Omega} <1\\
\frac{\lambda^{d}}{ n} + \frac{\log
  a\log^{2}(\frac{\lambda}{\Omega}a)}{n}I_{\frac{\lambda}{\Omega}\notin
\mathbb{Z}} + \frac{\log^{2}(a)[\log
\lambda+\log(\frac{\lambda}{\Omega}a)]}{\Omega^{}} + \frac{\lambda^{d}a^{d}}{n^{2}} &
\frac{\lambda}{\Omega}\geq 1 \\
\end{array}
\right..
\end{eqnarray}
\end{theorem}
PROOF. The proof is similar to the proof of Theorem
\ref{theorem:variance-nongaussian}.  We focus on the case
$d=1$. Expanding out (for both the case $\lambda/\Omega<1$ and
$\lambda/\Omega\geq 1$) we have  
\begin{eqnarray*}
&&\cov\left[\widetilde{Q}_{a,\Omega,\lambda}(g;r_{1}),\widetilde{Q}_{a,\Omega,\lambda}(g;r_{2})\right] \\
&=& \frac{1}{\Omega}A_{1}(r_{1},r_{2}) +
\frac{1}{\Omega}A_{2}(r_{1},r_{2}) + \widetilde{B}_{1}(r_{1},r_{2}) + 
\widetilde{B}_{2}(r_{1},r_{2}) +
O\left(\widetilde{\ell}_{a,\Omega,\lambda}\left[\frac{1}{\lambda}I_{\frac{\lambda}{\Omega}<1}
  + \frac{1}{\Omega}I_{\frac{\lambda}{\Omega}\geq 1}\right] \right),
\end{eqnarray*}
where $A_{1}(r_{1},r_{2})$ and $A_{2}(r_{1},r_{2})$ are defined in
Lemma \ref{lemma:covariance} and 
\begin{eqnarray*}
\widetilde{B}_{1}(r_{1},r_{2}) &=& \frac{c_{4}\lambda^{2}}{\Omega^{2}}
\sum_{k_{1},k_{2}=-a}^{a}
g(\omega_{\Omega,k_{1}})\overline{g(\omega_{\Omega,k_{2}})}\\
&&\Ex\bigg[\kappa_{4}(s_{2}-s_{1},s_{3}-s_{1},s_{4}-s_{1})e^{is_{1}\omega_{\Omega,k_{1}}}e^{-is_{2}\omega_{\Omega,k_{1}+r_{1}}} 
 e^{-is_{3}\omega_{\Omega,k_{2}}}e^{is_{4}\omega_{\Omega,k_{2}+r_{2}}} \bigg] \\
\widetilde{B}_{2}(r_{1},r_{2}) &=& \frac{\lambda^{2}}{\Omega^{2}n^{4}}\sum_{j_{1},\ldots,j_{4}\in \mathcal{D}_{3}}\sum_{k_{1},k_{2}=-a}^{a}
g(\omega_{\Omega,k_{1}})\overline{g(\omega_{\Omega,k_{2}})}\times \\
 &&\Ex\bigg[\kappa_{4}(s_{j_2}-s_{j_1},s_{j_3}-s_{j_1},s_{j_4}-s_{j_1})e^{is_{j_1}\omega_{\Omega,k_{1}}}e^{-is_{j_2}\omega_{\Omega,k_{1}+r_{1}}}
e^{-is_{j_3}\omega_{\Omega,k_{2}}}e^{is_{j_4}\omega_{\Omega,k_{2}+r_{2}}} \bigg]
\end{eqnarray*}
The focus in this proof will be on $\widetilde{B}_{1}(r_1,r_2)$ and
$\widetilde{B}_{2}(r_1,r_2)$. First we consider the ``leading term''
$\widetilde{B}_{1}(r_{1},r_{2})$ 
\begin{eqnarray*}
\widetilde{B}_{1}(r_{1},r_{2})&=& \frac{c_{4}\lambda^{2}}{(2\pi)^{3}\Omega^{2}\lambda^{4}}\sum_{k_{1},k_{2}=-a}^{a}g(\omega_{\Omega,k_{1}})\overline{g(\omega_{\Omega,k_{2}})}
\int_{\mathbb{R}^{3}}
f_{4}(\omega_{1},\omega_{2},\omega_{3})\int_{[-\lambda/2,\lambda/2]^{4}}
 e^{is_{1}(\omega_{1}+\omega_{2}+\omega_{3}+\omega_{k_{1}})} \\
&&e^{-is_{2}(\omega_{1}+\omega_{\Omega,k_{1}+r_{1}})}e^{-is_{3}(\omega_{2}+\omega_{\Omega,k_{2}})} 
e^{is_{4}(-\omega_{3}+\omega_{\Omega,k_{2}+r_{2}})}ds_{1}ds_{2}ds_{3}ds_{4}d\omega_{1}d\omega_{2}d\omega_{3}
\\
&=& \frac{c_{4}\lambda^{2}}{(2\pi)^{3}\Omega^{2}}\sum_{k_{1},k_{2}=-a}^{a}g(\omega_{k_{1}})\overline{g(\omega_{k_{2}})}
\int_{\mathbb{R}^{3}}
f_{4}(\omega_{1},\omega_{2},\omega_{3})
\sinc\left(\frac{\lambda}{2}\left(\omega_{1}+\omega_{2}+\omega_{3}+\omega_{\Omega,k_{1}}\right)\right) \\
&&\times\sinc\left(\frac{\lambda}{2}\left(\omega_{1} +
  \omega_{\Omega,k_{1}+r_{1}}\right)\right)
\sinc\left(\frac{\lambda}{2}\left(\omega_{2}+\omega_{\Omega,k_2}\right)\right)
\sinc\left(\frac{\lambda}{2}\left(\omega_{3}-\omega_{\Omega,k_{2}+r_{2}}\right)\right)
d\omega_{1}d\omega_{2}d\omega_{3}. 
\end{eqnarray*}
Making a change of variables we have 
\begin{eqnarray*}
\widetilde{B}_{1}(r_{1},r_{2}) 
 &=&   \frac{c_{4}}{\pi^{3}\Omega^{2}\lambda}\sum_{k_{1},k_{2}=-a}^{a}\int_{\mathbb{R}^{3}}g(\omega_{\Omega,k_{1}})\overline{g(\omega_{\Omega,k_{2}})}
f_{4}\left(\frac{2u_{1}}{\lambda}-\omega_{\Omega,k_{1}+r_{1}},\frac{2u_{2}}{\lambda}-\omega_{\Omega,k_2},\frac{2u_{3}}{\lambda}+
\omega_{\Omega,k_2 + r_2}\right)
 \times \nonumber\\
&&\times\sinc\left(u_{1}+u_{2}+u_{3}+\frac{\lambda}{\Omega}(r_{2}-r_{1})\pi\right)\sinc(u_{1})
\sinc(u_{2})\sinc(u_{3})du_{1}du_{2}du_{3},\label{eq:tildeB3}
\end{eqnarray*}
thus we observe that $\lambda
\widetilde{B}_{1}(r_1,r_2)=B_{1}(r_{1},r_{2})$, defined in
(\ref{eq:B1r1r2}).

We now show that $\widetilde{B}_{2}(r_{1},r_{2})$ is of lower order. To do so,
just as in the proof of Theorem \ref{theorem:variance-nongaussian} we
have 
\begin{eqnarray}
\label{eq:B2bound}
\widetilde{B}_{2}(r_{1},r_{2})
=\sum_{j=1}^{4}\widetilde{B}_{2,(3,j)}(r_{1},r_{2})+
\sum_{j=1}^{2}\widetilde{B}_{2,(2,j)}(r_{1},r_{2})
\end{eqnarray}
 where  
 \begin{eqnarray*}
 \widetilde{B}_{2,(3,1)}(r_{1},r_{2}) 
&=& \frac{|\mathcal{D}_{3,1}|\lambda^{2}}{n^{4}\Omega^{2}}\sum_{k_{1},k_{2}=-a}^{a}
g(\omega_{\Omega,k_{1}})\overline{g(\omega_{\Omega,k_{2}})}\times \\
 &&\Ex\bigg[\kappa_{4}(s_{j_2}-s_{j_1},0,s_{j_4}-s_{j_1})e^{is_{j_1}\omega_{\Omega,k_{1}}}e^{-is_{j_2}\omega_{\Omega,k_{1}+r_{1}}}
e^{-is_{j_1}\omega_{\Omega,k_{2}}}e^{is_{j_4}\omega_{\Omega,k_{2}+r_{2}}} \bigg] \\
\end{eqnarray*}
for $j=2,3,4$, $\widetilde{B}_{2,(3,j)}(r_{1},r_{2})$  are defined similarly and 
\begin{eqnarray*}
\widetilde{B}_{2,(2,1)}(r_{1},r_{2}) &=& \frac{|\mathcal{D}_{2,1}|\lambda^{2}}{n^{4}\Omega^{2}}\sum_{k_{1},k_{2}=-a}^{a}
g(\omega_{\Omega,k_{1}})\overline{g(\omega_{\Omega,k_{2}})}\times \\
 &&\Ex\bigg[\kappa_{4}(s_{j_2}-s_{j_1},0,s_{j_2}-s_{j_1})e^{is_{j_1}\omega_{\Omega,k_{1}}}e^{-is_{j_2}\omega_{\Omega,k_{1}+r_{1}}}
e^{-is_{j_1}\omega_{\Omega,k_{2}}}e^{is_{j_2}\omega_{\Omega,k_{2}+r_{2}}} \bigg].
\end{eqnarray*}
Integrating out the locations in $\widetilde{B}_{2,(3,1)}(r_{1},r_{2}) $ and
making a change of variables we have 
\begin{eqnarray*}
&& \widetilde{B}_{2,(3,1)}(r_{1},r_{2}) \\
&=& \frac{C\lambda^{2}|\mathcal{D}_{3,1}|}{n^{4}\Omega^{2}(2\pi)^{3}}\sum_{k_{1},k_{2}=-a}^{a}
g(\omega_{\Omega,k_{1}})\overline{g(\omega_{\Omega,k_{2}})} \int_{\mathbb{R}^{3}}
f_{4}(\omega_{1},\omega_{2},\omega_{3})
\sinc\left(\frac{\lambda}{2}\left(\omega_{1}+\omega_{3}+\omega_{\Omega,k_{2}-k_{1}}\right)\right)\times \\
&&\sinc\left(\frac{\lambda}{2} \left(\omega_{1} - \omega_{\Omega,k_{1}+r_{1}}\right)\right)
\times \sinc\left(\frac{\lambda}{2}\left( \omega_{3}+\omega_{\Omega,k_{2}+r_{2}}\right)\right)
d\omega_{1}d\omega_{2}d\omega_{3} \\
&=& \frac{C\lambda^{2}|\mathcal{D}_{3,1}|}{n^{4}\Omega^{2}(2\pi)^{3}}\sum_{k_{1},k_{2}=-a}^{a}
g(\omega_{\Omega,k_{1}})\overline{g(\omega_{\Omega,k_{2}})} \int_{\mathbb{R}^{3}}
f_{4}\left(\frac{2u_{1}}{\lambda}+\omega_{\Omega,k_1+r_1},\omega_{2},\frac{2u_{3}}{\lambda}-\omega_{\Omega,k_{2}+r_{2}}\right)\\
&&\times 
\sinc\left(u_{1}+u_{3} + \frac{\lambda}{\Omega}(r_{1}-r_{2})\right)\sinc(u_{1})\sinc(u_{3})
du_{1}d\omega_{2}du_{3}.
\end{eqnarray*}
This gives $\widetilde{B}_{2,(3,1)}(r_{1},r_{2}) = O(1/n)$. Next we consider
$\widetilde{B}_{2,(2,1)}(r_{1},r_{2})$, again by integrating out the locations and
changing variables we
have 
\begin{eqnarray*}
\widetilde{B}_{2,(2,1)}(r_{1},r_{2})
&=&  \frac{\lambda^{2}}{(2\pi)^{3}n^{2}\Omega^{2}}\sum_{k_{1},k_{2}=-a}^{a}
g(\omega_{\Omega,k_{1}})\overline{g(\omega_{\Omega,k_{2}})}\int_{\mathbb{R}^{3}}
f_{4}(\omega_{1},\omega_{2},\omega_{3})\sinc\left(\frac{\lambda}{2}\left(\omega_{1}+\omega_{2}
    + \omega_{\Omega,k_{2}-k_{1}}\right)\right)\times \\
&&\sinc\left(\frac{\lambda}{2}\left(\omega_{1}+\omega_{2} + \omega_{\Omega,k_{2}-k_{1}+r_{2}-r_{1}}\right)\right)
d\omega_{1}d\omega_{2} d\omega_{3}\\
&=&  \frac{2\lambda^{}}{(2\pi)^{3}n^{2}\Omega^{2}}\sum_{k_{1},k_{2}=-a}^{a}
g(\omega_{\Omega,k_{1}})\overline{g(\omega_{\Omega,k_{2}})}\int_{\mathbb{R}^{3}}
f_{4}\left(\frac{2u_{1}}{\lambda}-u_{2}-\omega_{\Omega,k_2-k_1},u_{2},u_{3}\right)
\sinc\left(u_{1}\right)\times \\
&&\sinc\left(u_{1}+\frac{\lambda}{\Omega}(r_{2}-r_{1})\pi)\right)
du_{1}du_{2} du_{3}.
\end{eqnarray*}
Thus $\widetilde{B}_{2,(2,1)}(r_{1},r_{2})
=O(\lambda^{}a/(n^{2}\Omega))$. Note this bound is the same
regardless of the sampling scheme on frequencies used, thus
using (\ref{eq:B2bound}) we have
\begin{eqnarray*}
|\widetilde{B}_{2}(r_1,r_2)| = O\left( \frac{1}{n} + \frac{\lambda
    a}{n^{2}\Omega} \right). 
\end{eqnarray*}
Thus (\ref{eq:tildeB3}) and the above prove (i) and (ii). 
\hfill $\Box$
\vspace{3mm}

We see from the above lemma that
$B_{1}(r_1,r_2)$ and
$B_3(r_{1},r_{2})$  are the leading higher order cumulant terms.
We now obtain some approximations for these rather complex terms. 
Let 
\begin{eqnarray}
&&D_{1}\left(\frac{a}{\Omega};\rb_1,\rb_2\right) \nonumber\\
&=& \frac{1}{(2\pi)^{2d}}
\int_{[-2\pi a/\Omega,2\pi a/\Omega]^{2d}}g(\ob_{1})\overline{g(\ob_{2})}
f_{4}\left(-\ob_{1}-\ob_{\Omega,\rb_{1}},-\ob_{2},\ob_{2}+
\ob_{\Omega,\rb_2}\right)\sinc\left(\frac{\lambda}{\Omega}(\rb_{2}-\rb_{1})\pi\right)\nonumber\\
&&D_{2}\left(\frac{a}{\Omega};\rb_1,\rb_2\right) \nonumber\\
&=& \frac{c_{4}}{(2\pi)^{2d}}
\int_{[-2\pi a/\Omega,2\pi a/\Omega]^{2d}}g(\ob_{1})g(\ob_{2})
f_{4}\left(-\ob_{1}-\ob_{\Omega,\rb_{1}},\ob_{2},-\ob_{2}-
\ob_{\Omega,r_2}\right)\sinc\left(\frac{\lambda}{\Omega}(\rb_{2}+\rb_{1})\pi\right). \nonumber\\
 && \label{eq:D1D2def}
\end{eqnarray}
\begin{theorem}
Let us suppose that  $\{Z(\ub);\ub \in\mathbb{R}^{d}\}$ is a fourth order stationary spatial 
random field that satisfies Assumption \ref{assum:S}(i),  
\ref{assum:uniform}, \ref{assum:G}, \ref{assum:nonGaussian} and
\ref{assum:GG}(a,c) or \ref{assum:GG}(b,c)  are satisfied. 
Let $D_{1}(\cdot)$ and $D_{2}(\cdot)$ be defined as in (\ref{eq:D1D2def}).
Then for $\|\rb_1\|_1,\|\rb_2\|_1\leq a$ 
\begin{eqnarray}
B_{1}(\rb_{1},\rb_{2}) 
 &=&   D_{1}\left(\frac{a}{\Omega};\rb_1,\rb_2\right)+ 
O\left( \frac{1}{\Omega} + 
\frac{\log^{2}(\lambda)}{\lambda} +
\sum_{i=1}^{d}\frac{I(\frac{\lambda}{\Omega}(r_{1i}-r_{2i})>e)\log(\frac{\lambda}{\Omega}|r_{i1}-r_{i2}|)}{\lambda}
+ \frac{1}{n}\right)\qquad \label{eq:B1r1r2approx}
\end{eqnarray}
and 
\begin{eqnarray}
B_{3}(\rb_{1},\rb_{2}) 
 &=&   D_{2}\left(\frac{a}{\Omega};\rb_1,\rb_2\right)+
O\left( \frac{1}{\Omega} + 
\frac{\log^{2}(\lambda)}{\lambda}+\sum_{i=1}^{d}\frac{I(\frac{\lambda}{\Omega}(r_{1i}-r_{2i})>e)\log(\frac{\lambda}{\Omega}|r_{i1}-r_{i2}|)}{\lambda}\right).\qquad \label{eq:B1r1r2approx2}
\end{eqnarray}
\end{theorem}
PROOF. We prove (\ref{eq:B1r1r2approx}) for $d=1$, the same proof applies to 
(\ref{eq:B1r1r2approx2}) and $d>1$. We define the series of approximations 
\begin{eqnarray*}
B_{1}^{(1)}(r_{1},r_{2}) 
 &=&   \frac{c_{4}}{\pi^{3}(2\pi)^{3}}\int_{[-2\pi a/\Omega,2\pi a/\Omega]^{2}}g(\omega_{1})\overline{g(\omega_{2})}
\int_{\mathbb{R}^{3}}
f_{4}\left(\frac{2u_{1}}{\lambda}-\omega_{1}-\omega_{\Omega,r_{1}},\frac{2u_{2}}{\lambda}-\omega_{2},\frac{2u_{3}}{\lambda}+
\omega_{2}+\omega_{\Omega,r_2}\right)
 \times \\
&&\times\sinc\left(u_{1}+u_{2}+u_{3}+\frac{\lambda}{\Omega}(r_{2}-r_{1})\pi\right)\sinc(u_{1})
\sinc(u_{2})\sinc(u_{3})du_{1}du_{2}du_{3},
\end{eqnarray*}
\begin{eqnarray*}
B_{1}^{(2)}(r_{1},r_{2}) 
 &=&  \frac{c_{4}}{(2\pi)^{2}}\int_{[-2\pi a/\Omega,2\pi a/\Omega]^{2}}
g(\omega_{1})\overline{g(\omega_{2})}
f_{4}\left(-\omega_{1}-\omega_{\Omega,r_{1}},-\omega_{2},\omega_{2}+
\omega_{\Omega,r_2}\right)\sinc\left(\frac{\lambda}{\Omega}(r_{2}-r_{1})\pi\right). 
\end{eqnarray*}
By taking differences (see the proof of Theorem \ref{theorem:variance-nongaussian}) we have 
\begin{eqnarray*}
\left|B_{1}(r_{1},r_{2}) - B_{1}^{(1)}(r_{1},r_{2}) \right| = 
O\left(\frac{\ell_{3}(\frac{\lambda}{\Omega}(r_{1}-r_{2}))}{\Omega}\right)
\end{eqnarray*}
and 
\begin{eqnarray*}
\left|B_{1}^{(1)}(r_{1},r_{2}) - B_{1}^{(2)}(r_{1},r_{2}) \right| = 
O\left(\frac{\log^{2}(\lambda)}{\lambda}\right),
\end{eqnarray*}
thus giving the required result. \hfill $\Box$

\vspace{3mm}
\noindent
The above result implies that in the case the process is non-Gaussian 
and $\frac{\lambda}{\Omega}<1$  then 
\begin{eqnarray*}
&&\lambda^{d}\var[\widetilde{Q}_{a,\Omega,\lambda}(g;\rb)] \\
&=&  C_{11}\left(\frac{a}{\lambda},\ob_{\Omega,\rb}\right)
\frac{\lambda^{d}}{\Omega^{d}}\sum_{\mb=-2a}^{2a}
\Sinc^{2}\left( \frac{\lambda}{\Omega}\mb\pi\right) +C_{12}\left(\frac{a}{\lambda},\ob_{\Omega,\rb}\right)
\frac{\lambda^{d}}{\Omega^{d}}\sum_{m=-2a}^{2a}\Sinc^{2}\left( \frac{\lambda}{\Omega}(\mb+\rb)
  \pi\right) \\
&& D_{1}\left(\frac{a}{\lambda};\rb_1,\rb_2\right) + 
O\left(\widetilde{\ell}_{a,\Omega,\lambda}^{(2)}+ \sum_{i=1}^{d}\frac{I(\frac{\lambda}{\Omega}(r_{1i}-r_{2i})>e)\log(\frac{\lambda}{\Omega}|r_{i1}-r_{i2}|)}{\lambda}\right).
\end{eqnarray*}
On the other hand if a coarse frequency grid is used with $\lambda
\geq \Omega$ then the rate of convergence is worse with  
\begin{eqnarray*}
&&\Omega^{d}\var[\widetilde{Q}_{a,\Omega,\lambda}(g;\rb)] \\
&=&  C_{11}\left(\frac{a}{\lambda},\ob_{\Omega,\rb}\right)
\sum_{\mb=-2a}^{2a}
\Sinc^{2}\left( \frac{\lambda}{\Omega}\mb\pi\right) +C_{12}\left(\frac{a}{\lambda},\ob_{\Omega,\rb}\right)
\sum_{m=-2a}^{2a}\Sinc^{2}\left( \frac{\lambda}{\Omega}(\mb+\rb)
  \pi\right) \\
&&  +\frac{\Omega^{d}}{\lambda^{d}}D_{1}\left(\frac{a}{\lambda};\rb_1,\rb_2\right)+
O\left(\widetilde{\ell}_{a,\Omega,\lambda}^{(2)}+\sum_{i=1}^{d}\frac{I(\frac{\lambda}{\Omega}(r_{1i}-r_{2i})>e)\log(\frac{\lambda}{\Omega}|r_{i1}-r_{i2}|)}{\lambda}\right).
\end{eqnarray*}

%We observe that if the sampling scheme is such that
%$\lambda/\Omega\leq 1$, then
%$\var[\widetilde{Q}_{a,\Omega,\lambda}(g;r)] = O(\lambda^{-1})$ and
%the $C(\frac{a}{\lambda},\omega_{\Omega,r})$ and
%$D(\frac{a}{\lambda})$ are of the same order. 
%However, just as in the Gaussian case for 
%the coarser sampling scheme the rate decreases when
%$\lambda/\Omega\geq 1$ and the variance depends on the sampling scheme with
%$\var[\widetilde{Q}_{a,\Omega,\lambda}(g;r)] =
%O(\Omega^{-1})$. However, we observe that when 
%$\Omega/\lambda \rightarrow \infty$ then the term involving the
% fourth order cumulant, $D_{1}(\frac{a}{\lambda};\omega_{\Omega,r})$, becomes asymptotically negligible  
%and is dominated by $C(\frac{a}{\lambda};\omega_{\Omega,r})$ which
%simply a function the spectral density function. 

%% file: B_6_fixed_domainv2.tex
\section{Fixed domain asymptotics}\label{sec:appendixfixed}

In this section our aim is to investigate the sampling properties of 
$Q_{a,\lambda,\Omega}(g;0)$ in the case that the domain, $\lambda$, over which is
the spatial process is defined in kept fixed but the number of locations 
that are sampled grows ($n\rightarrow \infty$). 

In the following theorem we evaluate an expression for the
covariance between the Fourier transforms when the domain is fixed.

\begin{theorem}[Fixed domain]\label{theorem:fixedDFT1}
%Suppose the domain is kept fixed, for convenience we set $\lambda
%=1$, and the locations $\{\sbb_{j}^{(n)};j=1,\ldots,n\}$ satisfy
%Assumption \ref{assum:fixedLocations}. 
Suppose that $\{Z(\sbb)\}$ is a second order stationary. Then 
\begin{itemize} 
\item[(i)] Under Assumption \ref{assum:nonuniform} (general random design of locations) we have 
\begin{eqnarray*}
&&\cov\left[J_{n}\left(\frac{2\pi\kb_1}{\Omega}\right),J_{n}\left(\frac{2\pi\kb_2}{\Omega}\right)
\right] = 
\left(\frac{\lambda}{2\pi}\right)^{d}\sum_{\jb_{1},\jb_{2}\in \mathbb{Z}^{d}}\gamma_{\jb_1}\gamma_{\jb_2}
\int_{\mathbb{R}^{d}}f(\ob)\Sinc\left(\frac{\lambda}{2}\left[\ob + \frac{2\pi \jb_{1}}{\lambda}
+\frac{2\pi \kb_{1}}{\Omega}\right]\right) \\
&& \times\Sinc\left(\frac{\lambda}{2}
\left[\ob - \frac{2\pi \jb_{2}}{\lambda} + \frac{2\pi
    \kb_{2}}{\Omega}\right]\right)d\ob 
+ O\left(\frac{\lambda^{d}}{n}\right) 
\end{eqnarray*}
\item[(ii)] Under Assumptions \ref{assum:fixedDFT}  and \ref{assum:GG}(e) we have
\begin{eqnarray*}
&&\cov\left[J_{n}\left(\frac{2\pi\kb_1}{\Omega}\right),J_{n}\left(\frac{2\pi\kb_2}{\Omega}\right)
\right] \\
&=&\left(\frac{2}{\lambda\pi}\right)^{d}
\int_{\mathbb{R}^{d}}f(\ob)\Sinc\left(\frac{\lambda}{2}\left[\ob+\frac{2\pi \kb_{1}}{\Omega}\right]\right)
\Sinc\left(\frac{\lambda}{2}
\left[\ob + \frac{2\pi
    \kb_{2}}{\Omega}\right]\right)d\ob +\\
&& O\left(\frac{\lambda}{n^{1/d}}\left(\frac{\|\kb_1\|_1 +
      \|\kb_2\|_{1}}{\Omega} + 1\right)\right)
\end{eqnarray*}
\end{itemize}
%where 
%\begin{eqnarray*}
%\frac{1}{(2\pi)^{d}}\int_{\mathbb{R}^{d}}f(\ob)\Sinc\left(\frac{\lambda}{2}\ob\right)\Sinc\left(\frac{1}{2}
%\left[\ob + \frac{2\pi \kb_{1}}{\Omega}\right]\right)
%\Sinc\left(\frac{1}{2}
%\left[\ob + \frac{2\pi \kb_{2}}{\Omega}\right]\right)
%d\ob +
%O\left( \Omega^{d}\prod_{\ell=1}^{d}\gamma(k_{\ell,1},k_{\ell,2})\right)
%\end{eqnarray*}
%and 
%\begin{eqnarray*}
%&&\gamma(k_{\ell,1},k_{\ell,2}) = 
%\left(\frac{I_{k_{\ell,1}\neq 0}}{|k_{\ell,1}|} + I_{k_{\ell,1}= 0}\right) \left(\frac{I_{k_{\ell,2}\neq 0}}{|k_{\ell,1}|} + I_{k_{\ell,2}= 0}\right).
%\end{eqnarray*}
\end{theorem}
{\bf PROOF}  The proof of (i) is straightforward and follows from the techniques
used to prove (ii), below.

We use $d=1$ and $\lambda = \Omega =1$ to  prove (ii).
Expanding $\cov\left(
J_{n}(\frac{2\pi k_1}{\Omega}),J_{n}(\frac{2\pi k_2}{\Omega})\right)$ gives 
\begin{eqnarray}
&&\cov\left[ J_{n}\left(\frac{2\pi k_1}{\Omega}\right),J\left(\frac{2\pi k_2}{\Omega}\right)\right]  =
\frac{\lambda}{n^{2}}\sum_{j_1,j_2=1}^{n}
c\left(s_{j_1}-s_{j_1}\right)\exp\left(\frac{2\pi i}{\Omega}\left[s_{j_1}k_{1}-s_{j_2}k_{2}\right]\right)
\nonumber\\
&=&
\frac{1}{2\pi \lambda}\int_{-\infty}^{\infty}f(\omega)\left[\frac{\lambda}{n}\sum_{j_{1}=1}^{n}\exp\left(is_{j_1}\left[\omega
  + \frac{2\pi k_{1}}{\Omega}\right]\right) \right]\left[\frac{\lambda}{n}\sum_{j_{2}=1}^{n}\exp\left(-is_{j_2}\left[\omega
  + \frac{2\pi k_{2}}{\Omega}\right]\right) \right]d\omega.\qquad \label{eq:JJf}.
\end{eqnarray} 
To simplify notation, in the following we let $s_{(j)} = s_{n,(j)}$.
Replacing the inner summand with an integral and using the mean value
theorem (on the integral) we have
%We note that under (\ref{eq:sj-sj-1})
\begin{eqnarray*}
&&\frac{\lambda}{n}\sum_{j_{}=1}^{n}\exp\left(is_{j}\left[\omega
  + \frac{2\pi k_{}}{\Omega}\right]\right) -
\int_{-\lambda/2}^{\lambda/2}\exp\left(is\left[\frac{2\pi k}{\Omega}+\omega\right]\right)ds \\
&=& \frac{\lambda\exp(is_{(n)}[
   \frac{2\pi k_{}}{\Omega}+\omega])}{n}+\sum_{j_{}=1}^{n-1}\left(\frac{\lambda}{n}\exp\left(is_{(j)}\left[
   \frac{2\pi k_{}}{\Omega}+\omega\right]\right) - \int_{-s_{(j)}}^{s_{(j+1)}}\exp\left( is\left[\frac{2\pi
      k}{\Omega}+\omega\right]\right)ds\right) \\
&=& \frac{\lambda\exp(is_{(n)}[
  \frac{2\pi k_{}}{\Omega} + \omega])}{n}+\sum_{j_{}=1}^{n-1}\left\{\frac{\lambda}{n}\exp\left(is_{(j)}\left[
   \frac{2\pi k_{}}{\Omega}+\omega\right]\right) - (s_{(j+1)}-s_{(j)})\exp\left( i\tilde{s}_{j}\left[\frac{2\pi
      k}{\Omega}+\omega\right]\right)\right\} \\
 &=& \frac{\lambda\exp(is_{(n)}[
  2\pi k_{}+\omega])}{n}+\sum_{j_{}=1}^{n-1}(s_{(j+1)}-s_{(j)})\left\{\exp\left(is_{(j)}\left[
  \frac{2\pi k_{}}{\Omega}+\omega\right]\right) - \exp\left( i\tilde{s}_{(j)}\left[\frac{2\pi k}{\Omega}+\omega\right]\right)\right\} \\
&& + \sum_{j_{}=1}^{n-1}\left(\frac{\lambda}{n} - (s_{(j+1)}-s_{(j)}) \right)\exp\left(is_{(j)}\left[
   \frac{2\pi k_{}}{\Omega}+\omega\right]\right). 
\end{eqnarray*}
By using the mean value theorem again we have 
\begin{eqnarray*}
&&\frac{\lambda}{n}\sum_{j_{}=1}^{n}\exp\left(is_{j}\left[\omega
  + \frac{2\pi k_{}}{\Omega}\right]\right) -
\int_{-\lambda/2}^{\lambda/2}\exp\left(is\left[\frac{2\pi k}{\Omega}+\omega\right]\right)ds \\
&=& \frac{\lambda\exp(is_{(n)}[\omega
  + 2\pi k_{}])}{n}+\sum_{j_{}=1}^{n-1}(s_{(j+1)}-s_{(j)}) (s_{(j)}-\tilde{s}_{(j)})\left[\omega
  + \frac{2\pi k_{}}{\Omega}\right]
\exp\left(i\hat{s}_{(j)}\left[\omega
  + \frac{2\pi k_{}}{\Omega}\right]\right)\\
&& + \sum_{j_{}=1}^{n-1}\left(\frac{\lambda}{n} - (s_{(j+1)}-s_{(j)}) \right)\exp\left(is_{(j)}\left[\omega
  + \frac{2\pi k_{}}{\Omega}\right]\right),
% + O\left(\frac{k+\omega}{n} \right) \\
\end{eqnarray*}
where $\widetilde{s}_{(j)},\widehat{s}_{(j)}\in
[s_{(j)},s_{(j+1)}]$. 
Thus we have 
\begin{eqnarray*}
&&\left|\frac{\lambda}{n}\sum_{j_{}=1}^{n}\exp\left(is_{j}\left[\omega
  + \frac{2\pi k_{}}{\Omega}\right]\right) -
\int_{-\lambda/2}^{\lambda/2}\exp\left(is\left[\frac{2\pi k}{\Omega}+\omega\right]\right)ds\right| \\
&\leq & \frac{\lambda}{n}+\left|\omega
  + \frac{2\pi k_{}}{\Omega}\right|\sum_{j_{}=1}^{n-1}(s_{(j+1)}-s_{(j)})^{2}
 + \sum_{j_{}=1}^{n-1}\left|\frac{\lambda}{n} - (s_{(j+1)}-s_{(j)})
 \right| = O\left(\frac{\lambda}{n}\left[\omega
  + \frac{2\pi k_{}}{\Omega}\right]\right)
% + O\left(\frac{k+\omega}{n} \right) \\
\end{eqnarray*}
By using (\ref{eq:sj-sj-1}) and
the mean value theorem (again) we have 
\begin{eqnarray*}
&&\frac{\lambda}{n}\sum_{j_{}=1}^{n}\exp(is_{j}[\omega
  + 2\pi k_{}] = \int_{-\lambda/2}^{\lambda/2}\exp\left(2\pi is\left[2\pi
      k+\omega\right]\right)ds + O\left(\frac{\lambda|\frac{2\pi k}{\Omega}+\omega|}{n}\right).
\end{eqnarray*}
Systematically replacing the summand with the integral in
(\ref{eq:JJf}), twice, and using that
$\int_{\mathbb{R}}f(\omega)|\omega|d\omega <\infty$
\begin{eqnarray*}
&&\cov\left[ J_{n}\left(\frac{2\pi k_1}{\Omega}\right),J\left(\frac{2\pi k_2}{\Omega}\right)\right]  =
\frac{\lambda}{n^{2}}\sum_{j_1,j_2=1}^{n}
c\left(s_{j_1}-s_{j_1}\right)\exp\left(\frac{2\pi i}{\Omega}\left[s_{j_1}k_{1}-s_{j_2}k_{2}\right]\right)
\nonumber\\
&=&
\frac{1}{2\pi \lambda}\int_{-\infty}^{\infty}f(\omega)\left[\frac{\lambda}{n}\sum_{j_{1}=1}^{n}\exp\left(is_{j_1}\left[\omega
  + \frac{2\pi k_{1}}{\Omega}\right]\right) \right]\left[\frac{\lambda}{n}\sum_{j_{2}=1}^{n}\exp\left(-is_{j_2}\left[\omega
  + \frac{2\pi k_{2}}{\Omega}\right]\right) \right]d\omega \\
&=& \frac{1}{2\pi
  \lambda}\int_{-\infty}^{\infty}f(\omega)\int_{-\lambda/2}^{\lambda/2}\int_{-\lambda/2}^{\lambda/2}
\exp\left( is_1\left[\frac{2\pi
      k_1}{\Omega}+\omega\right]\right)\exp\left(- is_2\left[\frac{2\pi
      k_2}{\Omega}+\omega\right]\right)  ds_1ds_2 +\\
&& O\left(\frac{\lambda}{n}\left(\frac{|k_1|+|k_2|}{\Omega}+1\right)\right) \\
&=&
\frac{2}{\lambda\pi}\int_{-\infty}^{\infty}f(\omega)\sinc\left(\frac{\lambda}{2}\left[\omega+\frac{2\pi
  k_1}{\Omega}\right]\right)\sinc\left(\frac{\lambda}{2}\left[\omega+\frac{2\pi
  k_2}{\Omega}\right]\right)d\omega + O\left(\frac{\lambda}{n}\left(\frac{|k_1|+|k_2|}{\Omega}+1\right)\right),
\end{eqnarray*} 
thus proving (ii).
%The final part of the proof involves showing that 
%\begin{eqnarray*}
%&&\frac{2}{\lambda\pi}\int_{-\infty}^{\infty}f(\omega)\sinc\left(\frac{\lambda}{2}\left[\omega+\frac{2\pi
%  k_1}{\Omega}\right]\right)\sinc\left(\frac{\lambda}{2}\left[\omega+\frac{2\pi
%  k_2}{\Omega}\right]\right)d\omega \\
%&=& \left\{
%\begin{array}{cc}
%A_{\lambda}\left(\frac{k}{\Omega}\right) & k_1 = k_2 (=k) \\
%\frac{(-1)^{k_1 - k_2 +1}}{\pi(k_1 -
 % k_2)}\left(B_{\lambda}(\frac{k_{1}}{\Omega}) - B_{\lambda}(\frac{k_{2}}{\Omega}) \right) & k_1\leq
%k_2 \\
%\end{array}
%\right.
%\end{eqnarray*}
%which we do so below. 
\hfill $\Box$

\vspace{3mm}
\noindent
The results below are used to prove Section \ref{sec:fixed-asy} 
we focus on the case $d=1$.  
We first consider the Fourier transform of the continuous  time
anologue $\{Z(s);s\in [-\lambda/2,\lambda/2]\}$. It is straightforward
to show 
\begin{eqnarray}
\cov\left[\mathcal{J}_{\lambda}\left(\frac{k_1}{\Omega}\right),\mathcal{J}_{\lambda}\left(\frac{k_2}{\Omega}\right)\right] &=& 
\frac{2}{\lambda\pi}\int_{-\infty}^{\infty}f(\omega)\sinc\left(\frac{\lambda}{2}\left[\omega+\frac{2\pi
  k_1}{\Omega}\right]\right)\sinc\left(\frac{\lambda}{2}\left[\omega+\frac{2\pi
  k_2}{\Omega}\right]\right)d\omega \nonumber\\
&=&  \frac{1}{\lambda}\int_{-\lambda/2}^{\lambda/2}\int_{-\lambda/2}^{\lambda/2}c(s_1-s_2)
\exp\left(\frac{2\pi is_1 k_1}{\Omega}\right)\exp\left(-
\frac{2\pi is_2 k_2}{\Omega}\right)  ds_1ds_2 \qquad \label{eq:Jcont}
\end{eqnarray}
where 
\begin{eqnarray} 
\label{eq:mathcalJ}
\mathcal{J}_{\lambda}\left(\frac{k}{\Omega}\right) = \frac{1}{\lambda^{1/2}}\int_{-\lambda/2}^{\lambda/2}Z(s)\exp\left(\frac{2\pi i k
  s}{\Omega}\right)ds.
\end{eqnarray}
Thus comparing the above with Theorem \ref{theorem:fixedDFT1} we have 
\begin{eqnarray*}
&&\cov\left[J_{n}\left(\frac{2\pi\kb_1}{\Omega}\right),J_{n}\left(\frac{2\pi\kb_2}{\Omega}\right)
\right] \\
&=&  \frac{1}{\lambda}\int_{-\lambda/2}^{\lambda/2}\int_{-\lambda/2}^{\lambda/2}c(s_1-s_2)
\exp\left(\frac{2\pi is_1 k_1}{\Omega}\right)\exp\left(-
\frac{2\pi is_2 k_2}{\Omega}\right)  ds_1ds_2
+ o\left(1\right),  
\end{eqnarray*}
where the error in the approximation depends on the sampling scheme. 
We now obtain an explicit expression for the above, which helps us
quantify the dependence.

\begin{theorem}\label{theorem:dftfixed2}
Suppose $\{Z(u)\}$ is a spatial second order stationary process and
$\mathcal{J}_{\lambda}\left(\frac{k}{\Omega}\right)$ is defined in 
 (\ref{eq:mathcalJ}).
Then 
\begin{eqnarray*}
\cov\left[\mathcal{J}_{\lambda}\left(\frac{k_1}{\Omega}\right),\mathcal{J}_{\lambda}\left(\frac{k_2}{\Omega}\right) \right]
= 
\left\{
\begin{array}{cc}
A_{\lambda}\left(\frac{k}{\Omega}\right) & k_1 = k_2 (=k) \\
B_{\lambda}\left(\frac{k_1}{\Omega},\frac{k_2}{\Omega}\right)
& k_1\leq
k_2 \\
\end{array}
\right.
\end{eqnarray*}
where 
\begin{eqnarray*}
A_{\lambda}\left(\frac{k}{\Omega}\right) = 
\int_{-\lambda}^{\lambda} T\left(\frac{u}{\lambda}\right)c(u)\exp\left(\frac{2i\pi ku}{\Omega}\right)du
\end{eqnarray*}
and
\begin{eqnarray*}
&& B_{\lambda}\left(\frac{k_1}{\Omega},\frac{k_2}{\Omega}\right) = 
\frac{\Omega}{(k_1-k_2)\pi \lambda}
\sin\left(\frac{2\pi \lambda(k_1-k_2)}{\Omega}\right)
\int_{-\lambda}^{\lambda}
c(v)\exp\left(\frac{2\pi i k_1 v}{\Omega}\right)dv \\
&& +\frac{\Omega}{(k_1-k_2)\pi\lambda} \Im\left[
e^{-\lambda i\pi (k_1-k_2)/\Omega}\left(\int^{\lambda}_{0}c(v) e^{2\pi i k_2 v/\Omega}dv
 - \int^{\lambda}_{0}c(v) e^{2\pi i k_1 v/\Omega}dv
\right)\right].
\end{eqnarray*}
\end{theorem}
{\bf PROOF} For the case $k_1=k_2 (=k)$ and (\ref{eq:Jcont}) we have 
\begin{eqnarray*}
\var\left[ \mathcal{J}_{\lambda}\left(\frac{k}{\Omega}\right) \right]
&=& 
\frac{1}{\lambda}\int_{-\lambda/2}^{\lambda/2}\int_{-\lambda/2}^{\lambda/2}c(s_1-s_2)
\exp\left(\frac{2\pi i(s_1 - s_2)k}{\Omega}\right)ds_1ds_2 \\
&=& \frac{1}{\lambda}\int_{-\lambda/2}^{\lambda/2}\left(\lambda - |u|\right)c(u)
\exp\left(\frac{2\pi i u k}{\Omega}\right)du 
= 
\int_{-\lambda/2}^{\lambda/2}\left(1 - \frac{|u|}{\lambda}\right)c(u)
\exp\left(\frac{2\pi i u k}{\Omega}\right)du.
\end{eqnarray*}
This gives $A_{\lambda}(\frac{k}{\Omega})$. For $k_1\neq k_2$ we have 
\begin{eqnarray*}
&& \cov\left[\mathcal{J}_{\lambda}\left(\frac{k_1}{\Omega}\right)^{},\mathcal{J}_{\lambda}\left(\frac{k_2}{\Omega}\right)^{} \right]\\
&=& 
\frac{1}{\lambda}\int_{-\lambda/2}^{\lambda/2}\int_{-\lambda/2}^{\lambda/2}
c(u_1 - u_2)\exp\left(\frac{2\pi i
    k_1(u_1-u_2)}{\Omega}\right)\exp\left(\frac{2\pi i
  u_2(k_1-k_2)}{\Omega}\right)du_1du_2 \\
&=& \frac{1}{\lambda}\int_{-\lambda}^{\lambda}
c(v)\exp\left(\frac{2\pi i
    k_1 v}{\Omega}\right) 
\int_{\max(-\lambda/2,-\lambda/2-v)}^{\min(\lambda/2,\lambda/2-v)} 
\exp\left(\frac{2\pi i
  u_2(k_1-k_2)}{\Omega}\right) du_2 dv \\
&=&\frac{1}{\lambda}\int_{-\lambda}^{\lambda}
c(v)\exp\left(\frac{2\pi i
    k_1 v}{\Omega}\right) 
\int_{-\lambda/2}^{\lambda/2} 
\exp\left(\frac{2\pi i
  u_2(k_1-k_2)}{\Omega}\right) du_2 dv \\
&& - \frac{1}{\lambda}\int_{-\lambda}^{0}
c(v)\exp\left(\frac{2\pi i
    k_1 v}{\Omega}\right) 
\int_{-\lambda/2}^{-\lambda/2-v} 
\exp\left(\frac{2\pi i
  u_2(k_1-k_2)}{\Omega}\right) du_2 dv \\
&& - \frac{1}{\lambda}\int^{\lambda}_{0}
c(v)\exp\left(\frac{2\pi i
    k_1 v}{\Omega}\right) 
\int_{\lambda/2-v}^{\lambda/2} 
\exp\left(\frac{2\pi i
  u_2(k_1-k_2)}{\Omega}\right) du_2 dv \\
&=& \frac{\Omega}{\lambda(k_1-k_2)}
\sin\left(\frac{2\pi \lambda(k_1-k_2)}{\Omega}\right)
\int_{-\lambda}^{\lambda}
c(v)\exp\left(\frac{2\pi i k_1 v}{\Omega}\right)dv +
I_{1} + I_{2}.
\end{eqnarray*}
We obtain expressions for $I_{1}$ and $I_{2}$
\begin{eqnarray*}
I_{1}  &=& -\frac{1}{\lambda}\int_{-\lambda}^{0}c(v)\exp\left(\frac{2\pi i
    k_1 v}{\Omega}\right) \left[
\frac{e^{2\pi i u(k_1-k_2)/\Omega}}{2\pi i (k_1-k_2)/\Omega}
\right]_{u=-\lambda/2}^{-v-\lambda/2} \\
&=& \frac{-\Omega}{\lambda 2\pi i(k_1-k_2)}
\int_{-\lambda}^{0}c(v)e^{2\pi i k_1 v/\lambda}
\left(e^{-\frac{2\pi i
      (k_1-k_2)}{\Omega}(\frac{\lambda}{2}+v)} - e^{-i\pi \lambda (k_1-k_2)/\Omega}\right)\\
&=& \frac{-\Omega e^{-\lambda i\pi (k_1-k_2)/\Omega}}{\lambda 2\pi
  i(k_1-k_2)}\int_{-\lambda}^{0}c(v)
\left( e^{2\pi i k_2 v/\Omega} - e^{2\pi i k_1 v/\Omega}\right)dv.
\end{eqnarray*}
Similarly 
\begin{eqnarray*}
I_{2}  &=& -\frac{1}{\lambda}\int^{\lambda}_{0}c(v)\exp\left(\frac{2\pi i
    k_1 v}{\Omega}\right) \left[
\frac{e^{2\pi i u(k_1-k_2)/\Omega}}{2\pi i (k_1-k_2)/\Omega}
\right]_{u=\lambda/2-v}^{\lambda/2} \\
&=& \frac{-\Omega}{\lambda 2\pi i(k_1-k_2)}
\int^{\lambda}_{0}c(v)e^{2\pi i k_1 v/\Omega}
\left(e^{\lambda i\pi (k_1-k_2)/\Omega}-e^{\frac{2\pi i
      (k_1-k_2)}{\Omega}(\frac{\lambda}{2}-v)} \right)\\
&=& \frac{-\Omega e^{\lambda i\pi (k_1-k_2)/\Omega}}{\lambda 2\pi
  i(k_1-k_2)}\int^{\lambda}_{0}c(v)
\left( e^{2\pi i k_1 v/\Omega} - e^{2\pi i k_2 v/\Omega}\right)dv.
\end{eqnarray*}
This gives 
\begin{eqnarray*}
I_{1} + I_{2}
&=&
\frac{\Omega}{(k_1-k_2)\pi\lambda} \Im\left[
e^{-\lambda i\pi (k_1-k_2)/\Omega}\left(\int^{\lambda}_{0}c(v) e^{2\pi i k_2 v/\Omega}dv
 - \int^{\lambda}_{0}c(v) e^{2\pi i k_1 v/\Omega}dv
\right)\right].
\end{eqnarray*}
Thus 
\begin{eqnarray*}
&&\cov\left[\mathcal{J}_{\lambda}\left(\frac{k_1}{\Omega}\right)^{},
\mathcal{J}_{\lambda}\left(\frac{k_2}{\Omega}\right)^{} \right]
= \frac{\Omega}{\lambda(k_1-k_2)}
\sin\left(\frac{2\pi \lambda(k_1-k_2)}{\Omega}\right)
\int_{-\lambda}^{\lambda}
c(v)\exp\left(\frac{2\pi i k_1 v}{\Omega}\right)dv +\\
&& \frac{\Omega}{(k_1-k_2)\pi\lambda} \Im\left[
e^{-\lambda i\pi (k_1-k_2)/\Omega}\left(\int^{\lambda}_{0}c(v) e^{2\pi i k_2 v/\Omega}dv
 - \int^{\lambda}_{0}c(v) e^{2\pi i k_1 v/\Omega}dv
\right)\right].
\end{eqnarray*}
giving the required result. \hfill $\Box$

\vspace{3mm}
In the following lemma we focus on
$\mathcal{J}_{\lambda}\left(\frac{k}{\lambda}\right)$ ($\Omega = \lambda$).
\begin{theorem}
Suppose $\{Z(u)\}$ is a spatial second order stationary process and
$\mathcal{J}_{\lambda}\left(\frac{k}{\lambda}\right)$ is defined in 
 (\ref{eq:mathcalJ}).
Then 
\begin{eqnarray*}
\cov\left[\mathcal{J}_{\lambda}\left(\frac{k_1}{\lambda}\right),\mathcal{J}_{\lambda}\left(\frac{k_2}{\lambda}\right) \right]
= 
\left\{
\begin{array}{cc}
A_{\lambda}\left(\frac{k}{\lambda}\right) & k_1 = k_2 (=k) \\
\frac{(-1)^{k_1 - k_2 +1}}{\pi(k_1 -
  k_2)}\left[B_{\lambda}(\frac{k_{1}}{\lambda}) - B_{\lambda}(\frac{k_{2}}{\lambda}) \right] & k_1\leq
k_2 \\
\end{array}
\right.
\end{eqnarray*}
if in addition the process is Gaussian  
then
\begin{eqnarray*}
&&\cov\left[\mathcal{J}_{\lambda}\left(\frac{k_1}{\lambda}\right)^{2},\mathcal{J}_{\lambda}\left(\frac{k_2}{\lambda}\right)^{2} \right]\\
&=& 
\left\{
\begin{array}{cc}
A_{\lambda}\left(\frac{k}{\lambda}\right)^{2}  + 
\frac{\lambda^{2}}{\pi^{2}k^{2}\lambda^{2}}B_{\lambda}(\frac{k}{\lambda})^{2}
& k_1 = k_2 (=k) \\
\frac{\lambda^{2}}{\lambda^{2}\pi^{2}(k_1 -
  k_2)^{2}}\left[B_{\lambda}(\frac{k_{1}}{\lambda}) -
  B_{\lambda}(\frac{k_{2}}{\lambda}) \right]^{2}+
\frac{\lambda^{2}}{\lambda^{2}\pi^{2}(k_1 +
  k_2)^{2}}\left[B_{\lambda}(\frac{k_{1}}{\lambda}) +
  B_{\lambda}(\frac{k_{2}}{\lambda}) \right]^{2}
& k_1\neq k_2 \\ 
\end{array}
\right.
\end{eqnarray*}
where 
\begin{eqnarray*}
A_{\lambda}\left(\frac{k}{\lambda}\right) = 
\frac{1}{2^{d}}\int_{-\lambda}^{\lambda} T\left(\frac{u}{\lambda}\right)c(u)\exp\left(\frac{2i\pi ku}{\lambda}\right)du
\end{eqnarray*}
and
\begin{eqnarray*}
B_{\lambda}\left(\frac{k}{\lambda}\right)  = \frac{1}{\lambda}\int_{0}^{\lambda}c(u)\sin\left(\frac{2\pi
k u}{\lambda}\right)du,
\end{eqnarray*}
\end{theorem}
PROOF. The proof for $k_{1}\neq k_2$ follows from Theorem \ref{theorem:dftfixed2}
and using that 
$e^{i\pi (k_1-k_2)}=(-1)^{k_1-k_2}$ and $\sin\left(2\pi
  (k_1-k_2)\right)=0$, this gives
\begin{eqnarray*}
\cov\left[\mathcal{J}_{\lambda}\left(\frac{k_1}{\lambda}\right),\mathcal{J}_{\lambda}\left(\frac{k_2}{\lambda}\right) \right]
&=& \frac{(-1)^{(k_1-k_2-1)}}{\lambda 2\pi
  (k_1-k_2)}\left[
\int^{\lambda}_{0}c(v)\left( \sin\left(\frac{2\pi k_1}{\lambda}\right) 
-  \sin\left(\frac{2\pi k_2}{\lambda} \right)
\right)dv\right] \\
&=& \frac{(-1)^{k_1-k_2+1}}{\pi
  (k_1-k_2)}
\left[B_{\lambda}(\frac{k_1}{\lambda}) - B_{\lambda}(\frac{k_2}{\lambda}) \right].
\end{eqnarray*}
The result for
$\cov\left[\mathcal{J}_{\lambda}\left(\frac{k_1}{\lambda}\right)^{2},
\mathcal{J}_{\lambda}\left(\frac{k_2}{\lambda}\right)^{2} \right]$
simply follows from the above and the covariance expansion in terms of
cumulants.  \hfill $\Box$

%% file: B_4_Q.tex
\section{Sample properties of $Q_{a,\lambda}(g;\rb)$}\label{appendix:Q}

In this section we summarize the result for $Q_{a,\lambda}(g;\rb)$,
where $Q_{a,\lambda}(g;\rb) = \widetilde{Q}_{a,\lambda,\lambda}(g;\rb)$
(which is defined in (\ref{eq:QQQr})). To
simplify notation we state the results only for the case that the
locations are uniformly distributed.  

We recall that 
\begin{eqnarray}\label{eq:Qalambda}
Q_{a,\lambda}(g;\rb) = \widetilde{Q}_{a,\lambda}(g;\rb) + G_{\lambda}V_{\rb}
\end{eqnarray}
where
$V_{\rb}=\frac{1}{n}\sum_{j=1}^{n}Z(\ub_{j})^{2}\exp(-i\ob_{\rb}\ub_{j})$,
$G_{\lambda} = \frac{1}{n}\sum_{\kb=-a}^{a}g(\ob_{\Omega,\kb})$ and
$\widetilde{Q}_{a,\lambda}(g;\rb)$ is defined in (\ref{eq:tildeQQQ}).

\begin{theorem}\label{theorem:mean-Q}
Suppose Assumptions \ref{assum:S}(i), \ref{assum:uniform}, 
$b=b(\rb)$ denotes the number of zero elements in the vector $\rb\in
\mathbb{Z}^{d}$ and 
\begin{itemize}
\item[(i)] Assumptions \ref{assum:G}(i) and \ref{assum:GG}(a,c) hold. 
Then we have 
\begin{eqnarray*}
&&\Ex\left[Q_{a,\lambda}(g;\rb)\right]  \nonumber\\
&=&
\left\{
\begin{array}{cl}
O(\frac{1}{\lambda^{d-b}}) &   \rb\in \mathbf{Z}^{d}/\{0\} 
\\
\frac{1}{(2\pi)^{d}}\int_{\ob \in 2\pi[-C,C]^{d}}f(\ob)g(\ob)d\ob + O(\frac{1}{\lambda} + \frac{\lambda^{d}}{n}) & \rb = 0 \\
\end{array}
\right. \qquad
\end{eqnarray*}
\item[(ii)] Suppose the Assumptions \ref{assum:G}(ii) and Assumption \ref{assum:GG}(b,c) hold and 
$\{m_{1},\ldots,m_{d-b}\}$ is the subset of
non-zero values in $\rb =(r_{1},\ldots,r_{d})$, then we have  
\begin{eqnarray*}
&&\Ex\left[Q_{a,\lambda}(g;\rb)\right] \nonumber\\
&=&\left\{ 
\begin{array}{cl}
 O\left( \frac{1}{\lambda^{d-b}}
\prod_{j=1}^{d-b}\left(\log\lambda + \log|m_{j}|\right)\right)
& \rb\in \mathbf{Z}^{d}/\{0\} \\
\frac{1}{(2\pi)^{d}}\int_{\ob \in \mathbb{R}^{d}}f(\ob)g(\ob)d\ob + \frac{c(0)}{n}\sum_{\kb = -a}^{a}g(\ob_{\kb}) + 
O\big(\frac{\log\lambda + \log\|\rb\|_{1}}{\lambda}+\frac{1}{n}\big)
& \rb = 0
\end{array}
\right.. \qquad
\end{eqnarray*}
\end{itemize}
\end{theorem}
{\bf PROOF} The proof of (i) immediately follows from Theorem
\ref{lemma:nonuniformdft}. 

The proof of (ii) follows from 
writing $Q_{a,\lambda}(g;\rb)$ as a quadratic form  and taking expectations
\begin{eqnarray*}
&& \Ex\left[Q_{a,\lambda}(g;\rb)\right] \nonumber\\
&=& c_{2}\sum_{\kb=-a}^{a}g(\ob_{\kb})
\frac{1}{\lambda^{d}}\int_{[-\lambda/2,\lambda/2]^{d}}c(\ub_{1}-\ub_{2})\exp(i\ob_{k}^{\prime}
(\ub_{1}-\ub_{2})-i\ub_{2}^{\prime}\ob_{\rb})d\ub_{1}d\ub_{2} +W_{\rb},
\end{eqnarray*}
where $c_{2}=n(n-1)/n^{2}$ and
$W_{\rb}=\frac{c(0)I(\rb=0)}{n}\sum_{\kb=-a}^{a}g(\ob_{\kb})$
($I(\rb=0)=1$ if $\rb=0$ else $I(\rb=0)=01$). We then follow the same
proof used to prove Theorem \ref{lemma:mean-stat}. \hfill $\Box$

\vspace{3mm}
\begin{theorem}\label{theorem:varianceQQ}[Asymptotic expression for variance]
Suppose Assumptions \ref{assum:S}, \ref{assum:uniform}, \ref{assum:G}
\ref{assum:GG}(b,c) hold. 
Then we have 
\begin{eqnarray*}
\lambda^{d}\cov\left[Q_{a,\lambda}(g;\rb_{1}),Q_{a,\lambda}(g;\rb_{2})\right]=
\left\{
\begin{array}{cl}
C_{1}(\ob_{\rb}) + 2\Re [C_{3}(\ob_{\rb}) G_{\lambda}]
+ 2f_{2}(\ob_{\rb})|G_{\lambda}|^{2} +  O(\ell_{\lambda,a,n}) & \rb_{1}=\rb_{2}(=\rb) \\
O(\ell_{\lambda,a,n}) & \rb_{1}\neq\rb_{2}\\
\end{array}
\right.
\end{eqnarray*}
\begin{eqnarray*}
\lambda^{d}\cov\left[Q_{a,\lambda}(g;\rb_{1}),\overline{Q_{a,\lambda}(g;\rb_{2})}\right]=
\left\{
\begin{array}{cl}
C_{2}(\ob_{\rb}) + 2G_{\lambda}C_{3}(\ob_{\rb}) + 2f_{2}(\ob_{\rb})G_{\lambda}^{2} + O(\ell_{\lambda,a,n}) & \rb_{1}=-\rb_{2}(=\rb) \\
O(\ell_{\lambda,a,n}) & \rb_{1}\neq-\rb_{2}\\
\end{array}
\right.,
\end{eqnarray*}
where $C_{1}(\ob_{\rb})$ and $C_{2}(\ob_{\rb})$ are defined in
(\ref{eq:Cr}), $G_{\lambda} = \frac{1}{n}\sum_{\kb
= -a}^{a}g(\ob_{\kb})$, $a^{d}=O(n)$,
\begin{eqnarray}
\label{eq:C3r}
C_{3}(\ob_{\rb}) =
\frac{2}{(2\pi)^{d}}\int_{2\pi[-a/\lambda,a/\lambda]^{d}}g(\ob)f(-\ob)f(\ob + \ob_{\rb})d\ob 
\quad \textrm{ and }\quad
f_{2}(\ob_{\rb}) =
\int_{\mathbb{R}^{d}}f(\boldsymbol{\lambda})f(\ob_{\rb}- \boldsymbol{\lambda})d\boldsymbol{\lambda}.
\end{eqnarray}
\end{theorem}
{\bf PROOF} To prove the result we use (\ref{eq:Qalambda}). 
Using this expansion we rewrite $\cov[
  Q_{a,\lambda}(g;\rb_{1}),Q_{a,\lambda}(g;\rb_{2})]$ in terms of 
$\widetilde{Q}_{a,\lambda}(g;\rb)$ and $V_{\rb}$
\begin{eqnarray*}
&&\cov\left[
  Q_{a,\lambda}(g;\rb_{1}),Q_{a,\lambda}(g;\rb_{2})\right] =\\
&& 
\cov\left[
  \widetilde{Q}_{a,\lambda}(g;\rb_{1}),\widetilde{Q}_{a,\lambda}(g;\rb_{2})\right]
+ |G_{\lambda}|^{2}\cov\left[V_{\rb_{1}},V_{\rb_{2}}\right] + 
\overline{G_{\lambda}}\cov\left[\widetilde{Q}_{a,\lambda}(g;\rb_{1}),V_{\rb_{2}}\right]
+ G_{\lambda}\overline{\cov\left[\widetilde{Q}_{a,\lambda}(g;\rb_{2}),V_{\rb_{1}}\right]}.
\end{eqnarray*}
In Theorem \ref{theorem:asymptotic} we have evaluated an asymptotic
expression for $\cov\left[
  \widetilde{Q}_{a,\lambda}(g;\rb_{1}),\widetilde{Q}_{a,\lambda}(g;\rb_{2})\right]$. We
now evaluate similar expressions for
$\cov\left[V_{\rb_{1}},V_{\rb_{2}}\right]$ and $\cov\left[\widetilde{Q}_{a,\lambda}(g;\rb_{1}),V_{\rb_{2}}\right]$.
It is straightforward to show that 
\begin{eqnarray}
\label{eq:covVV}
\lambda^{d} \cov[V_{\rb_{1}},V_{\rb_{2}}] &=& 
\left\{
\begin{array}{cc}
O(\frac{1}{\lambda}) & \rb_{1}\neq \rb_{2} \\
 2f_{2}(\ob_{\rb_{1}}) + O\left(\frac{1}{\lambda} + \frac{\lambda^{d}}{n}\right) & \rb_{1}= \rb_{2} \\
\end{array}
\right.
\end{eqnarray}
To evaluate an expression for
$\cov\left[\widetilde{Q}_{a,\lambda}(g;\rb_{1}),V_{\rb_{2}}\right]$ we
consider the case $d=1$. Using similar arguments to those in the proof
of Lemma \ref{lemma:covariance} we can show that 
\begin{eqnarray*}
&&\lambda\cov\left[\widetilde{Q}_{a,\lambda}(g;\rb_{1}),V_{\rb_{2}}\right]\\
&=& 2\lambda c_{3}\sum_{k=-a}^{a}g(\omega_{k})\Ex\left[
c(s_{1}-s_{3})c(s_{2}-s_{3})e^{is_{1}\omega_{k}}e^{is_{2}(y-\omega_{k}-\omega_{r})}e^{is_{3}(-x-y-\omega_{r})}
\right] + O\left(\frac{\lambda}{n^{2}}\right) \\
&=&
\frac{2c_{3}\lambda}{(2\pi)^{2}}\sum_{k=-a}^{a}g(\omega_{k})\int_{\mathbb{R}^{3}}
f(x)f(y)\sinc\left(\frac{\lambda x}{2}+k\pi\right)
\sinc\left(\frac{\lambda x}{2}-(k+r_{1})\pi\right)
\sinc\left(\frac{\lambda(x+y)}{2}-r_{2}\pi\right)dxdy \\
&&+O\left(\frac{\lambda^{}}{n^{2}}\right) \\
&=&\frac{2c_{3}}{\lambda^{}\pi^{2}}\sum_{k=-a}^{a}g(\omega_{k})\int_{\mathbb{R}^{3}}
f\left(\frac{2u}{\lambda}-\omega_{k}\right)f\left(\frac{2v}{\lambda}+\omega_{k+r}\right)\sinc\left(u\right)
\sinc\left(v\right)
\sinc\left(u+v+(r_{1}-r_{2})\pi\right)dudv \\
&&+O\left(\frac{\lambda^{}}{n^{2}}\right),
\end{eqnarray*}
where $c_{3}=n(n-1)(n-2)/n^{3}$. Finally by replacing sum with an
integral, using Lemma \ref{lemma:sum-integral}
and replacing
$f\left(\frac{2u}{\lambda}-\omega\right)f\left(\frac{2v}{\lambda}+\omega+\omega_{r}\right)$
with $f\left(-\omega\right)f\left(\omega_{}+\omega_{r}\right)$ and
 by using (\ref{eq:cum4bound2}) we have 
\begin{eqnarray}
\label{eq:covQV}
&&\lambda\cov\left[\widetilde{Q}_{a,\lambda}(g;r_{1}),V_{r_{2}}\right]\nonumber\\
&=&\frac{2c_{3}}{2\pi\pi^{2}}\int_{-2\pi a/\lambda}^{2\pi
  a/\lambda}g(\omega) f\left(-\omega\right)f\left(\omega+\omega_{r_{1}}\right)d\omega
\int_{\mathbb{R}^{3}}
\sinc\left(u\right)\sinc\left(v\right)
\sinc\left(u+v+(r_{1}-r_{2})\pi\right)dudv \nonumber\\
&&+O\left(\frac{\lambda^{}}{n^{2}}+\frac{\log^{2}\lambda}{\lambda}\right)
\nonumber\\
&=& \left\{
\begin{array}{cc}
C_{3}(\omega_{r_{1}})
+O\left(\frac{\lambda^{}}{n^{2}}+\frac{\log^{2}\lambda}{\lambda}+\frac{1}{n}\right)
& r_{1}=r_{2} \\
O\left(\frac{\lambda^{d}}{n^{2}}+\frac{\log^{2}\lambda}{\lambda}+\frac{1}{n}\right)
& r_{1}\neq r_{2}\\
\end{array}
\right.
\end{eqnarray}
Altogether, by using Theorem \ref{theorem:asymptotic},
(\ref{eq:covQV}) and (\ref{eq:covVV}) we obtain the result. \hfill
$\Box$

\vspace{3mm}

In the following lemma we show that we further simplify the expression
for the asymptotic variance if we keep $\rb$ fixed. 
\begin{corollary}\label{cor:QQ}
Suppose Assumption \ref{assum:G}, \ref{assum:GG}(a,c) or \ref{assum:GG}(b,c)
holds, and $\rb$ is fixed. 
 Let $C_{3}(\ob)$ and $f_{2}(\ob)$ be defined in  (\ref{eq:C3r}).
Then we have 
\begin{eqnarray*}
C_{3}(\ob_{\rb}) = C_{3} + O\left(\frac{\|\rb\|_{1}}{\lambda}\right),
\end{eqnarray*}
and $f_{2}(\ob_{\rb}) = f_{2} + O\left(\frac{\|\rb\|_{1}}{\lambda}\right)$,
where
\begin{eqnarray*}
C_{3}  &=& \frac{2}{(2\pi)^{d}}\int_{2\pi[-a/\lambda,a/\lambda]^{d}}g(\ob)f(\ob)^{2}d\ob 
\end{eqnarray*}
and $f_{2}=f_{2}(0)$ (note, if $\overline{g(\ob)} = g(-\ob)$, then $C_{1} = C_{2}$).
\end{corollary}
{\bf PROOF} The proof is the same as the proof of Corollary \ref{cor:C}. \hfill $\Box$

\begin{theorem}\label{theorem:CLTQ}[CLT on real and imaginary parts]
Suppose Assumptions \ref{assum:S}, \ref{assum:uniform},
\ref{assum:G} and \ref{assum:GG}(b,c)  hold. Let $C_{1}$, 
$C_{2}$, $C_{3}$ and $f_{2}$ be
defined as in Corollary \ref{cor:QQ}.
We define the $m$-dimension complex random vectors 
$\boldsymbol{Q}_{m} =
(Q_{a,\lambda}(g,\rb_{1}),\ldots,Q_{a,\lambda}(g,\rb_{m}))$,
where $\rb_{1},\ldots,\rb_{m}$ are such that $\rb_{i}\neq -\rb_{j}$ and $\rb_{i}\neq
0$. %then for any $i,j\in \{1,\ldots,m\}$ 
Under these conditions we have 
\begin{eqnarray}
\label{eq:cltQ}
\frac{2\lambda^{d/2}}{E_{1}}\left(
%\begin{array}{ccc}
\frac{E_{1}}{E_{1}+\Re E_{2}}\Re  Q_{a,\lambda}(g,0),  
\Re \boldsymbol{Q}_{m},
\Im \boldsymbol{Q}_{m}
%\end{array}
\right)\Pcon
 \mathcal{N}\big(0,I_{2m+1}\big),
\end{eqnarray}
where $E_{1} = C_{1} + 2\Re [G_{\lambda}C_{3}]
+ 2f_{2}|G_{\lambda}|^{2}$ and 
$E_{2} = C_{2} + 2G_{\lambda}C_{3} + 2f_{2}G_{\lambda}^{2}$
with $\frac{\lambda^{d}}{n} \rightarrow 0$ and $\frac{\log^{2}(a)}{\lambda^{1/2}}\rightarrow
0$  as $\lambda \rightarrow \infty$, $n\rightarrow \infty$ and
$a\rightarrow \infty$.
\end{theorem}
{\bf PROOF} Using the same method used to prove 
Lemma \ref{lemma:cumulantq}, and 
analogous results can be derived for the cumulants of
$Q_{a,\lambda}(g;\rb)$. Asymptotic normality follows from this.
We omit the details. \hfill $\Box$

\subsubsection{Application to nonparametric covariance estimator}\label{sec:nonparametric-prop}

In this section we apply the results to the nonparametric estimator
considered in Section \ref{sec:nonparametric}. We
define
\begin{eqnarray*}
\widetilde{c}_{\Omega,n}(\ubb) = 
\frac{1}{\Omega^{d}}\sum_{\kb=-a}^{a}|J_{n}(\ob_{\Omega,\kb})|^{2}
\exp(i\ubb^{\prime}\ob_{\Omega,\kb})
\textrm{ and }
\widehat{c}_{,\Omega,n}(\ubb) = 
 T\left(\frac{\ubb}{\widehat{\Omega}}\right)\widetilde{c}_{\Omega,n}(\ubb) 
 \end{eqnarray*}
where $T$ is the $d$-dimensional triangle kernel. 
It is clear the asymptotic sampling properties of
$\widehat{c}_{\Omega,n}(\ubb)$ are determined by
$\widetilde{c}_{\Omega,n}(\ubb)$. Therefore, we first derive the
asymptotic sampling properties of $\widetilde{c}_{\Omega,n}(\ubb)$. We
observe that 
$\widetilde{c}_{\Omega,n}(\ubb) = Q_{a,\Omega,\lambda}(e^{i{\boldsymbol v}^{\prime}\cdot};0)$, thus we
use the results in Section \ref{sec:uniform} to derive the asymptotic
sampling properties of $\widetilde{c}_{\Omega,n}(\ubb)$. 
By using Theorem \ref{lemma:summaryE} and under Assumptions \ref{assum:S},
\ref{assum:uniform} and \ref{assum:G}(ii) we have 
\begin{eqnarray}
\label{eq:Etildec}
\Ex[\widetilde{c}_{\Omega,n}(\ubb)] &=& \frac{1}{(2\pi)^{d}}\int_{2\pi[-a/\Omega,a/\Omega]^{d}}
f(\ob)\exp(i\ubb^{\prime}\ob)d\ob + O\left(\frac{\log\lambda}{\lambda}\right)\nonumber\\
&=& c(\ubb) +
O\left(\left(\frac{\Omega}{a} \right)^{\delta}+ \frac{\log \lambda}{\lambda}\right),
\end{eqnarray}
for $\ubb \in [-\lambda/2,\lambda/2]^{d}$ (if $\Omega = \lambda$) else
for $\ubb \in [-\lambda,\lambda]^{d}$ if $\Omega \geq 2\lambda$. 
We recall that $\delta$ is such that it satisfies Assumption
\ref{assum:G}(ii)(a).  By adapting the proof of Theorem
\ref{theorem:varianceQQ}, to fine grids (with $\lambda \leq \Omega$) we have  
$\lambda^{d}\var\left[\widetilde{Q}_{a,\lambda,\Omega}(g;0)\right] = 
\Sigma(\frac{a}{\Omega})+O\left(\widetilde{\ell}_{a,\lambda,\Omega}\right)$
where 
\begin{eqnarray*}
\Sigma\left(\frac{a}{\Omega}\right) &=&  
\left[\frac{1}{(2\pi)^{d}}\int_{-[2\pi a/\Omega,2\pi a/\Omega]^{d}} 
f(\ob_{})^{2}\left[1 + \exp(2i\ubb^{\prime}\ob)\right]d\ob\right]\left[
\frac{\lambda^{d}}{\Omega^{d}}\sum_{\kb=-2a}^{2a}\Sinc^{2}\left(
  \frac{\lambda}{\Omega}\kb 
  \pi\right)\right] + \\
&& 2f_{2}(0)
\left|\frac{\lambda^{d}}{\Omega^{d}n}\sum_{\kb=-a}^{a}\exp(i\ubb^{\prime}\ob_{\Omega,\kb})\right|^{2}
+ \frac{4}{(2\pi)^{d}}\int_{[-2\pi a/\Omega,2\pi
  a/\Omega]^{d}}\exp(i\ubb^{\prime}\ob)f(\ob)^{2}d\ob.
\end{eqnarray*}
Therefore, if $\Omega^{\delta}\lambda^{d/2}/a^{\delta}\rightarrow 0$ as
$a\rightarrow \infty$ and $\lambda\rightarrow \infty$, then
by using Theorem \ref{theorem:CLTQ} we have 
\begin{eqnarray}
\label{eq:chatCLT}
\lambda^{d/2}\left[\widetilde{c}_{\Omega,n}(\ubb) -c(\ubb)\right]
\Dcon \mathcal{N}\left(0,\Sigma\left(\frac{a}{\lambda}\right)\right),
\end{eqnarray}
$\frac{\lambda^{d}}{n} \rightarrow 0$ and 
$\frac{\log^{2}(a)}{\lambda^{1/2}}\rightarrow 0$ as $n\rightarrow \infty$,
$a\rightarrow \infty$ and $\lambda\rightarrow \infty$.
By using (\ref{eq:Etildec}) and (\ref{eq:chatCLT}) we have 
\begin{eqnarray*}
\Ex[\widehat{c}_{\Omega,n}(\ubb)] = T\left(\frac{\ubb}{\widehat{\Omega}}\right)c(\ubb) +
O\left(\widetilde{\ell}_{a,\lambda,\Omega}\right).
\end{eqnarray*}

\subsubsection*{Application to parameter estimation using an $L_{2}$ criterion}

In this section we consider the asymptotic sampling properties of
$\widehat{\theta}_{n} = \arg\min_{\theta \in \Theta} L_{n}(\theta)$, where $L_{n}(\cdot)$
is defined in (\ref{eq:L2estimator}) and $\Theta$ is a compact set. 
We will assume that there exists a $\theta_{0}\in \Theta$, such that
for all $\ob\in \mathbb{R}^{d}$, $f_{\theta_{0}}(\ob) = f(\ob)$ and there does
not exist another $\theta\in \Theta$ such that for all $\ob\in
\mathbb{R}^{d}$ $f_{\theta_{0}}(\ob)= f_{\theta}(\ob)$ and 
in addition
$\int_{\mathbb{R}^{d}}\|\nabla_{\theta}f(\ob;\theta_{0})\|_{1}^{2}d\ob <\infty$. Furthermore, we
will assume that $\widehat{\theta}_{n}\Pcon \theta_{0}$ as $\lambda
\rightarrow \infty$.

Making the usual Taylor expansion we have 
$\lambda^{d/2}(\widehat{\theta}_{n} - \theta_{0}) =
A^{-1}\lambda^{d/2}\nabla L_{n}(\theta_{0}) + o_{p}(1)$, 
where
\begin{eqnarray}
\label{eq:Wl2}
A=\frac{2}{(2\pi)^{d}}\int_{\mathbb{R}^{d}}[\nabla_{\theta}f(\ob;\theta_{0})]^{}[\nabla_{\theta}f(\ob;\theta_{0})]^{\prime}d\ob,
\end{eqnarray}
and  it is clear  the asymptotic sampling
properties of $\widehat{\theta}_{n}$ are determined by
$\nabla_{\theta}L_{n}(\theta_{0})$, which we see from (\ref{eq:DLn})
can be written as 
$\nabla_{\theta}L_{n}(\theta_{0})=Q_{a,\lambda}(-2 
\nabla_{\theta}f_{\theta_{0}}(\cdot);0) +
\frac{1}{\lambda^{d}}\sum_{\kb=-a}^{a}\nabla_{\theta}f_{\theta_{0}}(\ob_{\kb})^{2}$.

Thus by using Theorem
\ref{theorem:CLTQ} we have
$\lambda^{d/2}\nabla_{\theta}L_{n}(\theta_{0})\Dcon \mathcal{N}(0,B)$, where
\begin{eqnarray*}
B &=& \frac{4}{(2\pi)^{d}}\int_{\mathbb{R}^{d}}f(\ob)^{2}\left[
  \nabla_{\theta}f_{\theta}(\ob_{})
\right]\left[\nabla_{\theta}f_{\theta}(\ob_{})\right]^{\prime}\rfloor_{\theta
= \theta_{0}}d\ob
\\
&&
-\frac{4\Re G_{\lambda}}{(2\pi)^{d}}\int_{\mathbb{R}^{d}}f(\ob)^{2}\nabla_{\theta}f_{\theta}(\ob_{})\rfloor_{\theta=\theta_{0}}d\ob
+ 2|G_{\lambda}|^{2}\int_{\mathbb{R}^{d}}f(\ob)^{2}d\ob
\end{eqnarray*}
and $G_{\lambda} = \frac{2}{n}\sum_{\kb
= -a}^{a}\nabla_{\theta}f_{\theta}(\ob_{\kb})\rfloor_{\theta=\theta_{0}}$.
Therefore, by using the above we have 
\begin{eqnarray*}
\lambda^{d/2}(\widehat{\theta}_{n} - \theta_{0})\Pcon \mathcal{N}(0,A^{-1}BA^{-1})
\end{eqnarray*}
with $\frac{\lambda^{d}}{n} \rightarrow 0$ and 
$\frac{\log^{2}(a)}{\lambda^{1/2}}\rightarrow 0$ as
$a\rightarrow \infty$ and $\lambda\rightarrow \infty$.

% We recall  that for the exponential covariance  in dimension $d=2$, that
%$\delta = 1$ (see Remark \ref{remark:exp}). Thus a large number of frequencies, $a$, need to be used
%in the construction of $\widehat{c}_{n}({\boldsymbol v})$ in order to
%reduce the bias. For covariance whose spectral densities have thinner
%tails (such as the Whittle), using so many frequencies is not as necessary.  

%% file: B_7_simulations.tex
\section{Simulations}\label{sec:simulations}

In this section we illustrate the performance of the nonparametric non-negative definite estimator of the
  spatial covariance defined in  Section \ref{sec:nonparametric}. We
  compare our method to nonparametric estimator proposed in
  \citeA{p:hal-fis-hof-94}. We conduct all the simulations for $d=1$
  and the observations are observed over the spatial domain $[-20,20]$
  ($\lambda = 40$).

We use the estimator 
\begin{eqnarray*}
\widehat{c}_{a,\Omega,n}(u) &=&  
T\left(\frac{u}{\widehat{\Omega}}\right)
\widetilde{c}_{\Omega,n}(u),
\end{eqnarray*}
where 
\begin{eqnarray*}
\widetilde{c}_{a,\Omega, n}(u) = 
\frac{1}{\Omega^{}}\sum_{k  = -a}^{a}
\left|J_{n}\left(\omega_{k,\Omega}\right)\right|^{2}\exp(iu\omega_{k,\Omega}). 
\end{eqnarray*}
Since $\lambda = 40$ we use $\Omega = 80$ and
$\widehat{\Omega} = 40$ to construct the estimator. 

To evaluate the covariance estimator proposed in \citeA{p:hal-fis-hof-94} (from now on referred to as the HFH
estimator) we use the kernel method proposed in
\citeA{p:hal-fis-hof-94} to estimate the covariance $c(u)$ at 
$u=0,0.5,1,1.5,\ldots,30$, we denote this estimator as
$\widetilde{c}_{HRH}(u)$. As suggested by
\citeA{p:hal-fis-hof-94} for $u=30,30.5\ldots,35$ we taper the
covariance to zero and let
$\widetilde{c}_{HRH}(u) = (1-\frac{|u|-30}{5})\widehat{c}_{HRH}(30)$ for
$u\in[35,40]$ and let $\widetilde{c}_{HRH}(u)=0$ for $u\in [35,40]$. To
ensure the estimator is non-negative definite we evaluate the Fourier
transform of $\{\widetilde{c}_{HRH}(u);u=0,0.5,1,\ldots,40\}$ and set
negative values to zero and invert the Fourier transform. The result
is an estimator of the spatial covariance sequence which is a non-negative
definite function. We denote this estimator as $\widehat{c}_{HRH}(u)$. 

\begin{figure}
\centering
\includegraphics[scale = 0.3]{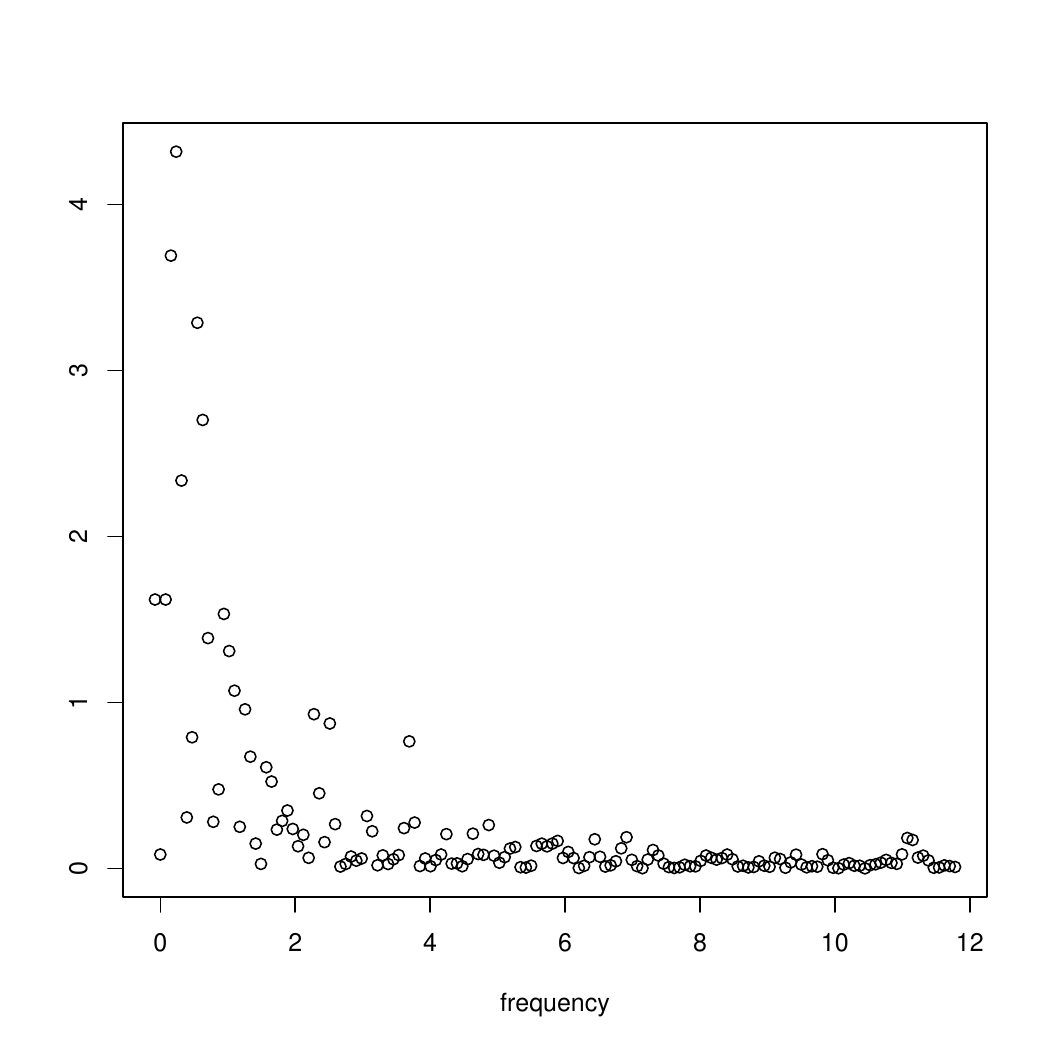}
\includegraphics[scale = 0.3]{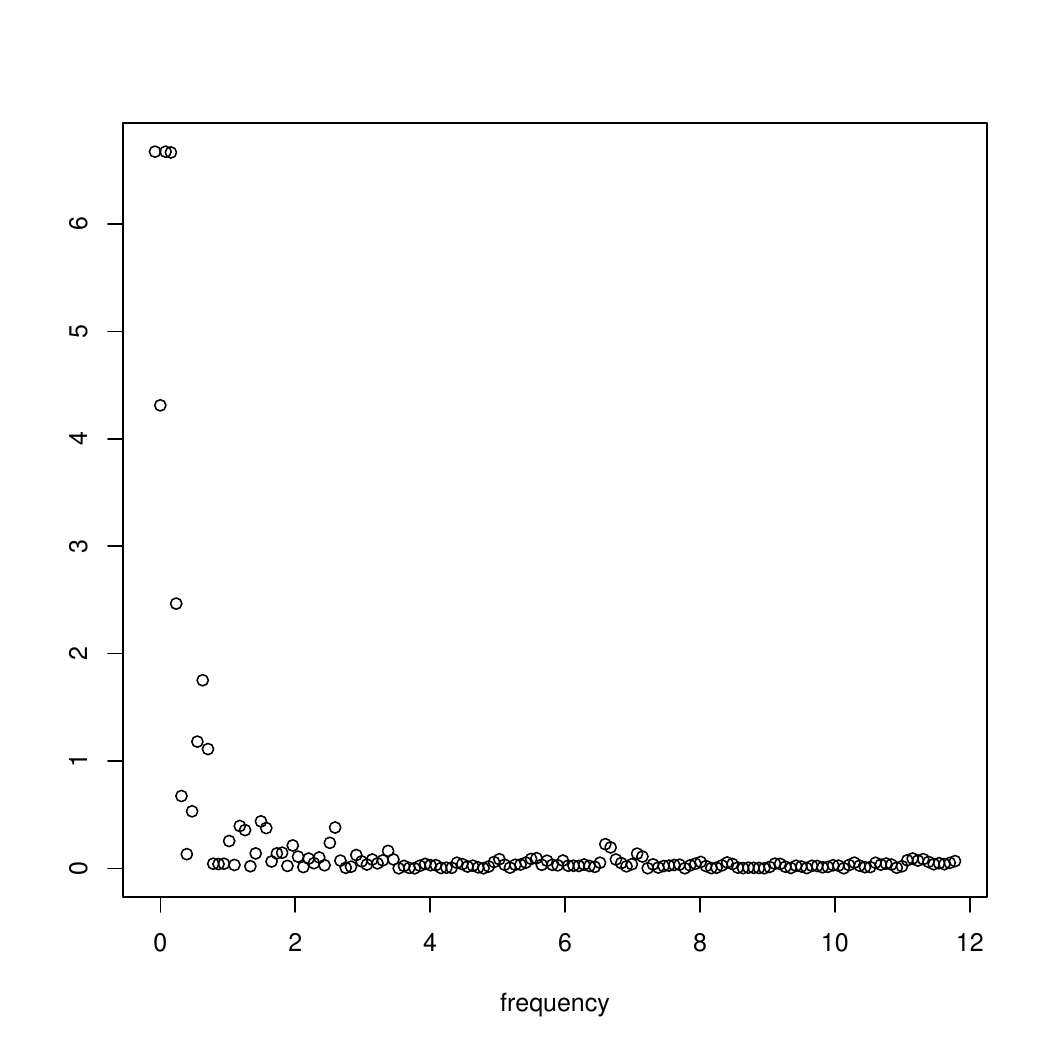}
\includegraphics[scale = 0.3]{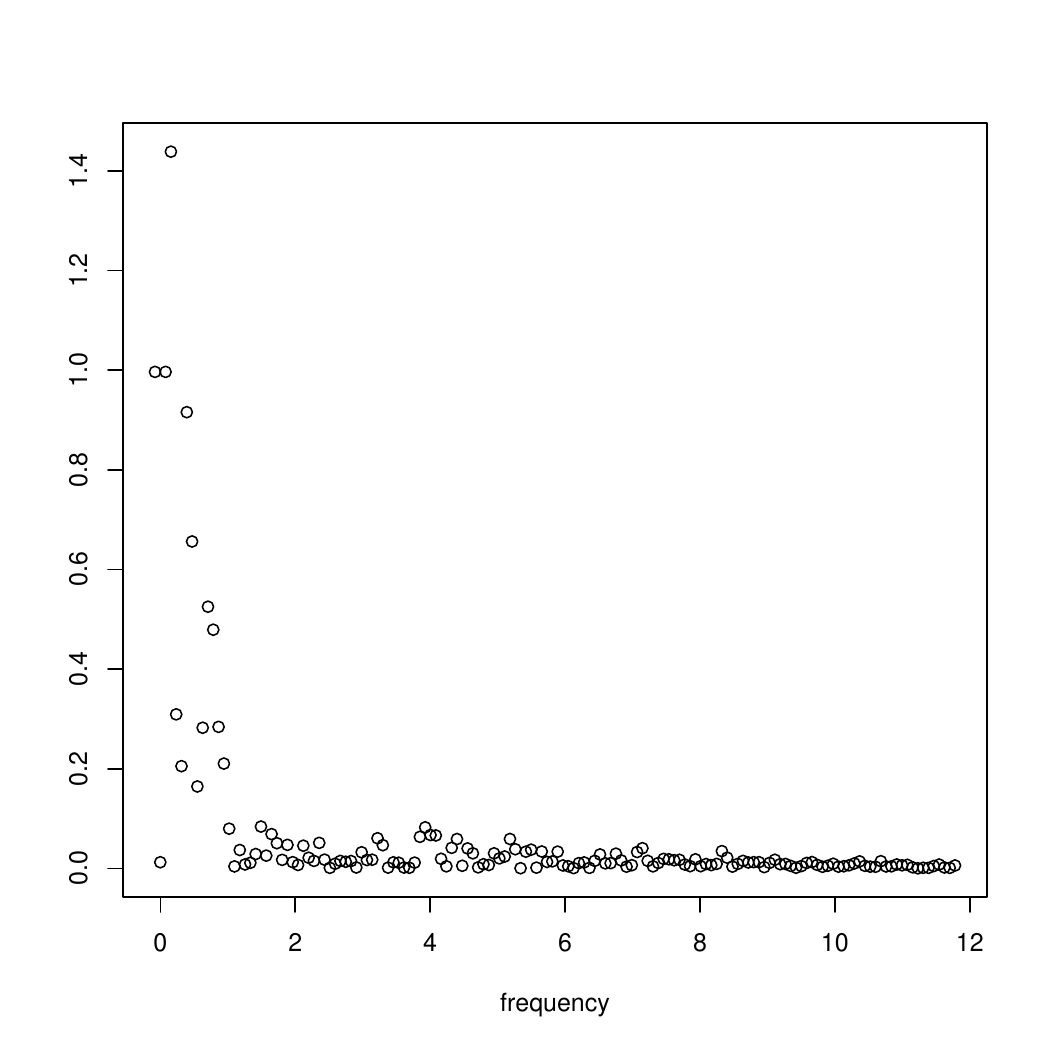}

\includegraphics[scale = 0.3]{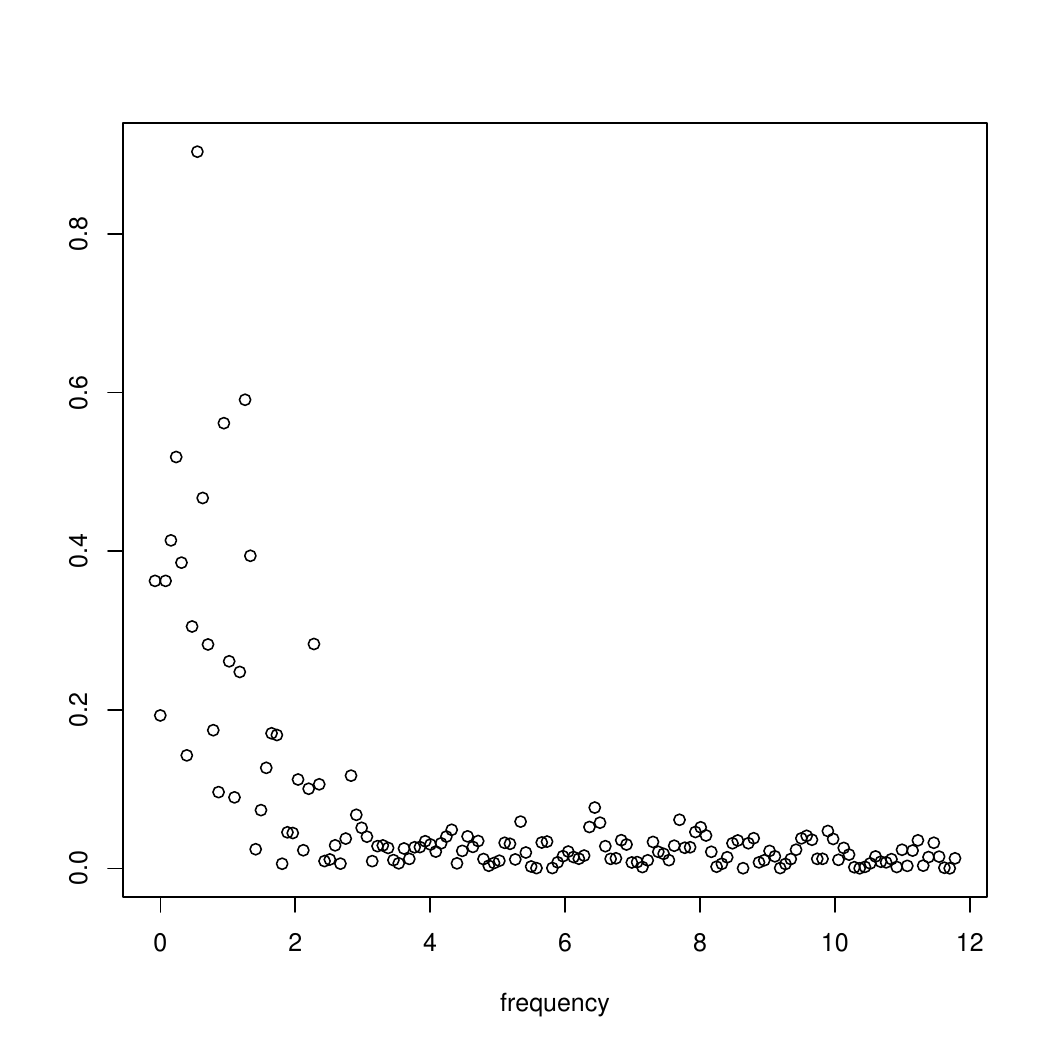}
\includegraphics[scale = 0.3]{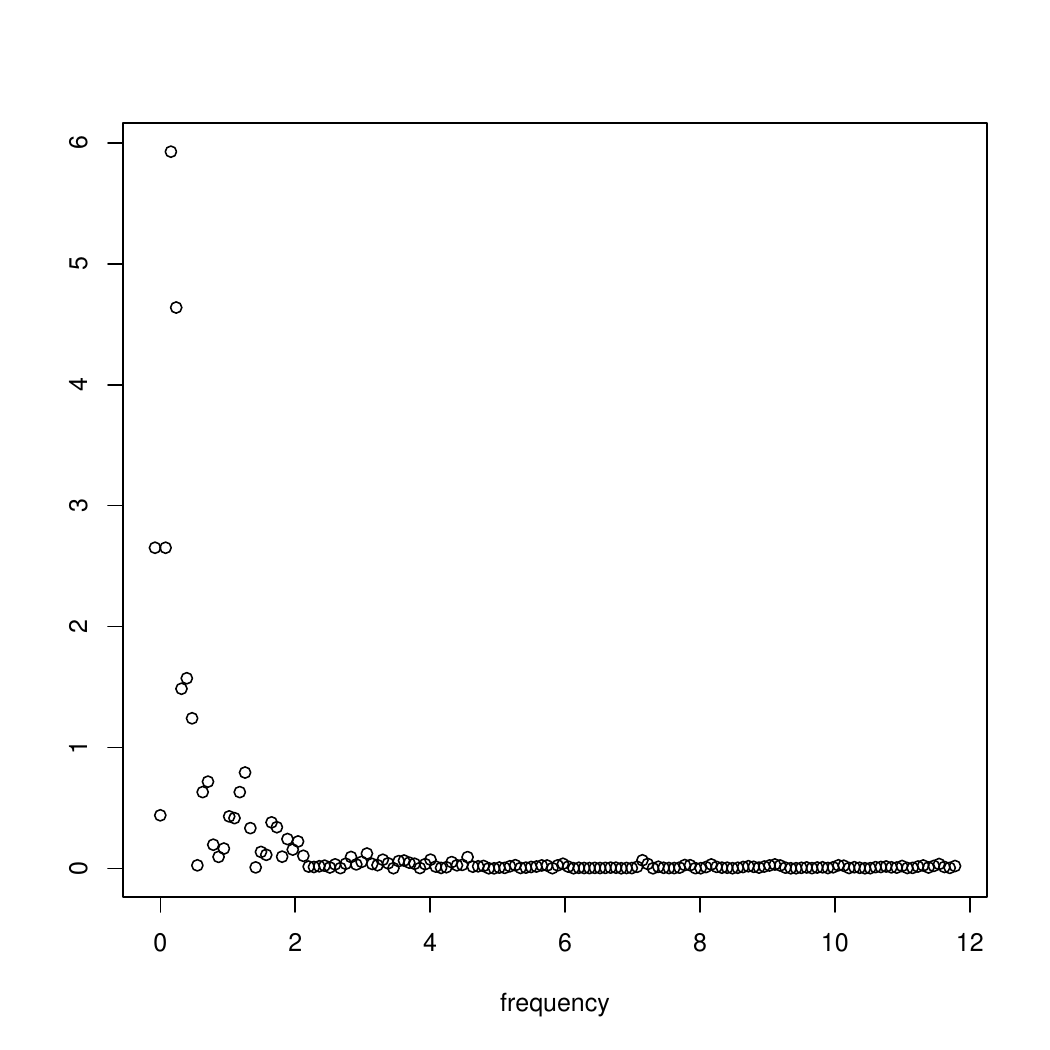}
\includegraphics[scale = 0.3]{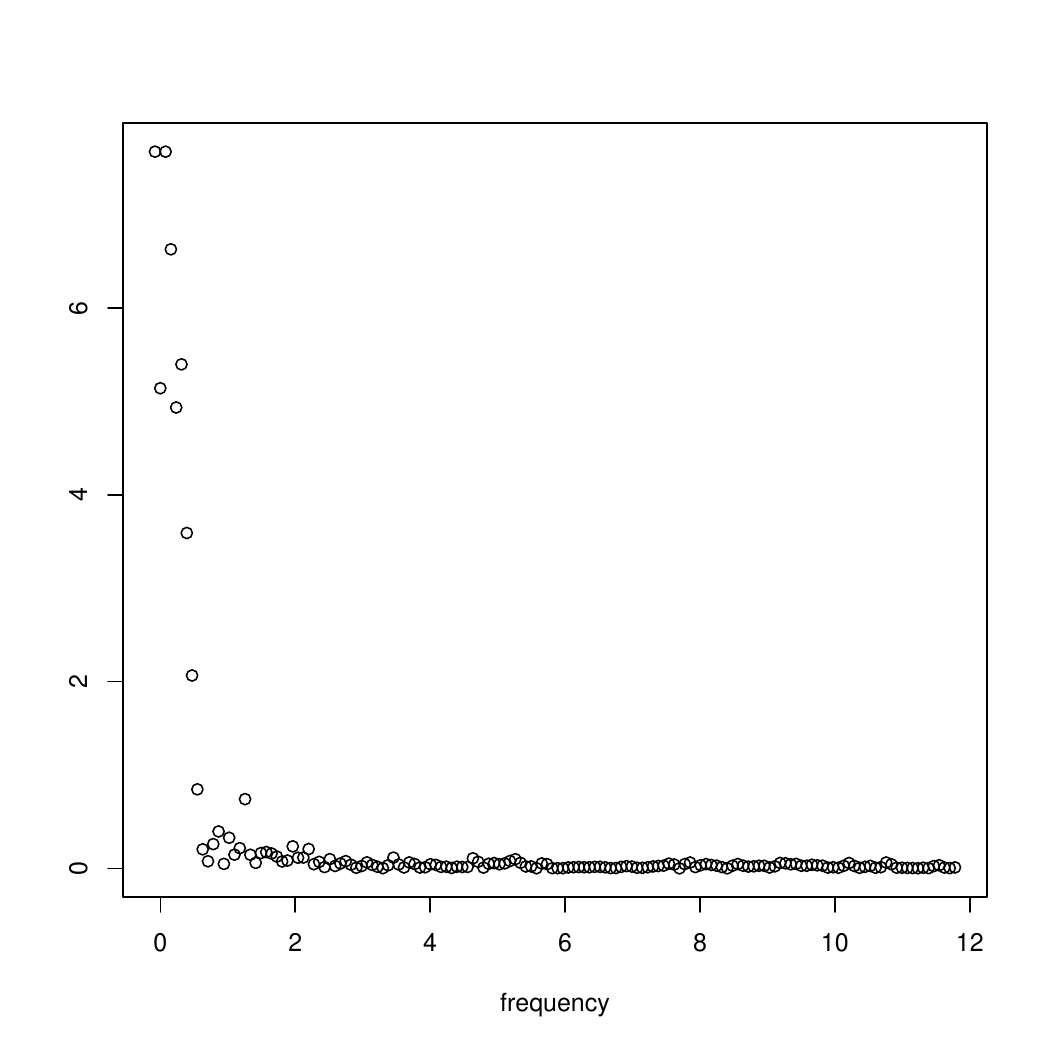}
\caption{Plot of $|J_{n}(\omega_{80,k})|^{2}$. Top for $n=1000$ and
  Left to Right $R=2,5$ and $10$. Bottom for $n=2000$ and
  Left to Right $R=2,5$ and $10$. 
\label{fig:1}}
\end{figure}

\begin{figure}
\centering
\includegraphics[scale = 0.3]{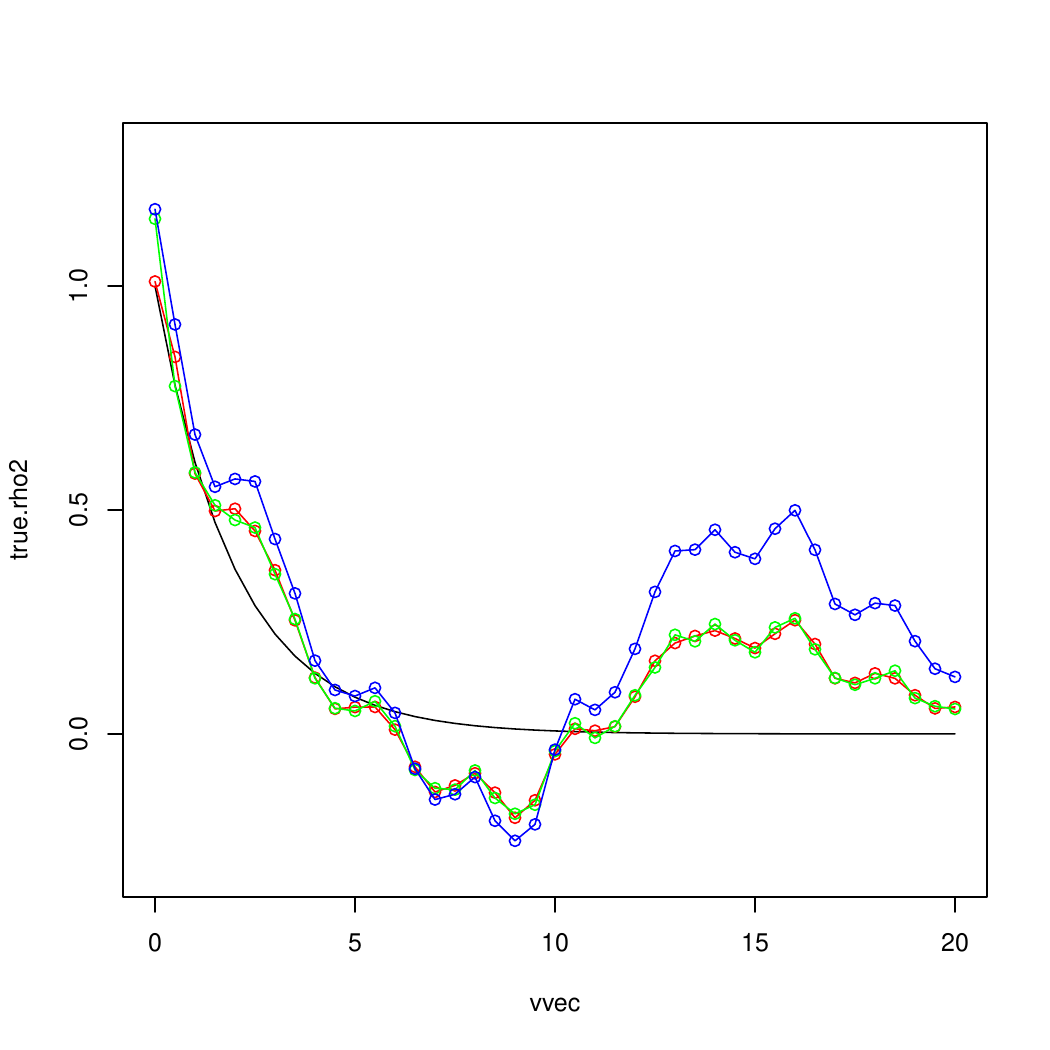}
\includegraphics[scale = 0.3]{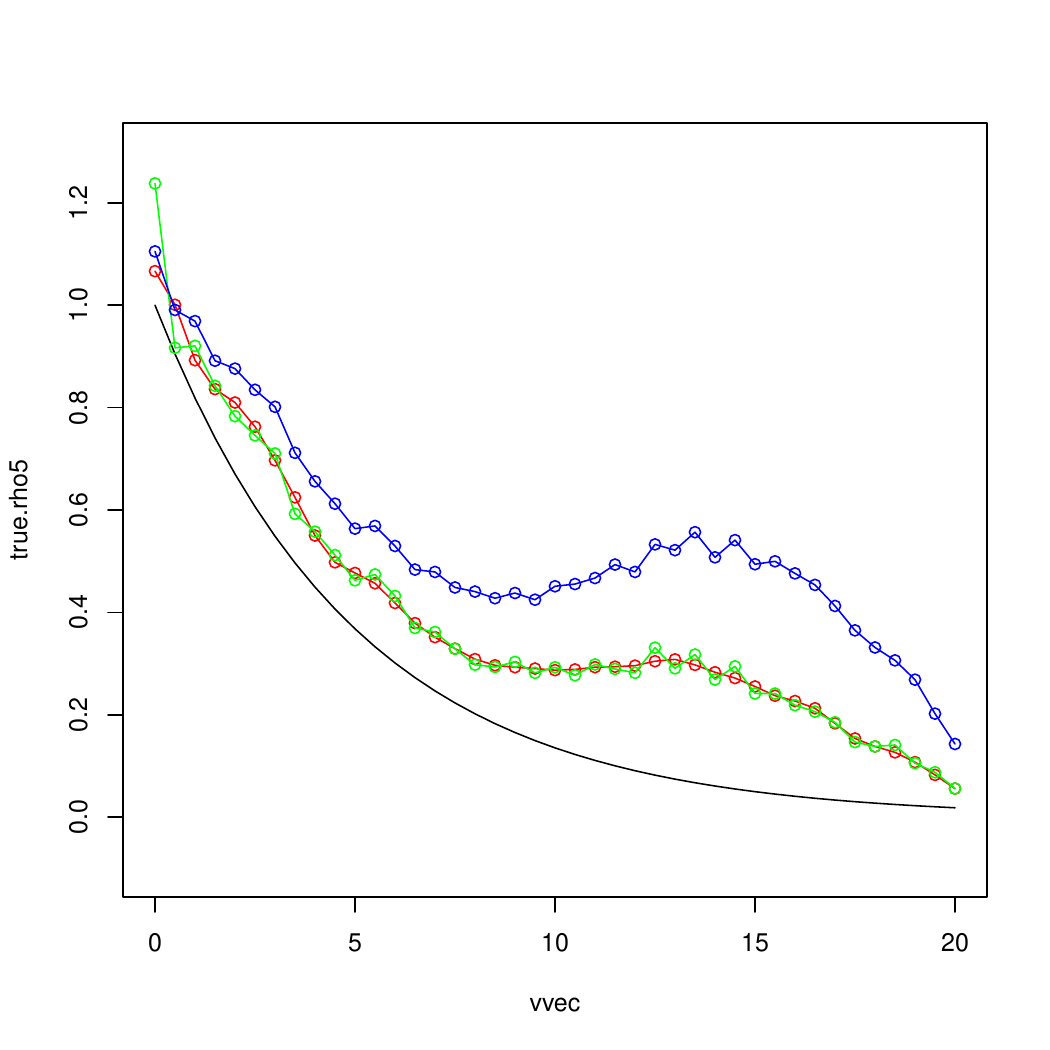}
\includegraphics[scale = 0.3]{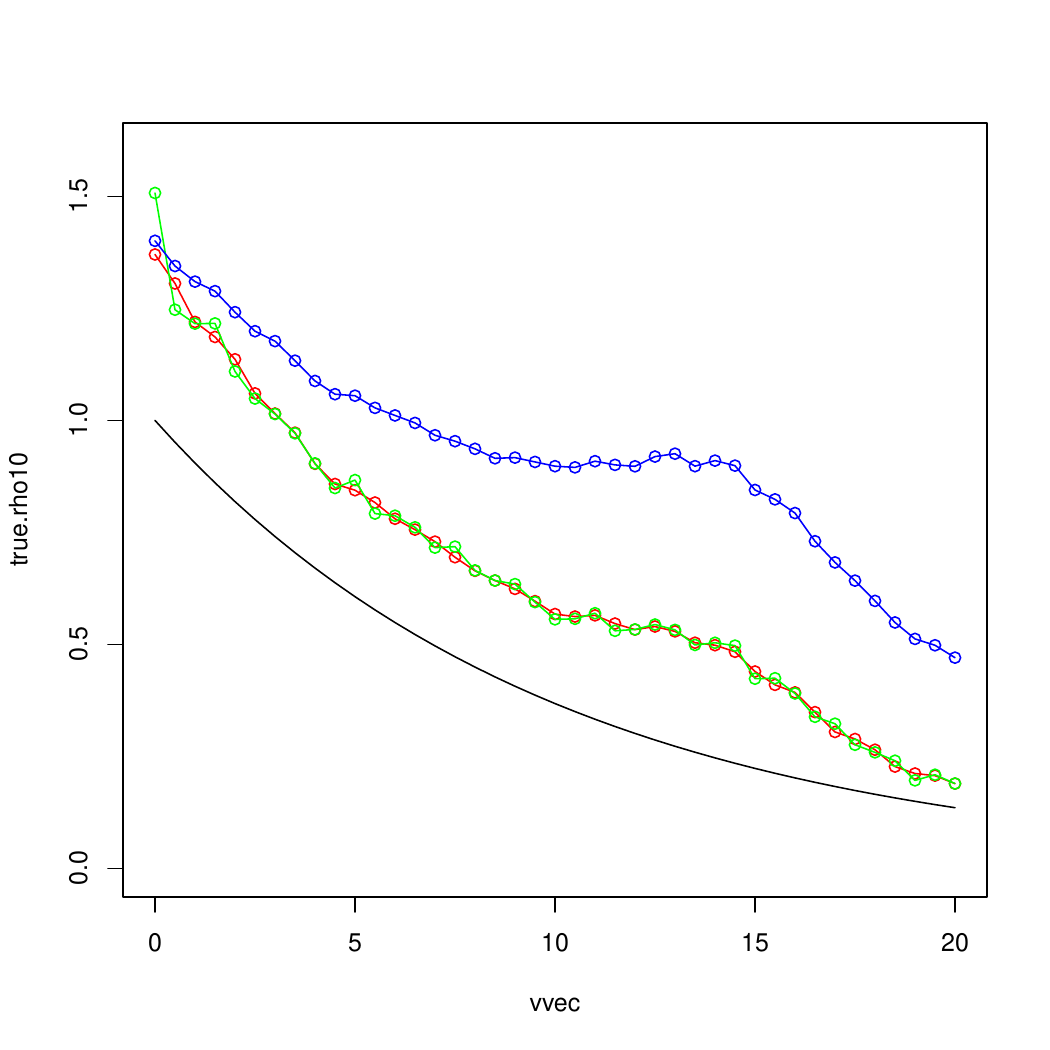}

\includegraphics[scale = 0.3]{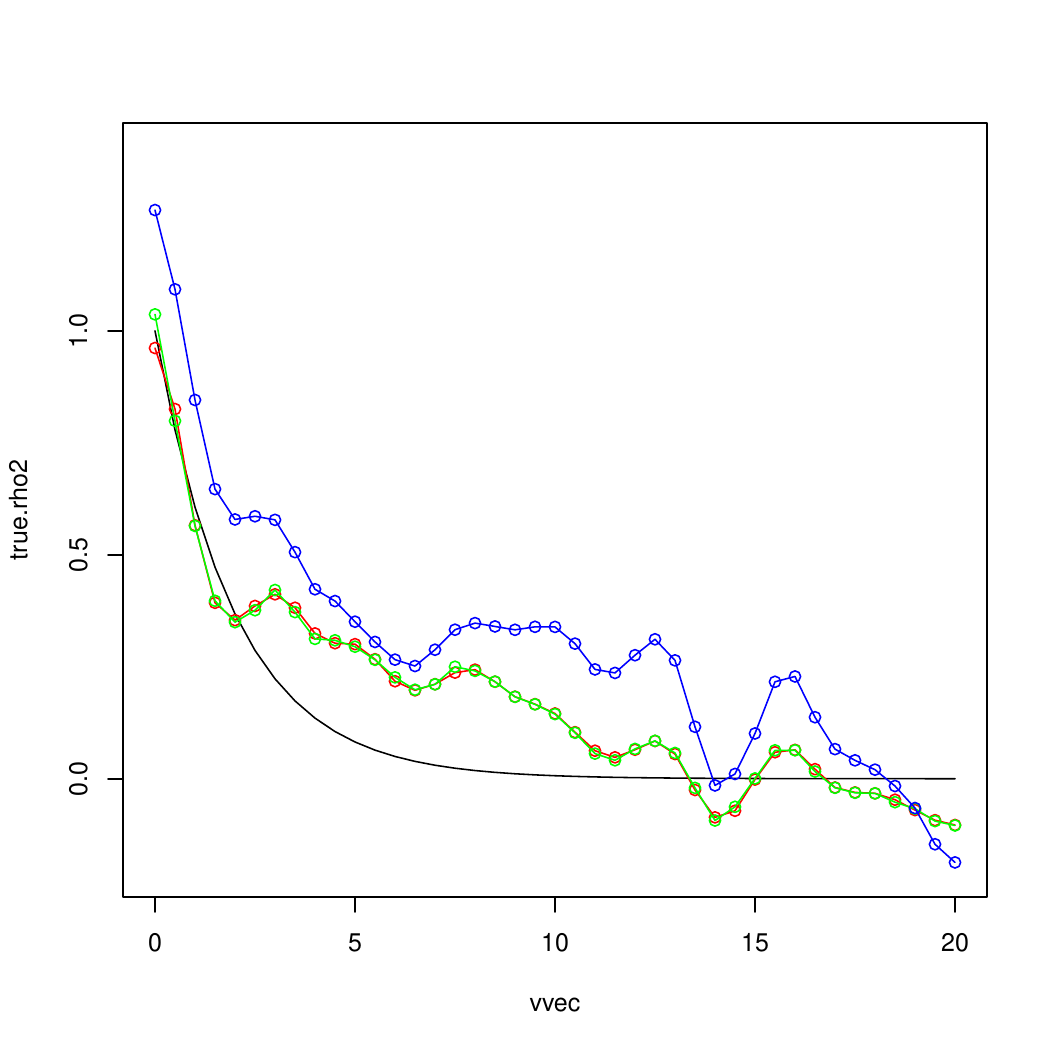}
\includegraphics[scale = 0.3]{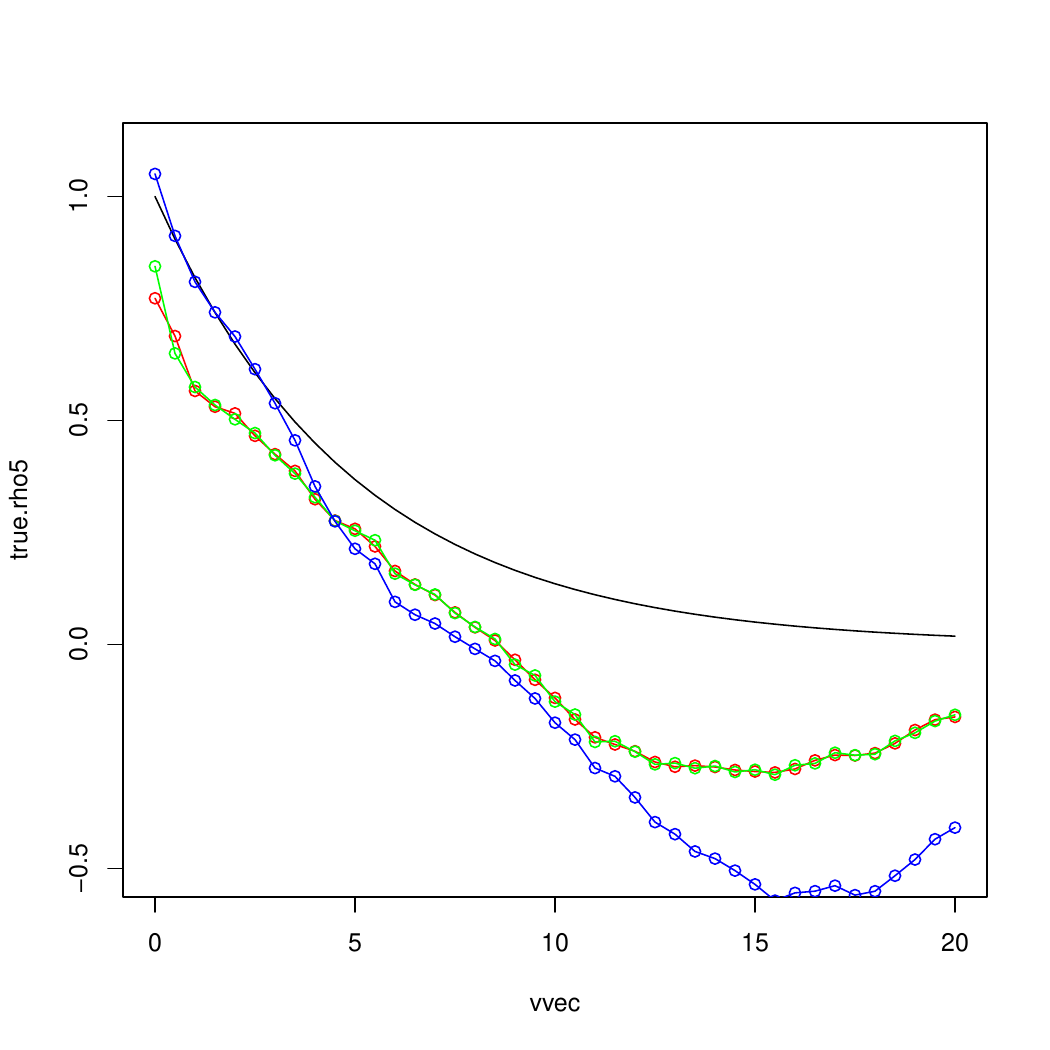}
\includegraphics[scale = 0.3]{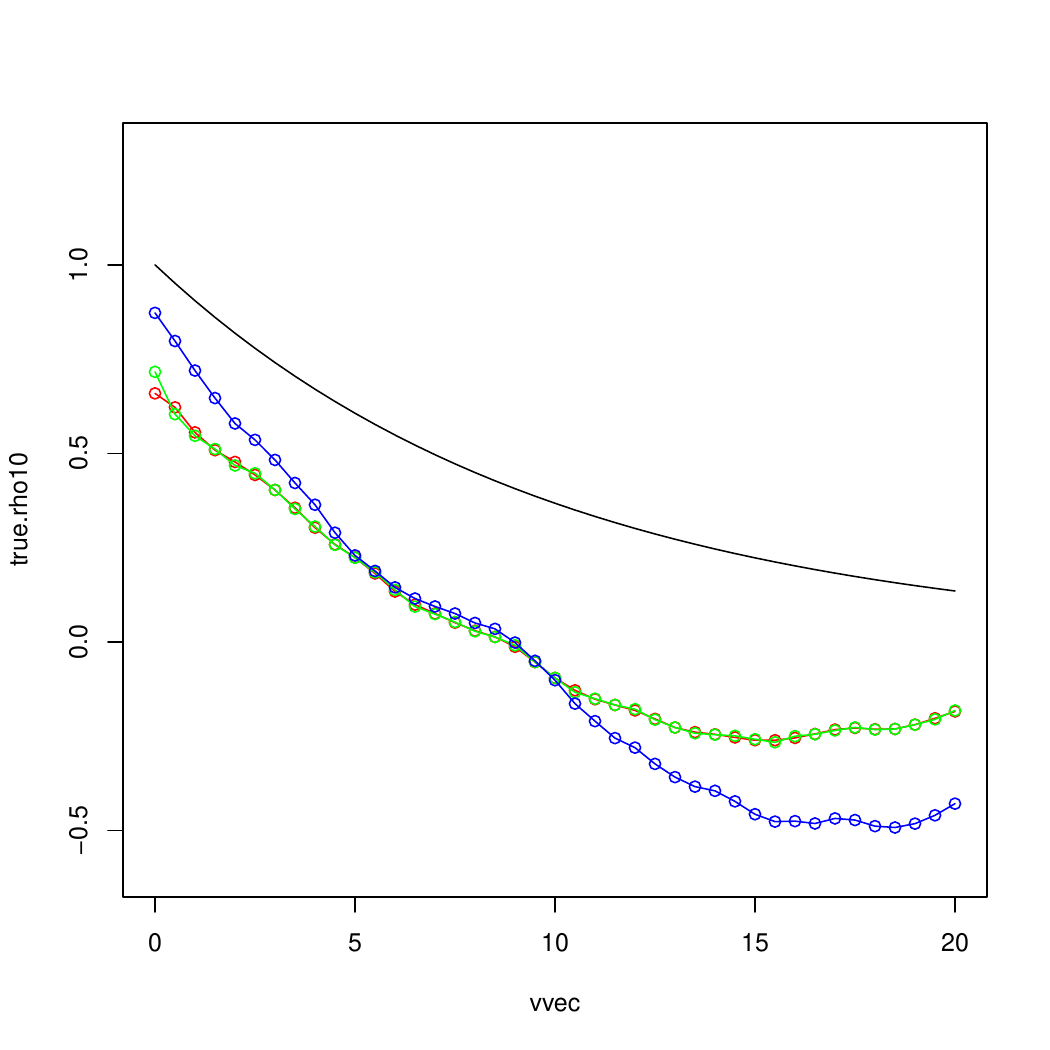}
\caption{Plot of $\widehat{c}_{a,80,n}(u)$ and $\widehat{c}_{HFH}(u)$
  evaluated at $u=0,0.5,1,\ldots,20$. The Black line is the true
  autocovariance $\exp(-|u|/R)$, Red line is $\widehat{c}_{50,80,n}(u)$
  ($a=50$). Green line is $\widehat{c}_{150,80,n}(u)$
  ($a=150$). Blue line is $\widehat{c}_{HRH}(u)$. 
  Left to Right $R=2,5$ and $10$. Top for $n=1000$ and bottom for $n=2000$.
\label{fig:2}}
\end{figure}

In the simulations we simulate from a spatial Gaussian random field
  with autocovariance $c(u) = \exp(-|u|/R)$ ($R$ denotes the range
  parameter) and use $R=2$ (range is $5\%$ of the field), $R=5$
  (range is $R=12.5\%$ of the field)  and $R=10$ (range is
  $25\%$ of the field). We also compare sample sizes $n=1000$ and
  $n=2000$. In Figure \ref{fig:1} we make a plot of
  $\{|J_{n}(\omega_{80,k})|^{2}\}_{k=0}^{150}$ against $\{\omega_{80,k}\}_{k=0}^{150}$ for $R=2,5,10$
  and $n=1000,2000$. The locations $\{u_{j}\}_{j=1}^{n}$ are sampled from a uniform distribution.
We observe that since $J_{n}(\omega)$ is being
  sampled on  a fine frequency grid ($2\pi/(2\times 40)$ compared
  with $2\pi /40$) adjacent values
  of $J_{n}(\omega_{80,k})$ are highly correlated. Furthermore, as
  expected, when $R=10$ the periodogram drops to zero ``faster''
  than when $R=2$. However, in all cases we see that the amplitude of 
$|J_{n}(\omega_{80,k})|^{2}$ drops to nearly zero when $a$ is larger than
$100$. A plot of the periodogram can be used to determine $a$ and 
in Figure \ref{fig:2} we give a plot of the estimator
$\widehat{c}_{a,80,k}(u)$ for $a=50$ and $a=150$  together with the 
\citeA{p:hal-fis-hof-94} estimator, $\widehat{c}_{HRH}(u)$. The plots
are for $R=2,5,10$ and $n=1000$ and $n=2000$. All the estimators are
evaluated at $u=0,0.5,\ldots,20$.

In order to compare the estimators we conduct a
simulation study using the specifications given above, where 500 replications are made for each $(R,n)$-pair. 
To understand how the choice of $a$ influences
the estimator we evaluate $\widehat{c}_{a,\Omega,n}(u)$ for  
$a=50,100,150,200$. We also evaluate $\widehat{c}_{HFH}(u)$
(defined in \citeA{p:hal-fis-hof-94}). The simulations are done over
500 replications. For each simulation (range parameter and sample size)
we calculate the square root average squared error
(SASE), the average (Ave) and the (g)lobal SASE (for
$u=0,0.5,\ldots,20$). They are defined as 
\begin{eqnarray*}
\textrm{SASE}(u) =
\sqrt{\frac{1}{N}\sum_{j=1}^{N}\left\{\widehat{c}^{(j)}(u)-c(u)\right\}^{2}},
 \textrm{Ave}(u) = \frac{1}{N}\sum_{j=1}^{N}\widehat{c}^{(j)}(u)
\textrm{ and } 
\textrm{gSASE} = \frac{1}{41}\sum_{i=0}^{40}\textrm{ASE}\left(\frac{i}{2}\right),
\end{eqnarray*}
where $\widehat{c}^{(j)}(u)$ denotes the estimator based on the $j$th
replication. 

The results of the simulations are reported in
Tables \ref{tab:1}-\ref{tab:6}. We observe that for most values of $a$
the estimator $\widehat{c}_{a,\Omega,n}$ seems to perform better than
the HRH estimator. The sensitivity of the estimator to the choice of
$a$ depends on the location the covariance is estimating. We observe
that when $n=2000$ the SASE is to roughly the same for all choices of
$a$ and all range parameters. However, when $n=1000$ for the estimator
at $u=0$ the SASE is larger for larger choices of $a$ (this is \emph{not}
seen for estimators at other values of $u$, except when the range
parameter is $R=10$, when something similar is seen at $u=2$). An
explanation can be found from the way in which
$\widetilde{c}_{a,\Omega,n}(u)$ is defined. We recall that
\begin{eqnarray*}
\widetilde{c}_{a,\Omega,n}(u)
= \widetilde{Q}_{a,\Omega,\lambda}(e^{iu\cdot}) +
\left(\frac{1}{n}\sum_{j=1}^{n}Z(\ub_{j})^{2}\right)
\times \left( \frac{1}{n}\sum_{k=-a}^{a}e^{iu2\pi k/\Omega}\right).
\end{eqnarray*}
The second term on the right hand side of the above is the ``so
called" finite bias its expectation is approximately equal to 
\begin{eqnarray*}
\frac{\sigma^{2}\Omega}{2\pi n}\int_{-2\pi a/\lambda}^{2\pi
  a/\lambda}e^{i u \omega}d\omega.
\end{eqnarray*}
If $n$ is large or $u\neq 0$ the above bias will be close to zero, however for
small $n$ and $u=0$ then the above may be quite large which may
explain the effect that we see. Despite this the estimator is
not overly sensitive to the choice of $a$, since we are comparing the
performance of the estimator over a very wide range of $a$ ($a=50-200$).  It seems
that the best way to choose $a$ is simply to make a plot of the
periodogram (similar to Figure \ref{fig:1}) and select the $a$ where
most of the amplitudes drop close to zero. 

To understand the effect sample size, $n$, has on the estimation
scheme simulation were conducted for $R=2$ and
$n=1000,2000,4000,6000$. We focussed on $a=50$ and $100$ and also
evaluated the HFH estimator. 100 replications were done for each
$(R,n)$-pair. The results are reported in Table \ref{tab:7}. As the
sample size increases there does not seem to any real change in the
SASE. This observations is supported by the theory developed in this
paper.

%% file: tables.tex
\begin{table}
\centering
\scalebox{0.8}{
\begin{tabular}{|c|c|c|c|c|c|c|c|c|}
\hline
\multicolumn{2}{|c|}{Method} & $c(0)=1$ & $c(2)=0.3678$ & $c(4)=0.1353$ & $c(6)=0.0497$ & $c(8)=0.0183$ & $c(10)=0.0067$ &
gSASE \\
\hline \hline
\multirow{2}{*}{$a = 50$} & SASE & 0.365 & 0.274 & 0.243 & 0.217 & 0.186 & 0.151 & 0.189 \\
         & Ave      &  0.964 & 0.366 & 0.124 & 0.032 & 0.009 & 0.005 & 
                \\
\hline 
\multirow{2}{*}{$a=100$} & SASE & 0.388 & 0.274 & 0.244 & 0.218 & 0.188 & 0.154 & 0.190 \\
  & Ave& 1.056 & 0.350 & 0.126 & 0.036 & 0.008 & 0.001 & \\
\hline 
\multirow{2}{*}{$a=150$} & SASE & 0.423 & 0.273 & 0.244 & 0.220 & 0.188 & 0.153 & 0.191 \\
 & Ave  & 1.119 & 0.342 & 0.125 & 0.039 & 0.008 & 0.001 & \\
\hline 
\multirow{2}{*}{$a=200$} & SASE & 0.463 & 0.275 & 0.246 & 0.220 & 0.188 & 0.153 & 0.192 \\
 &   Ave   &  1.177 & 0.350 & 0.126 & 0.037 & 0.007 & 0.003 & \\
\hline
 \multirow{2}{*}{HFH} & SASE &0.465&0.350 & 0.301 & 0.286 & 0.254 & 0.223 & 0.280 \\
 & Ave & 1.149 & 0.437 & 0.165 & 0.055 & 0.027 & 0.027 & \\
\hline
\end{tabular}}
\caption{$n=1000$, $R = 2$. Estimates evaluated at $u=0,2,\ldots,10$
  and globally.\label{tab:1}}
\end{table}

\begin{table}
\centering
\scalebox{0.8}{
\begin{tabular}{|c|c|c|c|c|c|c|c|c|}
\hline
\multicolumn{2}{|c|}{Method} & $c(0)=1$ & $c(2)=0.3678$ & $c(4)=0.1353$ & $c(6)=0.0497$ & $c(8)=0.0183$ & $c(10)=0.0067$ &
gSASE \\
\hline \hline
\multirow{2}{*}{$a = 50$} & SASE & 0.330 & 0.224 & 0.197 & 0.164 & 0.168 & 0.164 & 0.171 \\
     &  Ave   & 0.950 & 0.362 & 0.123 & 0.043 & 0.030 & 0.019 & \\
\hline
\multirow{2}{*}{$a=100$} & SASE & 0.333 & 0.225 & 0.197 & 0.164 & 0.169 & 0.166 & 0.171 \\
   &Ave& 1.017 & 0.349 & 0.124 & 0.047 & 0.029 & 0.016 & \\
\hline 
\multirow{2}{*}{$a=150$} & SASE & 0.343 &0.224 & 0.199 & 0.166 & 0.168 & 0.166 & 0.172 \\
      &Ave  & 1.055 & 0.346 & 0.123 & 0.049 & 0.029 & 0.016 & \\
\hline
\multirow{2}{*}{$a=200$}& SASE & 0.357 & 0.226 & 0.200 & 0.165 & 0.168 & 0.166 & 0.172 \\
       &Ave& 1.087 & 0.350 & 0.124 & 0.048 & 0.029 & 0.017 & \\
\hline 
\multirow{2}{*}{HFH} & SASE & 0.406 & 0.281 & 0.241 & 0.202 & 0.226 & 0.241 & 0.253 \\
   &Ave& 1.141 & 0.436 & 0.158 & 0.062 & 0.047 & 0.041 & \\ 
\hline
 \end{tabular}}
\caption{$n=2000$, $R = 2$. Estimates evaluated at $u=0,2,\ldots,10$
  and globally. \label{tab:2}}
\end{table}

\begin{table}
\centering
\scalebox{0.8}{
\begin{tabular}{|c|c|c|c|c|c|c|c|c|}
\hline
\multicolumn{2}{|c|}{Method} & $c(0)=1$ & $c(2)=0.670$ & $c(4)=0.449$ & $c(6)=0.301$ & $c(8)=0.202$ & $c(10)=0.135$ &
gSASE \\
\hline \hline
\multirow{2}{*}{$a=50$} & SASE & 0.564 & 0.475 & 0.417 & 0.376 & 0.314 & 0.266 & 0.312 \\
 &Ave& 1.045 & 0.655 & 0.394 & 0.237 & 0.159 & 0.106 &\\
\hline
\multirow{2}{*}{$a=100$} & SASE & 0.575 & 0.443 & 0.401 & 0.364 & 0.306 & 0.266 & 0.303 \\
  &Ave&  1.114 & 0.645 & 0.395 & 0.240 & 0.160 & 0.103 & \\
\hline
\multirow{2}{*}{$a=150$} & SASE & 0.617 & 0.441 & 0.401 & 0.366 & 0.304 & 0.265 & 0.304 \\
  &Ave&  1.171 & 0.637 & 0.395 & 0.242 & 0.160 & 0.102 & \\
\hline 
\multirow{2}{*}{$a=200$} & SASE & 0.660 & 0.444 & 0.402 & 0.368 & 0.305 & 0.266 & 0.306 \\
   &Ave& 1.227 & 0.645 & 0.396 & 0.241 & 0.160 & 0.104 & \\
\hline
\multirow{2}{*}{HFH} & SASE &  0.747 & 0.609 & 0.519 & 0.471 & 0.422 & 0.384 & 0.456 \\
   &Ave& 1.236 & 0.790 & 0.490 & 0.308 & 0.227 & 0.176 & \\
\hline
 \end{tabular}}
\caption{$n=1000$, $R = 5$. Estimates evaluated at $u=0,2,\ldots,10$
  and globally. \label{tab:3}}
\end{table}

\begin{table}
\centering
\scalebox{0.8}{
\begin{tabular}{|c|c|c|c|c|c|c|c|c|}
\hline
\multicolumn{2}{|c|}{Method} & $c(0)=1$ & $c(2)=0.670$ & $c(4)=0.449$ & $c(6)=0.301$ & $c(8)=0.202$ & $c(10)=0.135$ &
gSASE \\
\hline \hline
\multirow{2}{*}{$a=50$} & SASE &  0.492 & 0.445 & 0.403 & 0.363 & 0.327 & 0.299 & 0.308 \\
  &Ave& 0.965 & 0.605 & 0.363 & 0.217 & 0.127 & 0.074 & \\
\hline
\multirow{2}{*}{$a=100$} & SASE & 0.508 & 0.445 & 0.402 & 0.364 & 0.328 & 0.297 & 0.308 \\
  &Ave& 1.007 & 0.598 & 0.364 & 0.219 & 0.127 & 0.073 & \\ 
\hline
\multirow{2}{*}{$a=150$} & SASE & 0.525 & 0.443 & 0.402 & 0.364 & 0.327 & 0.297 & 0.308 \\
  &Ave& 1.038 & 0.594 & 0.363 & 0.220 & 0.127 & 0.073 & \\ 
\hline
\multirow{2}{*}{$a=200$} & SASE & 0.540 & 0.443 & 0.402 & 0.364 & 0.327 & 0.298 & 0.309 \\
   &Ave& 1.065 & 0.598 & 0.363 & 0.219 & 0.127 & 0.074 & \\
\hline 
\multirow{2}{*}{HFH} & SASE & 0.630 & 0.547 & 0.483 & 0.450 & 0.430 & 0.421 & 0.445 \\
       &Ave& 1.163 & 0.741 & 0.455 & 0.282 & 0.187 & 0.136    & \\
\hline
 \end{tabular}}
\caption{$n=2000$, $R = 5$. Estimates evaluated at $u=0,2,\ldots,10$
  and globally. \label{tab:4}}
\end{table}

\begin{table}
\centering
\scalebox{0.8}{
\begin{tabular}{|c|c|c|c|c|c|c|c|c|}
\hline
\multicolumn{2}{|c|}{Method} & $c(0)=1$ & $c(2)=0.819$ & $c(4)=0.670$ & $c(6)=0.549$ & $c(8)=0.449$ & $c(10)=0.368$ &
gSASE \\
\hline \hline
\multirow{2}{*}{$a=50$} & SASE &  0.605 & 0.540 & 0.495 & 0.457 & 0.424  &0.385 & 0.397 \\
   &Ave& 1.056 & 0.789 & 0.594 & 0.447 & 0.340 & 0.253 & \\
\hline
\multirow{2}{*}{$a=100$} & SASE & 0.634 & 0.539 & 0.495 & 0.456 & 0.423 & 0.386 & 0.397 \\
  &Ave& 1.117 & 0.780 & 0.596 & 0.450 & 0.340 & 0.251 & \\ 
\hline
\multirow{2}{*}{$a=150$} & SASE &  0.668 & 0.535  &0.495 & 0.457 & 0.424 & 0.386 & 0.398 \\
   &Ave& 1.172 & 0.774 & 0.594 & 0.452 & 0.339 & 0.250 & \\
\hline
\multirow{2}{*}{$a=200$} & SASE & 0.705 & 0.536 & 0.495 & 0.457 & 0.423 & 0.386 & 0.399 \\
   &Ave& 1.227 & 0.781 & 0.596 & 0.451 & 0.339 & 0.251 & \\
\hline 
\multirow{2}{*}{HFH} & SASE & 0.828 & 0.696 & 0.601 & 0.550 & 0.529 & 0.496 & 0.549\\
 &Ave & 1.275 & 0.964 & 0.728 & 0.561 & 0.462 & 0.388 & \\
\hline
 \end{tabular}}
\caption{$n=1000$, $R = 10$. Estimates evaluated at $u=0,2,\ldots,10$
  and globally. \label{tab:5}}
\end{table}

\begin{table}
\centering
\scalebox{0.8}{
\begin{tabular}{|c|c|c|c|c|c|c|c|c|}
\hline
\multicolumn{2}{|c|}{Method} & $c(0)=1$ & $c(2)=0.819$ & $c(4)=0.670$ & $c(6)=0.549$ & $c(8)=0.449$ & $c(10)=0.368$ &
gSASE \\
\hline \hline
\multirow{2}{*}{$a=50$} & SASE &  0.737 & 0.674 & 0.608 & 0.551 & 0.502 & 0.458 & 0.471 \\
   &Ave& 0.987 & 0.742 & 0.551 & 0.408 & 0.308 & 0.229 & \\
\hline
\multirow{2}{*}{$a=100$} & SASE & 0.757 & 0.673 & 0.609 & 0.553 & 0.504 & 0.460 & 0.473 \\
   &Ave& 1.021 & 0.737 & 0.552 & 0.409 & 0.308 & 0.227 & \\
\hline
\multirow{2}{*}{$a=150$} & SASE & 0.779 & 0.669  & 0.609 & 0.554 & 0.503 & 0.461 & 0.473 \\
   &Ave& 1.050 & 0.733 & 0.552 & 0.410 & 0.308& 0.227 & \\ 
\hline 
\multirow{2}{*}{$a=200$} & SASE & 0.803 & 0.673 & 0.608 & 0.554 & 0.504 &0.459 & 0.475 \\
  &Ave& 1.076 & 0.737 & 0.552 & 0.409 & 0.308 & 0.228 & \\
\hline 
\multirow{2}{*}{HFH} & SASE &  0.964 & 0.853 & 0.751 & 0.692 & 0.663 & 0.646 & 0.676 \\
 &Ave& 1.206 & 0.910 & 0.671 & 0.504 & 0.416 & 0.353 & \\
\hline
 \end{tabular}}
\caption{$n=2000$, $R = 10$. Estimates evaluated at $u=0,2,\ldots,10$
  and globally. \label{tab:6}}
\end{table}

\begin{table}
\centering
\scalebox{0.8}{
\begin{tabular}{|c|c|c|c|c|c|c|c|c|c|}
\hline
\multicolumn{2}{|c|}{Method} & $c(0)=1$ & $c(2)=0.3678$ & $c(4)=0.1353$ & $c(6)=0.0497$ & $c(8)=0.0183$ & $c(10)=0.0067$ &
gSASE \\
\hline \hline
\multicolumn{8}{|c|}{n=1000} \\
\hline
\multirow{2}{*}{$a = 50$} & SASE & 0.334 & 0.249 & 0.203 & 0.184 & 0.188 & 0.170 & 0.180 \\
      &Ave& 0.968 &  0.355 &  0.100 &  0.001  &0.011  & 0.011  & \\
\hline 
\multirow{2}{*}{$a=100$} & SASE & 0.388 & 0.274 & 0.244 & 0.218 & 0.188 & 0.154 & 0.190 \\
  &Ave& 1.056 & 0.350 & 0.126 & 0.036 & 0.008 & 0.001 & \\
\hline 
\multirow{2}{*}{HFH} & SASE&0.465&0.350 & 0.301 & 0.286 & 0.254 & 0.223 & 0.280 \\
   &Ave& 1.149 & 0.437 & 0.165 & 0.055 & 0.027 & 0.027 & \\
\hline
\hline
\multicolumn{8}{|c|}{n=2000} \\
\hline
\multirow{2}{*}{$a=50$} & SASE& 0.321 & 0.235 & 0.209 & 0.212 & 0.227 &  0.209 & 0.192 \\
  &Ave&  0.943  & 0.356  & 0.118 &  0.026 & -0.004 & -0.004 &  \\
\hline
\multirow{2}{*}{$a=100$} & SASE &  0.323 & 0.236 & 0.208 & 0.212 & 0.226 & 0.210 & 0.192 \\
  &Ave&  1.010 &  0.343 &  0.119 &  0.031 & -0.005 & -0.005 &  \\
\hline 
 \multirow{2}{*}{HFH} & SASE & 0.400 & 0.295 & 0.259  & 0.267 & 0.305 & 0.313 & 0.287 \\
  &A  & 1.132 & 0.431 & 0.153 & 0.039 & 0.001 & 0.006 & \\
\hline\hline
\multicolumn{8}{|c|}{n=4000} \\
\hline
\multirow{2}{*}{$a=50$} & SASE & 0.323 & 0.242 & 0.202 & 0.197 & 0.184 & 0.144 & 0.174 \\
   &Ave& 0.899 & 0.317 & 0.095 & 0.033 & 0.014 & 0.027 & \\
\hline
\multirow{2}{*}{$a=100$} & SASE & 0.313 & 0.245 & 0.202 & 0.197 & 0.185 & 0.144 & 0.174 \\
   &Ave& 0.950 & 0.307 & 0.097 & 0.034 & 0.014 & 0.025 & \\
\hline 
 \multirow{2}{*}{HFH} & SASE & 0.368 & 0.288 & 0.239 & 0.240 & 0.241 & 0.217 & 0.251 \\
   &Ave& 1.076 & 0.382 & 0.128 & 0.049 & 0.027 & 0.053 & \\
\hline \hline 
\multicolumn{8}{|c|}{n=6000} \\
\hline
\multirow{2}{*}{$a=50$} &SASE & 0.322 & 0.235 & 0.211 & 0.195 & 0.180 & 0.173 & 0.181 \\
   & A&  0.951  & 0.363 &  0.128 &  0.045 &  0.002 & -0.001 & \\
\hline 
\multirow{2}{*}{$a=100$} & SASE& 0.319 & 0.236 & 0.211 & 0.195 & 0.180 & 0.172 & 0.181 \\
   &Ave&  1.001 &  0.352 &  0.128 & 0.049 &  0.001 & -0.003 & \\
\hline 
\multirow{2}{*}{HFH} & SASE& 0.429 & 0.312 & 0.259 & 0.237 & 0.234 & 0.242 & 0.273 \\
  &Ave& 1.150 & 0.443 & 0.174 & 0.075 & 0.023 & 0.018 & \\
\hline
\end{tabular}}
\caption{$n=1000$, $2000$, $4000$ and $6000$; $R = 2$. Estimates evaluated at $u=0,2,\ldots,10$
  and globally. Simulations conducted over $100$ replications.\label{tab:7}}
\end{table}